\documentclass[11pt]{amsart}
\allowdisplaybreaks
\usepackage{microtype}
\usepackage{amsmath,amsfonts,amsbsy,amsgen,amscd,mathrsfs,amssymb,amsthm}
\usepackage{enumerate,mathtools}
\usepackage{tabularx}
\usepackage{array}
\newcolumntype{P}[1]{>{\centering\arraybackslash}p{#1}}
\usepackage{makecell}

\usepackage{stmaryrd}
\usepackage{graphicx}
\usepackage{fullpage}
\usepackage{xurl}
\usepackage{mathrsfs}
\usepackage{float}
\usepackage{bbm}
\usepackage{quiver}
\usepackage{multirow}
\usepackage[colorlinks=true,allcolors=blue]{hyperref}

\usepackage[
backend=biber,
style=alphabetic,
maxalphanames=5,
maxnames=5
]{biblatex}
\makeatletter
\renewcommand*{\eqref}[1]{%
  \hyperref[{#1}]{\textup{\tagform@{\ref*{#1}}}}%
}
\makeatother

\usepackage{tikz}
\usetikzlibrary{nfold}
\usetikzlibrary{shadings}
\usetikzlibrary{arrows.meta,calc}

\numberwithin{equation}{section}

\newtheorem{theorem}{Theorem}[section]
\newtheorem{proposition}[theorem]{Proposition}

\newtheorem{conjecture}[theorem]{Conjecture}
\newtheorem{corollary}[theorem]{Corollary}
\newtheorem{lemma}[theorem]{Lemma}
\theoremstyle{definition}
\newtheorem{definition}[theorem]{Definition}
\newtheorem{claim}[theorem]{Claim}

\newtheorem{example}[theorem]{Example}
\newtheorem{remark}[theorem]{Remark}

\newtheorem*{claim*}{Claim}

\newcommand{\A}{\mathbb{A}}
\newcommand{\C}{\mathbb{C}}

\newcommand{\E}{\mathbb{E}}
\newcommand{\F}{\mathbb{F}}
\newcommand{\G}{\mathbb{G}}

\newcommand{\N}{\mathbb{N}}

\newcommand{\Q}{\mathbb{Q}}
\newcommand{\R}{\mathbb{R}}

\newcommand{\Z}{\mathbb{Z}}

\DeclareMathOperator{\Aff}{Aff}

\DeclareMathOperator{\Aut}{Aut}

\DeclareMathOperator{\BiMod}{BiMod}

\DeclareMathOperator{\Ch}{Ch}

\DeclareMathOperator{\CHur}{CHur}
\DeclareMathOperator{\Cl}{Cl}

\DeclareMathOperator{\cofib}{cofib}

\DeclareMathOperator{\Conf}{Conf}

\DeclareMathOperator{\diag}{diag}
\DeclareMathOperator{\disc}{disc}

\DeclareMathOperator{\ev}{ev}

\DeclareMathOperator{\Frob}{Frob}
\DeclareMathOperator{\Fun}{Fun}

\DeclareMathOperator{\Gr}{Gr}
\DeclareMathOperator{\Hom}{Hom}
\DeclareMathOperator{\Hur}{Hur}

\DeclareMathOperator{\id}{id}
\DeclareMathOperator{\im}{im}

\DeclareMathOperator{\Jac}{Jac}

\DeclareMathOperator{\LMod}{LMod}
\DeclareMathOperator{\Map}{Map}

\DeclareMathOperator{\MHur}{MHur}
\DeclareMathOperator{\Mod}{Mod}

\DeclareMathOperator{\Norm}{Norm}

\DeclareMathOperator{\ord}{ord}

\DeclareMathOperator{\Pic}{Pic}

\DeclareMathOperator{\pr}{pr}

\DeclareMathOperator{\rk}{rk}
\DeclareMathOperator{\RMod}{RMod}

\DeclareMathOperator{\Sh}{Sh}

\DeclareMathOperator{\Spec}{Spec}

\DeclareMathOperator{\Surj}{Surj}

\DeclareMathOperator{\Tot}{Tot}
\DeclareMathOperator{\totcofib}{totcofib}
\DeclareMathOperator{\Tor}{Tor}
\DeclareMathOperator{\tr}{tr}

\newcommand{\GL}{\operatorname{GL}}

\newcommand{\SL}{\operatorname{SL}}

\newcommand{\Sp}{\operatorname{Sp}}

\usepackage[shortlabels]{enumitem}

\usepackage[OT2,T1]{fontenc}
\DeclareSymbolFont{cyrletters}{OT2}{wncyr}{m}{n}
\DeclareMathSymbol{\Sha}{\mathalpha}{cyrletters}{"58}

\usepackage{array, booktabs}

\title{Optimal homological stability: bounding secondary terms of Cohen-Lenstra moments and $S_d$-fields over $\F_q[t]$}

\author{Zhao Yu Ma}
\address{Princeton University, Princeton, NJ, USA}
\email{zm5336@princeton.edu}
\begin{document}
\begin{abstract}
We prove an explicit stability line for the homology $H_i(\Hur_n^c,\Q)$ of any Hurwitz space associated to a finite rack $c$ depending only on its homology up to some finite $n$. Moreover, as the range of $n$ increases, the slope of the resulting stability line converges to the optimal slope. Thus, for any Hurwitz space, a finite computation will yield near-optimal bounds on secondary terms for the corresponding arithmetic statistics problem. Our theorem also applies more generally to higher order stability, Frobenius weights, and braided vector spaces. 

As a special case of our theorem, we can deduce stability lines depending only on when $H_0(\Hur^c_n,\Q)$ stabilizes, which are sharp in some cases. We apply this to two families of Hurwitz spaces to obtain bounds on the secondary term for (1) the number of $S_d$-degree $d$ function fields with squarefree discriminant, and (2) Cohen-Lenstra moments.
\end{abstract}

\maketitle

\hypersetup{linktocpage=true}
\tableofcontents

\newcommand{\drawnaxis}[1]{
  \draw[->] (0,0) -- (#1+0.5,0) node[right] {$n$};
  \foreach \x in {1,...,#1}
    \node[below] at (\x-0.5,0) {\x};
}

\newcommand{\drawiaxis}[1]{
  \draw[->] (0,0) -- (0,#1+1.5) node[above] {$i$};
  \foreach \y in {1,...,\numexpr#1+1\relax}
    \pgfmathtruncatemacro{\yy}{\y-1}
    \node[left] at (0,\y-0.5) {\yy};
}

\newcommand{\drawblock}[6]{%
  \foreach \n in {\numexpr#1-1\relax,...,\numexpr#2-1\relax}{
    \foreach \i in {#3,...,#4}{
      \fill[#5] (\n,\i) rectangle (\n+1,\i+1);
      \node at (\n+0.5,\i+0.5) {$#6$};
    }
  }
}

\newcommand{\drawstaircase}[5]{
  \foreach \n in {\numexpr #1-1\relax,...,\numexpr#2-1\relax}{
    \foreach \i in {#3,...,\numexpr #3+\n+1-#1\relax}{
      \fill[#4] (\n,\i) rectangle (\n+1,\i+1);
      \node at (\n+0.5,\i +0.5) {$#5$};
    }
  }
}

\newcommand{\cell}[4]{
    \fill[#4] (#1,#2) rectangle (#1-1,#2+1);
    \node at (#1-0.5,#2+0.5) {$#3$};
  }
  
\newcommand{\smallcell}[4]{
    \fill[#4] (#1,#2) rectangle (#1-1,#2+1);
    \node at (#1-0.5,#2+0.5) {\fontsize{9}{11}\selectfont $#3$};
  }

\newcommand{\tinycell}[4]{
    \fill[#4] (#1,#2) rectangle (#1-1,#2+1);
    \node at (#1-0.5,#2+0.5) {\fontsize{6}{11}\selectfont $#3$};
  }

\newcommand{\verytinycell}[4]{
    \fill[#4] (#1,#2) rectangle (#1-1,#2+1);
    \node at (#1-0.5,#2+0.5) {\scalebox{0.8}{\fontsize{3}{1}\selectfont $#3$}};
  }
  
\section{Introduction}
Following Ellenberg-Venkatesh-Westerland \cite{EVW16}, many problems in arithmetic statistics over $\F_q[t]$ can be reduced to understanding the homology of the corresponding Hurwitz spaces via the Grothendieck-Lefschetz trace formula and an étale-to-complex comparison theorem. In particular, proving homological stability for the Hurwitz space tells us the leading term of the counting problem when $q\rightarrow \infty$. \cite{EVW16} proved this for $\Hur^c_n$ where $c\subseteq G$ is a conjugacy class which satisfies the non-splitting condition, leading to results on the Cohen-Lenstra heuristics. Here, the non-splitting condition is that for any subgroup $H\le G$, $c\cap H$ is either empty or a single $H$-conjugacy class.

Furthermore, knowing the stable homology can tell us that the leading term dominates not only in the limit $q \rightarrow \infty$, but for all sufficiently large $q$. This was done recently in a series of papers of Landesman-Levy \cite{landesman2025homologicalstabilityhurwitzspaces,landesman2025cohenlenstramomentsfunctionfields,landesman2025stablehomologyhurwitzmodules} where the stable homology of Hurwitz space was shown to be isomorphic to (multicolored) configuration space for all Hurwitz spaces $\Hur^c_n$ associated to a general finite rack $c$. Importantly, this was also done without the non-splitting condition, and led to the proof of various arithmetic statistics conjectures. There has also been a lot of other work related to Hurwitz spaces, such as \cite{ETW17,Randal-Williams_2020,ellenberglandesman,Liu_Wood_Zureick-Brown_2024,ma2026weightfiltrationhurwitzspaces,sawin2026distributionsunramifiedextensionsglobal}. 

However, there is still much more to be done. Apart from the leading term, most arithmetic statistical problems are expected to have a \textit{secondary term} of smaller exponent, but this is currently only known for cubic and quartic fields. For the classical Davenport-Heilbronn problem of counting the number of $S_3$-cubic fields \cite{Davenport_Heilbronn_1969,Davenport_Heilbronn_1971}, Roberts \cite{Roberts_2000} conjectured the existence of a secondary term, and an asymptotic count of $c_1X+c_2X^{5/6}+O(X^{2/3+\epsilon})$ was later proven by work of  \cite{Bhargava_Shankar_Tsimerman_2012,Taniguchi_Thorne_2013,Bhargava_Taniguchi_Thorne_2023}. The analogous result for cubic function fields was proven in \cite{zhaophd,kural2025geometrysecondarytermsarithmetic,ahlquist2025countingfunctioncubicfunction}. For $S_4$-quartic fields, the smoothed counting function was shown to be $c_1X+c_2X^{5/6}+O(X^{13/16+\epsilon})$ in \cite{shankar2025secondarytermscountingfunctions,shankar2025secondarytermscountingfunctions2}, and a similar result for quartic algebras over $\F_q[t]$ was proven in \cite{chang2024hurwitzspacesnicholsalgebras}. These proofs rely on the explicit parametrization of cubic and quartic rings in terms of prehomogeneous vector spaces due to Bhargava \cite{Bhargava_2004a,Bhargava_2004b,Bhargava_2004c,Bhargava_2008}, and hence they are difficult to generalize. 

There are only conjectural secondary terms in other cases. For example, data from \cite{malle2024computationaldatamathfraksnextensions} gives numerical evidence for a secondary term for $S_d$-degree $d$ extensions when $d\ge 5$; for example, it seems that the secondary term for $S_5$-quintic fields is roughly $X^{5/6}$. On the other hand, the best upper bound for the secondary term proven for $S_5$-quintic fields over number fields so far is only $c_1X+O(X^{199/200+\epsilon})$ in \cite{SHANKAR_TSIMERMAN_2014}, which is far from what can be predicted from data. Another example is the case of Cohen-Lenstra moments, where for odd $k$, \cite{Hough_2019} predicts that the \textit{average} of $k$-torsion in class groups over quadratic fields has a secondary term on the order of $X^{-\frac 1 2+
\frac 1 k}$, but apart from the $k=3$ case, no bound close to this was proven in either the number field or function field case. 

We introduce a new method to prove upper bounds on the exponent of secondary terms for \textit{any} counting problem over $\F_q[t]$ that is parametrized by some Hurwitz space. By bounding secondary terms, we mean bounding the full remainder after subtracting the main term, including any secondary or higher-order terms as well as the remaining error. For example, we show for large $q$ that:
\begin{itemize}
\item The secondary term of the number of $S_5$-quintic function fields with squarefree discriminant has exponent $\le 0.886$. Similarly, in degree $6$ we obtain the bound $\le 0.917$, in degree $7$ we have $\le 0.947$, and we also establish bounds for all degrees $d$.
\item For all odd $k$, the secondary term of the Cohen-Lenstra $k$-moment has exponent $\le -1/10+\epsilon$. For small $k$, we can obtain better bounds like $\le -1/4+\epsilon$ for $k=7$. We also obtain bounds for moments of odd abelian groups $H$ whose exponent depends only on the rank $r$ of $H$. 
\end{itemize}

We will state these more precisely later in Section \ref{sec: intro arithmetic applications}, together with various other arithmetic results. The reader who is only interested in arithmetic applications should skip to this section. 

More importantly, this method can, in theory, given sufficient computational power, prove upper bounds on the exponent that are \textit{optimal up to an arbitrarily small error}, assuming that there is no unexpected cancellation in Frobenius eigenvalues. This optimality result is a direct consequence of our main theorem, Theorem \ref{thm: main theorem}, which is a homological stability ``extrapolation'' theorem that given knowledge of homological stability for a finite range of $n$ lets us extrapolate this to a stability line of some slope. We remark that prior homological stability results also do prove a stability line with nontrivial slope; however, the slope is often very far from the optimal slope. We discuss the optimality of our technique in Section \ref{sec: intro optimality}.

Our results also work for higher order stability, provided we know that the higher order stabilization elements come from the original chains $A$ and not the quotient $A/(U_1,\ldots, U_k)$ where we quotient by the primary stabilization elements. \textit{If} this were true and \textit{secondary stability were known} for some Hurwitz space, then we prove that one can also obtain a secondary stability line with near-optimal slope with sufficient computation, and this would give near-optimal upper bounds for tertiary terms which appear after the secondary term. Unfortunately, no such secondary stability is known yet. However, this does give a hint at how one might obtain a ``power series'' expansion for counting functions over function fields, at least computationally, by proving homological stability and computing the stable homology as in \cite{EVW16,landesman2025homologicalstabilityhurwitzspaces,landesman2025cohenlenstramomentsfunctionfields,landesman2025stablehomologyhurwitzmodules}, improving the stability slope to the optimal value, and repeating with secondary stability and so on.

Another important application of our main theorem is that one can prove a non-trivial stability line given only knowledge of when $H_0$ stabilizes. This is especially useful when proving results in families as direct computation for the homology of each individual Hurwitz space is impossible, but one could potentially deduce a closed form for $H_0$ across the family. We discuss this in Section \ref{sec: intro stab from R}, and eventually apply this to the two examples of Cohen-Lenstra moments and $S_d$-degree $d$ function fields.

This paper is the sequel to the author's previous work \cite{prequel}, where an analogue of the main theorem was proved in the special case of homological vanishing using the formalism of bar-complexes (and the dual notion of cobar-complexes) developed in \cite{ETW17}. Although the present proof builds on the ideas of \cite{prequel}, the passage from vanishing to stability has nontrivial difficulties: one must work with the cofiber of the stabilization map rather than only with the Hurwitz spaces themselves. In terms of cobar-complexes, this means working with the cone of a map between cobar-complexes. It was not clear a priori how to adapt the multi-block filtration of the dual bar-complex constructed in \cite[Section 3]{prequel} to this setting.

The key technical innovation is a rather specific way to overcome this obstruction, which we briefly sketch here. We introduce a two-sided cone of the cobar-complex, involving both the left and right stabilization maps, and construct a two-block filtration on it whose graded pieces are tensor products involving one-sided cones of smaller cobar-complexes. Homological vanishing for these smaller one-sided cones then implies homological vanishing for the larger two-sided cone. A long exact sequence subsequently gives vanishing for the larger one-sided cone. Finally, iterating this process proves the result. A further difficulty arises in higher-order homological stability, where one must handle multiple iterated cofibers. Here, we use the theory of higher algebra to keep track of the resulting higher homotopies.

\subsection{Main homological theorem}\label{sec: intro main theorem} We will state a simplified version of our main theorem here so as to avoid notation from higher algebra. The general version applies more generally to braided vector spaces with stabilization maps by $\E_2$-central elements, as well as to higher order stability with iterated quotients. There are also versions which apply to racks with multiple components (multiple conjugacy classes) with stability in a chosen direction. We will present these in Section \ref{sec: general results}. 

For the purposes of the introduction, we consider the situation where $c$ is a finite rack and $\Hur^c_n$ is the corresponding Hurwitz space over $\C$ with $n$ points, and $\Hur^c\coloneqq \sqcup_{n\ge 0} \Hur^c_n$ is the union over all $n$. Just for the introduction, we let $c$ be non-splitting to avoid the case where $\Hur^c_n$ has polynomial stability as that requires iterated quotients.

Let $R\coloneqq H_0(\Hur^c,\Q)$ be the ring of connected components, which acts on the homology $H_*(\Hur^c,\Q)$ by the multiplicative structure of Hurwitz space. Let $U\in H_0(\Hur^c_a,\Q)$ be an element in $R$ where $a>0$, so that we can write $U=\sum b_i\gamma_i$ as a linear combination of components $\gamma_i\in \pi_0\Hur^c_a$ which are central in $\pi_0\Hur^c_a$. This also implies that $U$ is central in $R$. Then, \cite[Lemma 4.3.1]{landesman2025cohenlenstramomentsfunctionfields} tells us that for all $n,i\ge 0$, both left and right multiplication by $U$ give homotopic maps 
\begin{equation}\label{eqn: stab map}
U\colon H_i(\Hur^c_n,\Q)\rightarrow H_i(\Hur^c_{n+a},\Q).
\end{equation}

We say that $H_*(\Hur^c,\Q)$ is \textit{homologically stable} with respect to the stabilization map $U$ if for each $i$, the map $U$ induces isomorphisms of homology $H_i(\Hur^c_n,\Q)\xrightarrow \sim H_i(\Hur^c_{n+a},\Q)$ for all sufficiently large $n$. Furthermore, we say that it is homologically stable with \textit{slope} $s$ if there is some constant $t$ where these are isomorphisms when $sn+t\ge i$, and we call this equation the stability line. We remark that we will only be concerned with slopes in the range $0<s\le 1$, where the upper bound is because the homology will vanish for $i\ge n$. For example, when $c\subseteq G$ is a (non-splitting) conjugacy class, \cite{EVW16} proved that $H_*(\Hur^c,\Q)$ is homologically stable with respect to $U_D=\sum_{g\in c}g^{|g|D}$ for some sufficiently large $D$, with some explicitly computable but non-optimal slope. It was later shown in \cite[Proposition A.3.1]{ellenberglandesman} that we can take any $D\ge 1$.

We now state our (simplified) main theorem, in analogy to \cite[Theorem 1.2]{prequel}.
\begin{definition}\label{def: intro stabilize in staircase}
Let $N\ge 1$ and $I \ge 0$ be integers with $I<N$. We say that $H_*(\Hur^c,\Q)$ stabilizes in an \textit{$(N,I)$ staircase} with respect to $U$ if for all $n\le N$, we have that the stabilization map \eqref{eqn: stab map} is an isomorphism for $i\le I-N+n$ and a surjection for $i=I-N+n+1$.
\end{definition}
\begin{theorem}\label{thm: main theorem}
Suppose that $H_*(\Hur^c,\Q)$ has homological stability with respect to $U=\sum b_i\gamma_i$ where $\gamma_i\in \pi_0\Hur^c_a$ are central. Also suppose that $H_*(\Hur^c,\Q)$ stabilizes in an $(N,I)$ staircase. Then, it is homologically stable with slope $\frac{I+1}{N+1}$. In particular, for all $n\ge N$, the stabilization map \eqref{eqn: stab map} is an isomorphism for all $i$ satisfying
\begin{equation*}
i\le (I+1)\left\lfloor \frac{n+1}{N+1}\right\rfloor+ \max\left(I+\left\{ \frac{n+1}{N+1}\right\}-N\ ,\ 0\right)-1. 
\end{equation*}
and a surjection for all $i$ satisfying
\begin{equation*}
i=(I+1)\left\lfloor \frac{n+1}{N+1}\right\rfloor+ \max\left(I+\left\{ \frac{n+1}{N+1}\right\}-N\ ,\ 0\right). 
\end{equation*}
Here, we denote $\{\frac a b\}$ to be the remainder when $a$ is divided by $b$.
\end{theorem}
We remark that this simplified main theorem is in some sense orthogonal to the recent progress in Landesman-Levy. Indeed, it can be proven using only the methods developed in \cite{EVW16,ETW17}, without higher algebra and the knowledge of the stable homology. Nevertheless, knowing the stable homology is essential for arithmetic applications of our theorem, and is convenient in generalizations to iterated quotients, and other corollaries.

Let us demonstrate this in the classical Davenport-Heilbronn situation where $c=c^{\text{trn}}_{S_3}$ is the conjugacy class of transpositions in $S_3$, which counts $S_3$-cubic covers with squarefree discriminant. We choose $D=1$ with $U\coloneqq U_1=\sum_{g\in c}g^{|g|}$ as defined above, and take the iterated map $U^e$ as our stabilization map for some sufficiently large $e$. It is well-known that $H_0(\Hur^c_n,\Q)$ stabilizes after $n=3$ (this is also a special case of Theorem \ref{thm: cohen lenstra stab R}), so the stabilization map \eqref{eqn: stab map} is surjective for $(n,i)=(3,0)$ and is an isomorphism for $(n,i)=(4,0)$. Furthermore, if $e$ is sufficiently large, one can show that the stabilization map is surjective for $(n,i)=(4,1)$, as we will see from the discussion after Corollary \ref{cor: stab from R}. Hence, we can apply Theorem \ref{thm: main theorem} for $(N,I)=(4,0)$ to get a stability slope of $1/5$, corresponding to an upper bound of $9/10$ on the exponent of the secondary term.

\begin{figure}[!htbp]
    \centering
\begin{tikzpicture}[x=1cm,y=1cm,scale=0.65]
\drawnaxis{10}
\drawiaxis{9}
\drawstaircase{1}{10}{0}{black!10}{0}
\drawblock{5}{10}{0}{0}{blue!50}{6}
\drawblock{9}{10}{1}{1}{blue!50}{6}

\foreach \n/\i/\txt in {
3/0/6,
4/1/11
 }{
    \cell{\n}{\i}{\txt}{red!20}
  }
\foreach \n/\i/\txt in {
4/0/6
 }{
    \cell{\n}{\i}{\txt}{red!50}
  }
  \foreach \n/\i/\txt in {
5/1/12,
6/1/11,
7/1/9,
8/1/7,
9/2/0,
10/2/0
 }{
    \cell{\n}{\i}{\txt}{blue!20}
  }
\foreach \n/\i/\txt in {
1/0/3,
2/0/5,
2/1/5,
3/1/9,
3/2/3,
4/2/11,
4/3/6,
5/2/12,
5/3/15,
5/4/9,
6/2/10,
6/3/17,
6/4/41,
6/5/29,
7/2/6,
7/3/12,
7/4/33,
7/5/72,
7/6/48,
8/2/1,
8/3/5,
8/4/21,
8/5/75,
8/6/167,
8/7/108,
9/3/3,
9/4/9,
9/5/30,
9/6/120,
9/7/300,
9/8/204,
10/3/5,
10/4/10,
10/5/16,
10/6/43,
10/7/198,
10/8/718,
10/9/552
 }{
    \cell{\n}{\i}{\txt}{black!10}
  }

\end{tikzpicture}

    \caption{Theorem \ref{thm: main theorem} for $c^{\text{trn}}_{S_3}$ with $(N,I)=(4,0)$.}
    \label{fig: (4,0) S_3 trans}
\end{figure}

We illustrate this in Figure \ref{fig: (4,0) S_3 trans}, where the red cell denotes a known isomorphism of the stabilization map \eqref{eqn: stab map} \textit{starting from} that cell, and the light red cell denotes a known surjection from that cell. Theorem \ref{thm: main theorem} then gives us isomorphisms of \eqref{eqn: stab map} from the blue cells, and surjections from the light blue cells. This application of our main theorem tells us that $H_1(\Hur^c_n,\Q)$ stabilizes at $n=9$, which turns out to be exactly correct from our computations. Hence, we see that Theorem \ref{thm: main theorem} is \textit{tight} in this sense.

One can obtain progressively better stability slopes as one computes more homology. For example, computing the homology up to $n=10$ tells us that the homology stabilizes in a $(10,2)$ staircase. Applying Theorem \ref{thm: main theorem} gives a stability slope of $3/11$ and an upper bound on the exponent of $19/22$, improving the previous bounds obtained. We illustrate this in Figure \ref{fig: (10,2) S_3 trans}. To check stability in a $(10,2)$ staircase, we note that the stabilization map $U^e$ is an isomorphism for $i=0$ and $n\ge 3$ as mentioned earlier, and for $i=1$ it is a surjection for $n=8$ and an isomorphism for $n=9,10$ from the computed homological dimensions in the figure and the argument after Corollary \ref{cor: stab from R}. For $i=2,3$, using the fact that $H_2(\Hur^c_n,\Q)=H_3(\Hur^c_n,\Q)=0$ for sufficiently large $n$, we can conclude that $U^e$ is a surjection for $(n,i)=(10,3)$ and an isomorphism for $(n,i)=(9,2),(10,2)$ if $e$ is sufficiently large. 

\begin{figure}[!htbp]
    \centering
 \begin{tikzpicture}[x=1cm,y=1cm,scale=0.65]
\drawnaxis{22}
\drawiaxis{9}
\drawstaircase{1}{10}{0}{black!10}{0}
\drawblock{3}{6}{0}{0}{black!10}{6}
\drawblock{8}{10}{0}{0}{red!50}{6}
\drawblock{9}{10}{1}{1}{red!50}{6}
\drawblock{10}{10}{2}{2}{red!50}{0}
\drawblock{11}{22}{0}{9}{black!10}{*}
\drawblock{11}{22}{0}{1}{blue!50}{6}
\drawblock{11}{22}{2}{2}{blue!50}{0}
\drawblock{11}{22}{3}{3}{blue!20}{*}
\drawstaircase{19}{21}{3}{blue!50}{0}
\drawblock{22}{22}{3}{5}{blue!50}{0}

\foreach \n/\i/\txt in {
7/0/6,
8/1/7,
9/2/0,
10/3/5
 }{
    \cell{\n}{\i}{\txt}{red!20}
  }
\foreach \n/\i/\txt in {
19/4/*,
20/5/*,
21/6/*,
22/6/*
 }{
    \cell{\n}{\i}{\txt}{blue!20}
  }
\foreach \n/\i/\txt in {
1/0/3,
2/0/5,
2/1/5,
3/1/9,
3/2/3,
4/2/11,
4/3/6,
5/1/12,
5/2/12,
5/3/15,
5/4/9,
6/1/11,
6/2/10,
6/3/17,
6/4/41,
6/5/29,
7/1/9,
7/2/6,
7/3/12,
7/4/33,
7/5/72,
7/6/48,
8/2/1,
8/3/5,
8/4/21,
8/5/75,
8/6/167,
8/7/108,
9/3/3,
9/4/9,
9/5/30,
9/6/120,
9/7/300,
9/8/204,
10/4/10,
10/5/16,
10/6/43,
10/7/198,
10/8/718,
10/9/552
 }{
    \cell{\n}{\i}{\txt}{black!10}
  }

\end{tikzpicture}
    \caption{Theorem \ref{thm: main theorem} for $c^{\text{trn}}_{S_3}$ with $(N,I)=(10,2)$.}
    \label{fig: (10,2) S_3 trans}
\end{figure}

It turns out that one can push the homological computations even further, and we can verify computationally that the homology stabilizes in a $(28,8)$ staircase in this case, see Figure \ref{fig: homology computation S_d}. This proves a stability slope of $\frac 9 {29}\approx 0.310$, and establishes that the exponent is not more than $\frac{49}{58}\approx 0.845$. This is not too far off from the expected stability slope of $1/3\approx 0.333$ and secondary exponent of $\frac 5 6 \approx 0.833$. We also observe that the stability staircase seems to repeat under the shift $(n,i)\mapsto (n+6,i+2)$, which is consistent with the expected $1/3$ stability slope. If this pattern continues, we see that the upper bound on the exponent should approach the true exponent of $5/6$ as more computation is done. We will discuss the homological computations in more detail in Section \ref{sec: computations}.
\subsubsection{Optimality of slope, weights and exponent}\label{sec: intro optimality}
In the Davenport-Heilbronn example above, we saw that computing more homology gives progressively better stability slopes which should converge to the true optimal slope. The next corollary shows that this is indeed true for all Hurwitz spaces. This is in analogy to \cite[Corollary 1.3]{prequel}, and follows immediately from the definition of slope. 
\begin{corollary}\label{cor: optimal bound}
Suppose that $H_*(\Hur^c,\Q)$ is homologically stable with slope $s$ (with respect to $U=\sum b_i\gamma_i$ where $\gamma_i\in \pi_0\Hur^c_a$ are central). Then, for any $\epsilon>0$, there exist some $N>0$ and $I\ge 0$ such that $H_*(\Hur^c,\Q)$ stabilizes in an $(N,I)$ staircase and $\frac{I+1}{N+1}>s-\epsilon$. In particular, computing $H_*(\Hur^c,\Q)$ up to sufficiently large $n$ and applying Theorem \ref{thm: main theorem} will prove homological stability with slope $>s-\epsilon$.
\end{corollary}

Next, we discuss the implications of the optimal slope on Frobenius weights and the exponent of secondary terms. Before that, let us quickly recall how the stability slope $s$ leads to an upper bound of $1-s/2$ on the exponent for the secondary term. Note that the exponent in the Cohen-Lenstra case is instead $-s/2$ because we take the expectation over the quadratic fields. \textit{For the rest of this subsection,} let us restrict to the case where $c\subseteq G$ is a non-splitting conjugacy class in the group $G$, and we denote the Hurwitz space by $\Hur^{G,c}_n$. We also let $p$ be a prime satisfying $p\nmid |G|$ with $c$ closed under $p$-th powering. Here, we set this up so that $\Hur^{G,c}_n$ is a scheme defined over $\Z_p$. The Grothendieck-Lefschetz trace formula says that for $q$ a power of $p$, we have
$$\#\Hur^{G,c}_n(\F_q)=q^n\sum_{i=0}^{2n}(-1)^i\tr\left(\Frob_q,H_i(\Hur^{G,c}_{n,\overline\F_p},\overline \Q_l)\right).$$

Deligne's Weil bounds tell us that the Frobenius weights $w$ of $H_i(\Hur^{G,c}_{n,\overline\F_p},\overline \Q_l)$ are at most $-i$, so the eigenvalues of $\Frob_q$ have absolute value $\le q^{-i/2}$. We know by \cite[Theorem 1.4.2]{landesman2025cohenlenstramomentsfunctionfields} that the stable homology of each  component of $\Hur^{G,c}_n$ is $\Conf_n$. Letting $C_{q,n}$ be the number of geometric connected components fixed by Frobenius, the main term from the Grothendieck-Lefschetz trace formula is then $C_{q,n}\cdot\#\Conf_n(\F_q)$ which is easily computable. Suppose we have a stability line of $sn+t\ge i$, then for each $n$, the difference in the trace formulas for $\# \Hur^{G,c}_n(\F_q)-C_{q,n}\cdot\#\Conf_n(\F_q)$ only comes from terms in the unstable range where $i> sn+t$, and each eigenvalue of $\Frob_q$ there has absolute value $< q^{-(sn+t)/2}$. Hence, up to a Betti number factor and the constant factor $q^{-t/2}$, the growth of the secondary term is bounded by $(q^n)^{1-s/2}$. We also remark that we often look at the quotient of $[\Hur^{G,c}_n/G]$, and essentially the same argument will work.

One can then ask whether the bounds on the secondary exponent obtained this way converge to the optimal exponent  (modulo both the Betti number and constant terms). As explained in \cite[Section 1.0.2]{prequel}, there can be two reasons why this doesn't happen: (1) Deligne's Weil bounds are not tight, i.e. the Frobenius weights $w$ are strictly less than the homological degree $-i$, or (2) that the largest Frobenius eigenvalue has absolute value $q^{w/2}$, but the sum of eigenvalues in the trace formula cancels and has absolute value less than $q^{w/2}$. The first could indeed be a possibility, like in the case of $c=c^{\text{trn}}_{S_3}$ where \cite[Theorem 1.3]{ma2026weightfiltrationhurwitzspaces} showed that the Frobenius weights are \textit{all} strictly smaller than expected from homological degree when $i>n/3$ (although in this case it will only apply in the secondary stability range). On the other hand, we think that the second is unlikely to happen for most Hurwitz spaces.

Nevertheless, our methods in this paper, combined with the weight filtration developed in \cite{ma2026weightfiltrationhurwitzspaces}, can address the first factor. Indeed, one can always show, from computation, a stability line $-w\le sn+t$ in terms of Frobenius weights $w$ (instead of the homological degree $i$), whose \textit{weight} slope $s$ converges to the optimal value. Hence, we can get the optimal secondary exponent \textit{assuming} no unexpected cancellation of Frobenius eigenvalues. 

Let us state this more precisely. There exists a weight decomposition
$$H_i(\Hur^{G,c}_{n,\overline\F_p},\overline \Q_l)=\bigoplus_{w\le -i}\Gr_w H_i(\Hur^{G,c}_{n,\overline\F_p},\overline \Q_l)  $$
where $\Frob_q$ acts on the $w$-th graded piece with absolute value $q^{w/2}$. Then, if $U\in R$ is a central element, we will show that it preserves the weight grading $w$, so the stabilization map restricts to
\begin{equation}\label{eqn: WEIGHT stab map}
U\colon \Gr_w H_i(\Hur^{G,c}_{n,\overline\F_p},\overline \Q_l) \rightarrow \Gr_wH_i(\Hur^{G,c}_{n+a,\overline\F_p},\overline \Q_l).\end{equation}
Then, we have the following analogous theorem where we replace $i$ with $-w$ everywhere.
\begin{definition}
Let $N\ge 1$ and $W\ge 0$ be integers satisfying $W<N$. Say that $H_*(\Hur^{G,c}_{n,\overline\F_p},\overline \Q_l)$ stabilizes in an $(N,W)$-\textit{weight} staircase if for $n\le N$, the stabilization map \eqref{eqn: WEIGHT stab map} is an isomorphism for $-w\le W-N+n$ and a surjection for $-w=W-N+n+1$ over all homological degrees. 
\end{definition}
\begin{theorem}\label{thm: weight} Suppose that $H_*(\Hur^{G,c}_{n,\overline\F_p},\overline \Q_l)$ has homological stability with respect to $U=\sum b_i\gamma_i$ where $\gamma_i\in \pi_0\Hur^c_a$ are central, and stabilizes in an $(N,W)$-weight staircase. Then, we have homological stability with weight slope $\frac{W+1}{N+1}$, with the precise bounds analogous to those in Theorem \ref{thm: main theorem}. Furthermore, in analogy to Corollary \ref{cor: optimal bound}, computing the homology up to sufficiently large $n$ will prove homological stability with weight slope $s$ arbitrarily close to the optimal weight slope, as long as the optimal weight slope is $\le 1$.
\end{theorem}

However, we remark that we avoid using weight data for the explicit computations in this paper as it makes the computations much more complicated. We leave this as a theoretical optimality result for now.

\subsubsection{Stability line from $R$}\label{sec: intro stab from R} Naïve applications of our main theorem require explicit computation of the homology of Hurwitz space, and thus cannot give us results for families of Hurwitz spaces as it is impossible to carry out computations for each individual object in the infinite family. However, there is a useful special case of our theorem which does allow us to obtain results in families. Specifically, one would only need knowledge of when the ring of connected components $R=H_0(\Hur^c,\Q)$ stabilizes to deduce a non-trivial stability slope, and we state this in the following corollary. This is very useful as it is often possible to deduce when $R$ stabilizes in families of Hurwitz spaces.

\begin{corollary}\label{cor: stab from R}
Suppose that $H_*(\Hur^c,\Q)$ has homological stability with respect to $U=\sum b_i\gamma_i$ where $\gamma_i\in \pi_0\Hur^c_a$ are central, and there is some $N\ge 2$ where the stabilization map \eqref{eqn: stab map} is a surjection when $i=0$, $n=N-1$, and an isomorphism when $i=0$, $n\ge N$. Then, it is homologically stable with slope $\frac{1}{N+1}$, in particular, \eqref{eqn: stab map} is an isomorphism when $n\ge i(N+1)+N$ for all $i\ge 0$.
\end{corollary}

We now show that this is a consequence of Theorem \ref{thm: main theorem} when $I=0$. The condition of stabilization in an $(N,0)$ staircase is by definition \eqref{eqn: stab map} being an isomorphism for $(n,i)=(N,0)$ and a surjection for $(n,i)=(N,1),(N-1,0)$. One can deduce the $(N,0)$ and $(N-1,0)$ conditions from the assumptions of the corollary, but it seems like on top of that, we need information on $H_1(\Hur^c,\Q)$ to deduce the surjection for $(N,1)$.

However, there is a simple trick that lets us avoid the need for any computation of $H_1$. We use the fact that the homology of $\Hur^c_n$ stabilizes to a multiple of $\Conf_n$, one for each eventual connected component of $\Hur^c_n$ (see Theorem \ref{thm: stable hom conf all directions}). Let us consider the iterated stability map $$U^e\colon H_i(\Hur^c_N,\Q)\rightarrow H_i(\Hur^c_{N+ea},\Q)$$
for any sufficiently large integer $e>1$ so that $H_1$ has stabilized for $n\ge N+ea$. Then, we apply Theorem \ref{thm: main theorem} to $U^e$ instead of $U$. The assumptions of the corollary likewise guarantee the $(N,0)$ and $(N-1,0)$ conditions for $U^e$. For the $(N,1)$ condition, we consider the following commutative diagram.
\[\begin{tikzcd}
	{H_1(\Hur^c_N,\Q)} & {H_1(\Hur^c_{N+ea},\Q)} \\
	{H_1(\sqcup\Conf_N,\Q)} & {H_1(\sqcup\Conf_{N+ea},\Q)}
	\arrow["{U^e}", from=1-1, to=1-2]
	\arrow[two heads, from=1-1, to=2-1]
	\arrow["\sim", from=1-2, to=2-2]
	\arrow["\sim", from=2-1, to=2-2]
\end{tikzcd}\]
Note that since the $H_0$ has stabilized by $N$, the number of connected components on the bottom terms are the same, so the bottom map is an isomorphism. The right map is an isomorphism as $H_1$ has already stabilized, and the left is a surjection as $\Hur^c_n$ is a finite étale cover of $\sqcup \Conf_n$. This guarantees that $U^e$ is a surjection. 

Hence, Theorem \ref{thm: main theorem} tells us that the stability map $U^e$ is an isomorphism for the range stated in Corollary \ref{cor: stab from R}. It is then straightforward to conclude the same for the map $U$ by the following diagram, since we can simply apply the isomorphisms $U^e$ to both $H_i(\Hur^{G,c}_n)$ and $H_i(\Hur^{G,c}_{n+a})$ until we reach an isomorphism $H_i(\Hur^{G,c}_{n+kea})\xrightarrow\sim H_i(\Hur^{G,c}_{n+(ke+1)a})$. We remark that this trick of considering $U^e$ instead of $U$ is also useful for other applications of Theorem \ref{thm: main theorem}, as we will see later in the paper.

\[\begin{tikzcd}
	{H_i(\Hur^{G,c}_n)} & {H_i(\Hur^{G,c}_{n+a})} & \cdots & \cdots & {H_i(\Hur^{G,c}_{n+kea})} & {H_i(\Hur^{G,c}_{n+(ke+1)a})}
	\arrow["U"', from=1-1, to=1-2]
	\arrow["\sim", curve={height=-24pt}, from=1-1, to=1-3]
	\arrow["U"', from=1-2, to=1-3]
	\arrow["\sim"', curve={height=24pt}, from=1-2, to=1-4]
	\arrow["\sim", curve={height=-24pt}, from=1-3, to=1-5]
	\arrow["U"', from=1-4, to=1-5]
	\arrow["\sim"', curve={height=24pt}, from=1-4, to=1-6]
	\arrow["U"', from=1-5, to=1-6]
	\arrow["\sim", draw=none, from=1-5, to=1-6]
\end{tikzcd}\]

This type of result is not new; in fact, the original proof of homological stability in \cite{EVW16,Randal-Williams_2020} already gave a stability line in terms of when $R$ stabilizes. Again, the novelty in our theorem lies in the tightness of the stability line. Let us compare both results, using only knowledge of when $R$ stabilizes, in the case of $c=c^{\text{trn}}_{S_3}$. Taking $N=4$ and applying Corollary \ref{cor: stab from R} gives a stability line of $n\ge 5i+4$, which is sharp for $i=1$ as we saw earlier. On the other hand, as explained in \cite[Remark 5.3.2]{landesman2025cohenlenstramomentsfunctionfields}, the stability line that \cite{Randal-Williams_2020} obtains in this case is $n\ge 14i+12$.
\subsection{Arithmetic applications}\label{sec: intro arithmetic applications} In this subsection we will discuss the arithmetic consequences of our results on homological stability. In Section \ref{sec: intro cohen lenstra} we discuss results on Cohen-Lenstra moments over function fields, and in Section \ref{sec: intro S_d fields} we look at the number of $S_d$-degree $d$ fields with squarefree discriminant. We also apply our theorems to other Hurwitz spaces in Section \ref{sec: other hurwitz spaces}.
\subsubsection{Cohen-Lenstra moments}\label{sec: intro cohen lenstra}

We will consider a version of the Cohen-Lenstra heuristics with roots of unity, as in \cite[Section 1.2]{landesman2025cohenlenstramomentsfunctionfields}. Let $q$ be an odd prime power, and let $\mathcal {MH}_{n,q}$ denote the set of function fields $K$ of monic smooth hyperelliptic curves of the form $y^2=f(x)$, where $f(x)$ is a monic squarefree degree $n$ polynomial with coefficients in $\F_q$. Note that this is the analogue of quadratic number fields in the function field case. Let $\mathcal O_K$ be the normalization of $\F_q[t]$ in $K$, and let $\Cl(\mathcal O_K)$ denote its class group.

Let $H$ be a finite abelian group of odd order such that $|H|$ is coprime to $q$, and let $h\coloneqq \gcd(|H|,q-1)$. Then, \cite[Theorem 1.2.1, Theorem 5.4.1]{landesman2025cohenlenstramomentsfunctionfields} proved that there is a constant $c_H>0$ such that for all sufficiently large $q$, we have
$$\E_{K\in \mathcal{MH}_{n,q}} \left[\# \Surj(\Cl(\mathcal O_K),H)\right]=\alpha_n+o\left((q^n)^{-c_H}\right),$$
where we let
$$\alpha_n=\alpha_{n,H}=\begin{cases}
\left|\bigwedge^2 H[h]\right| & \text{if } n \text{ is odd}\\
\frac{\left|\bigwedge^2 H[h]\right|}{|H|} &\text{if } n \text{ is even}.
\end{cases}$$

The left hand side is often referred to as the $H$-moment of the class group $\Cl(\mathcal O_K)$. The $o((q^n)^{-c_H})$ term on the right contains all possible secondary, higher order and error terms, and the power-savings exponent $c_H$ is explicitly computable but not at all close to optimal and also not stated in \cite{landesman2025cohenlenstramomentsfunctionfields}. 

We give the following improved exponent which depends only on the rank of $H$.
\begin{theorem}\label{thm: cohen-lenstra H}
Fix an odd abelian group $H$, and let $r=\rk(H)$ be its rank. Then, for any $\epsilon>0$, we have for all sufficiently large $q$ that
$$\E_{K\in \mathcal{MH}_{n,q}} \left[\# \Surj(\Cl(\mathcal O_K),H)\right]=\alpha_n+O\left((q^n)^{-\frac{1}{4r+6}+\epsilon}\right).$$
\end{theorem}
Here, the $\epsilon$ is introduced solely to account for the Betti number factor. When the rank $r=1$, i.e. when $H$ is cyclic, this implies that the Cohen-Lenstra $k$-moment has secondary term $\le (q^n)^{-1/10+\epsilon}$ as promised earlier. This theorem will follow from Corollary \ref{cor: stab from R} and an analysis of when $R$ stabilizes for the relevant Hurwitz space $\Hur^{G,c}_n$ where $G=H\rtimes \Z/2\Z$ and $c=c^{\text{rfl}}_{G}$ is the conjugacy class of reflections in $G$. We will carry out this analysis in Section \ref{sec: cohen lenstra}.

One can improve the exponent of $-\frac{1}{4r+6}$ for individual cases of $H$ by computing the homology of the relevant Hurwitz space using the philosophy of Corollary \ref{cor: optimal bound}. We present these homological computations in Section \ref{sec: compute cohen lenstra}. These computations lead to the improvements in bounds for the stability slope and secondary exponent in Figures \ref{fig: table H without 3-torsion} and \ref{fig: table H with 3-torsion}. In all such tables, we omit the $+\epsilon$ in the error exponents, where $\epsilon>0$ is arbitrary. Note that as the size of $H$ increases, homology computations take longer which limits the amount of homology we can compute, and this results in weaker bounds on slopes and exponents. 

\begin{figure}[!htbp]
\centering

\begin{tabular}{
>{\centering\arraybackslash}m{1cm}
>{\centering\arraybackslash}m{5.5cm}
>{\centering\arraybackslash}m{2.5cm}
>{\centering\arraybackslash}m{2.5cm}
>{\centering\arraybackslash}m{1cm}
}
\toprule
\multicolumn{2}{c}{$\mathbf{H}$} &
\textbf{Slope} $\mathbf{\ge}$ & \textbf{Exponent} $\mathbf{\le}$ & $\mathbf{(N,I)}$ \\
\midrule

\multirow{5}{*}{$\Z/k\Z$}
& $k=5$ & $\frac{7}{15}\approx 0.467$ & $-\frac{7}{30}\approx \mathbf{-0.233}$ & $(14,6)$ \\
\cmidrule{2-5}
& $k=7$ & $\frac 1 2=0.5$ & $-\frac 1 4 =\mathbf{-0.25}$ & $(9,4)$ \\
\cmidrule{2-5}
& $k=11,13$ & $\frac 4 9\approx 0.444$ & $-\frac 2 9 \approx \mathbf{-0.222}$& $(8,3)$ \\
\cmidrule{2-5}
& $k=17,19,23,25$ & $\frac 3 8 = 0.375$& $-\frac{3}{16}\approx \mathbf{-0.187}$ & $(7,2)$ \\
\cmidrule{2-5}
& $k=29,31,35,37,41,43,47,49$ & $\frac 1 3 \approx  0.333$& $-\frac 1 6\approx \mathbf{-0.167}$ & $(5,1)$ \\
\midrule
\multicolumn{2}{c}{$(\Z/5\Z)^2$} & $\frac 1 3 \approx 0.333$ & $-\frac 1 6\approx \mathbf{-0.167}$ & $(8,2)$ \\
\midrule
\multicolumn{2}{c}{$(\Z/7\Z)^2$} & $\frac 1 4=0.25$ & $-\frac 1 8= \mathbf{-0.125}$ & $(7,1)$ \\
\bottomrule
\end{tabular}

\caption{Bounds on slope and exponent for various $H$ without $3$-torsion.}\label{fig: table H without 3-torsion}
\end{figure}

\begin{figure}[!htbp]
\centering

\begin{tabular}{
>{\centering\arraybackslash}m{1cm}
>{\centering\arraybackslash}m{5.5cm}
>{\centering\arraybackslash}m{2.5cm}
>{\centering\arraybackslash}m{2.5cm}
>{\centering\arraybackslash}m{1cm}
}
\toprule
\multicolumn{2}{c}{$\mathbf{H}$} &
\textbf{Slope} $\mathbf{\ge}$ & \textbf{Exponent} $\mathbf{\le}$ & $\mathbf{(N,I)}$ \\
\midrule

\multirow{3}{*}{$\Z/k\Z$}
& $k=3$ & $\frac 9 {29}\approx 0.310$ & $-\frac {9}{58}\approx\mathbf{-0.155}$ & $(28,8)$ \\
\cmidrule{2-5}
& $k=9$ & $\frac 1 4=0.25$ & $-\frac 1 8 =\mathbf{-0.125}$ & $(15,3)$ \\
\cmidrule{2-5}
& $k=15,21,27$ & $\frac{3}{11}\approx 0.273$ & $-\frac {3}{22}\approx \mathbf{-0.136}$& $(10,2)$ \\
\midrule
\multicolumn{2}{c}{$(\Z/3\Z)^2$} & $\frac 1 4=0.25$ & $-\frac 1 8 =\mathbf{-0.125}$ & $(15,3)$ \\
\midrule
\multicolumn{2}{c}{$\Z/3\Z\times \Z/9\Z$} & $\frac{3}{11}\approx 0.273$ & $-\frac {3}{22}\approx \mathbf{-0.136}$& $(10,2)$ \\
\bottomrule
\end{tabular}

\caption{Bounds on slope and exponent for various $H$ with $3$-torsion.}\label{fig: table H with 3-torsion}
\end{figure}

Here, we distinguish two cases according to whether $H$ has $3$-torsion as we observe significant differing stability patterns in these two cases. In the non-$3$-torsion case, we observe a stability slope of approximately $1/2$, and this can be seen most clearly in the $H=\Z/5\Z$ and $H=\Z/7\Z$ cases in Figure \ref{fig: homology computation non 3 torsion}. As seen in Figure \ref{fig: table H without 3-torsion}, our computations prove lower bounds on stability slopes that are slightly smaller than $1/2$, leading to an upper bound on the exponent that is slightly above $-1/4$. 

On the other hand, in the $3$-torsion case, our computations suggest that the stability staircase eventually repeats under the shift $(n,i)\mapsto (n+6,i+2)$ so that the stability slope is $1/3$, as we see in Figure \ref{fig: homology computation 3-torsion}. This is consistent with the lower bounds on stability slopes we proved via computation in Figure \ref{fig: table H with 3-torsion}, which are all smaller than $1/3$, unlike what we saw in Figure \ref{fig: table H without 3-torsion}. 

We previously observed a similar stability pattern for the conjugacy class $c^{\text{trn}}_{S_3}$ of transpositions in $S_3$. Indeed, this conjugacy class is isomorphic as a rack to the conjugacy class $c^{\text{rfl}}_{D_3}$ of reflections in $D_3$. When $H$ has $3$-torsion, then $D_3$ is a subgroup of $G=H\rtimes \Z/2\Z$ and so the relevant Hurwitz space $\Hur^{G,c}_n$ has connected components arising from $\Hur^{D_3,c^{\text{rfl}}_{D_3}}_n$. Moreover, just like the $c^{\text{trn}}_{S_3}$ case, our computations also exhibit evidence of secondary stability.

Motivated by this observation, we attempt to remove the contribution of the $S_3$ transposition case to the homology of $\Hur^{G,c}_n$. For each $n$, we consider the local system on $\Conf_n$ associated with the covering $\Hur^{G,c}_n\rightarrow \Conf_n$, decompose the corresponding $B_n$-representation, and remove all direct summands isomorphic to the reduced local system associated with $\CHur^{S_3,c^{\text{trn}}_{S_3}}_n$. The homology remaining after this removal is plotted in Figure \ref{fig: homology computation 3-torsion with S3 removed}. Remarkably, we observe an approximate stability slope of $1/2$, matching the pattern seen in the non-$3$-torsion case. This suggests that the smaller stability slope of $1/3$ in the $3$-torsion case may be entirely accounted for by the contributions arising from the $S_3$ transposition case.

Based on these observations and experiments, we make the following bold conjecture. Note that we do not claim that the exponent $-n/4$ in the error term is sharp, and we do not know whether our error term contains more secondary and higher order terms.
\begin{conjecture}
Fix an odd abelian group $H$. If $H$ has no $3$-torsion, then for all sufficiently large $q$ we have
$$\E_{K\in \mathcal{MH}_{n,q}} \left[\# \Surj(\Cl(\mathcal O_K),H)\right]=\alpha_n+O\left(q^{-n/4}\right).$$
If $H$ has $3$-torsion, then for all sufficiently large $q$, there exist constants $\beta_n=\beta_{n,q,H}$ which are periodic modulo $6$ in $n$ so that
$$\E_{K\in \mathcal{MH}_{n,q}} \left[\# \Surj(\Cl(\mathcal O_K),H)\right]=\alpha_n+\beta_n q^{-n/6}+O\left(q^{-n/4}\right).$$

We also conjecture the number field analogue: let $\mathcal K^+(X), \mathcal K^-(X)$ be the sets of real or imaginary quadratic fields with absolute discriminant at most $X$, respectively. Then we have
$$\E_{K\in \mathcal K^\pm(X)}[\#\Surj(\Cl(\mathcal O_K),H)]=\begin{cases}
\alpha^\pm +O\left(X^{-1/4}\right) & \text{if }H \text{ has no }3\text{-torsion}\\
\alpha^\pm +\beta^{\pm}X^{-1/6}+O\left(X^{-1/4}\right) & \text{if }H \text{ has }3\text{-torsion}
\end{cases}$$
where $\alpha^+=\frac{1}{|H|}$, $\alpha^-=1$, and $\beta^{\pm}=\beta^\pm_H$ are constants.
\end{conjecture}

However, this conjecture appears to be well beyond the reach of current techniques. We instead focus our attention to $H_1$ in the case of dihedral groups $G=D_k$ with $c=c^{\text{rfl}}_{D_k}$, where meaningful progress seems more attainable. From our computations, we also notice two differing patterns on $H_1(\Hur^{D_k,c}_n,\Q)$ depending on whether $3$ divides $k$. A priori, we can deduce from the stabilization of $R$ that $H_1(\Hur^{D_k,c}_n,\Q)$ always stabilizes by $n=9$, see Theorem \ref{thm: cohen lenstra stab slope}, and this seems to be sharp in the case when $k$ is a multiple of $3$. On the other hand, when $3\nmid k$, we observed that $H_1(\Hur^{D_k,c}_n,\Q)$ stabilizes much earlier at $n=5$. We formulate this in the following conjecture.
\begin{conjecture}\label{conj: H1}
Let $k$ be an odd integer with $3\nmid k$. Then we have that $H_1(\Hur^{D_k,c}_n,\Q)$ stabilizes for $n\ge 5$. Arithmetically, this implies that for any $\epsilon >0$, we have for all sufficiently large $q$ that $$\E_{K\in \mathcal{MH}_{n,q}} \left[\# \Surj(\Cl(\mathcal O_K),\Z/k\Z)\right]=\alpha_n+O\left((q^n)^{-\frac{1}{6}+\epsilon}\right).$$
\end{conjecture}

In Section \ref{sec: H1 conjectures}, we give a heuristic on why this conjecture should be true, at least in the case where $k$ is a prime, coming from the geometry of the Siegel modular threefold $\mathcal A_2$. In particular, this heuristic gives a reason why the homology would stabilize at $n=5$. The key idea is that $\Hur^{D_k,c}_5$ is the data of the hyper-elliptic curve of genus $g=2$ with the $5$ branch points along with extra data on a labelled point at $\infty$ as well as a choice of $k$-torsion in the Jacobian. Since the moduli space of hyperelliptic curves $\mathcal H_2\cong \mathcal M_2$ embeds into $\mathcal A_2$ via the Torelli map, we can view a quotient of $\Hur^c_5$ as an open subset of a finite cover of $\mathcal A_2$. This gives an interesting bridge between homological stability (of Hurwitz spaces) and the geometry of Hermitian symmetric spaces, which could then possibly be related to other things like theta lifts and automorphic forms.
\subsubsection{Counting $S_d$-fields with squarefree discriminant}\label{sec: intro S_d fields} Arguably the most classical problem in arithmetic statistics is to count the number of degree-$d$ number fields with discriminant bounded by $X$. Often, one also prescribes the Galois group $G\subseteq S_d$ of the Galois closure. In this paper we are concerned with the generic case of $G=S_d$. 

Here, we count the function field analogue of this, where for simplicity we also impose that the field has squarefree discriminant. This additional restriction only changes the counting function by a constant when $d\le 5$ by the squarefree sieve, but it is not clear whether this will also be true for larger $d$. We set this up as follows.

Let $q$ be a prime power. A finite separable degree $d$ field extension $K/\F_q(t)$ is equivalent to a degree $d$ branched cover $f\colon C\rightarrow \mathbb P^1_{\F_q}$, where we require $C$ to be a smooth irreducible curve and $f$ to be generically étale. The global discriminant $\Delta=\disc_{K/\F_q(t)}$ is a divisor on $\mathbb P^1$ which is the branch divisor of the corresponding cover. By Riemann-Hurwitz, the degree of the branch divisor must be an \textit{even} integer $2n$, and in this case $\Norm(\Delta)=q^{2n}$. Having squarefree discriminant $\Delta$ corresponds to $f$ being a simply-branched cover. Furthermore, the group $G$ corresponds to the arithmetic monodromy group $G_{\text{arith}}$ of the cover. In our case of $G=S_d$ with squarefree discriminant, the geometric monodromy group $G_{\text{geom}}$ is also $S_d$ because it is generated by geometric local monodromies which are transpositions and must be a normal subgroup of $G_{\text{arith}}=S_d$.

Let $N^{\text{sqf}}_{S_d}(2n)$ be the number of finite separable $S_d$-degree $d$ field extensions of $\F_q(t)$ with squarefree discriminant $\Delta$ of norm $q^{2n}$. Then, we have the following result across all $d$.

\begin{theorem}\label{thm: S_d count}
Let $d\ge 3$ be an integer. Then for any $\epsilon>0$, we have for all sufficiently large $q$ satisfying $\gcd(q,d!)=1$ that
$$N^{\text{sqf}}_{S_d}(2n)=(1-q^{-2})q^{2n}+O\left((q^{2n})^{c_d+\epsilon}\right)$$
with the exponent $c_d=1-\frac{1}{2a_d}$ where $a_d=\Theta(3^{d/3})$ is given explicitly in Equations \eqref{eqn: a_d values} and \eqref{eqn: small a_d values}.
\end{theorem}

Here, the relevant Hurwitz space $\Hur^{G,c}_n$ for this problem is when $c=c^{\text{trn}}_{S_d}$ is the conjugacy class of transpositions in $G=S_d$. Much like Theorem \ref{thm: cohen-lenstra H} in the Cohen-Lenstra case, we prove this theorem by deducing a stability line from only information of $R$ since that is easily computable in families. However, while Theorem \ref{thm: cohen-lenstra H} followed from a direct application of Corollary \ref{cor: stab from R}, the conjugacy class $c$ in this case is splitting, which means that it is no longer homologically stable with respect to a single element, bringing more technical difficulties. We carry this out in Section \ref{sec: S_d slope}. 

We can improve the exponents $c_d$ significantly for small $d$ by computing the homology of the respective Hurwitz spaces which we carry out in Section \ref{sec: compute S_d}, and then applying Corollary \ref{cor: general optimal} which is the general version of Corollary \ref{cor: optimal bound}. We present these improvements in Figure \ref{fig: improved c_d}, omitting the $+\epsilon$ in the exponents. Further computations should yield improved exponents for each of the $3\le d\le 7$ cases in the table, as well as extend improved bounds for higher $d$.

\begin{figure}[!htbp]
\centering
\begin{tabular}{>{\centering\arraybackslash}m{2.5cm}
>{\centering\arraybackslash}m{3.5cm}
>{\centering\arraybackslash}m{3.5cm}
>{\centering\arraybackslash}m{2.5cm}}
\toprule
$\mathbf{d}$ & \textbf{Slope} $\mathbf{\ge}$ & \textbf{Exponent} $\mathbf{\le}$ & $\mathbf{(N,I)}$ \\
\midrule
$3$ & $\frac 9 {29}\approx 0.310$ & $\frac{49}{58}\approx \mathbf{0.845}$ & $(28,8)$ \\
\midrule
$4$ & $\frac 5 {19}\approx 0.263$ & $\frac{33}{38}\approx \mathbf{0.868}$ & $(18,4)$ \\
\midrule
$5$ & $\frac 5 {22}\approx 0.227$ & $\frac{39}{44}\approx \mathbf{0.886}$ & $(21,4)$ \\
\midrule
$6$ & $\frac 1 {6}\approx 0.167$ & $\frac{11}{12}\approx \mathbf{0.917}$ & $(17,2)$ \\
\midrule
$7$ & $\frac 2 {19}\approx 0.105$ & $\frac{18}{19}\approx  \mathbf{0.947}$ & $(18,1)$ \\
\bottomrule
\end{tabular}
\caption{Bounds on slope and exponent for $3\le d\le 7$.}\label{fig: improved c_d}
\end{figure}

For $d=3,4$, the exponents we obtain are slightly weaker than what is currently known; however, they are the best bounds coming from the method of homological stability. For $d=5$, our upper bound on the exponent of $0.886$ is much smaller than the exponent of $199/200=0.995$ in \cite{SHANKAR_TSIMERMAN_2014} over number fields. For $d\ge 6$ there are no results over number fields, and over function fields the best exponent proven so far is approximately $1- 1/\left(\binom d 22^{\binom d 2}\right)$ as one can see from \cite[Remark 1.2.4]{landesman2025homologicalstabilityhurwitzspaces}. For $d=6$ and $d=7$, our exponents of $0.917$ and $0.947$ are much better than their exponents of $0.999998$ and $0.99999998$ obtained by substituting into the formula above respectively. Even as $d\rightarrow \infty$, the denominator in our exponent grows as $\Theta(3^{d/3})$, which is asymptotically smaller than $\binom d 22^{\binom d 2}$.  

We now give a brief explanation of the growth rate $\Theta(3^{d/3})$ obtained in Theorem \ref{thm: S_d count}. This grows slightly faster than $2^{d/2}$, and we first explain why a $2^{d/2}$ factor could appear. It turns out that this is due to the splitting of $c$. Consider a tower of subgroups $S_2\times\cdots\times S_2\leq S_4\times S_2\times \cdots \times S_2\leq S_{2\lfloor \frac{d}{2}\rfloor-1}\times S_2\leq S_d$, the conjugacy class of transpositions $c$ splits each time we go down the tower, and each time it splits we incur a factor of $2$. The reason why this factor appears is technical, but we try to give a rough idea here. When a conjugacy class splits after restricting it to a smaller subgroup, say into two different classes $c_1,c_2$, we now need to take quotients in both directions. This means that in order for there to be homological stability in both directions, both $n_1$ and $n_2$ have to be large. If say $n_1$ is small, then we find that we require $n_2$ to be twice as large as needed previously, so that the stronger homological stability can compensate for the lack of homological stability in the direction. The reason why our growth rate is $3^{d/3}$ instead of $2^{d/2}$ is due to another slight inefficiency, this time coming from the filtration relating the cohomology of Hurwitz space and maximal Hurwitz spaces, where we need a degree shift when we use the long exact sequence of homology coming from this filtration.

We also observe a striking pattern from our homological computations of $\CHur^{G,c}_n$ in Figure \ref{fig: homology computation S_d}, which is that the stability staircase always repeats under the shift $(n,i)\mapsto(n+6,i+2)$. We formulate this as the following conjecture.
\begin{conjecture}\label{conj: S_d}
Fix some $d\ge 3$, and let $G=S_d$ and $c=c^{\text{trn}}_{S_d}$. For all $i\ge 0$, the homology $H_i(\CHur^{G,c}_n,\Q)$ stabilizes when $n\ge 2d-3+6\lceil \frac i 2\rceil$. In particular, the rational homology of $\CHur^{G,c}$ has stability slope $1/3$. Moreover, there is secondary stability with respect to some element in degrees $(n,i)=(6,2)$ with slope $>1/3$. 

Arithmetically, this suggests that there exists a constant $e_d<5/6$, where for all sufficiently large $q$ there exist constants $\beta_n=\beta_{n,d,q}$ which are periodic modulo $3$ in $n$ so that 
$$N^{\text{sqf}}_{S_d}(2n)=(1-q^{-2})q^{2n}+\beta_n(q^{2n})^{5/6}+
O\left((q^{2n})^{e_d}\right).$$

We also conjecture the number field analogue: let $N^+_d(X),N^-_d(X)$ be the number of $S_d$-degree $d$ fields with positive and negative squarefree discriminant respectively with absolute discriminant at most $X$. Then there exist constants $\alpha^\pm_d,\beta^\pm_d$ and $e_d<5/6$ such that
$$N^\pm_d(X)=\alpha^\pm_d X+\beta^\pm_d X^{5/6}+O(X^{e_d}).$$
\end{conjecture}

We remark that the same stability slope and range should be true for the full Hurwitz space $\Hur^{S_d,c^{\text{trn}}_{S_d}}_n$ as well, but in this case we would need to quotient by multiple elements as explained in Section \ref{sec: general results}.

The conjecture also suggests a connection with homological stability for the moduli space of curves $\mathcal M_g$. For odd $n$, consider the locus with transposition boundary monodromy. The associated degree-$d$ cover has $n$ simple branch points over $\A^1$ and one over $\infty$, so Riemann-Hurwitz gives $n=2g+2d-3$. Forgetting the covering map gives a map to $\mathcal M_g$, and the conjectural stability bound of $n\ge (2d-3)+6\lceil \frac i 2 \rceil$ becomes $g\ge 3\lceil \frac i 2 \rceil$. For even $n$, the locus with trivial boundary monodromy also yields the same genus bound. The conjectural slope of $1/3$ in algebraic degree therefore becomes $2/3$ in genus. Combining classical homological stability \cite{Boldsen_2010,Wah13} with the rational refinement in \cite[Theorem B(iii)]{GKRW19} gives stability for $H_i(\mathcal M_g,\mathbb Q)$ in the exact same range $g\ge 3\lceil \frac i 2 \rceil$. We also note that this is the best possible general genus bound for rational homological stability of $\mathcal M_g$ by work of \cite{Morita_2003}, see \cite[Introduction]{Boldsen_2010} for an explanation. These numerical agreements suggest that the instability of Hurwitz spaces is closely related to that of moduli spaces of curves.

Finally, without going into too much detail, we mention that our results also have implications for the asymptotic Picard rank conjecture as discussed in \cite[Section 1.2]{landesman2025homologicalstabilityhurwitzspaces}. Indeed, the relevant Hurwitz space is the same, and an improvement in the stability line would directly result in a better lower bound for $n$ in terms of $d$. Here, $n$ is the number of branch points, with $n=2g+2d-2$. 
\begin{theorem}\label{thm: picard rank}
The Picard rank conjecture holds whenever $n\ge 10\cdot 3^{d/3}$. Moreover, Conjecture \ref{conj: S_d}, together with the results of \cite{Deopurkar_Patel_2015,Mullane_2023}, implies the full Picard rank conjecture.
\end{theorem}

\subsubsection{Other Hurwitz spaces}\label{sec: other hurwitz spaces} 
Finally, we will apply our main theorem to obtain stability slopes for various other Hurwitz spaces $\Hur^{G,c}_n$. We proceed in the same way as we did in the previous sections for $G=H\rtimes \Z/2\Z,c=c^{\text{rfl}}_G$ and $G=S_d,c=c^{\text{trn}}_G$, but now for other pairs of $G,c$. To be precise, we will do homological computations to prove increasingly optimal bounds on slope and the secondary exponent of $\#\Hur^{G,c}_n(\F_q)$ (see the discussion on the Grothendieck-Lefschetz trace formula in Section \ref{sec: intro optimality}), in the spirit of Corollary \ref{cor: optimal bound} when $c$ is non-splitting, and the generalization Corollary \ref{cor: general optimal} in the general case. We summarize here the results of the computations, which are carried out in Section \ref{sec: compute other hurwitz spaces}.

As computations take longer the larger $c$ is, we focus on cases where $c$ is small. In particular, we consider all other connected racks with size $\le 10$ which appear as a conjugacy class of some group. We first note that the conjugacy class of $4$-cycles in $S_4$ has the same rational homology as the transpositions in $S_4$ by \cite[Example 5.17(3), Lemma 5.18(3)]{Andruskiewitsch_Grana_2003}, so we omit this from the discussion. Figure \ref{fig: table G A_4} details the improved bound when $c$ is the conjugacy class of $(123)$ in $G=A_4$. All other examples happen to come from semidirect products, which we present in Figure \ref{fig: table G semidirect}. We again omit the $+\epsilon$ term for all exponent bounds. Conveniently, all of these $c$ are non-splitting, which makes our analysis later easier.

\begin{figure}[!htbp]
\centering
\begin{tabular}{>{\centering\arraybackslash}m{0.5cm}
>{\centering\arraybackslash}m{3cm}
>{\centering\arraybackslash}m{4cm}
>{\centering\arraybackslash}m{2.25cm}
>{\centering\arraybackslash}m{2.25cm}
>{\centering\arraybackslash}m{1cm}}
\toprule
$\mathbf{|c|}$ &
$\mathbf{G}$ & $\mathbf{c}$ & \textbf{Slope} $\mathbf{\ge}$ & \textbf{Exponent} $\mathbf{\le}$ & $\mathbf{(N,I)}$ \\
\midrule
$4$ & $A_4$ & Conjugates of $(123)$ & $\frac 8 {29}\approx 0.276$ & $\frac{25}{29}\approx \mathbf{0.862}$ & $(28,7)$ \\
\bottomrule
\end{tabular}
\caption{Bounds on slope and exponent for the conjugacy class of $(123)$ in $A_4$.}\label{fig: table G A_4}
\end{figure}

\begin{figure}[!htbp]
\centering
\begin{tabular}{>{\centering\arraybackslash}m{0.5cm}
>{\centering\arraybackslash}m{3cm}
>{\centering\arraybackslash}m{4cm}
>{\centering\arraybackslash}m{2.25cm}
>{\centering\arraybackslash}m{2.25cm}
>{\centering\arraybackslash}m{1cm}}
\toprule
$\mathbf{|c|}$ &
$\mathbf{G}$ & \textbf{Generating automorphism} & \textbf{Slope} $\mathbf{\ge}$ & \textbf{Exponent} $\mathbf{\le}$ & $\mathbf{(N,I)}$ \\
\midrule
$5$ & $\Z/5\Z \rtimes \Z/4\Z$ & $x\mapsto 2x$ & $\frac 2 5 \approx 0.4$ & $\frac{4}{5}= \mathbf{0.8}$& $(14,5)$ \\
\midrule
$7$ & $\Z/7\Z \rtimes \Z/3\Z$ & $x\mapsto 2x$ & $\frac 5 {11}\approx 0.455$ & $\frac{17}{22}\approx \mathbf{0.773}$ & $(10,4)$ \\
\midrule
$7$ & $\Z/7\Z\rtimes \Z/6\Z$ & $x\mapsto 3x$ & $\frac 4 {15} \approx 0.267$ & $\frac {13}{15} \approx  \mathbf{0.867}$ & $(14,3)$ \\
\midrule
$8$ & $(\Z/2\Z)^3\rtimes \Z/7\Z$ & $(x,y,z)\mapsto (z,x+z,y)$ & $\frac 4 7 \approx 0.571$ & $\frac 5 7\approx \mathbf{0.714}$ & $(6,3)$ \\
\midrule
$9$ & $\Z/9\Z\rtimes \Z/6\Z$ & $x\mapsto 2x$ & $\frac 3 {11}\approx 0.273$& $\frac {19}{22} \approx \mathbf{0.864}$ & $(10,2)$\\
\midrule
$9$ & $(\Z/3\Z)^2\rtimes \Z/4\Z$ & $(x,y)\mapsto (y,-x)$ & $\frac 4 {11} \approx 0.364$ & $\frac 9 {11}\approx \mathbf{0.818}$& $(10,3)$ \\
\midrule
$9$ & $(\Z/3\Z)^2\rtimes \Z/6\Z$ & $(x,y)\mapsto(-x+y,-y)$ & $\frac 1 4=0.25$& $\frac 7 8=\mathbf{0.875}$  & $(15,3)$\\
\midrule
$9$ & $(\Z/3\Z)^2\rtimes \Z/8\Z$ & $(x,y)\mapsto (x-y,x+y)$ & $\frac 1 2=0.5$ & $\frac 3 4 = \mathbf{0.75}$ & $(7,3)$\\
\bottomrule
\end{tabular}
\caption{Bounds on slope and exponent for conjugacy classes in various semidirect products $G=A\rtimes B$. Here, all conjugacy classes are taken to be $c=\{(x,1)\mid x\in A\}$.}\label{fig: table G semidirect}
\end{figure}

We also observe some interesting patterns from our computations in Figure \ref{fig: computations other hurwitz}. For example, a number of these exhibit clear stability staircases. For example, the conjugacy class of $(123)$ in $A_4$ has a stability staircase with step $(\Delta n,\Delta i)=(6,2)$. Moreover, the conjugacy classes coming from the groups $\Z/5\Z\rtimes \Z/4\Z,\Z/7\Z\rtimes \Z/3\Z, \Z/7\Z\rtimes \Z/6\Z$ and $\Z/9\Z\rtimes \Z/6\Z$ have stability staircases with steps $(\Delta n,\Delta i)=(4,2), (3,2), (6,2)$ and $(6,2)$ respectively. The other examples do not have a clear stability pattern. 

On top of this, the conjugacy class of $(123)$ in $A_4$ seems to exhibit secondary stability. Furthermore, the cases $\Z/5\Z\rtimes \Z/4\Z$ and $\Z/7\Z\rtimes \Z/6\Z$ seem to exhibit some variant of secondary stability, where in the first case the last nonzero homologies in the $i=2,4,6$ rows have dimensions in the ratio $1:3:5$, and in the second case the last nonzero homologies in the $i=2,4$ rows have the ratio $3:5$. It would be interesting to account for these observations in general.

\subsection{Acknowledgments} I would like to thank my advisor Will Sawin for providing me with much valuable guidance along the way. I would also like to thank Jordan Ellenberg, Aaron Landesman, Ishan Levy and Mark Shusterman for helpful discussions and feedback. 

ChatGPT was used to develop code for the homological computations in Section \ref{sec: computations} and the recurrence computations in Section \ref{sec: explicit slope computations}, and to assist with copyediting. I take full responsibility for the correctness of the code and its outputs. The mathematical results and proofs are my own.

\section{Preliminaries}\label{sec: preliminaries}
\subsection{Racks and Hurwitz spaces} Following \cite{landesman2025homologicalstabilityhurwitzspaces,landesman2025stablehomologyhurwitzmodules}, we define racks as well as Hurwitz spaces $\Hur^c_n$ associated to a finite rack $c$ over $\C$. This generalizes the more classical definition of Hurwitz spaces $\Hur^{G,c}_n$ associated to a union of conjugacy classes $c$ of a group $G$. These are more directly relevant to arithmetic applications, and we will recall the construction of these over $\Z[\frac{1}{|G|}]$. We do, however, need the generalized version for racks as some of the stability theorems involve quotient racks.

\subsubsection{Racks} We first define racks, quotient racks and components of racks.
\begin{definition}\label{def: rack}
A \textit{rack} is a set $c$ with a binary operation $(a,b)\mapsto b^a$ satisfying (i) $(c^a)^{b^a} = (c^b)^a$ and (ii) for each $a,b\in c$ there is a unique element in $c$ which we denote ${}^ab$ such that $({}^ab)^a=b$. 

One can also replace both conditions (i) and (ii) with the following equivalent condition: for all $n\ge 1$ and $1\le i\le n-1$, the  operation $\sigma_i\colon c^n\rightarrow c^n$ sending
$$(a_1,\ldots,a_i,a_{i+1},\ldots a_n)\mapsto (a_1,\ldots, a_{i+1},a_i^{a_{i+1}},\ldots , a_n)$$
where $a_i,a_{i+1}$ are swapped with $a_{i+1}$ acting on $a_i$, is an action of the Artin braid group on $n$ strands
$$B_n\coloneqq \langle \sigma_1,\ldots, \sigma_{n-1}\mid \sigma_i\sigma_{i+1}\sigma_i=\sigma_{i+1}\sigma_i\sigma_{i+1},\ \sigma_i\sigma_j=\sigma_j\sigma_i \text{ if }|i-j|>1\rangle$$
on $c^n$.
\end{definition}

Now, we discuss some properties and constructions of racks, in analogy to group theory.
\begin{definition}
Let $c$ be a rack and $c'\subseteq c$ be a subrack. The \textit{normalizer} $N_c(c')$ is defined to be the set of $a\in c$ such that $b^a\in c'$ for all $b\in c'$. Say that $c'\subseteq c$ is \textit{normal} if $N_c(c')=c$. 

If $c'\subseteq c$ is normal, then we can construct the \textit{quotient rack} $c/c'$ whose elements are the equivalence classes of elements $\overline b$ generated by equivalences of the form $b\sim b^a$ for $a\in c'$. One can check that the rack operation descends to $c/c'$ (see \cite[Lemma 2.3.1]{landesman2025stablehomologyhurwitzmodules}), so we define ${\overline b}^{\overline a}\coloneqq \overline{b^a}.$

As a special case when $c'=c$, we call the elements of $c/c$ the \textit{components} of the rack $c$. There is a natural projection of racks $c\rightarrow c/c$, and we call each fiber of this map a \textit{component} of $c$.
\end{definition}
Lastly, we define the reduced structure group, which basically takes the place of the group $G$ for a general rack $c$. We will need this to define connected Hurwitz spaces later.
\begin{definition}
The \textit{reduced structure group} $G^0_c$ of a rack $c$ is the subgroup of $\Aut(c)$ generated by the operations $(-)^a\colon b\mapsto b^a$ for $a\in c$. We also write $\ord_c(a)$ to be the order of the operator $(-)^a$ in $G^0_c\subseteq \Aut(c)$.
\end{definition}

\begin{example}
As mentioned earlier, the most important example of a rack $c$ for us is when $c=c_1\sqcup\cdots \sqcup c_k$ is a union of conjugacy classes of a group $G$. Then $c$ is a rack under the conjugation operation $g^h=h^{-1}gh$. Let us assume that $c$ generates $G$; otherwise we can just replace $G$ with the subgroup $\langle c\rangle$ generated by $c$. Then, the components of $c$ are simply the individual conjugacy classes $c_1,\ldots, c_k$, and the reduced structure group is isomorphic to $G/Z(G)$ where $Z(G)$ is the center of $G$.
\end{example}
\subsubsection{Topological Hurwitz spaces $\Hur^c_n$}
Now, we define the Hurwitz space associated to a rack $c$ as a complex variety. We refer the reader to \cite[Section 2.2]{prequel}, \cite[Section 2.1.1]{ma2026weightfiltrationhurwitzspaces}, and \cite[Section 2.2]{landesman2025homologicalstabilityhurwitzspaces} for more details. First, we define configuration space $\Conf_n$ and multicolored configuration space $\Conf_{n_1,\ldots, n_k}$ over $\Spec \Z$. 
\begin{definition}
Let $\Conf_n$ be the configuration space of $n$ distinct unordered points on $\A^1$. This can be defined as the open subscheme of $\A^n$ where the discriminant is invertible, by identifying the $n$ points with the coefficients of the monic polynomial with these $n$ roots. Alternatively, we can consider the open subset $U\subseteq \A^n$ where no two coordinates are equal, and define $\Conf_n\coloneqq U/S_n$ where $S_n$ permutes the coordinates. 

Likewise, for non-negative integers $n_1,\ldots, n_k$ with $n_1+\cdots +n_k=n$, we let $\Conf_{n_1,\ldots, n_k}$ be the multicolored configuration space of $n$ distinct points on $\A^1$ with $n_i$ points of each color that are indistinguishable. We can define this as the open subscheme of $\A^{n_1}\times \cdots \times \A^{n_k}$ of $k$ monic polynomials of degree $n_i$ where the discriminant of the product is invertible, or as the quotient $U/S_{n_1}\times \cdots \times S_{n_k}$.
\end{definition}

The complex analytification of $\Conf_n$ has fundamental group $\pi_1(\Conf_n)\cong B_n$; in fact, it is the Eilenberg-MacLane space $K(B_n,1)$. Here, the usual generators $\sigma_i$ of the braid group correspond to the half Dehn twist swapping the $i$ and $i+1$ points by moving them counterclockwise around each other. Similarly, we have $\pi_1(\Conf_{n_1,\ldots, n_k})\cong B_{n_1,\ldots, n_k}$, which is the subgroup of $B_n$ that preserves the first $n_1$ strands, next $n_2$ strands and so on. For notational simplicity, we will refer to both complex analytifications also as $\Conf_n$ and $\Conf_{n_1,\ldots, n_k}$. 

We are now ready to define Hurwitz spaces and some of their variants over $\C$.
\begin{definition}\label{def: hurwitz space}
For a rack $c$, let the \textit{Hurwitz space} $\Hur_n^c$ over $\C$ be the finite étale cover of $\Conf_n$ over $\C$ with fibers $c^n$ and monodromy action given by the $B_n$ action on $c^n$. More generally, given a rack $c=c_1\sqcup \cdots \sqcup c_k$ where each $c_i$ is a union of components, we can define the multicolored analogues as follows. The \textit{multicolored Hurwitz space} $\Hur_{n_1,\ldots, n_k}^{c_1,\ldots, c_k}$ over $\C$ is the finite étale cover of $\Conf_{n_1,\ldots, n_k}$ over $\C$ with fibers $c_1^{n_1}\times \cdots \times c_k^{n_k}$ with the monodromy action of $B_{n_1,\ldots, n_k}$. Since $\Conf_{n_1,\ldots, n_k}$ is itself a finite étale cover of $\Conf_{n}$, we can also describe $\Hur_{n_1,\ldots, n_k}^{c_1,\ldots, c_k}$ as the finite étale cover of $\Conf_n$ with fibers $c^{n_1,\ldots, n_k}\subseteq c^n$ which is the subset of $n$-tuples with $n_i$ elements in $c_i$ for each $i$, with the action of $B_n$.

Furthermore, we define the \textit{connected Hurwitz space} $\CHur^c_n$ over $\C$ to be the finite étale cover of $\Conf_n$ over $\C$ with fibers $(c^n)^\circ$, which is defined to be the subset of $c^n$ consisting of $n$-tuples of elements in $c$ whose actions generate the reduced structure group $G^0_c$, with the same monodromy action of $B_n$. Likewise, the \textit{connected multicolored Hurwitz space }$\CHur^{c_1,\ldots, c_k}_{n_1,\ldots, n_k}$ can be defined as the finite étale cover of $\Conf_n$ with fibers $(c^{n_1,\ldots, n_k})^\circ\coloneqq (c^n)^\circ \cap c^{n_1,\ldots, n_k}$. 

Finally, we define the \textit{maximal Hurwitz space} $\MHur^c_n$ over $\C$ which will be introduced in later versions of \cite{landesman2025homologicalstabilityhurwitzspaces,landesman2025stablehomologyhurwitzmodules}. It is defined to be the finite étale cover of $\Conf_n$ over $\C$ with fibers $(c^n)^{\max}\coloneqq c^n\setminus \left(\bigcup_{c'\subsetneq c} c'^n\right)$. The \textit{maximal multicolored Hurwitz space} $\MHur^{c_1,\ldots, c_k}_{n_1,\ldots, n_k}$ is then the finite étale cover of $\Conf_n$ with fibers $(c^{n_1,\ldots, n_k})^{\max}\coloneqq (c^n)^{\max}\cap c^{n_1,\ldots, n_k}$. 

We can summarize these in the following six equations

\begin{alignat*}{2}
\Hur^c_n
&\coloneqq \widetilde{\Conf_n}\times_{B_n}c^n,
&\qquad \Hur^{c_1,\ldots,c_k}_{n_1,\ldots,n_k}
&\coloneqq \widetilde{\Conf_n}\times_{B_n}c^{n_1,\ldots,n_k},\\
\CHur^c_n
&\coloneqq \widetilde{\Conf_n}\times_{B_n}(c^n)^\circ,
&\qquad \CHur^{c_1,\ldots,c_k}_{n_1,\ldots,n_k}
&\coloneqq \widetilde{\Conf_n}\times_{B_n}(c^{n_1,\ldots,n_k})^\circ,\\
\MHur^c_n
&\coloneqq \widetilde{\Conf_n}\times_{B_n}(c^n)^{\max},
&\qquad \MHur^{c_1,\ldots,c_k}_{n_1,\ldots,n_k}
&\coloneqq \widetilde{\Conf_n}\times_{B_n}(c^{n_1,\ldots,n_k})^{\max}.
\end{alignat*}
where $\widetilde{\Conf_n}$ is the universal cover of $\Conf_n$.
\end{definition}

One can think of these Hurwitz spaces as configuration space with extra data where each point is labeled by an element of $c$, and moving these points changes the labels according to the group action. 

We note that both $\CHur^c_n$ and $\MHur^c_n$ are open and closed subsets (i.e. unions of connected components) of $\Hur^c_n$. In fact, by definition, $\MHur^c_n$ is simply the union of the connected components which do not appear in $\Hur^{c'}_n$ for a proper subrack $c'\subsetneq c$, and this property turns out to be very useful later. 

For a general rack $c$, we have that $\MHur^c_n$ is an open and closed subset of $\CHur^c_n$. In fact, it turns out that these are very closely related. Indeed, one can check that a component of $\CHur^c_n$ belongs to $\MHur^c_n$ if and only if the corresponding tuple contains an element of every component of the rack $c$. In particular, if $c$ is connected, then $\CHur^c_n=\MHur^c_n$, and if $c=c_1\sqcup \cdots \sqcup c_k$ where each $c_i$ is a single connected component, then $\CHur^{c_1,\ldots, c_k}_{n_1,\ldots, n_k}=\MHur^{c_1,\ldots, c_k}_{n_1,\ldots, n_k}$ when all $n_i\ge 1$. 

In the special case when $c\subseteq G$ is a union of conjugacy classes, we write the Hurwitz space as $\Hur^{G,c}_n\coloneqq \Hur^c_n$, as we already saw in the introduction. In this case, $\Hur^{G,c}_n$ has the following moduli interpretation, as given in \cite[Section 2.3]{EVW16}. It is the moduli space of $n$-branched $G$-covers of a disk $D$ that have a
marked point on the boundary of the cover, such that the local monodromy around each branch
point lies in $c$. We also write $\CHur^{G,c}_n\coloneqq \CHur^c_n$, $\MHur^{G,c}_n\coloneqq \MHur^c_n$, and similarly for the multicolored versions. We note that if $c$ is a single conjugacy class generating $G$, then $\CHur^{G,c}_n$ parametrizes exactly those $G$-covers which are connected, as suggested by the name.
\subsubsection{Algebraic Hurwitz spaces $\Hur^{G,c}_n$} For our arithmetic applications over $\F_q[t]$, we need to define Hurwitz spaces over $\Spec \Z_p$ rather than just over the complex numbers, so that we can take the fiber over $\F_q$. This also allows us to compare the homology of the fibers over $\F_q$ and $\C$. It is not known how to construct Hurwitz spaces associated to a general rack $c$, so we restrict our attention to the special case of $\Hur^{G,c}_n$ where this is possible. 

We use the construction given in \cite[Section 2.1]{landesman2025cohenlenstramomentsfunctionfields}, which we summarize in the following definition.
\begin{definition}\label{def: algebraic hurwitz space}
Let $c$ be a union of conjugacy classes in a group $G$, where $p\nmid |G|$, and assume that $c$ is closed under $p$-th powering. Then, one can construct the Hurwitz space $\Hur^{G,c}_n$ over $\Spec\Z_p$ as shown in \cite[Section 2.1]{landesman2025cohenlenstramomentsfunctionfields}. Roughly, the $T$-points of this are finite $G$-covers over $\mathbb P^1_T$ where the local monodromy of a point over $\A^1_T$ lies in $c$. Furthermore, we mark an unramified point in the cover ``near infinity'' like for topological Hurwitz spaces above, which can be done easily if $\infty$ is not a branch point. Otherwise, if $\infty$ is a branch point with inertia group of order $e$, one uses a root stack of order $e$ at $\infty$ as defined in \cite[Definition 2.2.4]{Cadman_2007}, which lets us ``mark $1/e$ of a ramification point''.
\end{definition}
\subsection{Cobar-complexes and the homology of Hurwitz spaces} To apply the Grothendieck-Lefschetz trace formula, we need to compute the rational homology of these Hurwitz spaces over $\C$. One way to do this is to look at the cellular decomposition of Hurwitz space, and the corresponding cellular chain complex is what we call the (co)bar-complex. This was carried out in \cite{ETW17}, and in fact, they did this in the more general context of braided vector spaces $V$, which allows one to compute the homology of local systems like $V^{\otimes n}$ on $\Conf_n$ as well, which turns out to be useful for other arithmetic sums like in \cite{ellenberglandesman,bergstrom2023hyperellipticcurvesscanningmap,ES26,prequel}. Hence, we choose to work in this generalized setting of braided vector spaces, and most of our theorems will work here with no additional effort.

We now introduce braided vector spaces and (co)bar-complexes following \cite{ETW17,KS20,prequel,ma2026weightfiltrationhurwitzspaces}. However, our convention in this paper will be \textit{dual} to that in the above papers. This is done so that we look at chain complexes and homology instead of cochain complexes and cohomology, to better match the notation in \cite{landesman2025homologicalstabilityhurwitzspaces,landesman2025cohenlenstramomentsfunctionfields,landesman2025stablehomologyhurwitzmodules}. 

\subsubsection{Braided vector spaces}
Throughout this paper, we let $k$ be a characteristic zero field which we will use as our coefficient field.

\begin{definition}
A braided vector space $V$ over $k$ is a finite-dimensional $k$-vector space equipped with an invertible braiding $R\colon V\otimes V\xrightarrow \sim V \otimes V$ which satisfies the Yang-Baxter equation
$$(R\otimes \id)\circ (\id \otimes R)\circ (R\otimes \id)=(\id \otimes R)\circ (R\otimes \id)\circ (\id \otimes R).$$
This equation is exactly equivalent to the condition that $V^{\otimes n}$ is a representation of $B_n$ for all $n$, where the generator $\sigma_i$ acts on the $i$-th and $i+1$-st copies of $V$ by $R$. 
\end{definition}
\begin{example}
The following are the main constructions of braided vector spaces that we will use.
\begin{enumerate}
    \item We can construct the braided vector space $k c$ coming from a rack $c$. This has basis $c$ with braiding $R(a\otimes b)=b\otimes a^b$.
    \item We can twist any braided vector space $V$ by a constant $\zeta\in k^{\times}$ (which is usually a root of unity), to obtain the twisted braided vector space $V_\zeta$ with braiding $R_{V_\zeta}(a\otimes b)=\zeta\cdot R_V(a\otimes b)$. 
\end{enumerate}
\end{example}

There are other more general constructions; for example, we can twist $k c$ by an arbitrary cocycle, we can take the direct sum of two braided vector spaces, or even consider a braided monoidal category. However, we will not use these, and we refer the reader to \cite[Section 5]{ES26} and \cite[Section 2.2]{ma2026weightfiltrationhurwitzspaces} for more details.

\subsubsection{Tensor algebras}
From now on, let $V$ be a braided vector space over $k$. Next, we define braiding operators on $V^{\otimes n}$.
\begin{definition}
Given a permutation $\sigma\in S_n$, we define the \textit{braiding operator} $R_\sigma\colon V^{\otimes n}\xrightarrow \sim V^{\otimes n}$ as follows. There is a Matsumoto lift $\widetilde{\sigma}\in B_n$ given by expressing $\sigma$ as a minimum length (i.e. reduced) word in terms of transpositions $s_i=(i,i+1)$ and replacing these with the braid $\sigma_i$. This is well defined as any two reduced words are related by the braid relations. Then, $R_\sigma$ is defined to be the action of $\widetilde{\sigma}$ on the $B_n$-representation $V^{\otimes n}$.
\end{definition}

From this, we can define shuffle products and the dual notion of unshuffle coproducts.
\begin{definition}
The \textit{shuffle product} $\tau_{p,q}\colon V^{\otimes p}\otimes V^{\otimes q}\rightarrow V^{\otimes p+q}$ is defined as
$$\tau_{p,q}(u,v)=\sum_{\sigma\in \Sh(p,q)}R_\sigma(u\otimes v)$$
where we define the set of $(p,q)$-shuffles $\Sh(p,q)\subseteq S_{p+q}$ to be those permutations that preserve the order of the first $p$ elements and the last $q$ elements (i.e. $\sigma(1)<\cdots <\sigma(p)$ and $\sigma(p+1)<\cdots<\sigma(p+q)$).

The \textit{unshuffle coproduct} $\Delta_{p,q}\colon V^{\otimes{p+q}}\rightarrow V^{\otimes p}\otimes V^{\otimes q}$ is defined as
$$\Delta_{p,q}(w)=\sum_{\sigma\in \Sh(p,q)} R_{\sigma^{-1}}w$$
where we sum over the $(p,q)$-unshuffles which are the inverses of the $(p,q)$-shuffles.
\end{definition}

This allows us to define the tensor algebra as a graded bialgebra. One could also define a quantum shuffle algebra with trivial comultiplication and shuffle product as multiplication dually like in \cite{prequel,ma2026weightfiltrationhurwitzspaces}, but we will not need this here.
\begin{definition}
The \textit{tensor algebra} $T(V)=\bigoplus_{n=0}^{\infty} V^{\otimes n}$ is a graded bialgebra with free multiplication $\mu$, i.e. it has components $\mu_{p,q} \colon V^{\otimes p}\otimes V^{\otimes q}\xrightarrow \sim V^{\otimes p+q}$ given by the identity map, and comultiplication $\Delta$ with components $\Delta_{p,q}$ given by the unshuffle coproduct. The comultiplication map is also the unique algebra map that extends $\Delta\mid_{T^1(V)}=1\otimes \id + \id \otimes 1$. We sometimes write $T=T(V)$ and $T_n=T(V)_n=V^{\otimes n}$ as the $n$-th graded piece.
\end{definition}

\subsubsection{Cobar-complexes}
Finally, we define the cobar-cube and the cobar-complex of a tensor algebra.
\begin{definition}
The $n$-th \textit{cobar-cube} associated to the tensor algebra $T=T(V)$ is the commutative hypercube diagram whose vertices are labeled by the $2^{n-1}$ ordered partitions of $[n]$, with the vertex $\lambda=(\lambda_1,\ldots, \lambda_k)$ having the term $T_{\lambda_1}\otimes \cdots \otimes T_{\lambda_k}$, where the edges are given by the unshuffle coproducts. For example, for $n=3$ we have
\[\begin{tikzcd}
	{T_3} & {T_2\otimes T_1} \\
	{T_1\otimes T_2} & {T_1\otimes T_1\otimes T_1.}
	\arrow["{\Delta_{2,1}}", from=1-1, to=1-2]
	\arrow["{\Delta_{1,2}}"', from=1-1, to=2-1]
	\arrow["{\Delta_{1,1}\otimes \id}", from=1-2, to=2-2]
	\arrow["{\id \otimes \Delta_{1,1}}"', from=2-1, to=2-2]
\end{tikzcd}\]
From this, we can form the $n$-th \textit{cobar-complex} 
$$B_n^\dagger (T)=\left\{T_n\rightarrow \bigoplus_{\substack{p+q=n\\ 0<p,q<n}}T_p\otimes T_q\rightarrow \bigoplus_{\substack{p+q+r=n\\ 0<p,q,r<n}}T_p\otimes T_q\otimes T_r \rightarrow \cdots \rightarrow T_1^{\otimes n}\right\}$$
by collapsing the maps in the cobar-cube and adding a Koszul sign twist so that the differentials $d$ of the complex satisfy $d^2=0$. Here, we place $T_n$ in degree $i=n-1$ and $T_1^{\otimes n}$ in degree $i=0$. This convention differs by a shift of $1$ from that used in some other papers, where the complex is concentrated in degrees $n$ to $1$ instead. One can also interpret this homological degree convention as placing $T_n$ in homological degree $n-1$.

One can define the full cobar-complex to be
$$B^\dagger(T)=\bigoplus_{n=0}^\infty B_n^\dagger(T)$$
where we define $B^\dagger_0(T)\coloneqq k$.
\end{definition}
Given a basis $v_1,\ldots, v_k$ of $V$, it is common to consider the basis of $T_{\lambda_1}\otimes \cdots T_{\lambda_k}$ written in the form of $v_{i^1_1}\cdots v_{i^1_{\lambda_1}}| \cdots| v_{i^k_1}\cdots v_{i^k_{\lambda_k}}$, and every term in the cobar-complex can be written as a linear combination of these basis vectors.

We chose this grading convention so that the cobar-complex computes the homology of $V^{\otimes n}$:
\begin{equation}\label{eqn: H_i(B_n,V^n)}
H_i(B_n,V^{\otimes n})=H_i(B_n^\dagger(T(V_{-1}))).
\end{equation}
This follows directly from taking the dual of \cite[Equation (2.1)]{prequel}. 

Lastly, we return to the setting of Hurwitz spaces. One can consider the pushforward of the constant sheaf $k$ on $\Hur^c_n$ to $\Conf_n$. This gives us a local system of $\Conf_n$, and since $\pi_1(\Conf_n)\cong B_n$, one can view this as a representation of $B_n$. It is not hard to see that this representation is given by $(k c)^{\otimes n}$. Hence, we have
$$H_i(\Hur^c_n,k)=H_i(B_n,(k c)^{\otimes n})=H_i(B_n^\dagger (T((k c)_{-1}))).$$

For multicolored Hurwitz space $\Hur^{c_1,\ldots, c_k}_{n_1,\ldots, n_k}$, the corresponding local system is the sub-representation $k c^{n_1,\ldots, n_k}\subseteq (k c)^{\otimes n}$, and to compute the homology from the cobar-complex, we simply look at the part of the cobar-complex with the appropriate grading of $n_1,\ldots, n_k$. To be precise, the decomposition $k c=\bigoplus_{i=1}^k (k c_i)$ induces an $\N^k$-multigrading $n_1,\ldots, n_k$ on the tensor algebra where the original algebraic grading $n=n_1+\cdots +n_k$ is the sum of the new gradings. This also induces a multigrading on the cobar-complex, and we denote $B_{n_1,\ldots, n_k}^\dagger(T((k c)_{-1}))\subseteq B_n^\dagger(T((k c)_{-1}))$ to be the subcomplex with grading $n_1,\ldots, n_k$. Then, we have
$$H_i(\Hur^{c_1,\ldots, c_k}_{n_1,\ldots, n_k},k)=H_i(B_n,k c^{n_1,\ldots, n_k})=H_i(B_{n_1,\ldots, n_k}^\dagger (T((k c)_{-1}))).$$
\subsection{Higher algebra} We will need the formalism of $\infty$-categories to keep track of higher homotopies which are required for iterated quotients which appear in the statement of higher order homological stability. Furthermore, it is essential that we quotient by $\E_2$-central elements, which we also define here.
\subsubsection{$\infty$-categories} We use the theory of $\infty$-categories developed by Joyal and Lurie, see \cite{HTT,HA,Kerodon}. We will work in the bounded-below derived $\infty$-category of $k$-vector spaces $\mathcal D^-(k)$, which is a stable symmetric monoidal $\infty$-category. This is also a subcategory of the unbounded derived $\infty$-category $\mathcal D(k)$ which is equivalent to the symmetric monoidal $\infty$-category of $Hk$ modules in spectra $\Mod_{Hk}(\Sp)$ where $Hk$ is the Eilenberg-Maclane spectrum of $k$. One could also choose to work with any commutative ring, or in the larger $\infty$-category of spectra like in \cite{landesman2025cohenlenstramomentsfunctionfields}.

The bounded-below derived $\infty$-category $\mathcal D^-(k)$ has the following concrete description, as we recall from \cite[Section 1.3.2]{HA} or \cite[Section 10.4]{Aaron_notes}. This will be needed later when we view iterated quotients from the perspective of cobar-complexes.
\begin{definition}\label{def: D-(k)}
The \textit{bounded-below derived $\infty$-category} $\mathcal D^-(k)$ is the differential graded nerve $N_{dg}(\Ch^-(k))$ of the differential graded category of bounded-below chain complexes of $k$-vector spaces $\Ch^-(k)$ whose objects are chain complexes $X_*$ where $X_n=0$ for $n\ll 0$. Concretely, $\mathcal D^-(k)$ is the simplicial set whose $n$-simplices are the set of ordered pairs $(\{X_i\}_{0\le i\le n},\{f_I\}_{I\subseteq [n],|I|\ge 2})$, where:
\begin{enumerate}
\item For $0\le i\le n$, $X_i$ is an object of $\Ch^-(k)$,
\item For every subset $I=\{i_0<i_1<\cdots<i_{m+1}\}\subseteq [n]$ with $m\geq 0$, $f_I$ is an element of $\Hom_m(X_{i_0},X_{i_{m+1}})$, satisfying the equation
$$df_I=\sum_{1\le j\le m} (-1)^j (f_{I\setminus\{i_j\}}-f_{\{i_j,\ldots,i_{m+1}\}}\circ f_{\{i_0,\ldots, i_j\}}).$$
Here, recall that given $f\in \Hom_m(X,Y)$, we define $df\coloneqq  d_Y\circ f + (-1)^{m+1} f \circ d_X\in \Hom_{m-1}(X,Y)$. 
\end{enumerate}
The maps between simplices are defined in the obvious way, which we omit here.
\end{definition}
To make sense of this, we see that $f_I$ records an $m$-homotopy between the chain complexes $X_{i_0}$ and $X_{i_{m+1}}$. For example, when $I=\{i_0,i_1\}$, $f_I$ is a chain map in $\Hom_0(X_0,X_1)$, and when $I=\{i_0,i_1,i_2\}$, $f_I$ is a chain homotopy living in $\Hom_1(X_0,X_2)$ between the chain maps $f_{\{i_1,i_2\}}\circ f_{\{i_0,i_1\}}$ and $f_{\{i_0,i_2\}}$. Hence, this definition allows us to keep track of all the homotopies at once.

Since we are essentially always working with chain complexes, we will simply use the notation $H_i(A)$ to denote the $i$-th homology group of an object $A\in \mathcal D^-(k)$. This is in contrast with the notation in \cite{landesman2025cohenlenstramomentsfunctionfields}, where stable homotopy groups $\pi_i(A)$ are used for spectra, but we mean exactly the same thing. Similarly, we will also use the notation of the shift operator $A[i]$ defined by $A[i]_n=A_{n-i}$, instead of the suspension operator $\Sigma^i A$ used for spectra. 
\subsubsection{Simplicial chains and $\E_1$-algebras} To any symmetric monoidal $\infty$-category there is a notion of an $\E_1$-algebra, which is an associative algebra in the $\infty$-categorical setting. We refer the reader to \cite[Section 4]{HA} for the full definition which requires the notion of $\infty$-operads.

Furthermore, to an $\E_1$-algebra $A$, one can define the $\infty$-categories of right modules $\RMod_A$, left modules $\LMod_A$ and bimodules $\BiMod_A$; these are stable $\infty$-categories by \cite[Proposition 4.3.3.8, Corollary 7.1.1.5]{HA}. For example, $A$ itself is an object in $\BiMod_A$. There are also forgetful functors $\BiMod_A\rightarrow \RMod_A$ and $\BiMod_A\rightarrow \LMod_A$, and there is a relative tensor product
$$-\otimes_A-\colon \BiMod_A \times \BiMod_A \rightarrow \BiMod_A$$
which is an $\infty$-functor, and also preserves cofibers separately in each variable \cite[Corollary 4.4.2.15]{HA}.

The main instance of an $\E_1$-algebra in $\mathcal D^-(k)$ that we will need is simplicial chains of local systems on configuration space coming from a braided vector space $V$ as follows. 
\begin{definition}\label{def: A_V}
For $n\ge 1$, define the chain complex $C_*(B_n,V^{\otimes n})\coloneqq B_n^\dagger(T(V_{-1}))$ as an object of $\mathcal D^-(k)$; this notation is justified by Equation \eqref{eqn: H_i(B_n,V^n)}. Then, we define
$$A_V \coloneqq \bigoplus_{n=0}^\infty C_*(B_n,V^{\otimes n})$$
to be the direct sum of $C_*(B_n,V^{\otimes n})$, where for $n=0$ we define $C_*(B_0,V^{\otimes 0})\coloneqq k$. Alternatively, we also have that $A_V=B^\dagger(T)$ is the full cobar-complex. This is an $\E_1$-algebra in $\mathcal D^-(k)$ with respect to the concatenation action; indeed, one can verify this on the level of cobar-complexes. 
\end{definition}

The following example specializes this construction to the setting of Hurwitz spaces, where we can also look at simplicial chains of connected Hurwitz spaces, and relate them via a filtration.
\begin{example}\label{eg: A_c, CA_c}
When we take $V=k c$ for a finite rack $c$, this construction gives the simplicial chains of Hurwitz spaces as defined in \cite[Notation 3.2.1]{landesman2025cohenlenstramomentsfunctionfields}, where
$$A_c\coloneqq C_*(\Hur^c,k)=\bigoplus_{n=0}^\infty C_*(\Hur^c_n,k).$$
Furthermore, we can do the same for connected and maximal Hurwitz spaces
$$CA_c\coloneqq C_*(\CHur^c,k)=\bigoplus_{n=0}^\infty C_*(\CHur^c_n,k)$$
$$MA_c\coloneqq C_*(\MHur^c,k)=\bigoplus_{n=0}^\infty C_*(\MHur^c_n,k)$$
which are both $A_c$-bimodules. 

One can relate the simplicial chains $A_c$ and maximal simplicial chains $MA_c$ using the following filtration of $A_c$-bimodules constructed in \cite[Construction 5.0.4, Lemma 5.0.5]{landesman2025homologicalstabilityhurwitzspaces}. Let $F_m\Hur^c$ be the union of components of $\Hur^c$ such that if $c'\subseteq c$ is the smallest subrack such that the component comes from $\Hur^{c'}$, then $|c'|\ge m$. Then, one can check that $F_mA_c\coloneqq C_*(F_m\Hur^c,k)$ is an $A_c$-bimodule, and there is an isomorphism
\begin{equation}\label{eqn: filtration A_c and MA_c}
F_mA_c/F_{m+1}A_c\cong \bigoplus_{c'\subseteq c,|c'|=m} MA_{c'} 
\end{equation}
of $A_c$-bimodules, where $MA_{c'}$ is given the structure of an $A_c$-bimodule by having elements of $c$ that are not in $c'$ act by zero. This will let us deduce properties of $A_c$ from $MA_c$ and vice versa, as we will see in Section \ref{sec: comparing A_c and MA_c}.
\end{example}
\subsubsection{Iterated quotients}
We first give the definition of $\E_2$-centrality with respect to an $\E_1$-algebra $A$, following \cite[Section 3.1.5]{landesman2025cohenlenstramomentsfunctionfields} or \cite[Section 5.3]{HA}. 
\begin{definition}
We say that $U\in H_d(A)$ is $\E_2$-central if there exists a map $U\colon A[d]\rightarrow A$ in $\BiMod_A$ sending $1$ to $U$. We will implicitly identify $U$ with one such choice of an $A$-bimodule map.
\end{definition}

\begin{example}
In the case when $A=A_c$ comes from Hurwitz space as in Example \ref{eg: A_c, CA_c}, \cite[Corollary 4.3.2]{landesman2025cohenlenstramomentsfunctionfields} showed that if a connected component $\gamma \in \pi_0\Hur^c$ is central, then the associated class $\gamma\in H_0(A)$ is $\E_2$-central. This means that a sum of such classes $\sum b_i\gamma_i$ is also $\E_2$-central. This will be our main source of examples for $\E_2$-central elements. We will also often require $\gamma_i$ to have the same algebraic degree.
\end{example}

We can now define quotients as well as iterated quotients by $\E_2$-central elements.
\begin{definition}\label{def: iterated quotients}
Let $U\in H_d(A)$ be an $\E_2$-central element. We define the \textit{quotient} $A/U\coloneqq \cofib(A[d]\xrightarrow U A)$ to be the cofiber of the map $U$ which is an object in $\BiMod_A$. In our setting of the derived $\infty$-category, this is also the cone of the map $U$. 

Let $U_1,\ldots, U_k$ be $\E_2$-central elements. We define the \textit{iterated quotient} $A/(U_i)=A/(U_1,\ldots, U_k)$ to be the relative tensor product $A/U_1\otimes_A \cdots \otimes_A A/U_k$ which is an object in $\BiMod_A$ as each quotient is in $\BiMod_A$. Note that the ordering of $U_i$ matters in this definition.

Furthermore, if $M$ is an $A$-bimodule, we will denote $M/(U_i)\coloneqq M\otimes_A A/(U_i)$.
\end{definition}
For any $A$-bimodule $M$, there is a cofiber long exact sequence
\begin{equation}\label{eqn: single quotient LES}
\cdots \rightarrow H_{i+1}(M/U)\rightarrow H_{i-d}(M)\xrightarrow U H_i(M)\rightarrow H_i(M/U)\rightarrow H_{i-d-1}(M)\rightarrow \cdots
\end{equation}
where $U\in H_d(A)$. Moreover, since $\otimes_A$ preserves cofibers in each variable, we also have the long exact sequence
\begin{equation}\label{eqn: iterated quotient LES}
\cdots \rightarrow H_{i-d_i}(M/(U_1,\ldots,\widehat{U_i},\ldots,U_k))\xrightarrow {U_i} H_i(M/(U_1,\ldots,\widehat{U_i},\ldots,U_k))\rightarrow H_i(M/(U_1,\ldots ,U_k))\rightarrow \cdots
\end{equation}
where $U_i\in H_{d_i}(A)$. 
\subsection{Homological stability from Landesman-Levy} We now state the key homological stability results proven by Landesman-Levy \cite{landesman2025homologicalstabilityhurwitzspaces,landesman2025cohenlenstramomentsfunctionfields,landesman2025stablehomologyhurwitzmodules}. We will not need these to prove our main theorem, Theorem \ref{thm: main theorem}, but it is sometimes convenient to use these when applying our main theorem, and we also need these for our arithmetic applications. Let $c$ be a finite rack whose components are $c_1,\ldots, c_k$. Let $\Hur^{c_1,\ldots, c_k}_{n_1,\ldots,n_k}$ be the relevant topological Hurwitz space, with $A_c=C_*(\Hur^c,k)$, $CA_c=C_*(\CHur^c,k)$ and $MA_c=C_*(\MHur^c,k)$ as in Example \ref{eg: A_c, CA_c}.

\subsubsection{Stable homology of Hurwitz spaces}
We first state the theorem showing that the stable homology of Hurwitz space is the homology of configuration space. We state the theorems with $k$-coefficients for convenience, although the actual theorems apply for some $\Z[\frac 1 N]$.
\begin{theorem}[{\cite[Theorem 1.4.2]{landesman2025homologicalstabilityhurwitzspaces}}]\label{thm: stable hom conf all directions} There exist constants $I$ and $J$ depending on $c$, so that for any $i\ge 0$, $n_1,\ldots, n_k>Ii+J$, and any connected component $Z\subseteq \CHur^{c_1,\ldots, c_k}_{n_1,\ldots, n_k}$, the projection map gives an isomorphism $$H_i(Z,k)\xrightarrow \sim H_i(\Conf_{n_1,\ldots, n_k},k).$$
\end{theorem}

When $k=1$, we recover that for any connected component $Z\subseteq \CHur^c_n$ and $n>Ii+J$ that $H_i(Z,k)\xrightarrow\sim H_i(\Conf_n,k)$. In particular, if $c$ is non-splitting, this also applies to all $Z\subseteq \Hur^c_n$ (as each connected component is part of $\CHur^{c'}_n$ for $c'\subseteq c$), and this is the content of \cite[Theorem 1.4.7]{landesman2025cohenlenstramomentsfunctionfields}.

We also have stability in a single direction as in the following theorem. Note that we can recover Theorem \ref{thm: stable hom conf all directions} above by applying this theorem iteratively for each component $c_1,\ldots, c_k$. 
\begin{theorem}[{\cite[Theorem 1.4.6]{landesman2025stablehomologyhurwitzmodules}}]\label{thm: stable hom conf one direction} There exist constants $I$ and $J$ depending on $c$, so that for any $i\ge 0$ and $n_1>Ii+J$, and any connected component $Z\subseteq \CHur^{c_1,\ldots, c_k}_{n_1,\ldots, n_k}$ which maps to a component $Z'\subseteq \CHur^{c_1/c_1,\ldots, c_k/c_1}_{n_1,\ldots, n_k}$ under the projection map induced by $c\rightarrow c/c_1$, this projection map gives an isomorphism
$$H_i(Z,k)\xrightarrow \sim H_i(Z',k).$$
\end{theorem}

\subsubsection{Stabilization maps for Hurwitz spaces}\label{sec: stab map for Hurwitz space}
We will also need information on the stabilization maps, which are the elements we take iterated quotients by. The following theorem says that for homological stability of $\MHur^c_n$, it is sufficient to take the quotient by a power of an element in each conjugacy class. Furthermore, let $c'=c_1'\sqcup \cdots \sqcup c_j'\subseteq c$ be a subrack consisting of $j$ components of $c$, and denote $n'=n_1+\cdots +n_j$ as the grading in $c'$, then if we only want to consider stability in the direction of $c'$, we only need to take quotients by those in $c_i'$.
\begin{theorem}[{\cite[Theorem 5.0.6]{landesman2025homologicalstabilityhurwitzspaces}}]\label{thm: LL CA_c iterated quotients}
There exist constants $I$ and $J$ depending on $c$, so that for any choice of elements $s_i\in c_i'$, we have that
$$H_i\left(\left(MA_c/(s_i^{\ord_c(s_i)})\right)_{n'}\right)=0$$
if $n'>Ii+J$, where $n'$ is the grading in $c'$.
\end{theorem}

Note that $s_i^{\ord_c(s_i)}$ is central in $\pi_0\Hur^c_n$, so the associated class is $\E_2$-central in $A_c$, and the iterated quotient is well defined. If $c$ only has a single component and we take $c'=c$, then we have that for any element $s\in c$ that $H_i(\CHur^c,k)=H_i(\MHur^c,k)$ is homologically stable with respect to $s^{\ord_c(s)}$. If $c$ has multiple components and we choose $c'=c_1$ to be one of them, then this shows that multiplication by $s_1^{\ord_c(s_1)}$ gives homological stability in the $n_1$ direction. 

There is also an analogous theorem for the full Hurwitz space $\Hur^c_n$, where $MA_c$ is replaced with $A_c$, but we need to take quotients by every element in $c'$. This theorem can also be deduced from the previous theorem via the filtration in Equation \eqref{eqn: filtration A_c and MA_c}.
\begin{theorem}[{\cite[Theorem 3.2.6]{landesman2025homologicalstabilityhurwitzspaces}}]\label{thm: LL A_c stability} There exist constants $I$ and $J$ depending on $c$, so that if we write $c'=\{s_1,\ldots, s_{|c'|}\}$ (here, we choose an arbitrary ordering), we have that
$$H_i\left(\left(A_c/(s_i^{\ord_c(s_i)})\right)_{n'}\right)=0$$
if $n'>Ii+J$, where $n'$ is the grading in $c'$.
\end{theorem}
\section{Statements of generalized homological theorems}\label{sec: general results} In this section, we discuss the general versions of our theorems as promised in Section \ref{sec: intro main theorem}, with the language introduced in Section \ref{sec: preliminaries}. We will do this in the setting of braided vector spaces with stability with respect to multiple $\E_2$-central elements. We will first discuss the setup and definitions in Section \ref{sec: setup and definition generalized homology}, then state the theorems in Section \ref{sec: subsec generalized theorem statements}. Finally, we will discuss some lemmas on extra stability elements in Section \ref{sec: lem extra elements}, which will be useful in our proofs later.
\subsection{Setup and definitions}\label{sec: setup and definition generalized homology}
Let $V$ be a braided vector space with simplicial chains $A=A_V$ as defined in Definition \ref{def: A_V} with algebraic grading $n$. Let $U_1,\ldots, U_k$ be $\E_2$-central elements, where $U_i$ is in homological degree $d_i$ and algebraic degree $a_i>0$. Note that since $d_i$ is not necessarily zero, multiplication by $U_i$ no longer goes only to the right (in the coordinates of $(n,i)$ like in Figure \ref{fig: (4,0) S_3 trans} or \ref{fig: (10,2) S_3 trans}), but goes up and right with slope $s_i=d_i/a_i$. Let $s_0=\max_is_i$ be the maximal slope which clearly must satisfy $0\le s_0<1$. We will need to account for this when we define what it means to be homologically stable with respect to these elements.
\begin{definition}\label{def: general hom stability}
We say that $A$ is \textit{homologically stable with respect to $U_1,\ldots, U_k$} if there is some slope $s>s_0$ (satisfying $0<s\le 1$) and some $t$ where $H_i\left((A/(U_1,\ldots, U_k))_n\right)$ vanishes whenever $i<sn+t$. In this case, we say that it is homologically stable \textit{with slope $s$}.
\end{definition}
The reason why we require $s>s_0$ is that multiplication by $U_i$ gives stable values for every line of slope $s_i$ and we want stability to be true for every such line. As secondary stability is not known in any case, we will actually only apply our homological theorems to arithmetic problems in the case when all $d_i=s_i=0$, but we believe it is useful in the future to carry this out in this generality.

We also have the following analogue of Definition \ref{def: intro stabilize in staircase}.  
\begin{definition}\label{def: staircase iterated version}
Let $N\ge 1$ and $I \ge 0$ be integers with $I<N$. We say that $A$ \textit{stabilizes in an $(N,I)$ staircase} with respect to $U_1,\ldots, U_k$ if for all integers $n\le N$ and $i\le I-N+n$ we have that 
\begin{equation}\label{eqn: vanishing homology iterated quotient}
H_{i+k+d_1+\cdots +d_k}\left((A/(U_1,\ldots, U_k))_{n+a_1+\cdots+a_k}\right)=0.
\end{equation}
Here, $n$ and $i$ can be taken to be negative integers. We also define stability for $A$-bimodules $M$ in the same way.
\end{definition}

Let us see how this specializes to Definition \ref{def: intro stabilize in staircase} when we take the case of Hurwitz space $A=A_c$, with $k=1$, $U=U_1$ with homological degree $d_1=0$ and algebraic degree $a=a_1$. In this case, the long exact sequence in \eqref{eqn: single quotient LES} tells us that $H_i((A/U)_{n+a})=0$ for $i\le I-N+n+1$ if and only if the stabilization map $U\colon H_i(A_n)\rightarrow H_i(A_{n+a})$ is an isomorphism for $i\le I-N+n$ and a surjection for $i=I-N+n+1$.

\begin{remark}\label{rmk: shifted}
It is important to remember that the algebraic grading and homological degree in Equation \eqref{eqn: vanishing homology iterated quotient} are shifted by $a$ and $k+d$ respectively. It will be clear why this is a rather natural choice from the lemmas we prove in Section \ref{sec: lem extra elements}.
\end{remark}

\subsubsection{Strong stability in a staircase}
We realize that the notion of stabilization in a staircase is rather weak in the following sense. For instance, consider the case of Hurwitz space with just one $U$ like in the introduction, and recall that stability in an $(N,I)$ staircase tells us that the stabilization map \eqref{eqn: stab map} is an isomorphism for $n\le N$ and $i\le I-N+n$. However, this does not guarantee that $H_i(\Hur^c_n,\Q)$ is actually the stable homology of $\Hur^c_n$. For that, we need the stronger statement that we have a chain of isomorphisms 
$$H_i(\Hur^c_n,\Q)\xrightarrow \sim H_i(\Hur^c_{n+a},\Q)\xrightarrow \sim H_i(\Hur^c_{n+2a},\Q)\xrightarrow \sim H_i(\Hur^c_{n+3a},\Q)\xrightarrow\sim\cdots $$
induced by the stabilization map $U$. This motivates us to introduce the following stronger version of Definition \ref{def: staircase iterated version}. 
\begin{definition}\label{def: strong staircase}
Let $N\ge 1$ and $I \ge 0$ be integers satisfying $I<N$. We say that $A$ \textit{stabilizes strongly in an $(N,I)$ staircase} with respect to $U_1,\ldots, U_k$ if for all $i,n$ satisfying both $i\le I-N+n$ and $i\le I+(n-N)s_0$ we have Equation \eqref{eqn: vanishing homology iterated quotient}. One can also define stability for $A$-bimodules $M$ in the same way.
\end{definition}

Again, it will become clear, for example, from the theorem statements in Section \ref{sec: subsec generalized theorem statements} and the lemmas in Section \ref{sec: lem extra elements}, why this stronger notion is natural.

\begin{figure}[!htbp]
    \centering
\begin{tikzpicture}[x=1cm,y=1cm,scale=0.5]
\drawnaxis{12}
\drawiaxis{6}
\drawstaircase{1}{7}{0}{black!10}{*}
\drawblock{8}{12}{0}{6}{black!10}{*}
\drawstaircase{3}{5}{0}{red!40}{0}
\drawblock{6}{12}{0}{2}{blue!40}{0}
\end{tikzpicture}
\begin{tikzpicture}[x=1cm,y=1cm,scale=0.5]
\drawnaxis{12}
\drawiaxis{6}
\drawstaircase{1}{7}{0}{black!10}{*}
\drawblock{8}{12}{0}{6}{black!10}{*}
\drawstaircase{3}{5}{0}{red!40}{0}
\drawblock{6}{12}{0}{2}{blue!40}{0}
\drawblock{8}{12}{3}{3}{blue!40}{0}
\drawblock{11}{12}{4}{4}{blue!40}{0}
\end{tikzpicture}
\caption{Illustration of usual stabilization and strong stabilization in an $(N,I)$ staircase for $(N,I)=(5,2)$, where the left and right panels have $s_0=0$ and $s_0=1/3$ respectively.}\label{fig: strong vs usual staircase}
\end{figure}
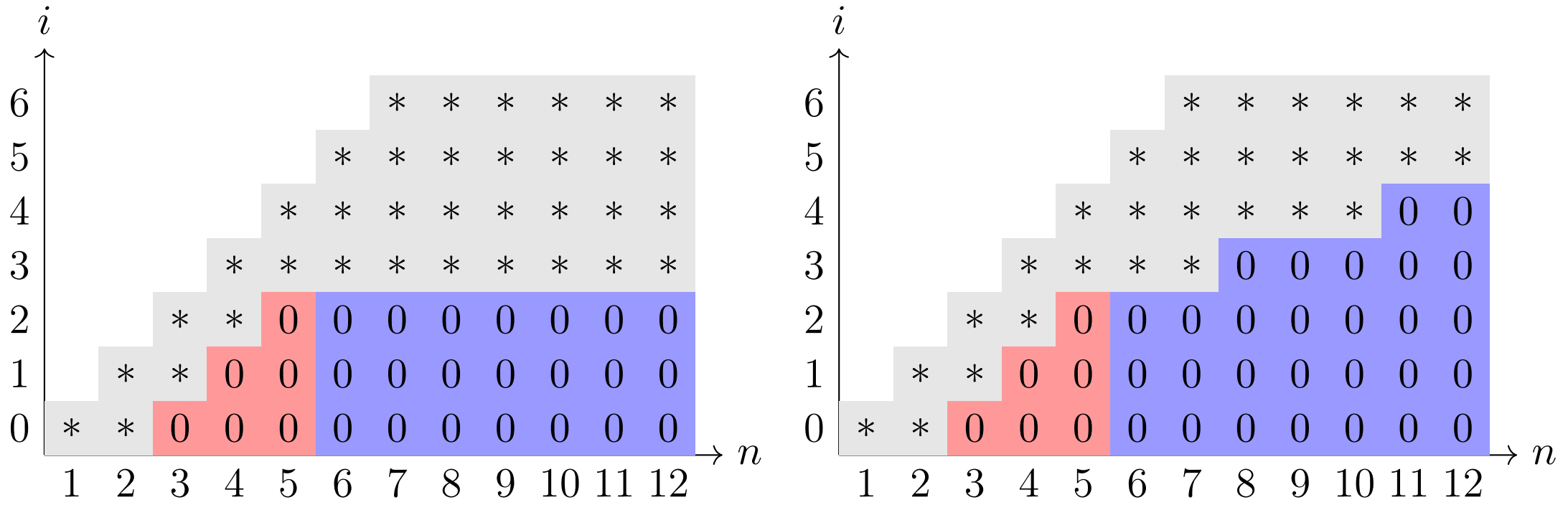

It is useful to illustrate and compare the regions of $(n,i)$ for which Equation \eqref{eqn: vanishing homology iterated quotient} is true in Definitions \ref{def: staircase iterated version} and \ref{def: strong staircase}, and we do this in Figure \ref{fig: strong vs usual staircase}. Here, usual stability in an $(N,I)$ staircase only guarantees vanishing of $(n,i)$ in the red staircase, while strong stability guarantees vanishing in both the red and blue regions. Going back to the one variable case, we see that the stronger definition guarantees the chain of isomorphisms needed to conclude that $H_i(\Hur^c_n,\Q)$ is indeed the stable value.

Somewhat surprisingly, we have that these two definitions are equivalent in the case when all $d_i=0$, i.e. when all $U_i\in H_0(A)$, which is what we are concerned with for our arithmetic applications. In particular, this also applies to the theorems in the introduction. Thus, we will usually not need to distinguish between these two notions.
\begin{theorem}\label{thm: usual vs strong stability}
Suppose that all $U_i\in H_0(A)$. Then $A$ stabilizes in an $(N,I)$ staircase with respect to $U_1,\ldots, U_k$ if and only if it stabilizes strongly in an $(N,I)$ staircase with respect to $U_1,\ldots, U_k$.
\end{theorem}
\subsubsection{Stability in a chosen direction}\label{sec: stability in chosen direction}
We can also modify the above setup and definitions for stability in a chosen direction. Much like how we considered a subrack $c'\subseteq c$ which is a union of components of $c$ in Section \ref{sec: stab map for Hurwitz space} with corresponding grading $n'$, we will consider a sub-braided vector space $V_1\subseteq V$ as follows. Suppose that we can write the vector space $V=V_1\oplus V_2$ such that the braiding $R_V$ restricts to 
\begin{equation}\label{eqn: V_1+V_2 condition}
R_V|_{V_i\otimes V_j} \colon V_i\otimes V_j\xrightarrow \sim V_j\otimes V_i
\end{equation}
for $i,j\in \{1,2\}$. This is equivalent to $V$ being the direct sum of the two braided vector spaces $V_1,V_2$ in the sense of \cite[Section 6.3]{prequel} or \cite[Section 5.2]{ES26}. Since the braiding preserves both $V_1$ and $V_2$, we can define a bigrading $(n_1,n_2)$ with $n_1+n_2=n$ on $V^{\otimes n}$, the tensor algebra, cobar-complex, and hence $A_V$. 

Let us redefine the properties of the $\E_2$-central elements $U_1,\ldots, U_k$ in this adjacent setup. Let $U_i$ have homological degree $d_i$ and algebraic grading $(a_i,b_i)$ so that each $a_i>0$. We will only be concerned about the algebraic grading $a_i$ in the $n_1$ direction. We now define the slope $s_i=d_i/a_i$ to be the slope in the $n_1$ direction, and likewise define $s_0=\max_is_i$ which is the maximal slope in the $n_1$ direction.

We can simply modify the definitions for stability in the $n_1$ direction, by only considering the algebraic gradings in that direction. For example, we say that $A$ is homologically stable with respect to $U_1,\ldots, U_k$ \textit{in the $n_1$ direction} if for each $i$ we have that $H_i\left((A/(U_1,\ldots, U_k))_{n_1,*}\right)$ vanishes whenever $i<sn_1+t$ for some directional slope $s>s_0$ and $t$. Likewise, we can also define (strong) stability in an $(N,I)$ staircase in the $n_1$ direction, by replacing the grading $n$ with the grading $n_1$, and so on. Furthermore, we have that the analogue of Theorem \ref{thm: usual vs strong stability} is true.
\begin{theorem}\label{thm: usual vs strong stability direction}
Suppose that all $U_i\in H_0(A)$. Then $A$ stabilizes in an $(N,I)$ staircase with respect to $U_1,\ldots, U_k$ in the $n_1$ direction if and only if it stabilizes strongly in an $(N,I)$ staircase with respect to $U_1,\ldots, U_k$ in the $n_1$ direction.
\end{theorem}

We remark that this setup for stability in a chosen direction generalizes the original setup, as we can simply take $V=V_1$ so that the grading $n_1$ is simply the grading $n$, allowing us to recover the original definitions. We will assume this setup whenever stability in the $n_1$ direction is mentioned, and assume the original setup otherwise.
\subsection{Generalized theorem statements}\label{sec: subsec generalized theorem statements}
We now present the generalization of Theorem \ref{thm: main theorem} in the multiple element directional context as in Section \ref{sec: stability in chosen direction}. 
\begin{theorem}\label{thm: general main theorem}
Suppose that $A$ is homologically stable with respect to the $\E_2$-central elements $U_1,\ldots, U_k$ in the $n_1$ direction, and that it also stabilizes strongly in an $(N,I)$ staircase with respect to these elements in the $n_1$ direction. Also suppose that $\frac{I+1}{N+1}>s_0$, then it is homologically stable with slope $\frac{I+1}{N+1}$ in the $n_1$ direction. In particular, $$H_{i+k+d_1+\cdots +d_k}\left((A/(U_1,\ldots, U_k))_{n_1+a_1+\cdots+a_k,*}\right)=0$$
for all $n_1\ge N$ and
\begin{equation*}
i\le (I+1)q+ s_0r+\max\left(I-N+(1-s_0)r\ ,\ 0\right)-1. 
\end{equation*}
where $q=\left\lfloor \frac{n_1+1}{N+1}\right\rfloor$ and $r=\left\{ \frac{n_1+1}{N+1}\right\}$ are the quotient and remainder respectively.
\end{theorem}

We can recover Theorem \ref{thm: main theorem} by removing the directionality and taking a single element with zero homological degree. We can see this from Theorem \ref{thm: usual vs strong stability} which tells us that in this case strong stability is the same as usual stability, and the discussion below Definition \ref{def: staircase iterated version} which tells us the two notions of stability in an $(N,I)$ staircase are the same. We also note that this equation specializes to the equation in Theorem \ref{thm: main theorem} when $s_0=0$. 

Now, we give a more general version of the optimality theorem presented in Corollary \ref{cor: optimal bound} with multiple elements. However, we would need to specialize to the case where the homological degrees $d_i=0$ so that strong stability is usual stability, and also remove the directionality context, so that the initial conditions required are finite and computable. In this case, the generalized optimality theorem will then follow directly from Theorem \ref{thm: general main theorem} as usual.
\begin{corollary}\label{cor: general optimal}
Suppose that $A$ is homologically stable of slope $s$ with respect to the $\E_2$-central elements $U_1,\ldots, U_k\in H_0(A)$. Then, for any $\epsilon>0$, there exist some $N>0$ and $I\ge 0$ such that $A$ stabilizes in an $(N,I)$ staircase with respect to $U_1,\ldots, U_k$ and $\frac{I+1}{N+1}>s-\epsilon$. In particular, computing $A$ up to sufficiently large $n$ and applying Theorem \ref{thm: general main theorem} will prove homological stability with slope $>s-\epsilon$.
\end{corollary}

\subsubsection{Stability from $R$ for general Hurwitz spaces}
Lastly, we comment on the generalization of Corollary \ref{cor: stab from R} where we deduced a stability slope from stability of $R=H_0(\Hur^c,\Q)$ in the case of $\Hur^c$ when $c$ is non-splitting. We cannot generalize this to a general braided vector space, as the argument in Section \ref{sec: intro stab from R} uses the fact that $\Hur^c_n$ stabilizes to a multiple of $\Conf_n$ as in Theorem \ref{thm: stable hom conf all directions}, which is not necessarily true for a general braided vector space. 

However, we can attempt to generalize this by removing the assumption that $c$ is non-splitting, and allowing for iterated quotients. We were able to do this in the case of $G=S_d$ and $c$ the conjugacy class of transpositions, where we concluded a nontrivial stability slope only from knowledge of $R$. We carry this out in Section \ref{sec: S_d slope}. However, it is not clear whether this method works in other cases, as we used several facts specific to this case. 

The argument in this case is much more complicated than the one for Corollary \ref{cor: stab from R}, and we give a very rough sketch of the argument. The main idea is to use the filtration in Equation \eqref{eqn: filtration A_c and MA_c} to go between the usual chains $A_c$ and maximal chains $MA_c$. This allows us to get a desired stability staircase for $A_c$ from stability of $MA_c$ which we deduce from $R$, together with stability of $MA_{c'}$ for proper subracks $c'\subsetneq c$, which we will deduce by induction on $d$. It also turns out that we need to define and use directional maximal Hurwitz spaces. A flowchart of the argument is presented in Figure \ref{fig: summary of S_d slope}.
\subsection{Lemmas on extra stability elements} \label{sec: lem extra elements} We record a few lemmas which allow us to add or remove extra stability elements to the iterated quotient, which will be helpful for the proofs later. These will basically follow from the definition we introduced previously together with the appropriate cofiber long exact sequence. 

Let us work in the directional setup as in Section \ref{sec: stability in chosen direction}, and suppose that we have the additional $\E_2$-central elements $V_1,\ldots, V_j$ (which may overlap with $U_1,\ldots, U_k$) with algebraic grading $(a_i',b_i')$ with $a_i'>0$, homological degree $d_i'$, slope $s_i'=d_i'/a_i'$, and maximal slope $s_0'=\max(s_i')$. 

Unlike the previous subsection, the lemmas here all apply to a general $A$-bimodule $M$, as the theorems here do not follow from the filtrations constructed in Section \ref{sec: filtrations} which only apply to $A$. Of course, we can recover the case for $A$ by taking $M=A$.

The first lemma tells us that adding more stability elements preserves homological stability, as long as their slopes are not too large.
\begin{lemma}\label{lem: homologically stable with extra elements}
Suppose that $M$ is homologically stable with respect to $U_1,\ldots, U_k$ in the $n_1$ direction with slope $s>s_0$. If we have $s>s_0'$, then $M$ is homologically stable with respect to $U_1,\ldots, U_k,V_1,\ldots, V_j$ in the $n_1$ direction with the same slope $s$. Furthermore, the $V_i$ may also be inserted arbitrarily while preserving the order of the $U_i$.
\end{lemma}
\begin{proof}
It suffices to prove the lemma for the simplest case of one element $V_1$ as we can obtain the statement with $V_1,\ldots, V_j$ by iteratively applying the one element version. This special case will then follow rather simply from the long exact sequence \eqref{eqn: iterated quotient LES} applied to $U_1,\ldots,  V_1,\ldots, U_k$ where $V_1$ is inserted arbitrarily in the sequence:
\begin{equation*}
\begin{split}  &\qquad \qquad \qquad \qquad \qquad  \cdots \xrightarrow{V_1} H_i((M/(U_1,\ldots U_k))_{n_1,n_2})\rightarrow H_i((M/(U_1,\ldots,V_1,\ldots,U_k))_{n_1,n_2}) \xrightarrow\delta
   \\ 
    \rightarrow & H_{i-d_1'-1}((M/(U_1,\ldots, U_k))_{n_1-a_1',n_2-b_1'}) \xrightarrow{V_1}  \cdots.
\end{split}
\end{equation*}
We have by definition that the first term $H_i((M/(U_1,\ldots,U_k))_{n_1,n_2})$ vanishes for all $i<sn_1+t$. By possibly changing $t$, we can also ensure that the last term vanishes for $i<sn_1+t$, which tells us that the middle term vanishes, proving the desired statement.
\end{proof}

The second lemma tells us that we can remove extra stability elements (as long as their slopes are not too large) while preserving the notion of strong stability in a staircase. This will be used in the proof of our main theorem in Section \ref{sec: pf general results}.
\begin{lemma}\label{lem: remove extra}
Suppose that $M$ is homologically stable with respect to $U_1,\ldots, U_k$ in the $n_1$ direction, and that $s_0\ge s_0'$. If $M$ stabilizes strongly in an $(N,I)$ staircase with respect to $U_1,\ldots, U_k,V_1,\ldots, V_j$ (or with $V_1,\ldots, V_j$ inserted arbitrarily in $U_1,\ldots, U_k$) in the $n_1$ direction, then $M$ stabilizes strongly in an $(N,I)$ staircase with respect to $U_1,\ldots, U_k$ in the $n_1$ direction.
\end{lemma}
\begin{proof}
By Lemma \ref{lem: homologically stable with extra elements}, our assumption tells us that $M$ is also homologically stable with respect to $U_1,\ldots, U_k,V_1,\ldots, V_j$; here we could also insert $V_1,\ldots, V_j$ arbitrarily or leave some of them out. Again, we only need to prove this where we only have one such $V_1$, as we can deduce the lemma by removing one $V_i$ at a time. The key is a slightly different part of the same long exact sequence \eqref{eqn: iterated quotient LES}:
\begin{equation*}
\begin{split}   \cdots \rightarrow  H_{i+d_1'+1}((M/(U_1,\ldots,V_1,\ldots, U_k))_{n_1+a_1',n_2+b_1'}) \xrightarrow\delta &
   \\ 
    \rightarrow  H_{i}((M/(U_1,\ldots, U_k))_{n_1,n_2}) \xrightarrow{V_1} H_{i+d_1'}((M/(U_1,\ldots, U_k))_{n_1+a_1',n_2+b_1'})\rightarrow \cdots . \qquad\qquad\qquad\qquad&
\end{split}
\end{equation*}

For notational convenience, let us denote $S=N+\sum_ia_i$ and $T=I+k+\sum_i d_i$. Then, by unwinding the definition of strong stability, the condition tells us that the first term vanishes when $(n_1,i)$ satisfies both $i\le T+(n_1-S)$ and $i\le T+(n_1-S)s_0$, and we want to show that the second term vanishes for such $(n_1,i)$. Let us now fix such an $(n_1,i)$, and it suffices to show that the last term vanishes.

We observe that the last term is exactly the middle term with the algebraic and homological degrees shifted up by $(a_1',b_1')$ and $d_1'$. Motivated by this, let us shift the algebraic and homological degrees of the exact sequence up the same amount. It is clear that the first term still vanishes as $(n_1+a_1',i+d_1')$ still satisfies both inequalities required for strong stability as the right side always increases more than the left since $d_1'=s_1'a_1'\le a_1's_0'\le a_1'$ (indeed, this is one of the reasons we defined strong stability that way), and so vanishing of the middle term of this new equation (which is the last term of the old equation) follows from the vanishing of the last term of the new equation. 

Iterating this, we see that it suffices to show that
$$H_{i+ed_1'}((M/(U_1,\ldots, U_k))_{n_1+ea_1',n_2+eb_1'})=0$$
for some sufficiently large $e$. However, this follows from the fact that $M$ is homologically stable with respect to $U_1,\ldots, U_k$ (with some slope $s>s_0$). Indeed, as $sa_1'>a_1's_0\ge a_1's_0'\ge a_1's_1'= d_1'$, there will be some $e$ where $i+ed_1'<s(n_1+ea_1')+t$ as the right side of the inequality grows faster than the left.
\end{proof}

There are some limitations to the above lemma. Firstly, it works only for strong stability (and indeed, this is the reason why Theorem \ref{thm: general main theorem} is stated using strong stability). This is not really a problem when we are taking $M=A$ as here strong stability is equivalent to usual stability by Theorem \ref{thm: usual vs strong stability direction} in the case when the stability elements all belong to $H_0(A)$. However, if we want to apply it to a general $A$-bimodule $M$, such as $MA_c$, the notion of strong stability for $M$ is not as useful because Theorem \ref{thm: usual vs strong stability direction} no longer holds for $M$.

It turns out that one can fix this problem by considering sufficiently large powers of stability elements. We first give a helpful lemma which tells us taking powers of the stability elements still preserves homological stability.

\begin{lemma}\label{lem: homological stability for powers}
Suppose that $M$ is homologically stable with respect to $U_1,\ldots, U_k$ in the $n_1$ direction. Then, for any positive integers $e_1,\ldots, e_k$, we have that $M$ is also homologically stable with respect to $U_1^{e_1},\ldots, U_k^{e_k}$ in the $n_1$ direction.
\end{lemma}
\begin{proof}
We show that we can replace $U_1$ by $U_1^{e_1}$ for any $e_1\ge 1$, and the rest of the proof will follow by repeatedly replacing elements by larger powers. Again, from the long exact sequence \eqref{eqn: iterated quotient LES}, we have that $H_i((M/(U_1,\ldots, U_k))_{n_1,n_2})$ vanishes for $i\le I$ if and only if 
$$H_{i-d_1}((M/(U_2,\ldots, U_k))_{n_1-a_1,n_2-b_1})\xrightarrow{U_1}H_i((M/(U_2,\ldots, U_k))_{n_1,n_2})$$
is a surjection for $i=I$ and an isomorphism for $i< I$. By chaining $e_1$ of these maps together, we can conclude that $H_i((M/(U_1^{e_1},\ldots, U_k))_{n_1,n_2})$ vanishes for $i\le I$ if we have that $H_{i-rd_1}((M/(U_1,\ldots, U_k))_{n_1-ra_1,n_2-rb_1})$ vanishes for $i\le I$ and $0\le r<e_1$. We can deduce the statement easily from this.
\end{proof}

By considering sufficiently large powers of stability elements, we can add and remove these powers without affecting homological stability in a staircase, which we state in the following lemma. This will be very useful for our applications in Section \ref{sec: S_d slope}.
\begin{lemma}\label{lem: add extra elements}
Suppose that $M$ is homologically stable with respect to $U_1,\ldots, U_k$, and that $s_0\ge s_0'$. Then, there exists a constant $C$ depending only on $U_1,\ldots, U_k,V_1,\ldots, V_j$ and $M,N,I$ (not on $e_i,f_i$) such that the following are equivalent.
\begin{enumerate}[(\alph*)]
\item $M$ stabilizes (strongly) in an $(N,I)$ staircase with respect to $U_1^{e_1},\ldots, U_k^{e_k}$ in the $n_1$ direction.
\item For all exponents $f_1,\ldots, f_j>C$, $M$ stabilizes (strongly) in an $(N,I)$ staircase with respect to $U_1^{e_1},\ldots, U_k^{e_k},$ $V_1^{f_1},\ldots, V_j^{f_j}$ in the $n_1$ direction.
\end{enumerate}
\end{lemma}
\begin{proof}
We first discuss the case of adding or removing one element. Writing $f=f_1$ for simplicity, the usual long exact sequence \eqref{eqn: iterated quotient LES} gives  
\begin{equation*}
\small
\begin{split} &  \xrightarrow{V_1^{f}} H_{i+fd_1'+1}((M/(U_1^{e_1},\ldots, U_k^{e_k}))_{n_1+fa_1',n_2+fb_1'})\rightarrow H_{i+fd_1'+1}((M/(U_1^{e_1},\ldots,V_1^{f},\ldots,U_k^{e_k}))_{n_1+fa_1',n_2+fb_1'}) \rightarrow
   \\ &
    \rightarrow  H_{i}((M/(U_1^{e_1},\ldots, U_k^{e_k}))_{n_1,n_2}) \xrightarrow{V_1^f}  H_{i+fd_1'}((M/(U_1^{e_1},\ldots, U_k^{e_k}))_{n_1+fa_1',n_2+fb_1'})\rightarrow \cdots.
\end{split}
\end{equation*}

Let us again write $S=N+\sum_i a_ie_i$ and $T=I+k+\sum_i d_ie_i$, then condition $(a)$ tells us exactly that the third term vanishes for $(n_1,i)$ satisfying $i\le T+(n_1-S)$ and $n_1\le S$, and condition $(b)$ tells us exactly that the second term vanishes for such $(n_1,i)$ (for all $f>C$). To show these two are equivalent, it suffices to show that the first and last terms vanish for all $f>C$. However, this is clear, because $M$ is homologically stable with respect to $U_1,\ldots,U_k$ and hence $U_1^{e_1},\ldots, U_k^{e_k}$. We also need to check that the bound $C$ does not depend on $e_1,\ldots, e_k$, and indeed, we check that the boundary $i<sn_1+t$ of the stability slope in Definition \ref{def: general hom stability} for $U_1^{e_1},\ldots, U_k^{e_k}$ compared to $U_1,\ldots, U_k$ is simply shifted by the sum of the homological and algebraic degrees of the stability elements by the proof of Lemma \ref{lem: homological stability for powers}, and this is already accounted for in the definition of $S$ and $T$. Another way of putting this is that the stability slope does not change in the shifted coordinates of Equation \eqref{eqn: vanishing homology iterated quotient}. Hence, the constant $C$, which is dependent on how far $(n_1,i)$ is away from the stability slope of $U_1^{e_1},\ldots, U_k^{e_k}$ as well as the slope of $V_1$, does not change as we change $e_i$.

For multiple elements, we iteratively apply the same argument. Note that during the intermediate steps, by Lemma \ref{lem: homologically stable with extra elements}, $M$ is still homologically stable when we add the elements $V_1^{f_1},\ldots, V_i^{f_i}$ for $1\le i\le j$. Furthermore, from the proof of Lemma \ref{lem: homologically stable with extra elements}, we can check that the stability line shifts by the homological and algebraic gradings of the powers of the new elements. Hence, by the same reasoning, the constant $C$ will not change based on the choice of $f_i$.
\end{proof}
\section{Iterated quotients as a homotopy-coherent Koszul complex}\label{sec: iterated quotients as koszul}

Having stated the general versions of our theorems, we now turn to the proof of our main theorem, which occupies this section and the next section. The goal of this section is to write iterated quotients that appear in the theorem statements in the form of an explicit Koszul-like cobar-complex. Then, in Section \ref{sec: filtrations}, we will define filtrations on this Koszul-like cobar-complex, which will help us understand quotients in larger algebraic degree from smaller algebraic degree, leading to a proof of our main theorem. In particular, we note that having a concrete chain-complex model for these iterated quotients is essential for our argument.

\subsection{Motivation and statements} We first discuss the usual notion of Koszul complexes as a motivating example, following \cite[Section 15.29]{Stacksproject} and \cite[Section 4.5]{Weibel_2011}. Let $A$ be a \textit{commutative} ring with elements $f_1,\ldots, f_k$, and we consider the derived quotient $A/_{\mathbf L}(f_1,\ldots, f_k)\coloneqq (\cdots(A/_{\mathbf L}f_1)\cdots/_{\mathbf L}f_k)$. The Koszul complex $K_A(f_1,\ldots, f_k)$ gives a chain-level model of this, where $K_A(f_i)=\{A\xrightarrow {f_i}A\}$ is a two-term chain-complex in degrees $1$ and $0$, and $K_A(f_1,\ldots, f_k)$ is the total complex associated to the $k$-dimensional cube $K_A(f_1)\otimes \cdots\otimes K_A(f_k)$. In other words, one can write
$$K_A(f_1,\ldots, f_k)\coloneqq A \otimes \bigwedge(e_1,\ldots, e_k)$$
where each $e_i$ is in homological degree $1$ and the differential sends
$$d(e_{i_1}\wedge\cdots \wedge e_{i_m})=\sum_{j=1}^m (-1)^{j-1} f_{i_j}e_{i_1}\wedge\cdots\wedge \widehat{e_{i_j}}\wedge\cdots \wedge e_{i_m}.$$

Now, we return to our situation where $A$ is an $\E_1$-algebra in $\mathcal D^-(k)$ and $U_1,\ldots, U_k$ are $\E_2$-central elements with homological degrees $d_1,\ldots, d_k$. Since our discussion will work for $\E_1$-algebras in general,  we choose not to specialize to the case of $A=A_V$ as in Definition \ref{def: A_V} for most of this section, although we will look at some examples in this special case later in Section \ref{sec: koszul example}.

We want a chain-complex model for our iterated quotient $A/(U_1,\ldots, U_k)$ as defined in Definition \ref{def: iterated quotients}. However, in our case, $A$ is no longer commutative, so the usual Koszul complex construction does not work. Nevertheless, it turns out that we can get such a model by modifying the Koszul complex to include higher homotopies. By abuse of notation, let us also denote the chain-complex representing $A$ by $A$.
\begin{definition}\label{def: right koszul}
We define the \textit{right homotopy-coherent Koszul complex} of the quotient $A/(U_1,\ldots, U_k)$ to be
$$\mathcal K^R_{A}(U_1,\ldots, U_k)\coloneqq A \otimes \bigwedge (e_1,\ldots, e_k),$$
where here $e_i$ is in homological degree $d_i+1$. We define $e_I=e_{i_1}\wedge\cdots\wedge e_{i_m}$ for $I=\{i_1<\cdots <i_m\}$, and we can also just view $\bigwedge (e_1,\ldots, e_k)=\bigoplus_{I\subseteq [k]}e_I$ as a direct sum of these formal variables. The differential is of the form
$$d=d_A+\sum_{\emptyset\neq I\subseteq [k]}d_I,$$
where the components of the differentials are defined to be
\begin{equation*}
\begin{split}
d_A(a\otimes e_J) &= da\otimes e_J,\\
d_I(a\otimes e_J) &= \begin{cases}
(-1)^{\kappa(a,I,J)}\cdot  a \cdot h_I \otimes e_{J\setminus I} & \text{if }I\subseteq J\\
0 & \text{otherwise}
\end{cases}
\end{split}
\end{equation*}
where $h_{\{i\}}=U_i$ for each $i\in [k]$, and $h_I\in A$ for $|I|\ge 2$ records higher homotopies which are given by right multiplication by the elements $h_I$ which are in homological degree $\sum_{i\in I}(d_i+1)-1$. Here, $(-1)^{\kappa(a,I,J)}$ is the usual Koszul sign twist given by 
$$\kappa(a,I,J)=|a|+\sum_{\substack{x\in I, y\in J\setminus I\\ y<x}}(d_x+1)(d_y+1).$$
\end{definition}

Furthermore, in the case that $A$ has a grading, we will let $e_i$ have the same grading as $U_i$. The first goal of this section is to prove that this is indeed a model for the iterated quotient.
\begin{proposition}\label{prop: right koszul}
$\mathcal K^R_A(U_1,\ldots,U_k)$ is a chain-complex representing $A/(U_1,\ldots,U_k)$. Furthermore, the elements $h_I$ arise from the higher homotopies of the image in $\LMod_A$ of the cube constructed below Definition \ref{def: k-cube}, whose total cofiber, as defined in Definition \ref{def: total cofiber}, is $A/(U_1,\ldots,U_k)$.
\end{proposition}

Furthermore, we define the \textit{left homotopy-coherent Koszul complex} as
$$\mathcal K^L_{A}(U_1,\ldots, U_k)\coloneqq \bigwedge (e_1,\ldots, e_k)\otimes A.$$
Compared with Definition \ref{def: right koszul}, the exterior generators are placed on the left and the differentials act by left multiplication. Note that the Koszul sign twists are different, and for the precise differential formulas, see the special case of Definition \ref{def: two-sided koszul} with no right-hand generators. The higher homotopies come from the same ordered cube constructed below Definition \ref{def: k-cube}, now viewed in $\RMod_A$, and may be represented by different elements. Then, the analogue of Proposition \ref{prop: right koszul} holds, and the resulting complex represents $A/(U_1,\ldots,U_k)$ as a right $A$-module.

Pushing this idea further, we can define a two-sided version as follows. Let $U_1,\ldots, U_j,V_1,\ldots, V_k$ be $\E_2$-central elements with homological degree $c_1,\ldots, c_j, d_1,\ldots, d_k$.
\begin{definition}\label{def: two-sided koszul}
We define the \textit{two-sided homotopy-coherent Koszul complex} of the quotient $A/(U_1,\ldots, U_j,V_1,\ldots, V_k)$ to be
$$\mathcal K_A(U_1,\ldots, U_j\mid V_1,\ldots, V_k)\coloneqq \bigwedge(e_1,\ldots, e_j)\otimes A\otimes \bigwedge(f_1,\ldots, f_k)$$
where $e_i$ is in degree $c_i+1$, and $f_i$ in degree $d_i+1$. The differentials are defined to be
$$d=d_A+\sum_{\emptyset \neq I \subseteq [j]}d^L_I+\sum_{\emptyset \neq K \subseteq [k]}d^R_K$$
where
\begin{equation*}
\begin{split}
d_A(e_J \otimes a \otimes f_L)&=(-1)^{\kappa_1(J)}e_J\otimes da\otimes f_L\\
d^L_I(e_J\otimes a\otimes f_L)&=\begin{cases}
(-1)^{\kappa_2(I,J)} e_{J\setminus I}\otimes g_I\cdot a \otimes f_L & \text{if } I\subseteq J\\
0 & \text{otherwise}
\end{cases}\\
d^R_K(e_J\otimes a\otimes f_L)&=\begin{cases}
(-1)^{\kappa_3(J,a,K,L)} e_J\otimes a \cdot h_K\otimes f_{L\setminus K} &\text{if } K\subseteq L\\
0 &\text{otherwise.}
\end{cases}
\end{split}
\end{equation*}
where $g_{\{i\}}=U_i$ and $h_{\{i\}}=V_i$, where $\kappa_i$ are the usual Koszul sign twists.
\end{definition}
\begin{remark}
For completeness, we record the Koszul sign twists here, which can be derived from the usual sign conventions.
\begin{equation*}
\begin{split}
\kappa_1(J)
&=
\sum_{x\in J}(c_x+1),\\
\kappa_2(I,J)
&=
\sum_{y\in J\setminus I}(c_y+1)
+
\sum_{\substack{x\in I,\ y\in J\setminus I\\ x<y}}
(c_x+1)(c_y+1),\\
\kappa_3(J,a,K,L)
&=
\sum_{x\in J}(c_x+1)
+
|a|
+
\sum_{\substack{x\in K,\ y\in L\setminus K\\ y<x}}
(d_x+1)(d_y+1).
\end{split}
\end{equation*}
It is not difficult to see that these signs are essentially forced (up to sign flips in $g_I$ and $f_L$) so that we have $d^2=0$. As these signs are not important for our use case and are somewhat tedious to keep track of, we will omit any further discussions on sign twists, including in the proof of the statement.
\end{remark}

Like just now, we will also prove that this complex is indeed a model for the iterated quotient.
\begin{proposition}\label{prop: two sided koszul}
$\mathcal K_A(U_1,\ldots, U_j\mid V_1,\ldots, V_k)$ is a chain-complex for $A/(U_1,\ldots, U_j,V_1,\ldots, V_k)$. Furthermore, the elements $g_I$ and $h_K$ are the same as those appearing in $\mathcal K^L_A(U_1,\ldots, U_j)$ and $\mathcal K^R_A(V_1,\ldots, V_k)$ respectively.
\end{proposition}
\subsection{Homotopy-coherent hypercubes} We will first need some theory on homotopy-coherent hypercubes, mainly following \cite{Goodwillie_1991,Dyckerhoff_Jasso_Walde_2019,Varolgunes_2021,ye2026infinitycategoricaltqftinstantons}. There have also been similar constructions of hypercubes of chain complexes in \cite{Manolescu_Ozsváth_2025,hendricks2025surgeryexacttrianglesinvolutive,guth2026satellitestelescopesconcordanceformula}. Most of these works, however, formulate homotopy-coherent hypercubes directly in the setting of chain complexes rather than in the language of $\infty$-categories. Hence, they are not directly applicable in our paper as we need to begin with the $\infty$-categorical notion of a hypercube and then identify its corresponding concrete model in chain complexes.

We first define $k$-cubes in $\infty$-categories, which we also call homotopy-coherent hypercubes.
\begin{definition}\label{def: k-cube}
Let $\mathcal C$ be an $\infty$-category. A \textit{$k$-cube} in $\mathcal C$ is a simplicial map $X\colon (\Delta^1)^k\rightarrow \mathcal C$, where $\Delta^1$ is the standard $1$-simplex. One can also view $(\Delta^1)^k$ as the nerve $\mathcal N(\mathcal P([k]))$ of the poset of subsets of $[k]$. It will be convenient to use the more general notion of an \textit{$S$-cube} for a finite set $S$, which is a simplicial map $X\colon (\Delta^1)^{S}\cong\mathcal N(\mathcal P(S))\rightarrow \mathcal C$.
\end{definition}

For disjoint finite sets $S,T$, since $$\Fun((\Delta^1)^{S\sqcup T},\mathcal C)=\Fun((\Delta^1)^S,\Fun((\Delta^1)^T,\mathcal C)),$$ 
one can view an $S\sqcup T$-cube as an $S$-cube in the functor category of $T$-cubes.

In our context, the iterated quotient $A/(U_1,\ldots, U_k)$ we consider has an associated homotopy-coherent hypercube as we explain next. Recall that since $U_i$ is $\E_2$-central, it defines a $\BiMod_A$-map $U_i\colon A[d_i]\rightarrow A$, which is a simplicial map $\Delta^1\rightarrow \mathcal \BiMod_A$. Since the relative tensor product $\otimes_A$ is an $\infty$-functor, tensoring the maps for $U_1,\ldots, U_k$ gives a $k$-cube $(\Delta^1)^k\rightarrow \BiMod_A$ with vertices $A[\sum_{i\not\in S}d_i]$ indexed by subsets $S\subseteq [k]$. We illustrate the $k=3$ case in the following diagram.
\[\begin{tikzcd}[sep=small]
	& {A[d_1+d_3]} && {A[d_1]} \\
	{A[d_1+d_2+d_3]} && {A[d_1+d_2]} \\
	& {A[d_3]} && A \\
	{A[d_2+d_3]} && {A[d_2]}
	\arrow[from=1-2, to=1-4]
	\arrow[from=1-2, to=3-2]
	\arrow[from=1-4, to=3-4]
	\arrow[from=2-1, to=1-2]
	\arrow[from=2-1, to=2-3]
	\arrow[from=2-1, to=4-1]
	\arrow[from=2-3, to=1-4]
	\arrow[from=2-3, to=4-3]
	\arrow[from=3-2, to=3-4]
	\arrow[from=4-1, to=3-2]
	\arrow[from=4-1, to=4-3]
	\arrow[from=4-3, to=3-4]
\end{tikzcd}\]

Given an $S$-cube, we can associate the total cofiber as introduced in \cite[Section 1]{Goodwillie_1991}. One of the many equivalent definitions of this is given by iteratively taking the cofiber which we describe here. 
\begin{definition}\label{def: total cofiber}
Let $\mathcal C$ be a stable $\infty$-category. We define the \textit{total cofiber} $\totcofib\colon \Fun((\Delta^1)^S,\mathcal C)\rightarrow \mathcal C$ inductively on the size of $S$ as follows. Let $X\colon (\Delta^1)^S\rightarrow \mathcal C$ be an $S$-cube in $\mathcal C$. If $S=\emptyset$, we let $\totcofib(X)\coloneqq X$. Otherwise, choose any $i\in S$, and we can view $X$ as a morphism ($1$-cube) $X^i_0\rightarrow X^i_1$ of $S\setminus \{i\}$-cubes. Then, we define 
$$\totcofib(X)\coloneqq \cofib(\totcofib(X^i_0)\rightarrow \totcofib(X^i_1))
$$
where the total cofibers of the smaller homotopy-coherent cubes $X^i_0,X^i_1$ are defined by induction. 
\end{definition}

It is known that the order in which we take cofibers does not matter so this is well defined. Furthermore, by defining cofibers functorially, one can also write $\totcofib\simeq \cofib^{|S|}$. We refer the reader to \cite[Appendix A2]{Dyckerhoff_Jasso_Walde_2019} for these two facts.

Next, we see that the iterated quotient $A/(U_1,\ldots, U_k)$ is simply the total cofiber of the associated $k$-cube described above. Indeed, since $\otimes_A$ preserves cofibers, we can write
\begin{equation*}
\begin{split}
A/(U_1,\ldots, U_k) &=\cofib( A[d_1]\xrightarrow{U_1} A)\otimes_A A/(U_2,\ldots, U_k) \\
&= \cofib(A/(U_2,\ldots, U_k)[d_1]\xrightarrow{U_1}A/(U_2,\ldots, U_k))
\end{split}   
\end{equation*}
which is simply the recursive definition of the total cofiber explained above.

We now specialize to the derived $\infty$-category, where \cite[Section 7.1]{ye2026infinitycategoricaltqftinstantons}, following \cite[Section 2.1]{Varolgunes_2021}, gives a chain-complex for the total cofiber of a $k$-cube in this context. This is stated in \cite[Definition 2.1, Construction 7.6]{ye2026infinitycategoricaltqftinstantons}; here we have made some minor changes to the statement.
\begin{lemma} \label{lem: fan ye} Let $X\colon (\Delta^1)^k\rightarrow \mathcal D^-(k)$ be a $k$-cube. Then, its total cofiber $\totcofib(X)$ is given by the chain complex
$$C=\bigoplus_{S\subseteq [k]}X(S)[k-|S|]$$
with differential $$d=\sum_{S\subseteq S'\subseteq [k]}d_{S,S'}$$
where the component $d_{S,S'}$ maps $X(S)[k-|S|]\rightarrow X(S')[k-|S'|]$.

Furthermore, for $S=S'$ we have $d_{S,S}=(-1)^{k-|S|}d_{X(S)}$, and otherwise for $S\subsetneq S'$ with $|S'\setminus S|=j+1$, $d_{S,S'}$ is the $j$-homotopy of the subcube with initial vertex $S$ and final vertex $S'$. 
\end{lemma}

This was stated without proof in \cite{ye2026infinitycategoricaltqftinstantons}, and we did not find a full proof in \cite{Varolgunes_2021} either as it was not done in the setting of $\infty$-categories. Furthermore, the setup there is also slightly different. While the result is natural and one would certainly expect this to be true, we still opt to give a self-contained proof of this by induction on $k$.
\begin{proof}[Proof of Lemma \ref{lem: fan ye}] 
By \cite[Remark 1.3.2.17]{HA}, the cofiber $\cofib(f)$ of a map $f\colon A\rightarrow B$ in the derived $\infty$-category is given by the mapping cone $C(f)\coloneqq A[1]\oplus B$ with differential $$d_{C(f)}=\begin{pmatrix}
-d_A& 0\\
f & d_B
\end{pmatrix}$$
in this basis. Since the total cofiber is just the cofiber when $k=1$, the lemma is proven in this case.

The cone satisfies the $\infty$-categorical universal property of the cofiber as the pushout $B\coprod_A 0$ given by
$$\Map(C(f),Z)\simeq\Map(B,Z)\times_{\Map(A,Z)}\Map(0,Z)$$
for every $Z$, and here we also have $\Map(0,Z)\simeq *$. On the level of points, this gives
$$
\Hom(C(f),Z)\simeq\Hom(B,Z)\times_{\Hom(A,Z)}0,$$
in other words, the data of a map $C(f)\rightarrow Z$ is equivalent to the data of a map $g\colon B\rightarrow Z$ together with a homotopy $h$ from $g\circ f$ to $0$, as shown in the following diagram.
\[\begin{tikzcd}
	A && Z \\
	\\
	B && {C(f)}
	\arrow[""{name=0, anchor=center, inner sep=0}, "0", from=1-1, to=1-3]
	\arrow["f"', from=1-1, to=3-1]
	\arrow[""{name=1, anchor=center, inner sep=0}, "g"', from=3-1, to=1-3]
	\arrow[from=3-1, to=3-3]
	\arrow[dashed, from=3-3, to=1-3]
	\arrow["h", shift left=4, between={0}{0.9}, Rightarrow, nfold, from=1, to=0]
\end{tikzcd}\]

This is not difficult to check, as a map $C(f)\rightarrow Z$ can be written in terms of its components $g\colon B\rightarrow Z$ and $h\colon A[1]\rightarrow Z$, and since $C(f)\rightarrow Z$ is a chain map we have
$$-h\circ d_A + g \circ f=d_Z\circ h.$$ Then, the map $h\colon A[1]\rightarrow Z$ is a homotopy between $g\circ f$ and $0$ as
$$dh=d_Z\circ h+h\circ d_A=g\circ f - 0$$
using the previous equation. The converse is also true in the same way. The equivalence we described also extends to the mapping complexes, where a degree $n$ map $\overline{\Hom}_n(C(f),Z)$ can be decomposed as a degree $n$ map $\overline{\Hom}(B,Z)_n$ and a degree $n+1$ map $\overline{\Hom}(A,Z)_{n+1}$.

Now, we move on to the case of $k=2$, where we have a square $X\colon (\Delta^1)^2\rightarrow \mathcal D^-(k)$ 
\[\begin{tikzcd}
	A & C \\
	B & D
	\arrow["y", from=1-1, to=1-2]
	\arrow["f"', from=1-1, to=2-1]
	\arrow["g", from=1-2, to=2-2]
	\arrow["h", between={0.1}{0.9}, Rightarrow, nfold, from=2-1, to=1-2]
	\arrow["z"', from=2-1, to=2-2]
\end{tikzcd}\]
where $h\colon A[1]\rightarrow D$ is the homotopy of the square between $z\circ f$ and $g\circ y$. The total cofiber is then
$$\totcofib(X)= \cofib(C(f)\rightarrow C(g)),$$
where we recall that $C(f)=A[1]\oplus B$ and $C(g)=C[1]\oplus D$ with the differential as written above. Next, we deduce the map $C(f)\rightarrow C(g)$ from the universal property of cones discussed earlier and the following diagram.
\[\begin{tikzcd}
	A & C & \\
	B & D \\
	& {C(f)} & {C(g)}
	\arrow["y", from=1-1, to=1-2]
	\arrow["f"', from=1-1, to=2-1]
	\arrow["g", from=1-2, to=2-2]
	\arrow[""{name=0, anchor=center, inner sep=0}, "0", from=1-2, to=3-3]
	\arrow["h", between={0.1}{0.9}, Rightarrow, nfold, from=2-1, to=1-2]
	\arrow["z"', from=2-1, to=2-2]
	\arrow[from=2-1, to=3-2]
	\arrow["\iota"', from=2-2, to=3-3]
	\arrow[dashed, from=3-2, to=3-3]
	\arrow["{h'}", shift right=2, between={0.2}{1}, Rightarrow, nfold, from=2-2, to=0]
\end{tikzcd}\]
We break the map into the components $e_1\colon B\rightarrow C(g)$ and $e_2\colon A[1] \rightarrow C(g)$. By the description above, we see that $e_1=\iota \circ z$ where $\iota\colon D\rightarrow C(g)=C[1]\oplus D$ is given by the identity map on $D$. On the other hand, $e_2$ is the homotopy between $\iota\circ z\circ f\colon A\rightarrow C(g)$ and the zero map $0\circ y$, and from the diagram we see that this is $\iota \circ h + h'\circ y$. Again, by the universal property of $C(g)$, we have that $h'\colon C[1] \rightarrow C(g)=C[1]\oplus D$ is given by the identity map. Hence, we see that the map between two cones in this setting has three components $h,y,z$. Looking at degree $n$ maps $\overline{\Hom}_n(C(f),C(g))$ in an analogous situation will also give us the analogous three components.

Combining all of these, we see that $\cofib(C(f)\rightarrow C(g))=A[2]\oplus B[1]\oplus C[1]\oplus D$ with differential
$$
\begin{pmatrix}
d_A & 0& 0 & 0\\
-f &-d_B & 0 & 0\\
y & 0 & -d_C & 0\\
h & z & g & d_D
\end{pmatrix}
$$
in this basis. In particular, this proves the lemma for $k=2$.

We will complete the proof by induction. Suppose for some $k\ge 2$ that the lemma is true for $k-1$ and we want to prove the lemma for $k$. Let $X\colon (\Delta^1)^{[k]}\rightarrow \mathcal D^-(k)$ be a $k$-cube. To get the total cofiber of $X$ we can first apply cofibers in the last direction to get a $(k-1)$-cube of cofibers $\cofib(X)\colon (\Delta^1)^{[k-1]}\rightarrow\mathcal D^-(k)$ and then apply the total cofiber, i.e. we have $\totcofib(X)=\totcofib(\cofib(X))$. By the induction hypothesis, this is given by $$C=\bigoplus_{S\subseteq [k-1]}\cofib(X)(S)[k-1-|S|],\qquad d=\sum_{S\subseteq S'\subseteq [k-1]}d_{S,S'}.$$
Since we know that $\cofib(X)(S)=X(S)[1]\oplus X(S\sqcup\{k\})$ we get the desired form of $C$.

It now remains to check that the differential $d$ is as claimed. For $S=S'$, the component $d_{S,S}=d_{C(X(S)\rightarrow X(S\sqcup \{k\}))}$ which by the definition of the mapping cone has further sub-components $d_{S,S},d_{S,S\sqcup \{k\}},d_{S\sqcup \{k\},S\sqcup \{k\}}$ (up to sign) as desired. For $S\subsetneq S'$ with $|S'\setminus S|=j+1$, the component $$d_{S,S'}\colon \cofib(X)(S)[k-1-|S|]\rightarrow \cofib(X)(S')[k-1-|S'|]$$ is the $j$-homotopy of the subcube with initial vertex $S$ and final vertex $S'$. Now, from the proof of the $k=2$ case applied analogously to the case of degree $j$ maps $\overline{\Hom}_j(C(f),C(g))$ between cones, we see that $d_{S,S'}$ has the desired three components $d_{S,S'},d_{S\sqcup\{k\},S'\sqcup\{k\}},d_{S,S'\sqcup\{k\}}$ given by homotopies between the initial and end vertices, completing our proof.
\end{proof}
\subsection{Proof of chain models}
We now prove Propositions \ref{prop: right koszul}, \ref{prop: two sided koszul} by applying Lemma \ref{lem: fan ye}. For the first proposition, the key is to realize that homotopies in $\LMod_A$ (resp. $\RMod_A$) are given by right (resp. left) multiplication. 
\begin{proof}[Proof of Proposition \ref{prop: right koszul}] Consider the $k$-cube $X:(\Delta^1)^k\rightarrow \BiMod_A$ constructed below Definition \ref{def: k-cube}. By applying the forgetful functor $\BiMod_A\rightarrow \LMod_A$, we get a $k$-cube in $\LMod_A$. By applying Lemma \ref{lem: fan ye} to this $k$-cube and using the fact that the iterated quotient is the total cofiber, it remains to check that a $j$-homotopy from $A$ to $A$ in $\LMod_A$ in Lemma \ref{lem: fan ye} is given by right multiplication by a degree $j$ element. Indeed, by free-forgetful adjunction we have $\Map_{\LMod_A}(A,A)\simeq \Map_{D^-(k)}(k,A)\simeq A$, so the $j$-simplices on the left are identified with $j$-simplices of $A$. Hence, $j$-homotopies in $\LMod_A$ are given by right multiplication by a degree $j$ element $h\in A_j$, and one can check that under the adjunction it corresponds to right multiplication by $h$.
\end{proof}

Intuitively, the second proposition should be true by choosing to view the homotopies coming from $V_1,\ldots, V_k$ as right multiplication, while viewing the homotopies coming from $U_1,\ldots, U_j$ as left multiplication, and there should be no cross terms as left multiplication commutes strictly with right multiplication. However, it seems difficult to make this argument rigorous, and we will instead prove it by considering the tensor product of the left and right models $\mathcal K^L_A(U_1,\ldots, U_j)\otimes_A \mathcal K^R_A(V_1,\ldots, V_k)$ as dg modules.
\begin{proof}[Proof of Proposition \ref{prop: two sided koszul}]
Recall from the definition of iterated quotients that $$A/(U_1,\ldots, U_j,V_1,\ldots, V_k)=A/(U_1,\ldots, U_j)\otimes_A A/(V_1,\ldots, V_k).$$By Proposition \ref{prop: right koszul} and the left multiplication analogue described after its statement, we know that $\mathcal K^L_A(U_1,\ldots, U_j)$ is a chain model for $A/(U_1,\ldots, U_j)$, and likewise $\mathcal K^R_A(V_1,\ldots, V_k)$ for $A/(V_1,\ldots, V_k)$. By \cite[Proposition 2.19]{Christ_2022}, we can view $\mathcal K^L_A(U_1,\ldots, U_j)$ as a dg right $A$-module, and $\mathcal K^R_A(V_1,\ldots, V_k)$ as a dg left $A$-module, and in this situation \cite[Proposition 2.20]{Christ_2022} then tells us that Lurie's relative tensor product $\otimes_A$ agrees with the derived tensor product $\otimes^{\mathbf L}_A$ of dg-modules. 

We will show using a notion of flatness that the derived tensor product $\mathcal K^L_A(U_1,\ldots, U_j)\otimes^{\mathbf L}_A \mathcal K^R_A(V_1,\ldots, V_k)$ is given by the usual chain complex tensor product. Indeed, one can check that $\mathcal K^L_A(U_1,\ldots, U_j)$ is semifree in the sense of \cite[Definition 11.4.3]{Yekutieli_2020} by considering the filtration by exterior degree $F_p\mathcal K^L_A(U_1,\ldots, U_j)= \bigoplus_{|I|\le p}e_IA$. Then, by \cite[Proposition 10.3.4, Theorem 11.4.14]{Yekutieli_2020}, it is $K$-flat, so \cite[Theorem 12.3.1]{Yekutieli_2020} tells us that the derived tensor product is the usual tensor product. This gives the desired form of $\mathcal K_A(U_1,\ldots, U_j|V_1,\ldots, V_k)$ as stated in the proposition.
\end{proof}
\subsection{Examples}\label{sec: koszul example} We now give some examples of these homotopy-coherent Koszul complexes in the setting of Hurwitz spaces to better illustrate these constructions. We will consider a simple example where $c=c^{\text{trn}}_{S_3}$ is the conjugacy class of transpositions in $S_3$, and let us denote the elements of $c$ by $x,y,z$. As in Section \ref{sec: preliminaries}, we construct the corresponding tensor algebra $T=T((kc)_{-1})=\bigoplus_{n=0}^\infty T_n$ where each $T_n\cong (kc)^{\otimes n}$ as vector spaces. Then, the simplicial chains of $\Hur^c$ are given by
$$A_c=B^\dagger(T),$$
and we will take this cobar-complex model from now on. The full cobar-complex has infinitely many terms, so in our examples below we will look at a single algebraic degree. For example, $(A_c)_3$ can be written explicitly as 
$$T_3\rightarrow (T_1\otimes T_2)\oplus (T_2\otimes T_1) \rightarrow T_1\otimes T_1\otimes T_1.$$

\begin{figure}[!htbp]
    \centering
\[\begin{tikzcd}
	{\text{homological degree}} & 3 & 2 & 1 & 0 \\
	{B^\dagger(T)_2\colon} && {T_2} & {T_1^{\otimes 2}} \\
	{B^\dagger(T)_4\colon} & {T_4} & \begin{array}{c} \substack{ (T_3\otimes T_1)\\\oplus (T_2\otimes T_2)\\ \oplus (T_1\otimes T_3)} \end{array} & \begin{array}{c} \substack{(T_2\otimes T_1\otimes T_1)\\\oplus (T_1\otimes T_2\otimes T_1)\\\oplus (T_1\otimes T_1\otimes T_2)} \end{array} & {T_1^{\otimes 4}}
	\arrow[from=2-3, to=2-4]
	\arrow["\bigoplus"{description}, draw=none, from=2-3, to=3-3]
	\arrow["{\cdot x|x}"{description}, from=2-3, to=3-4]
	\arrow["\bigoplus"{description}, draw=none, from=2-4, to=3-4]
	\arrow["{\cdot x|x}"{description}, from=2-4, to=3-5]
	\arrow[from=3-2, to=3-3]
	\arrow[from=3-3, to=3-4]
	\arrow[from=3-4, to=3-5]
\end{tikzcd}\]
\caption{Chain complex for $(A_c/U)_4\simeq \mathcal K^R_{A_c}(U)_4$.}
    \label{fig: K^R_4}
\end{figure}
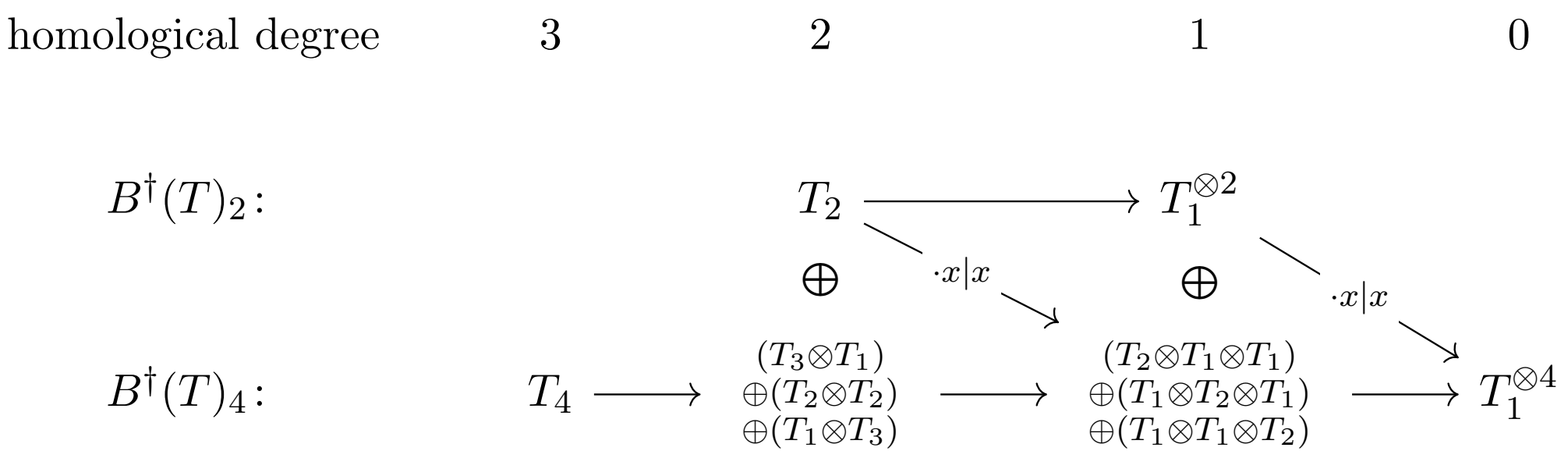

We first look at the case of a single quotient. Take $U=x| x\in H_0(A_c)$. Note that this is not a stabilizing map for $\Hur^c_n$, as it does not stabilize the components with only $y$s or $z$s; one can replace this with $U=x|x+y|y+z|z$ to fix this, but we take $U=x|x$ for simplicity. We showed that $\mathcal K^R_{A_c}(U)$ is a chain model for $A_c/U$. In Definition \ref{def: right koszul}, $e_1$ is defined to have the same algebraic grading as $U$ which is $n=2$. Say, for example, that we take the $n=4$ part. Then $\mathcal K^R_{A_c}(U)_4$ is the total complex $\Tot(B^\dagger(T)_2\xrightarrow{x\mid x}B^\dagger(T)_4)$ which we illustrate in Figure \ref{fig: K^R_4} ignoring the Koszul signs needed.

The left homotopy-coherent Koszul complex $\mathcal K^L_{A_c}(U)$ has an analogous description, which we omit. It is cumbersome to construct an explicit homotopy between multiplying by $U$ on the left and on the right, although such a homotopy must exist from geometry and the proof of $\E_2$-centrality.

\begin{figure}[!htbp]
    \centering

\[\begin{tikzcd}[row sep=scriptsize]
	\begin{array}{c} \substack{\text{homological}\\\text{degree}} \end{array} & 4 & 3 & 2 & 1 & 0 \\
	{B^\dagger(T)_1\colon} &&& {T_1} \\
	{B^\dagger(T)_3\colon} && {T_3} & \begin{array}{c} \substack{(T_2\otimes T_1)\\\oplus (T_1\otimes T_2)} \end{array} & {T_1^{\otimes 3}} \\
	{B^\dagger(T)_3\colon} && {T_3} & \begin{array}{c} \substack{(T_2\otimes T_1)\\\oplus (T_1\otimes T_2)} \end{array} & {T_1^{\otimes 3}} \\
	{B^\dagger(T)_5\colon} & {T_5} & \begin{array}{c} \substack{(T_4\otimes T_1)\\\oplus (T_3\otimes T_2)\\\oplus (T_2\otimes T_3)\\\oplus (T_1\otimes T_4)} \end{array} & \begin{array}{c} \substack{(T_3\otimes T_1\otimes T_1)\\\oplus (T_2\otimes T_2\otimes T_1)\\\oplus(T_2\otimes T_1\otimes T_2)\\\oplus\ \cdots} \end{array} & \begin{array}{c} \substack{(T_2\otimes T_1\otimes T_1\otimes T_1)\\\oplus (T_1\otimes T_2\otimes T_1\otimes T_1)\\\oplus (T_1\otimes T_1\otimes T_2\otimes T_1)\\\oplus (T_1\otimes T_1\otimes T_1\otimes T_2)} \end{array} & {T_1^{\otimes 5}}
	\arrow["\bigoplus"{description}, draw=none, from=2-4, to=3-4]
	\arrow["{\cdot  x|x}"{description}, from=2-4, to=3-5]
	\arrow["{y|y\cdot}"{description}, from=2-4, to=4-5]
	\arrow[from=3-3, to=3-4]
	\arrow["\bigoplus"{description}, draw=none, from=3-3, to=4-3]
	\arrow["{y|y\cdot}"{description}, from=3-3, to=5-4]
	\arrow[from=3-4, to=3-5]
	\arrow["\bigoplus"{description}, draw=none, from=3-4, to=4-4]
	\arrow["{y|y\cdot}"{description}, from=3-4, to=5-5]
	\arrow["\bigoplus"{description}, draw=none, from=3-5, to=4-5]
	\arrow["{y|y\cdot}"{description}, from=3-5, to=5-6]
	\arrow[from=4-3, to=4-4]
	\arrow["\bigoplus"{description}, draw=none, from=4-3, to=5-3]
	\arrow["{\cdot x|x}"{description}, from=4-3, to=5-4]
	\arrow[from=4-4, to=4-5]
	\arrow["\bigoplus"{description}, draw=none, from=4-4, to=5-4]
	\arrow["{\cdot x|x}"{description}, from=4-4, to=5-5]
	\arrow["\bigoplus"{description}, draw=none, from=4-5, to=5-5]
	\arrow["{\cdot x|x}"{description}, from=4-5, to=5-6]
	\arrow[from=5-2, to=5-3]
	\arrow[from=5-3, to=5-4]
	\arrow[from=5-4, to=5-5]
	\arrow[from=5-5, to=5-6]
\end{tikzcd}\]
\caption{Chain complex for $(A_c/(U,V))_5\simeq \mathcal K_{A_c}(U|V)_5$.}
    \label{fig: K_5}
\end{figure}
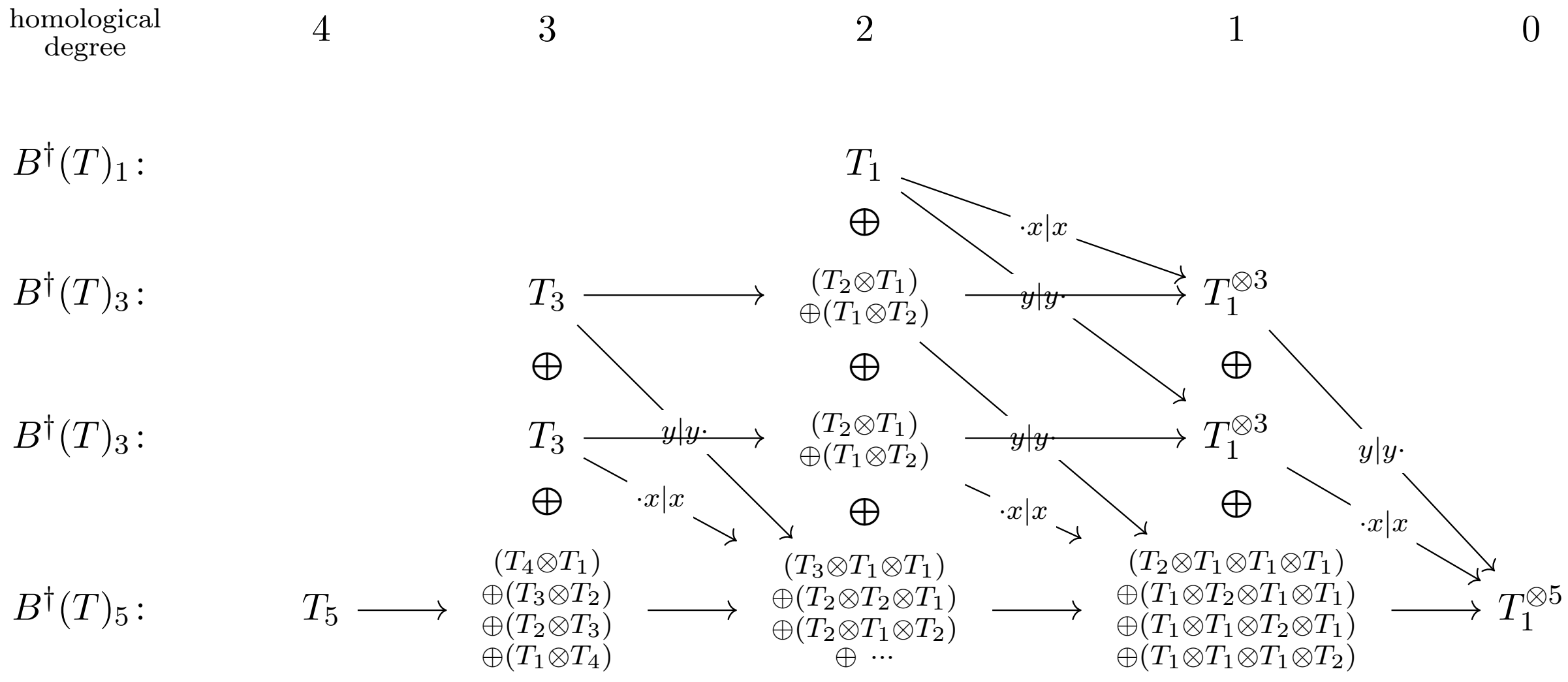

Next, let us look at quotients by two elements, say for example $U=y|y$ and $V=x|x$. There are a few models for the double quotient $A_c/(U,V)$, namely the two-sided complex $\mathcal K_{A_c}(U|V)$, and the one-sided complexes $\mathcal K^R_{A_c}(U,V)$ and $\mathcal K^L_{A_c}(U,V)$. We first look at $\mathcal K_{A_c}(U|V)$ in algebraic grading $n=5$. Since left and right multiplication commute, this is the total complex of the commuting square of chain complexes
\[\begin{tikzcd}
	{B^\dagger(T)_1} & {B^\dagger(T)_3} \\
	{B^\dagger(T)_3} & {B^\dagger(T)_5}
	\arrow["{\cdot x|x}", from=1-1, to=1-2]
	\arrow["{y|y\cdot}"', from=1-1, to=2-1]
	\arrow["{y|y\cdot}"', from=1-2, to=2-2]
	\arrow["{\cdot x|x}"', from=2-1, to=2-2]
\end{tikzcd}\]
which we illustrate in Figure \ref{fig: K_5}, again dropping the Koszul signs for convenience.

For the one-sided complex $\mathcal K^R_{A_c}(U,V)$, we have a diagram of chain complexes 
\[\begin{tikzcd}
	{B^\dagger(T)_1} & {B^\dagger(T)_3} \\
	{B^\dagger(T)_3} & {B^\dagger(T)_5}
	\arrow["{\cdot x|x}", from=1-1, to=1-2]
	\arrow["{\cdot y|y}"', from=1-1, to=2-1]
	\arrow["{\cdot H}", from=1-1, to=2-2]
	\arrow["{\cdot y|y}", from=1-2, to=2-2]
	\arrow["{\cdot x|x}"', from=2-1, to=2-2]
\end{tikzcd}\]
where the square no longer commutes, so we have to introduce a diagonal homotopy of homological degree $1$ given by right multiplication by some $H$, where we have $$dH=x|x|y|y-y|y|x|x$$
to account for the non-commutativity. Explicitly, by unwinding the description of $\E_2$-centrality in the proof of \cite[Corollary 4.3.2]{landesman2025cohenlenstramomentsfunctionfields}, $H$ is given geometrically by moving the two $y$s across the two $x$s, so we have
$$H=x|y|zy+x|xy|y+y|zy|x+xy|y|x.$$ 
Hence, Figure \ref{fig: K^R(U,V)_5} gives the chain-complex for $\mathcal K^R_{A_c}(U,V)$, where we note that compared to Figure \ref{fig: K_5} there is an additional map given by right multiplication by $H$.
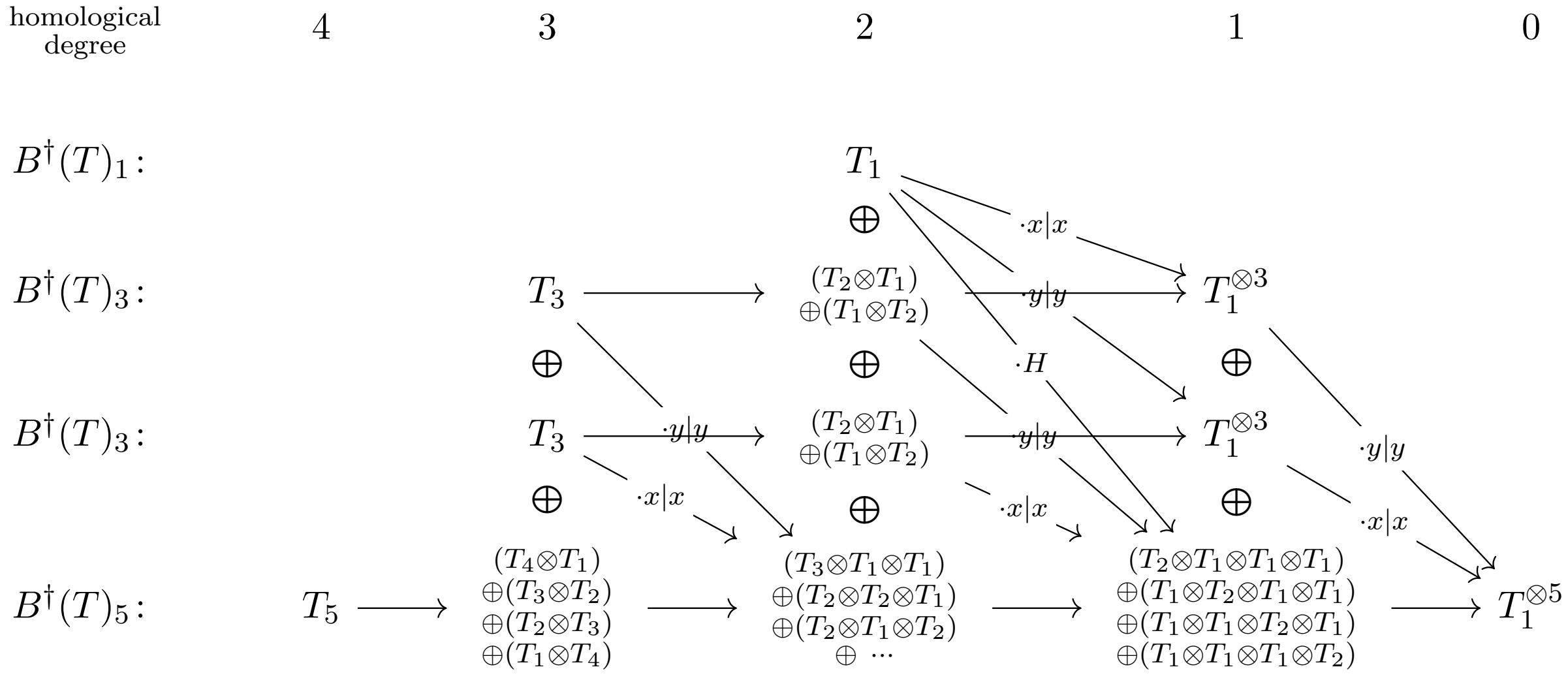
\begin{figure}[!htbp]
    \centering
\[\begin{tikzcd}[row sep=scriptsize]
	\begin{array}{c} \substack{\text{homological}\\\text{degree}} \end{array} & 4 & 3 & 2 & 1 & 0 \\
	{B^\dagger(T)_1\colon} &&& {T_1} \\
	{B^\dagger(T)_3\colon} && {T_3} & \begin{array}{c} \substack{(T_2\otimes T_1)\\\oplus (T_1\otimes T_2)} \end{array} & {T_1^{\otimes 3}} \\
	{B^\dagger(T)_3\colon} && {T_3} & \begin{array}{c} \substack{(T_2\otimes T_1)\\\oplus (T_1\otimes T_2)} \end{array} & {T_1^{\otimes 3}} \\
	{B^\dagger(T)_5\colon} & {T_5} & \begin{array}{c} \substack{(T_4\otimes T_1)\\\oplus (T_3\otimes T_2)\\\oplus (T_2\otimes T_3)\\\oplus (T_1\otimes T_4)} \end{array} & \begin{array}{c} \substack{(T_3\otimes T_1\otimes T_1)\\\oplus (T_2\otimes T_2\otimes T_1)\\\oplus(T_2\otimes T_1\otimes T_2)\\\oplus\ \cdots} \end{array} & \begin{array}{c} \substack{(T_2\otimes T_1\otimes T_1\otimes T_1)\\\oplus (T_1\otimes T_2\otimes T_1\otimes T_1)\\\oplus (T_1\otimes T_1\otimes T_2\otimes T_1)\\\oplus (T_1\otimes T_1\otimes T_1\otimes T_2)} \end{array} & {T_1^{\otimes 5}}
	\arrow["\bigoplus"{description}, draw=none, from=2-4, to=3-4]
	\arrow["{\cdot  x|x}"{description}, from=2-4, to=3-5]
	\arrow["{\cdot y|y}"{description}, from=2-4, to=4-5]
	\arrow["{\cdot H}"{description}, from=2-4, to=5-5]
	\arrow[from=3-3, to=3-4]
	\arrow["\bigoplus"{description}, draw=none, from=3-3, to=4-3]
	\arrow["{\cdot y|y}"{description}, from=3-3, to=5-4]
	\arrow[from=3-4, to=3-5]
	\arrow["\bigoplus"{description}, draw=none, from=3-4, to=4-4]
	\arrow["{\cdot y|y}"{description}, from=3-4, to=5-5]
	\arrow["\bigoplus"{description}, draw=none, from=3-5, to=4-5]
	\arrow["{\cdot y|y}"{description}, from=3-5, to=5-6]
	\arrow[from=4-3, to=4-4]
	\arrow["\bigoplus"{description}, draw=none, from=4-3, to=5-3]
	\arrow["{\cdot x|x}"{description}, from=4-3, to=5-4]
	\arrow[from=4-4, to=4-5]
	\arrow["\bigoplus"{description}, draw=none, from=4-4, to=5-4]
	\arrow["{\cdot x|x}"{description}, from=4-4, to=5-5]
	\arrow["\bigoplus"{description}, draw=none, from=4-5, to=5-5]
	\arrow["{\cdot x|x}"{description}, from=4-5, to=5-6]
	\arrow[from=5-2, to=5-3]
	\arrow[from=5-3, to=5-4]
	\arrow[from=5-4, to=5-5]
	\arrow[from=5-5, to=5-6]
\end{tikzcd}\]
    \caption{Chain complex for $(A_c/(U,V))_5\simeq \mathcal K^R_{A_c}(U,V)_5$.}
    \label{fig: K^R(U,V)_5}
\end{figure}

We give an important warning that not any choice of $H$ satisfying $dH=x|x|y|y-y|y|x|x$ works, and one does really have to trace through the given construction to deduce $H$. Indeed, the bimodule maps and relative tensor-product construction determine a specific homotopy-coherent square, and from this square $H$ will be defined uniquely up to homotopy, i.e. up to a boundary $dK$ for some element $K$ of homological degree $2$. To illustrate this further, in the case when $U=V=x|x$, where we require $dH=0$, the correct homotopy given by our construction is $H=x|x|xx+2x|xx|x+xx|x|x$, and one \textit{cannot} simply choose $H=0$. In fact, since the sum of coefficients differs, these two choices cannot be homotopic.

\section{Filtrations on cobar-complexes}\label{sec: filtrations}
In the previous section, we defined the homotopy-coherent Koszul complexes that we will be working with. In this section, we will specialize to the case when $A=A_V$ is the chain complex coming from a braided vector space $V$, which we have already seen some examples of in Section \ref{sec: koszul example}. In this case, $A$ can be written in terms of an explicit cobar-complex, and for $\E_2$-central $U_1,\ldots, U_k$, the quotient $A/(U_1,\ldots, U_k)$ can be written as a (homotopy-coherent) Koszul cobar-complex. 

We will construct filtrations on these Koszul cobar-complexes in a fixed algebraic degree whose graded pieces are made of similar Koszul cobar-complexes of smaller algebraic degree. Much like in \cite{prequel}, these filtrations will be combinatorial in nature. For each term of the Koszul cobar-complex, we will look at the partition of the algebraic grading determined by the term, and this will tell us which piece of the filtration the term is in. After introducing some notation on partitions in Section \ref{sec: algebraic grading}, we will carry out the construction of these filtrations in Section \ref{sec: construction filtration}. Finally, in Section \ref{sec: pf general results}, we will use these filtrations to deduce vanishing of iterated quotients $A/(U_1,\ldots, U_k)$ in large algebraic degree from vanishing in smaller degrees, giving the proof of our main theorems that we stated in Section \ref{sec: general results}. 

\subsection{Partitions on algebraic grading}\label{sec: algebraic grading} Before constructing the filtrations, we first discuss notations and conventions on associating a partition of the algebraic grading to terms in the cobar-complex. Throughout this section, we will work in the context of stability in a chosen direction, as in Section \ref{sec: stability in chosen direction}, which we recall briefly. Let $V=V_1\oplus V_2$ be a direct sum of braided vector spaces (i.e. it satisfies Equation \ref{eqn: V_1+V_2 condition}), and we define a corresponding algebraic bigrading $(n_1,n_2)$ such that $n_1+n_2=n$, and for our theorems we will be concerned with the grading $n_1$. This induces a bigrading on the tensor algebra $T=T(V)=\bigoplus_{n_1,n_2\ge 0}T_{n_1,n_2}$, and hence on the cobar-complex as well as the Koszul complex.

We discuss these partitions in the two-sided case $\mathcal K_A(U_1,\ldots, U_j|V_1,\ldots, V_k)$ as that will specialize to the left and right cases when we take $k=0$ or $j=0$ respectively. Let us take the $(n_1,n_2)$ graded part of the complex, which we write as $\mathcal K_A(U_1,\ldots, U_j|V_1,\ldots, V_k)_{n_1,n_2}$. Suppose that $U_i,V_i$ have algebraic gradings $(\rho_1^i,\rho_2^i)$ and $(\sigma_1^i, \sigma_2^i)$ respectively, as well as homological gradings $c_i,d_i$. Consider the term of the form 
\begin{equation}\label{eqn: term in cobar-complex}
e_J\otimes T_{\lambda_1^1,\lambda_2^1}\otimes\cdots \otimes  T_{\lambda_1^\ell,\lambda_2^\ell}\otimes f_L
\end{equation}
in the cobar-complex, where we must have $\lambda^i\coloneqq \lambda_1^i+\lambda_2^i>0$. One can compute that this term is in homological degree $$\sum_{i\in J}(c_i+1)+\sum_{i\in L}(d_i+1)+\sum_{i=1}^\ell (\lambda^i-1).$$

Let us now consider the algebraic bigrading of $e_J$ and $f_L$. By definition, $e_J$ has algebraic grading equal to the sum of the algebraic gradings of the $U_i$ for $i\in J$ which we write as $(\rho_1^J,\rho_2^J)\coloneqq \sum_{i\in J}(\rho_1^i, \rho_2^i)$. In the same way, $f_L$ has algebraic grading $(\sigma_1^L,\sigma_2^L)\coloneqq \sum_{i\in L}(\sigma_1^i,\sigma_2^i)$. For the algebraic grading to add up to $(n_1,n_2)$ we then must have $$n_i=\rho_i^J+\lambda_i^1+\cdots +\lambda_i^\ell+\sigma_i^L$$ 
for both gradings $i=1,2$.

\begin{figure}[!htbp]
\centering
\begin{tikzpicture}[
    x=0.8cm,
    y=0.8cm,
    >=Latex,
    path/.style={thick,-{Latex[length=2mm]}},
    comp/.style={dashed,thin},
    pt/.style={circle,fill=black,inner sep=2.4pt},
    endpt/.style={circle,draw=black,fill=white,inner sep=2.7pt},
    cornerlab/.style={font=\small},
    complab/.style={font=\scriptsize,fill=white,inner sep=1pt}
]

\draw[step=1, very thin, black!15] (0,0) grid (10,7);
\draw[semithick] (0,0) rectangle (10,7);

\node[cornerlab, below left]  at (0,0)  {$(0,0)$};
\node[cornerlab, below right] at (10,0) {$(n_1,0)$};
\node[cornerlab, above left]  at (0,7)  {$(0,n_2)$};
\node[cornerlab, above right] at (10,7) {$(n_1,n_2)$};

%

\coordinate (O)  at (0,0);
\coordinate (P0) at (2,1);
\coordinate (P1) at (4,1);
\coordinate (P2) at (5,3);
\coordinate (P3) at (5,4);
\coordinate (P4) at (8,5);
\coordinate (N)  at (10,7);

\coordinate (Q0) at (2,0);
\coordinate (Q2) at (5,1);
\coordinate (Q4) at (8,4);
\coordinate (Q5) at (10,5);

\draw[path] (O)  -- (P0);
\draw[path] (P0) -- (P1);
\draw[path] (P1) -- (P2);
\draw[path] (P2) -- (P3);
\draw[path] (P3) -- (P4);
\draw[path] (P4) -- (N);

\draw[comp]
    (O) -- node[below,complab] {$\rho^J_1$} (Q0)
        -- node[right,complab] {$\rho^J_2$} (P0);

\node[below,complab] at ($(P0)!0.5!(P1)$)
    {$\lambda^1_1$};

\draw[comp]
    (P1) -- node[below,complab] {$\lambda^2_1$} (Q2)
         -- node[right,complab] {$\lambda^2_2$} (P2);

\node[right,complab] at ($(P2)!0.5!(P3)$)
    {$\lambda^3_2$};

\draw[comp]
    (P3) -- node[below,complab] {$\lambda^4_1$} (Q4)
         -- node[right,complab] {$\lambda^4_2$} (P4);

\draw[comp]
    (P4) -- node[below,complab] {$\sigma^L_1$} (Q5)
         -- node[right,complab] {$\sigma^L_2$} (N);

\node[endpt] at (O) {};
\node[pt] at (P0) {};
\node[pt] at (P1) {};
\node[pt] at (P2) {};
\node[pt] at (P3) {};
\node[pt] at (P4) {};
\node[endpt] at (N) {};

\end{tikzpicture}

\caption{
Lattice path associated to a term in the Koszul cobar-complex, illustrated
in the case $l=4$. The filled points are the elements of the sequence $S$.
}
\label{fig: grading-path}
\end{figure}
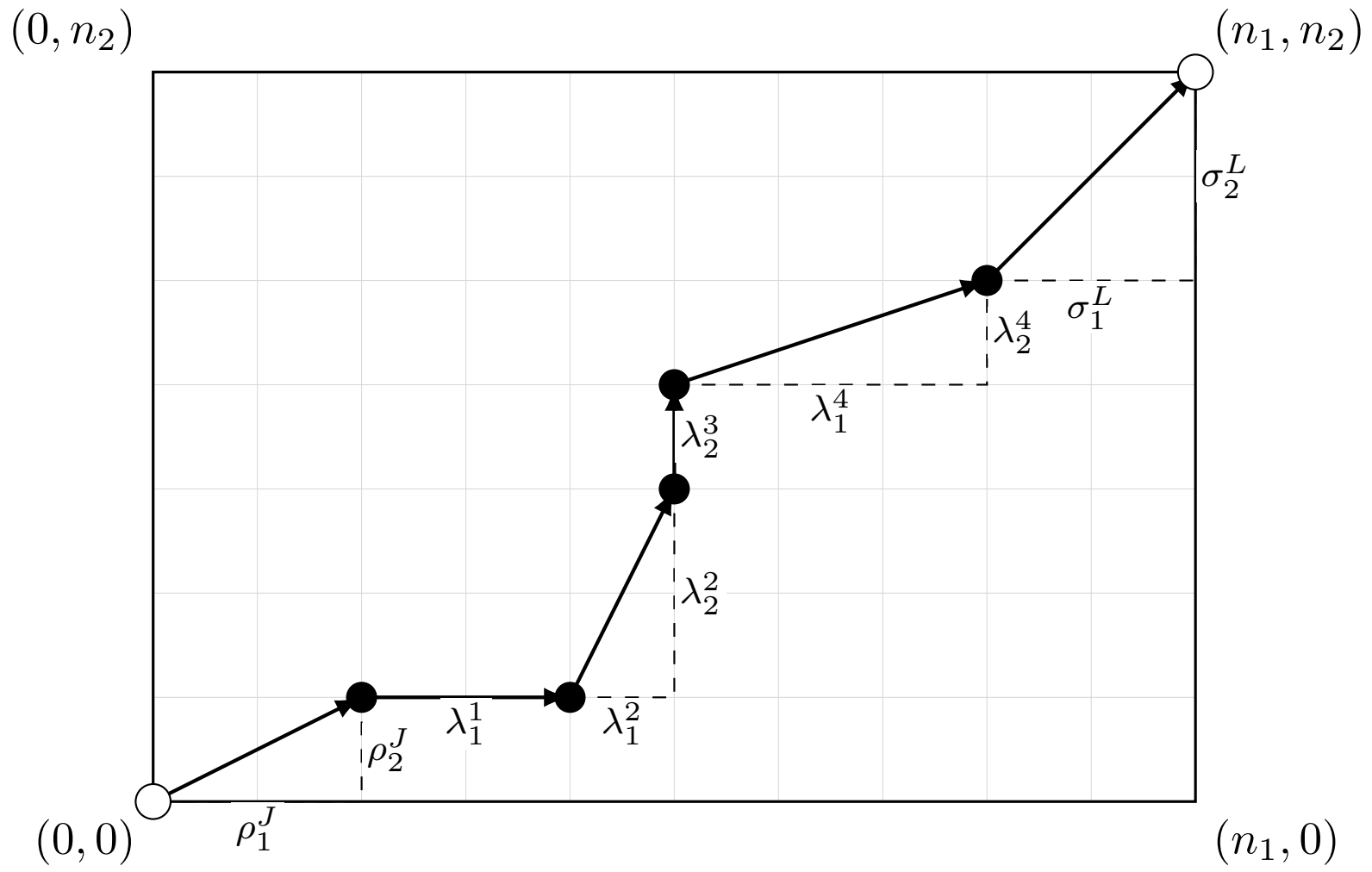

Then, to this term we associate the sequence of tuples
$$S=\left((\rho_1^J,\rho_2^J),(\rho_1^J+\lambda_1^1,\rho_2^J+\lambda_2^1),\ldots, (\rho_1^J+\lambda_1^1+\cdots +\lambda_1^\ell,\rho_2^J+\lambda_2^1+\cdots +\lambda_2^\ell)\right).$$
One should picture this as a sequence of lattice points on $[0,n_1]\times [0,n_2]$ which traces out a path going from $(0,0)$ to $(n_1,n_2)$ whose coordinates are non-decreasing, as illustrated in Figure \ref{fig: grading-path}. 

Note that no two tuples in $S$ are the same, and the successive differences are simply the algebraic grading of each tensor factor. Also note that the first and last points can be $(0,0)$ and $(n_1,n_2)$ when $J=\emptyset$ and $L=\emptyset$ respectively. We remark that this construction is analogous to the one-variable setting, in which one associates to a partition its set of separating points.

Lastly, we have the following lemma on how the differential of the chain complex changes the sequence $S$. Once again, this also applies to the left and right Koszul complexes as special cases.
\begin{lemma}\label{lem: sequence S}
Suppose we have an element $x$ in $\mathcal K_A(U_1,\ldots, U_j|V_1,\ldots, V_k)$ that lies in a term of the form \eqref{eqn: term in cobar-complex} with associated sequence $S$. Then, $dx$ is contained in a direct sum of terms \eqref{eqn: term in cobar-complex} whose associated sequences $S'$ contain $S$ as a subsequence.
\end{lemma}
\begin{proof}
By Definition \ref{def: two-sided koszul}, it suffices to check that $d_A, d^L_I$ and $d^R_K$ acting on $x$ satisfy the desired condition. The differential $d_A$ acts via the differential of the cobar-complex, so it is given by the sum over $i$ of comultiplication on each $T_{\lambda_1^i,\lambda_2^i}$. Each comultiplication splits $T_{\lambda_1^i,\lambda_2^i}$ into a direct sum of tensor products $T_{\alpha_1,\alpha_2}\otimes T_{\beta_1,\beta_2}$ where $\alpha_1+\beta_1=\lambda_1^i$ and $\alpha_2+\beta_2=\lambda_2^i$, which effectively adds a new tuple to the sequence $S$. 

For the differential $d^L_I$, $g_I$ must lie in the direct sum of terms of the form $T_{\alpha_1^1,\alpha_2^1}\otimes \cdots \otimes T_{\alpha_1^p,\alpha_2^p}$ as it is an element of $A$. Hence, $d^L_I x$ lies within the direct sum of terms of the form $e_{J\setminus I} \otimes T_{\alpha_1^1,\alpha_2^1}\otimes \cdots \otimes T_{\alpha_1^p,\alpha_2^p}\otimes T_{\lambda_1^1,\lambda_2^1}\otimes\cdots \otimes  T_{\lambda_1^\ell,\lambda_2^\ell}\otimes f_L$. This has the effect of adding several tuples to the beginning of the sequence. The same argument works for $d^R_K$, where tuples are added instead to the end of the sequence, finishing the proof.
\end{proof}

\subsection{Construction of filtrations}\label{sec: construction filtration}

We will first construct the general two-sided filtration on the two-sided Koszul cobar-complex $\mathcal K_A(U_1,\ldots, U_j|V_1,\ldots, V_k)_{n_1,n_2}$ with algebraic grading $(n_1,n_2)$, which will depend on a parameter $\rho_1\le c<n_1-\sigma_1$ which gives the ``cutoff line'' in the $n_1$ direction. Here, we let $\rho_1=\sum_{i=1}^j\rho_1^i$ be the sum of the $n_1$ gradings of the left generators $U_1,\ldots, U_j$, and likewise let $\sigma_1=\sum_{i=1}^k\sigma_1^i$ be the corresponding sum for the right generators.

\begin{proposition}\label{prop: two-sided filtration}
Let $c$ be an integer such that $\rho_1\le c<n_1-\sigma_1$. Then, there exists a filtration on $\mathcal K_A(U_1,\ldots, U_j|V_1,\ldots, V_k)_{n_1,n_2}$ with graded pieces of the form
\begin{equation}\label{eqn: graded piece form}
\mathcal K^L_A(U_1,\ldots, U_j)_{\alpha_1,\alpha_2}\otimes T_{\lambda_1,\lambda_2}\otimes \mathcal K^R_A(V_1,\ldots, V_k)_{\beta_1,\beta_2} 
\end{equation}
for all $\alpha_i,\lambda_i,\beta_i\ge 0$ satisfying 
\begin{equation*}
\begin{split}
\alpha_i+\lambda_i+\beta_i&=n_i \qquad 
\text{for } i=1,2,\\
\alpha_1&\le c, \\
\beta_1&\le n_1-c-1
\end{split}
\end{equation*}
Here, we take the usual non-derived tensor product of complexes, where the differential acts trivially (by zero) on the middle factor and acts as in Definition \ref{def: right koszul} on the left and right factors.
\end{proposition}
\begin{proof}
Recall from the previous section that the terms in the complex are given by Equation \eqref{eqn: term in cobar-complex} with bigrading $(n_1,n_2)$. Consider such a term, and  let us associate the corresponding sequence $S$. Let $(\alpha_1,\alpha_2)$ be the last tuple in the sequence where $\alpha_1\le c$, and let $(n_1-\beta_1,n_2-\beta_2)$ be the first tuple in the sequence where $n_1-\beta_1>c$. Both of these must exist because the condition tells us that $\rho_1\le c<n_1-\sigma_1$, and the very first term in $S$ is $(\rho^J_1,\rho^J_2)$ and $\rho^J_1\le \rho_1$ by definition, and likewise for the last term. It is clear by definition that these two tuples $(\alpha_1,\alpha_2),(n_1-\beta_1,n_2-\beta_2)$ appear next to each other in $S$; we call these the \textit{central points} of $S$. We then let $(\lambda_1,\lambda_2)$ be the difference between these two terms, and note that by construction $\alpha_i,\lambda_i,\beta_i$ satisfy the relations in the statement of the proposition. Let us also call $\lambda\coloneqq \lambda_1+\lambda_2$ the \textit{span} of the term; this is the taxicab distance between the two central points. We illustrate this in Figure \ref{fig: two sequences}.

\begin{figure}[!htbp]
\centering
\begin{tikzpicture}[
    x=0.8cm,
    y=0.8cm,
    >=Latex,
    redpath/.style={thick,draw=red!75!black,-{Latex[length=2mm]}},
    bluepath/.style={thick,draw=blue!75!black,-{Latex[length=2mm]}},
    redcomp/.style={dashed,thin,draw=red!75!black},
    bluecomp/.style={dashed,thin,draw=blue!75!black},
    redpt/.style={circle,fill=red!75!black,inner sep=2.4pt},
    bluept/.style={circle,fill=blue!75!black,inner sep=2.4pt},
    endpt/.style={circle,draw=black,fill=white,inner sep=2.7pt},
    centralred/.style={
        circle,draw=red!75!black,fill=white,
        line width=0.9pt,inner sep=4.0pt
    },
    centralblue/.style={
        circle,draw=blue!75!black,fill=white,
        line width=0.9pt,inner sep=4.0pt
    },
    cornerlab/.style={font=\small},
    complab/.style={font=\scriptsize,fill=white,inner sep=1pt},
    spanlab/.style={font=\small,fill=white,inner sep=1pt}
]

\draw[step=1, very thin, black!15] (0,0) grid (10,7);
\draw[semithick] (0,0) rectangle (10,7);

\node[cornerlab, below left]  at (0,0)  {$(0,0)$};
\node[cornerlab, below right] at (10,0) {$(n_1,0)$};
\node[cornerlab, above left]  at (0,7)  {$(0,n_2)$};
\node[cornerlab, above right] at (10,7) {$(n_1,n_2)$};

\draw[densely dotted,thick] (5.5,0) -- (5.5,7);

\node[cornerlab,anchor=base] at (5,-0.43) {$\mathstrut c$};
\node[cornerlab,anchor=base] at (6,-0.43) {$\mathstrut c+1$};

\coordinate (O) at (0,0);
\coordinate (N) at (10,7);


\coordinate (R0) at (2,1);
\coordinate (R1) at (4,2);
\coordinate (R2) at (6,3);
\coordinate (R3) at (8,4);
\coordinate (R4) at (9,5);
\coordinate (RQ) at (6,2);

\draw[redpath] (O)  -- (R0);
\draw[redpath] (R0) -- (R1);
\draw[redpath] (R1) -- (R2);
\draw[redpath] (R2) -- (R3);
\draw[redpath] (R3) -- (R4);
\draw[redpath] (R4) -- (N);

\node[redpt] at (R0) {};
\node[redpt] at (R1) {};
\node[redpt] at (R2) {};
\node[redpt] at (R3) {};
\node[redpt] at (R4) {};

\node[centralred] at (R1) {};
\node[centralred] at (R2) {};
\node[redpt] at (R1) {};
\node[redpt] at (R2) {};

\draw[redcomp]
    (R1) -- node[below,complab] {$\lambda^R_1$} (RQ)
         -- node[right,complab] {$\lambda^R_2$} (R2);

\node[spanlab,text=red!75!black]
    at ($(R1)!0.5!(R2)+(0,0.5)$)
    {$\lambda^R$};


\coordinate (B0) at (1,2);
\coordinate (B1) at (3,4);
\coordinate (B2) at (5,4);
\coordinate (B3) at (6,6);
\coordinate (B4) at (8,6);
\coordinate (BQ) at (6,4);

\draw[bluepath] (O)  -- (B0);
\draw[bluepath] (B0) -- (B1);
\draw[bluepath] (B1) -- (B2);
\draw[bluepath] (B2) -- (B3);
\draw[bluepath] (B3) -- (B4);
\draw[bluepath] (B4) -- (N);

\node[bluept] at (B0) {};
\node[bluept] at (B1) {};
\node[bluept] at (B2) {};
\node[bluept] at (B3) {};
\node[bluept] at (B4) {};

\node[centralblue] at (B2) {};
\node[centralblue] at (B3) {};
\node[bluept] at (B2) {};
\node[bluept] at (B3) {};
\node[endpt] at (O) {};
\node[endpt] at (N) {};

\draw[bluecomp]
    (B2) -- node[below,complab] {$\lambda^B_1$} (BQ)
         -- node[right,complab] {$\lambda^B_2$} (B3);

\node[spanlab,text=blue!75!black]
    at ($(B2)!0.5!(B3)+(-0.35,0.25)$)
    {$\lambda^B$};

\end{tikzpicture}

\caption{
Two sequences of equal span $\lambda^R=\lambda^B=3$. The dotted line separates points with first coordinate $\le c$ and $>c$. The circled vertices are the central points.
}
\label{fig: two sequences}
\end{figure}
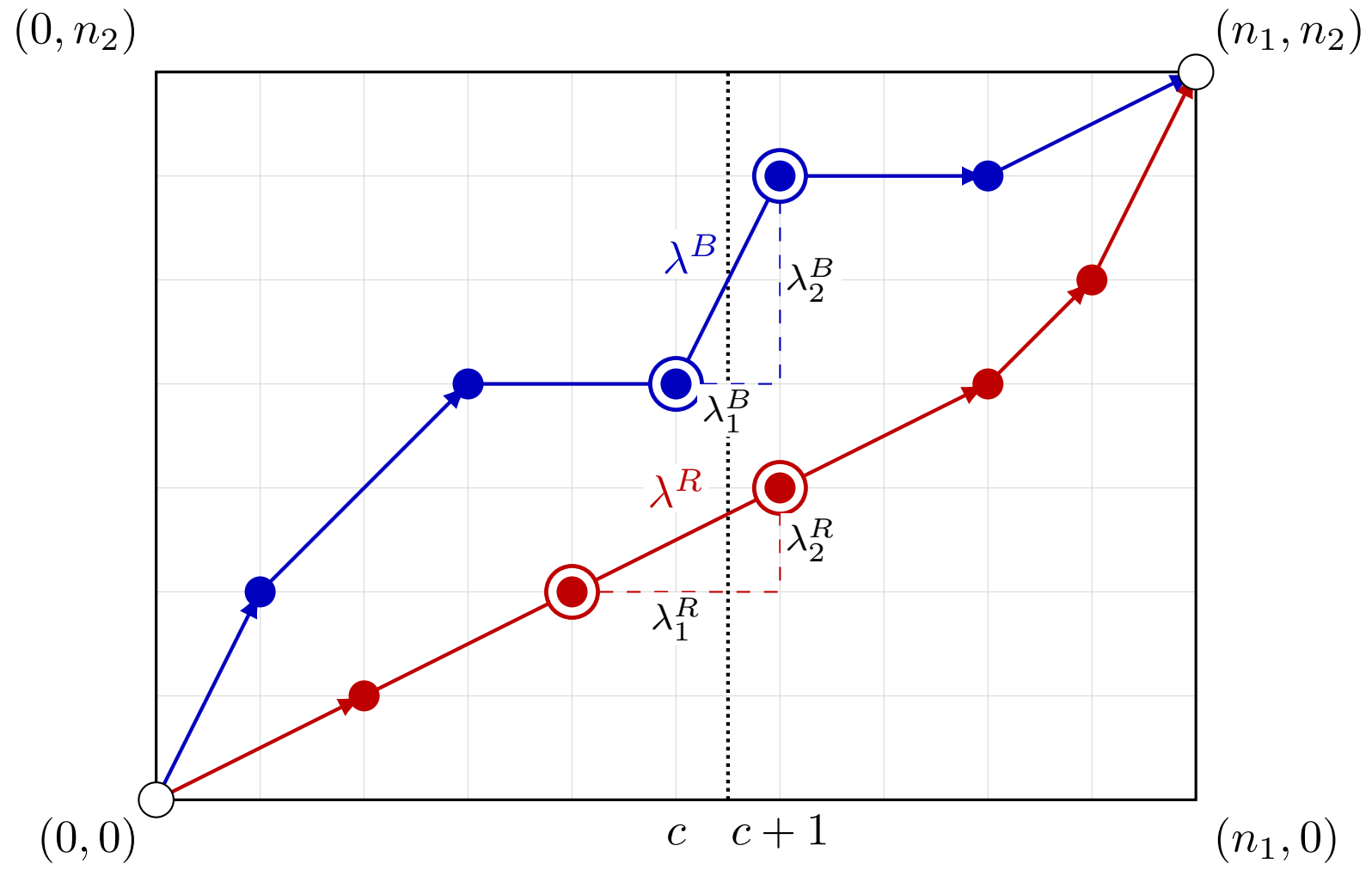

We define an increasing filtration
$$F_{p}\mathcal K_A(U_1,\ldots, U_j|V_1,\ldots, V_k),$$
where each filtered piece contains exactly the terms of the complex whose associated span satisfies $\lambda\le p$. To check that the differential preserves each filtered piece, note that by Lemma \ref{lem: sequence S}, since the differential can only add terms to the sequence $S$, the span $\lambda$ can only stay the same or decrease.

We now look at the $p$-th graded piece of the filtration $\Gr_p\mathcal K_A(U_1,\ldots, U_j|V_1,\ldots, V_k)$, which consists of terms with span $\lambda=p$. We claim that this decomposes into a direct sum of subcomplexes, where each complex consists of terms that have the same pair of central points $(\alpha_1,\alpha_2), (n_1-\beta_1,n_2-\beta_2)$. To prove this, we simply need to check that within this graded piece, the differential must send the term with central points to a term with the same central points. Indeed, this is true again by Lemma \ref{lem: sequence S} as the differential can only add points to the sequence $S$, and so it either decreases the span or keeps the central points unchanged. The first option sends it deeper into the filtration and thus acts by zero within the graded piece, leaving the second choice as desired. It may be helpful to look at Figure \ref{fig: two sequences} to convince oneself that the differential cannot send terms with the blue sequence to terms with the red sequence and vice versa.

By further refining this filtration to separate each graded piece, we obtain another filtration where each graded piece consists of exactly those terms with some central points $(\alpha_1,\alpha_2),(n_1-\beta_1,n_2-\beta_2)$. It is clear from construction that these terms are exactly the terms appearing in the complex \eqref{eqn: graded piece form}, so it remains to check that the differential agrees. The components $d^L_I,d^R_K$ of the differential act within the graded piece as they do not affect the term $T_{\lambda_1,\lambda_2}$. The other component $d_A$ is the sum of comultiplication of each $T_{\lambda_1^i,\lambda_2^i}$, which will stay in this graded piece unless it acts on the term $T_{\lambda_1,\lambda_2}$, in which case it sends it deeper into the filtration and thus sends it to zero on the graded piece. Hence, we see that $d^L_I$ and the comultiplications to the left of $T_{\lambda_1,\lambda_2}$ are exactly the components of the differentials for $\mathcal K^L_A(U_1,\ldots, U_j)$, and the analogue holds for $\mathcal K^R_A(V_1,\ldots, V_k)$ as well. Following the Koszul sign convention ensures that the signs agree, finishing the proof.
\end{proof}

We will now specialize this to a one-sided filtration on $\mathcal K^R_A(U_1,\ldots, U_k)_{n_1,n_2}$ as a corollary.
\begin{corollary}\label{cor: one-sided filtration}
Suppose that $n_1>\sigma_1$. Then, there exists a filtration on $\mathcal K^R_A(U_1,\ldots, U_k)_{n_1,n_2}$ with graded pieces of the form
\begin{equation}\label{eqn: graded piece one-sided}
A_{0,\alpha_2}\otimes T_{\lambda_1,\lambda_2}\otimes \mathcal K^R_A(U_1,\ldots, U_k)_{\beta_1,\beta_2} 
\end{equation}
for all $\alpha_i,\beta_i,\lambda_i$ satisfying $\lambda_1+\beta_1=n_1$, $\alpha_2+\lambda_2+\beta_2=n_2$ and $\lambda_1>0$. 
\end{corollary}
\begin{proof}
This is a special case of Proposition \ref{prop: two-sided filtration} applied to $\mathcal K_A(-|U_1,\ldots, U_k)_{n_1,n_2}$ with $c=0$. Note that in this case $\mathcal K^L_A(-)$ is simply $A$, and the inequalities there force $\alpha_1=0$ and $\lambda_1>0$.
\end{proof}
\subsection{Proof of homological theorems}\label{sec: pf general results}

Finally, we apply the filtrations on the Koszul cobar-complex that we constructed to prove our homological theorems in Section \ref{sec: general results} and Section \ref{sec: intro main theorem}. We will use the one-sided filtration to prove Theorem \ref{thm: usual vs strong stability direction}, and the two-sided filtration to prove Theorem \ref{thm: general main theorem}.

Let us recall some notation in the directional case. Let $A=A_V$ be the simplicial chains for some braided vector space $V$, and let $U_1,\ldots, U_k$ be $\E_2$-central elements with algebraic grading $(a_i,b_i)$ and homological degree $d_i$ such that $a_i>0$. Denote the sums $d=\sum_i d_i$, $a=\sum_i a_i$, and let $s_i=d_i/a_i$ be the directional slopes. Furthermore, let $s_0=\max_i s_i$ be the maximal slope.

We first prove Theorem \ref{thm: usual vs strong stability direction} that the notions of stability and strong stability in a staircase are equivalent when all $d_i=0$.
\begin{proof}[Proof of Theorem \ref{thm: usual vs strong stability direction}] For the proof we will use shifted coordinates (in the sense of Remark \ref{rmk: shifted}) to minimize confusion. Referring to the left of Figure \ref{fig: strong vs usual staircase}, we want to show that the blue region vanishes, i.e. that $H_i((A/(U_1,\ldots, U_k))_{n_1,n_2})$ vanishes for all $n_1>N+a$ and $i\le I+k+d$. The condition tells us that this homology vanishes for $n_1\le N+a$ and $i\le I+n_1-N-a+k+d$. 

In fact, it suffices to prove only for $n_1=N+a+1$ that $H_i((A/(U_1,\ldots, U_k))_{n_1,n_2})$ vanishes for $i\le I+k+d$, as that would imply that our assumption holds true when $N$ is replaced by $N+1$, and we can iterate this to prove the claim for all $n_1$. Put visually as in the left of Figure \ref{fig: strong vs usual staircase}, it suffices to prove vanishing for the leftmost vertical slice of the blue block as iterating that implies vanishing for each vertical slice of the blue block from left to right.

We now take $n_1=N+a+1>a$, so it satisfies the $n_1>\sigma_1$ condition in Corollary \ref{cor: one-sided filtration}. By Proposition \ref{prop: right koszul}, we use the right-sided Koszul cobar-complex $\mathcal K^R_A(U_1,\ldots, U_k)$ as a chain complex model for the iterated quotient $A/(U_1,\ldots, U_k)$. Then, by Corollary \ref{cor: one-sided filtration}, it suffices to show that each graded piece \eqref{eqn: graded piece one-sided} vanishes in homological degrees $\le I+k+d$. The first factor $A_{0,\alpha_2}$ has non-negative homological degree and the second $T_{\lambda_1,\lambda_2}$ has homological degree $\lambda_1+\lambda_2-1\ge \lambda_1-1$, where we also recall that $\lambda_1>0$. For the last factor $\mathcal K^R_A(U_1,\ldots, U_k)_{\beta_1,\beta_2}$, since $\beta_1\le N+a$, then by our assumption this vanishes in homological degree $\le I+\beta_1-N-a+k+d$. Putting these facts together, the graded piece vanishes in homological degree $\le I+\beta_1+\lambda_1-N-a-1+k+d=I+k+d$ as $\beta_1+\lambda_1=n_1=N+a+1$, which completes the proof.
\end{proof}

We remark that this proof also shows the same vanishing rectangular blue region for $(n_1,i)$ in Figure \ref{fig: strong vs usual staircase} even when some $d_i\neq 0$, although that is not the correct region for strong stability when $s_0>0$.

Now, we turn to the proof of Theorem \ref{thm: general main theorem}. We first show the following lemma, which tells us that we can deduce strong stability for a larger staircase in terms of two smaller staircases. In the proof, we will use the two-sided filtration of the complex $\mathcal K_A(U_1,\ldots, U_k|U_1,\ldots, U_k)$ which has graded pieces which are given in terms of the left and right complexes $\mathcal K^L_A(U_1,\ldots, U_k)$ and $\mathcal K^R_A(U_1,\ldots, U_k)$. 
\begin{lemma}\label{lem: staircase implies staircase}
Suppose that $A$ is homologically stable, stabilizes in an $(N,I)$ staircase, and stabilizes strongly in an $(N',I')$ staircase, all with respect to $U_1,\ldots, U_k$ in the $n_1$ direction. Then $A$ stabilizes strongly in an $(N+N'+1,I+I'+1)$ staircase with respect to $U_1,\ldots, U_k,U_1,\ldots, U_k$ in the $n_1$ direction. Hence, by Lemma \ref{lem: remove extra}, $A$ stabilizes strongly in an $(N+N'+1,I+I'+1)$ staircase with respect to $U_1,\ldots, U_k$ in the $n_1$ direction.
\end{lemma}

\begin{proof}
Let us use $\mathcal K_A(U_1,\ldots, U_k|U_1,\ldots, U_k)_{n_1,n_2}$ as a chain complex for the iterated quotient $(A/(U_1,\ldots, U_k,U_1,\ldots, U_k))_{n_1,n_2}$. Define the constants $S=N+N'+1+2a,T=I+I'+1+2k+2d$. We need to show that the $i$-th homology of this complex vanishes for all $(n_1,i)$ satisfying both $i\le T+(n_1-S)$ and $i\le T+(n_1-S)s_0$. By Proposition \ref{prop: two-sided filtration}, it suffices to show this for each graded piece \eqref{eqn: graded piece form} for some suitable value of $c$ which we choose to be $c=N+a$. This is valid in the case of $n_1>N+2a$ as here the cutoff $c$ will lie in the range given by Proposition \ref{prop: two-sided filtration}. We will deal with this case first, and consider the $n_1\le N+2a$ case at the end of the proof.

Let us now look at the homology of each tensor factor. Like in the previous proof, the middle factor $T_{\lambda_1,\lambda_2}$ has homological degree $\ge \lambda_1-1$. We now look at the left tensor factor $\mathcal K^L_A(U_1,\ldots, U_k)_{\alpha_1,\alpha_2}$. Since $\alpha_1\le c = N+a$, the assumption $A$ stabilizes in an $(N,I)$ staircase tells us by definition that this factor vanishes in homological degrees $\le I-N+\alpha_1-a+k+d$. Finally, for the last tensor factor $\mathcal K^R_A(U_1,\ldots, U_k)_{\beta_1,\beta_2}$, strong vanishing in an $(N',I')$ staircase tells us that this factor vanishes in homological degrees 
$$\le I'+k+d+  \min\left(\beta_1-a-N',(\beta_1-a-N')s_0\right).$$
Adding everything up, the graded piece vanishes in homological degrees less than or equal to
\begin{equation*}
(I-N+\alpha_1-a+k+d)+(\lambda_1-1)+I'+k+d+  \min\left(\beta_1-a-N',(\beta_1-a-N')s_0\right)+1
\end{equation*}

Substituting the first argument in $\min$, we obtain $T+n_1-S$. On the other hand, substituting the second argument, we obtain
$$T+(n_1-S)s_0 +(\alpha_1+\lambda_1-N-a-1)(1-s_0).$$
However, we have that $\alpha_1+\lambda_1=n_1-\beta_1\ge c+1=N+a+1$ by the second inequality in Proposition \ref{prop: two-sided filtration}, so the latter product is non-negative so this is lower bounded by $T+(n_1-S)s_0$. This gives the desired vanishing homological degrees, proving the lemma when $n_1>N+2a$.

Let us now deal with the edge case where $n_1\le N+2a$, where we need to supply an additional argument. We do not need this for our applications, as we are only interested in extrapolating to larger $n_1$, but we include it for completeness. Suppose again that $n_1>a$, and choose the cutoff $c=n_1-a-1$. Let $P$ be the subcomplex of terms \eqref{eqn: term in cobar-complex} where $a_J\le n_1-a-1$, where $a_J=\sum_{i\in J}a_i$. This forces the first term in the sequence $S$ to still have $n_1$ coordinate $\le c$, so the span filtration holds as usual and we have the desired vanishing from the same argument above. If $n_1\le a$, we simply take $P$ to be the zero complex. 

Now, we consider $K/P$, and filter the terms by left exterior cardinality $|J|$, so the graded pieces look like $e_J\otimes \mathcal K^R_A(U_1,\ldots, U_k)_{\beta_1,\beta_2}$, where $\beta_1\le a$ by the definition of $P$. Strong stability in the $(N',I')$ staircase tells us that the homology of the right term vanishes in degrees $\le I'+k+d+\beta_1-a-N'$. The homological grading of $e_J$ is $d_J+|J|=\sum_{j\in J}(d_j+1)$. Adding these together, we get that the graded piece vanishes in homological degrees 
$$\le I'+k+d-N'-a+n_1-a_J+d_J+|J|$$
and subtracting the desired bound of $T+n_1-S$, we obtain a difference of 
$$\ge N-I+(a-a_J)-(d+k-d_J-|J|)=N-I+\sum_{j\not\in J}(a_j-d_j-1)\ge 0$$
as we have $I<N$ and $d_j+1\le a_j$ by definition, from the setup.
\end{proof}

Now, by iteratively applying the lemma starting from one staircase, we immediately obtain the following corollary.
\begin{corollary}
Suppose that $A$ is homologically stable and stabilizes strongly in an $(N,I)$ staircase with respect to $U_1,\ldots, U_k$ in the $n_1$ direction. Then, for each integer $r\ge 1$, $A$ stabilizes strongly in an $(r(N+1)-1,r(I+1)-1)$ staircase with respect to $U_1,\ldots, U_k$ in the $n_1$ direction.
\end{corollary}
\begin{proof}[Proof of Theorem \ref{thm: general main theorem}] 
By applying the corollary, we get vanishing in the union of the regions $$i\le r(I+1)-1+\min\left(n_1-r(N+1)+1,(n_1-r(N+1)+1)s_0\right),$$
and one can check that the union of these is exactly the desired equation.
\end{proof}

Finally, we prove Theorem \ref{thm: weight}.
\begin{proof}[Proof of Theorem \ref{thm: weight}]
To prove this theorem, we simply need to replace the homological grading $i$ with the negative of the Frobenius weight grading $-w$ as defined in \cite{ma2026weightfiltrationhurwitzspaces} in the proofs. The arguments work as expected; for example, the span filtration respects this grading, the weights add in the tensor product pieces, and the factor $T_{\lambda_1,\lambda_2}$ has weight grading $-w\ge i\ge \lambda_1-1$.
\end{proof}
All other theorems in Sections \ref{sec: general results} and \ref{sec: intro main theorem} then follow from the ones proven here (or were already proven earlier).
\section{Applications to Cohen-Lenstra moments}\label{sec: cohen lenstra}

For the next few sections, we will use our homological theorems to prove the arithmetic applications as described in Section \ref{sec: intro arithmetic applications}. There will be two ways we do this. First, we can deduce a stability slope from knowledge of $R=H_0(\Hur^{G,c},\Q)$, using Corollary \ref{cor: stab from R} in the classical case when there is only one stabilization element, and using a generalization of that argument when there are multiple. Secondly, we can directly compute the homology for some specific Hurwitz space and apply the main theorem directly, which in theory gives us near-optimal secondary term bounds by Corollary \ref{cor: optimal bound} or Corollary \ref{cor: general optimal}.

In this section, we will carry these out in the case of Cohen-Lenstra moments. Here, the relevant Hurwitz space $\Hur^{G,c}_n$ is for the semidirect product $G=H\rtimes \Z/2\Z$ for an odd abelian group $H$ where the $\Z/2\Z$ acts by inversion and $c=G\setminus H$ is the union of conjugacy classes consisting of elements $(x,1)\in H\rtimes \Z/2$ i.e. reflections. We will use this choice of $G$ and $c$ throughout this section. To see why this Hurwitz space is relevant for Cohen-Lenstra moments, we refer the reader to \cite[Proposition 8.7]{EVW16}. 

In Section \ref{sec: cohen lenstra R}, we give an explicit parametrization of the connected components of $\Hur^{G,c}_n$. This allows us to prove the following theorem on the stability of $R$ in Section \ref{sec: cohen lenstra slope}.
\begin{theorem}\label{thm: cohen lenstra stab R}
Let $D\ge 1$ be any integer and $U_D=\sum_{g\in c} g^{2D}$. Then, for $n\ge 2r+1$ where $r=\rk(H)$, we have that
$$U_D\colon H_0(\Hur^{G,c}_n,\Q)\xrightarrow\sim H_0(\Hur^{G,c}_{n+2D},\Q)$$
is an isomorphism.
\end{theorem}

By applying Corollary \ref{cor: stab from R} with $N=2r+2$ to this theorem, while noting that we have homological stability with respect to $U_D$ by \cite[Proposition A.3.1]{ellenberglandesman}, we get a stability slope of $\frac{1}{2r+3}$ in the following theorem.
\begin{theorem}\label{thm: cohen lenstra stab slope}
Let $D\ge 1$ be any integer and $U_D = \sum_{g\in c}g^{2D}$. Then for $n\ge (2r+3)i+(2r+2)$ where $r=\rk(H)$, we have that $$U_D\colon H_i(\Hur^{G,c}_n,\Q)\xrightarrow\sim H_i(\Hur^{G,c}_{n+2D},\Q)$$
is an isomorphism.
\end{theorem}

This theorem allows us to deduce the arithmetic application in Theorem \ref{thm: cohen-lenstra H}, which we will explain in Section \ref{sec: cohen lenstra slope}. In Section \ref{sec: H1 conjectures}, we explain Conjecture \ref{conj: H1} on the $H_1(\Hur^{G,c},\Q)$ in the case of dihedral groups $D_k$ where $k$ is prime heuristically from the perspective of the Siegel threefold $\mathcal A_2$. 

\subsection{Parametrization of $R$}\label{sec: cohen lenstra R} There are many papers that tackle the Cohen-Lenstra heuristics over function fields using Hurwitz spaces, such as \cite{EVW16, liu2022nonabeliancohenlenstraheuristicspresence,Liu_Wood_Zureick-Brown_2024,landesman2025cohenlenstramomentsfunctionfields,Sawin_Wood_2026}, and part of these papers consists of computing how many components of $\Hur^{G,c}_n$ there are when $n$ is large. However, we need a slight improvement to this statement that we could not find in the literature, which is after which $n$ the number of components starts to stabilize. We remark that the dihedral case of $G=\Z/k\Z\rtimes \Z/2\Z$ was studied in \cite{Hou_2008,Sia_2009,Berger_2011}, but we will see that this case is significantly easier as it does not involve any additional lifting invariant.

We first give a quick exposition of the current literature roughly following the proof of \cite[Proposition 5.3.3]{landesman2025cohenlenstramomentsfunctionfields}, state our refined theorem for $G=H\rtimes \Z/2\Z$, then give an explicit self-contained proof of the theorem using the Burau representation.

Recall from \cite[Section 3.3]{EVW16} that the components of Hurwitz space $\Hur^{G,c}_n$ are parametrized by the Hurwitz orbits $c^n/B_n$ of tuples $(g_1,\ldots, g_n)$ under the braid action as described in Definition \ref{def: rack}. There are some rather clear combinatorial invariants of these orbits as described in \cite[Section 2.4]{EVW16}, such as the global monodromy $\langle g_1,\ldots, g_n\rangle$ which is the subgroup generated by the elements of the tuple and the boundary monodromy $g_1\ldots g_n$ which is the product of the elements. However, these two invariants turn out to be insufficient for the current context of Cohen-Lenstra moments. 

It turns out that for a general Hurwitz space one can define an algebraic lifting invariant $U(G,c)$ which completely characterizes the stable components of Hurwitz space. This is described in \cite{Wood_2021}, and is based on ideas in the retracted paper \cite{ellenberg2013homologicalstabilityhurwitzspaces}. Indeed, for a general Hurwitz space if we let $c=c_1\sqcup\cdots \sqcup c_k$ with grading $\Z^k$, then by \cite[Theorem 3.1]{Wood_2021}, there is a bijection between the Hurwitz orbits with full global monodromy $\langle g_1,\ldots, g_n\rangle = G$ and elements of $U(G,c)$ as long as each grading is sufficiently large. For those components without full global monodromy $G'\subsetneq G$, one can carry out the same argument by restricting to $G'$ and looking at $U(G',c\cap G')$ instead.

Let us now specialize to the case when $G=H\rtimes \Z/2\Z$ with $c=G\setminus H$ where here $H=(H,0)\leq G$. We first enumerate all possible subgroups $G'\leq G$ that intersect $c$ nontrivially. Let $H'=G'\cap H$ and take a reflection $(h_0,1)\in G'\cap c\neq \emptyset$, then it is clear that we must have $G'=H'\sqcup (H'+h_0,1)$. We remark that this $G'$ is abstractly isomorphic to $H'\rtimes \Z/2\Z$. Conversely, any $G'$ of this form is a valid subgroup. Hence, the subgroups $G'\leq G$ can be parametrized by the cosets $H'+h_0$ of $H$. 

Let us fix such a subgroup $G'=H'\sqcup (H'+h_0,1)$, as well as the algebraic grading $n$. Since $c$ is non-splitting, the restriction $c'=c\cap G'=(H'+h_0,1)$ is a single conjugacy class. We also fix the boundary monodromy $\gamma\in G'$, which must be in $H'$ if $n$ is even and $(H'+h_0,1)$ otherwise. Then, \cite[Theorem 2.5]{Wood_2021} tells us that $U(G',c')=S_{c'}\times_{G'^{\text{ab}}}\Z$, and since $H_2(G',c')\cong \ker(S_{c'}\rightarrow G')$ in the proof of \cite[Lemma 2.2]{Wood_2021}, the $n$-th graded part of $U(G',c')$ with boundary monodromy $\gamma$ is an $H_2(G',c')$-torsor. On the other hand, the proof of \cite[Theorem 3.1]{Sawin_Wood_2026} (which follows \cite{Liu_Wood_Zureick-Brown_2024}) tells us that $H_2(G',c')\cong (\bigwedge^2 H')^{\Z/2\Z}$, and this is simply $\bigwedge^2 H'$ since $\Z/2\Z$ acts by inversion on $H'$ and hence the identity on $\bigwedge^2H'$. 

This tells us that in sufficiently large algebraic degree $n$, the components of $\Hur^{G,c}_n$ have a one-to-one correspondence to compatible choices of the global monodromy $G'$ (which corresponds to cosets $H'+h_0\subseteq H$), boundary monodromy $\gamma$, and an element of the $\bigwedge^2H'$-torsor. In our paper, we make the ``sufficiently large'' part precise, as in the final sentence of the following proposition.
\begin{proposition}\label{prop: param cohen lenstra R} For each $n\ge 1$, there is an injective map
$$\pi_0\left(\Hur^{G,c}_n\right) \hookrightarrow \left\{(H'+h_0,\gamma,\omega) \middle|\ H'+h_0\subseteq H\text{ coset},\  \gamma\in \begin{cases}(H',0) & \text{for } n \text{ even}\\
(H'+h_0,1) & \text{for } n \text{ odd}\end{cases},  \ \omega\in \bigwedge^2H'\right\}$$
which sends a connected component to $(H'+h_0,\gamma,\omega)$ such that the global monodromy is $H'\sqcup (H'+h_0,1)$, boundary monodromy is $\gamma$, and we explain the construction of $\omega$ later. 

Furthermore, this is an isomorphism for $n\ge 2r+1$ where $r=\rk(H)$.
\end{proposition}

The goal of this subsection is to prove this proposition. We tackle this using a more explicit method which is specific to the $G=H\rtimes \Z/2\Z$ case, in contrast to the above general discussion. First, we note that \cite{Lipnowski_Tsimerman_2019} constructs a similar exterior square product, which gives us a geometric interpretation of the invariant $\omega\in \bigwedge^2H'$. They considered the exterior square $\bigwedge^2J[l^n](1)$ on the $l^n$-torsion points of the Jacobian $J$ associated to the hyperelliptic curve with corresponding $n$ branch points, together with an element $\omega$ which represents the Weil pairing. This comes with an action of the braid group $B_n$ from the action of the branch points of the hyperelliptic curve, and this is essentially the geometric definition of the Burau representation. 

We will basically use this point of view, but we find that it is more convenient to carry this out explicitly in terms of matrices. Let us first recall the definition of the Burau representation. 
\begin{definition}
The \textit{unreduced Burau representation} of the braid group $B_n$ is the representation $B_n\rightarrow \GL_n$ given by sending the generators
$$\sigma_i \mapsto \left( \begin{array}{c|cc|c} I_{i-1} & 0 & 0 & 0 \\ \hline 0 & 1-t & t & 0 \\ 0 & 1 & 0 & 0 \\ \hline 0 & 0 & 0 & I_{n-i-1} \end{array} \right)$$
for $1\le i\le n-1$.

The \textit{reduced Burau representation} of the braid group $B_n$ is a quotient of the unreduced Burau representation by the invariant diagonal subspace. Explicitly, it is the representation $B_n\rightarrow \GL_{n-1}$ given explicitly by sending the generators
\begin{equation*}
\sigma_1 \mapsto \left( \begin{array}{cc|c}-t & 1 & 0 \\ 0 & 1 & 0 \\ \hline 0 & 0 & I_{n-3} \end{array} \right), \quad
\sigma_i \mapsto \left( \begin{array}{c|ccc|c} I_{i-2} & 0 & 0 & 0 & 0 \\ \hline 0 & 1 & 0 & 0 & 0 \\ 0 & t & -t & 1 & 0 \\ 0 & 0 & 0 & 1 & 0 \\ \hline 0 & 0 & 0 & 0 & I_{n-i-2} \end{array} \right), \quad 
\sigma_{n-1} \mapsto \left( \begin{array}{c|cc} I_{n-3} & 0 & 0 \\ \hline 0 & 1 & 0 \\ 0 & t & -t \end{array} \right),
\end{equation*}
where the middle map holds for the indices $2\le i\le n-2$.

We remark that we will only need the specialization of these representations to $t=-1$, which are also known as the integral Burau representations, and in this case we denote the unreduced and reduced representations $\beta_n$ and $\widetilde \beta_n$ respectively.
\end{definition}

Let us now see how these representations arise in our context. Consider a tuple $(g_1,\ldots, g_n)\in c^n$ where $g_i=(a_i,1)\in H\rtimes \Z/2\Z$, so we can think of these as tuples $(a_1,\ldots, a_n)\in H^n$. Then, a generator $\sigma_i$ of the braid group $B_n$ would send the piece $(g_i,g_{i+1})\mapsto(g_{i+1},g_i^{g_{i+1}})$. In terms of $a_i$ coordinates, since we have $(a_{i+1},1)^{-1}(a_i,1)(a_{i+1},1)=(2a_{i+1}-a_i,1)$, we have the action 
\begin{equation}\label{eqn: braid on a_i}
\sigma_i\colon (a_1,\ldots,a_i,a_{i+1},\ldots, a_n)\mapsto (a_1,\ldots,a_{i+1},2a_{i+1}-a_i,\ldots, a_n)
\end{equation}
instead. One can check that this is the \textit{inverse} representation of $\beta_n$, the \textit{unreduced} Burau representation at $t=-1$. In fact, if we conjugate by $\diag(1,-1,1,\cdots)$, i.e. adding a sign change to the coordinates $a_2,a_4,\ldots,$ this is the dual (contragredient) representation $\beta_n^*$ of $\beta_n$. 

We now state the global and boundary monodromy in terms of these coordinates $a_i$. Recall that the global monodromy $G'=H'\sqcup (H'+h_0,1)$ is determined by the coset $H'+h_0$, and it is not hard to see that this coset is the affine space $a_1+\langle a_2-a_1,\ldots, a_n-a_1\rangle$ spanned by integer linear combinations of $a_i$ such that the sum of coefficients is $1$. On the other hand, the boundary monodromy is $\gamma=(a_1,1)\cdots (a_n,1)=(a_1-a_2+\cdots +(-1)^{n-1}a_n, \mathbf 1_{n \text{ odd}})$. 

For the rest of the analysis, we will need to split into two cases depending on the parity of $n$. The reason for this, as will soon become apparent, is that we need to consider symplectic representations from the Burau representations, and since symplectic representations are necessarily of even dimension, their construction will necessarily depend on the parity of $n$. From a geometric perspective, this essentially boils down to the fact that hyperelliptic curves that have an odd or even number $n$ of branch points on affine space are treated differently, since an odd number of affine branch points occurs if and only if $\infty$ is also a branch point.

\subsubsection{$n$ is odd} Suppose that $n=2g+1$ for some $g\ge 0$. Let us then consider the new coordinates $(x_1,\ldots,x_{n-1},d)$ defined by $$x_1=a_2-a_1,\ \ldots,\  x_{n-1}=a_n-a_{n-1} ,\ d=a_1-a_2+\cdots -a_{n-1}+a_n$$ 
so that the boundary monodromy is $\gamma=(d,1)$, and $\langle x_1,\ldots, x_{n-1}\rangle = H'$. The braid action thus preserves $d$, and one can check that on the rest of the coordinates $(x_1,\ldots, x_{n-1})$, the braid group acts by $\widetilde \beta_n^*$, the dual (contragredient) representation of the \textit{reduced} Burau representation at $t=-1$.

Put in another way, let $V=\Z^{2g}$ be the lattice with basis $c_1,\ldots, c_{2g}$, equipped with the action $\widetilde \beta_n\colon B_n\rightarrow \GL(V)$. Then, one can identify the $2g$-tuples $(x_1,\ldots, x_{n-1})\in H'^{2g}$ which generate $H'$ with surjective homomorphisms $f\in \Surj(V,H')$ by identifying $(x_1,\ldots, x_{n-1})$ with the images $(f(c_1),\ldots, f(c_{2g}))$, and now $\Surj(V,H')$ comes with a natural dual action $\widetilde \beta_n^*$.

We now collect some facts on the reduced Burau representation at $t=-1$ in our $n$ odd case.
\begin{lemma}\label{lem: symplectic odd}
For odd $n=2g+1$, $\widetilde \beta_n\colon B_n \rightarrow \GL(V)$ preserves a unimodular symplectic form on $\Z^{2g}$, and with respect to a symplectic basis, the image contains the principal level-2 congruence subgroup $\Sp_{2g}(\Z)[2]$. In particular, for odd $m$, the reduction of the image modulo $m$ induces a surjection to $\Sp_{2g}(\Z/m\Z)$.
\end{lemma}
\begin{proof}
The first sentence is due to \cite{A’Campo_1979}, where on pg 318 he mentions that the monodromy representation of $B_n$ on the rank $2g$ homology of the hyperelliptic fiber preserves the nondegenerate intersection pairing. It is known that the reduced Burau representation at $t=-1$ is exactly this monodromy representation, see for example \cite[pg 267]{Brendle_Margalit_Putman_2014}. Then, \cite[Theorem 1(2)]{A’Campo_1979} tells us that the image contains $\Sp_{2g}(\Z)[2]$. For the last sentence, it is well-known that $\Sp_{2g}(\Z)\twoheadrightarrow\Sp_{2g}(\Z/2m\Z)$ is surjective (see for example \cite[Lemma 5.16]{Deligne_Mumford_1969}), and by the Chinese remainder theorem we have $\Sp_{2g}(\Z/2m\Z)\cong \Sp_{2g}(\F_2)\times \Sp_{2g}(\Z/m\Z)$. Hence, for any $M\in \Sp_{2g}(\Z/m\Z)$, the element $(I,M)$ can be lifted to $\Sp_{2g}(\Z)$ and the lift is necessarily in the congruence subgroup $\Sp_{2g}(\Z)[2]$.
\end{proof}

It is worth stating the explicit symplectic pairing on $V$ that $\widetilde\beta_n$ preserves, 
$$\langle c_i,c_{i+1}\rangle = -1, \quad \langle c_i,c_j\rangle =0 \text{ for }|i-j|\neq 1,$$
which can be found in \cite[pg 279]{Goldschmidt_1990} and is originally due to \cite{Squier_1984}. This 
gives us one possible symplectic basis $p_1,q_1,\ldots,p_g,q_g$ where
$$p_i=c_1+c_3+\cdots +c_{2i-1}, \quad q_i=-c_{2i}.$$
The symplectic form can then be explicitly written as an element in $\bigwedge^2V$, $$\omega_0=\sum_{i=1}^g p_i\wedge q_i.$$

Let us finally explain how to obtain the form $\omega\in \bigwedge^2H'$ from the tuple $(x_1,\ldots, x_{n-1})\in H'^{2g}$ above. Recall that such a tuple corresponds to some $f\in \Surj(V,H')$, which will induce a map $\bigwedge^2f \colon \bigwedge^2V\rightarrow \bigwedge^2H'$. Then, we define the form $\omega$ to be the image of $\omega_0$ under this map, i.e.
$$\omega \coloneqq \left(\bigwedge^2f\right)(\omega_0)=\sum_{i=1}^g f(p_i)\wedge f(q_i).$$
Since the form $\omega_0$ is $B_n$-invariant by Lemma \ref{lem: symplectic odd}, this descends to a map $\omega\colon \Surj(V,H')/B_n\rightarrow \bigwedge^2H'$.

We can now prove the proposition in the odd case.
\begin{proof}[Proof of Proposition \ref{prop: param cohen lenstra R} for odd $n$]
It suffices to prove that the $B_n$-orbits on $(x_1,\ldots, x_{n-1},d)$ are parametrized as in the right side of the map in the proposition. Let us fix the choice of $(H'+h_0,\gamma)$. This corresponds exactly to those orbits where $\langle x_1,\ldots, x_{n-1}\rangle =H'$ with $d$ coming from $\gamma=(d,1)$. Here, we note that the choice of coset $H'+h_0$ is determined by $H'$ and $d$ as we need $d\in H'+h_0$. 

It then suffices to show that $\omega\colon \Surj(V,H')/B_n\rightarrow \bigwedge^2H'$ is an injection for all $n$, and an isomorphism when $n\ge 2r+1$. Let $m$ be the exponent of $H'$. We may take our base ring to be $A=\Z/m\Z$ and replace $V$ with $A^{2g}$. By Lemma \ref{lem: symplectic odd}, the image of $\widetilde\beta_n$ modulo $m$ is the full group $\Sp_{2g}(A)$, so
$$
\Surj(V,H')/B_n\cong \Surj(V,H')/\Sp_{2g}(A).$$
By the Chinese remainder theorem, this orbit set decomposes as a product over the primary components of $H'$, as does $\bigwedge^2H'$, and $\omega$ respects these decompositions. Hence it suffices to consider the case when $H'$ is an $\ell$-group for an odd prime $\ell$. Write $H'=\Z/\ell^{\alpha_1}\Z\times\cdots\times\Z/\ell^{\alpha_k}\Z$ with $\alpha_1\ge\cdots\ge\alpha_k$, where $k=\rk(H')\le r=\rk(H)$. In this case, $A=\Z/\ell^{\alpha_1}\Z$ and $V=A^{2g}$.

Now, \cite[Theorem 2.14]{Michael_2006} tells us that $\omega$ gives a bijection
$$\frac{\Surj(V,H')}{\Sp_{2g}(A)} \xrightarrow\sim \left\{\eta\in\bigwedge^2H' \middle| \ \frac 1 2 \rk (M_{\overline\eta})\ge k-g  \right\},$$
where here $M_{\overline \eta}$ is the alternating $k\times k$ matrix over $\F_\ell$ corresponding to the residue $\overline\eta$ of $\eta$ modulo $\ell$ by the projection $\bigwedge^2H'\twoheadrightarrow \bigwedge^2 \F_\ell^k$. Note this is stated differently in Theorem 2.14 of loc. cit. -- we obtained this from unwinding the notation in Assumption 2.10 of loc. cit. and the definition of the index above Lemma 2.2 of loc. cit.. In fact, the condition $\frac 1 2 \rk (M_{\overline\eta})\ge k-g  $ is equivalent to the condition that the maximal isotropic subspace has dimension $\le g$. 

This completes the proof because the above map is clearly an injection into $\bigwedge^2H'$ for all $n$, and for $n\ge 2r+1$, we have $g\ge r\ge k$, so the condition $\frac 1 2\rk(M_{\overline\eta})\ge 0 \ge  k-g$ is vacuously true. \end{proof}

We remark that this proof also tells us the precise image of the map in the proposition, and the same will be true for the even case as well.
\subsubsection{$n$ is even} Suppose that $n=2g+2$ for some $g\ge 0$. We will carry out an analogous argument as in the odd case. It will be more convenient here to work directly with the \textit{unreduced} Burau representation $\beta_n$ at $t=-1$ instead of the reduced version $\widetilde\beta_n$, as now $\beta_n$ has even dimension $2g+2$ and we will see that it is also symplectic.

Let us fix the global monodromy $G'=H'\sqcup (H'+h_0,1)$ by choosing the coset $H'+h_0$. To define $\omega\in \bigwedge^2H'$, we will need to choose an explicit representative $h_0$ of the coset $H'+h_0$; such a choice is unavoidable because the actual invariant lives in a $\bigwedge^2H'$-torsor as discussed above. Then, consider the new coordinates $(b_1,\ldots, b_n)$ where $b_i=a_i-h_0\in H'$. Note that $B_n$ still acts on this via $\beta_n$ via Equation \eqref{eqn: braid on a_i}.

Just like the odd case, it is convenient to regard these tuples as surjective homomorphisms $f\in \Surj(V,H')$ where $V=\Z^{2g+2}$ is the lattice with basis $c_i$ which comes with the action of $B_n$ via $\beta_n$. More precisely, let $f\in \Surj(V,H')$ correspond to the tuple $(b_1,\ldots, b_n)=(f(c_1),-f(c_2),\ldots, -f(c_n))$. By the discussion after Equation \eqref{eqn: braid on a_i}, we see the alternating signs are required for the action on $\Surj(V,H')$ to correspond to the action on the tuples.

Before we continue our discussion, we first state the analogous lemma to Lemma \ref{lem: symplectic odd} for $\beta_n$.
\begin{lemma}\label{lem: symplectic even}
For even $n=2g+2$, $\beta_n\colon B_n\rightarrow \GL(V)$ preserves the unimodular symplectic form given by $\langle c_i,c_j\rangle =(-1)^{j-i}$ for $i<j$, and fixes the vector $u=c_1+\cdots +c_n$. Furthermore, the image contains $\Sp(V)_u[2] \coloneqq \Sp(V)_u\cap \Sp(V)[2]$ which is the principal level-$2$ congruence subgroup of the stabilizer $\Sp(V)_u$ of $u$. In particular, for odd $m$, the reduction of the image modulo $m$ induces a surjection to $\Sp(V/mV)_{\overline u}$.
\end{lemma}
\begin{proof}
We can verify the first sentence from a simple direct computation. The second statement is proven in \cite[Theorem 3.3]{Brendle_Margalit_2015} for $g\ge 2$, extended to $g\ge 1$ in \cite[Proposition 3.3]{Bloomquist_Patzt_Scherich_2024} and the remaining case of $g=0$ follows from a direct computation. 

For the final sentence, we first show that reduction gives a surjection $\Sp(V)_u\twoheadrightarrow \Sp(V/mV)_{\overline u}$. To do this, pick a $v\in V$ where $\langle u,v\rangle = 1$, and let $W=u^\perp\cap v^\perp$ be a symplectic sublattice of rank $2g$. Then, $\Sp(V)_u$ is an extension of $\Sp_{2g}(W)$ by a free unipotent group. The same holds true modulo $m$, and reduction is surjective on both factors, so the reduction $\Sp(V)_u\twoheadrightarrow \Sp(V/mV)_{\overline u}$ is indeed a surjection. Now, using the Chinese remainder theorem gives us $\Sp(V/2mV)_{\overline u}\cong \Sp(V/2V)_{\overline u} \times \Sp(V/mV)_{\overline u}$, and lifting $(I,M)$ finishes the proof in the same way as in Lemma \ref{lem: symplectic odd}.
\end{proof}
An explicit symplectic basis $p_1,q_1,\ldots, p_{g+1},q_{g+1}$ is given by
$$p_i = c_{2i-1}+c_{2i}, \quad q_i=c_1+\cdots +c_{2i-1}.$$
Then, one can define the invariant $\omega=\sum_{i=1}^{g+1}f(p_i)\wedge f(q_i)$ in the same way.

Going back to the original discussion, we realize that the fixed vector $u=c_1+\cdots+c_n$ plays an important role in our context. Given some $f\in \Surj(V,H')$ corresponding to a tuple $(b_1,\ldots, b_n)$ above, the boundary monodromy is simply given by $f(u)=\gamma$. Furthermore, it is not difficult to check that $u^\perp$ is spanned by $c_1+c_2,c_2+c_3,\ldots, c_{n-1}+c_n$, so $f(u^\perp)=\langle b_2-b_1,\ldots, b_n-b_{n-1}\rangle =\langle a_2-a_1,\ldots, a_n-a_{n-1}\rangle = H'$. 

Hence, we need the condition $f(u^\perp)=H'$ which is stronger than surjectivity of $f$. To deal with this, let $\ev_u=\langle u,-\rangle\colon V\rightarrow \Z$ and consider the map $\phi=(f,\ev_u)\colon V\rightarrow H'\times \Z$ instead. We claim that $\phi$ is surjective if and only if $f(u^\perp)=H'$. Indeed, since $\ker(\ev_u)=u^\perp$, we have that $f(u^\perp)=H'$ if and only if $\im(\phi)\supseteq H'\times \{0\}$ which happens if and only if $\phi$ is surjective (since we can translate by $v$ satisfying $\ev_u(v)=\langle u,v\rangle = 1$). Thus, the space of tuples we are dealing with is
$$\{\phi \in \Surj(V,H'\times \Z)\mid \pr_2 \circ \phi = \ev_u\} \subseteq \Surj(V,H'\times \Z).$$

We can now tackle the proof of the proposition in the even case.
\begin{proof}[Proof of Proposition \ref{prop: param cohen lenstra R} for even $n$] Recall that in the above setup, after fixing the coset $H'+h_0$ and making a choice of $h_0$, we are left with the $B_n$-orbits of $\{\phi \in \Surj(V,H'\times \Z)\mid \phi=(f,\ev_u)\}$. Again assume that $H'$ is an $\ell$-group and write $H'=\Z/\ell^{\alpha_1}\Z\times \cdots \times \Z/\ell^{\alpha_k}\Z$ with $\alpha_1\ge \ldots \geq \alpha_k$ with $k\le r$. This allows us to take the base ring $A=\Z/\ell^{\alpha_1}\Z$ and to replace $V$ with $A^{2g+2}$ and $\Z$ with $A$ above. Furthermore, by Lemma \ref{lem: symplectic even}, $\beta_n$ generates the stabilizer $\Sp(V)_u$ (here we also use $u$ to denote the reduction of $u$ modulo $\ell^{\alpha_1}$), so we are concerned with the orbits
$$\frac{\{\phi \in \Surj(V,H'\times A)\mid \pr_2\circ \phi = \ev_u\}}{\Sp(V)_u}.$$

We claim that there is an isomorphism of orbits
$$\frac{\{\phi \in \Surj(V,H'\times A)\mid \phi=(f,\ev_u)\}}{\Sp(V)_u}\xrightarrow\sim \frac{\Surj(V,H'\times A)}{\Sp(V)}.$$
Indeed, consider the projection $\pr_2\colon \Surj(V,H'\times A) \twoheadrightarrow \Surj(V,A)$. Note that $\Sp(V)$ acts transitively on the target $\Surj(V,A)$ as we can view $\Surj(V,A)\subseteq \Hom(V,A) \cong V^\vee$ as the primitive covectors in $V^\vee$ and it is clear $\Sp(V)$ acts transitively on these primitive covectors. Now, the left side is simply the fiber $\pr_2^{-1}(\ev_u)$ quotiented by the stabilizer $\Sp(V)_u\subseteq \Sp(V)$ of $\ev_u\in \Surj(V,A)$, so the isomorphism must be true by transitivity.

Finally, we apply \cite[Theorem 2.14]{Michael_2006} again, which says that pushing forward the symplectic form $\omega_0$ on $V$ to $H'\times A$ induces a bijection 
$$\frac{\Surj(V,H'\times A)}{\Sp(V)} \xrightarrow\sim \left\{\eta\in\bigwedge^2 (H'\times A) \middle| \ \frac 1 2 \rk (M_{\overline\eta})\ge (k+1)-(g+1)  \right\}.$$
Here, note that since $\bigwedge^2 (H'\times A)\cong \bigwedge^2H' \times H'$, we can view $\eta=(\omega,\gamma)$ where $\omega\in \bigwedge^2H'$ and $\gamma \in H'$. It is not so difficult to check that since the map from $V\rightarrow A$ is given by $\ev_u$ and we are pushing forward the symplectic pairing, the $\gamma$ here must be $f(u)$ which is indeed the boundary monodromy. 

Hence, for $n\ge 2r+2$, we have $g\ge r\ge k$ so the condition on the right is vacuous and the orbits corresponding to the coset $H'+h_0$ are parametrized by $\omega\in \bigwedge^2H'$ and the boundary monodromy $\gamma \in H'$. Likewise, injectivity for all $n$ is clear.
\end{proof}

\subsection{Stability slope and arithmetic applications}\label{sec: cohen lenstra slope} We start by proving Theorem \ref{thm: cohen lenstra stab R}. Using the parametrization in Proposition \ref{prop: param cohen lenstra R}, we first prove the following lemma, from which we will deduce the theorem.
\begin{lemma}\label{lem: bijection cohen lenstra}
Fix integers $D\ge 1$, $n\ge 2r+1$, and a subgroup $G'=H'\sqcup (H'+h_0,1)$ of $G$. For any $g\in c\cap G'=(H'+h_0,1)$, we have that right (or left) multiplication by $g^{2D}$ induces a bijection between the connected components of $\Hur^{G,c}_n$ with global monodromy $G'$ and the connected components of $\Hur^{G,c}_{n+2D}$ with global monodromy $G'$. Furthermore, this bijection is the same for all choices of $g$.
\end{lemma}
\begin{proof}
Right multiplication by $g^{2D}$ is given by adding $2D$ copies of $g$ on the right of the original tuple, changing $(g_1,\ldots, g_n)$ to $(g_1,\ldots, g_n,g,\ldots, g)$. Since $g^2=1$, this does not change the boundary monodromy $\gamma$. We will also check that this does not change the invariant $\omega$, which will finish the proof by Proposition \ref{prop: param cohen lenstra R}. For the odd case where $n=2g+1$, then adding $2D$ copies of $g$ changes $\omega$ by $\sum_{i=g+1}^{g+D} f(p_i)\wedge f(q_i)$. However, we have $q_i=-c_{2i}$ so $f(q_i)=-x_{2i}=a_{2i}-a_{2i+1}$, and for $g+1\le i\le g+D$, we have $g_{2i}=g_{2i+1}=g$ so this is zero. Similarly, for the even case $n=2g+2$, we need to show that $\sum_{i=g+2}^{g+D+1}f(p_i)\wedge f(q_i)=0$, but in this range we have $f(p_i)=f(c_{2i-1}+c_{2i})=b_{2i-1}-b_{2i}=a_{2i-1}-a_{2i}=0$ for the same reason.
\end{proof}
\begin{proof}[Proof of Theorem \ref{thm: cohen lenstra stab R}]
We will view $U_D\colon H_0(\Hur^{G,c}_n,\Q)\rightarrow H_0(\Hur^{G,c}_{n+2D},\Q)$ as a matrix with basis indexed by connected components. Let us group the basis vectors according to global monodromy $G'$ of the connected component, and arrange these groups according to the size of $G'$. Then, $U_D$ is a block upper triangular matrix since multiplying by $g^{2D}$ can only make the global monodromy $G'$ larger. Hence, to show that this matrix is invertible, it suffices to check each block on the diagonal is invertible. Consider such a block which sends basis vectors with global monodromy $G'$ to $G'$. It is clear that this block is given by multiplication by $\sum_{g\in c\cap G'}g^{2D}$ since for those $g\not\in G'$ multiplication by $g^{2D}$ makes the global monodromy larger and is thus not within this block on the diagonal. However, by Lemma \ref{lem: bijection cohen lenstra}, since each $g^{2D}$ in the sum gives the same bijection between the basis vectors in the block, adding these up is necessarily an invertible map.
\end{proof}

As mentioned previously, from this we obtain Theorem \ref{thm: cohen lenstra stab slope} which gives a stability slope of $\frac{1}{2r+3}$. We will next prove Theorem \ref{thm: cohen-lenstra H} by just following the proof of \cite{landesman2025cohenlenstramomentsfunctionfields} with this new input.
\begin{proof}[Proof of Theorem \ref{thm: cohen-lenstra H}]
Our theorem will follow from \cite[Theorem 5.4.1]{landesman2025cohenlenstramomentsfunctionfields} once we show that the constant $I=2r+3$. Tracing through the proof, we see that the constants $I$ and $J$ there come from the homological stability statement of \cite[Theorem 1.4.7]{landesman2025cohenlenstramomentsfunctionfields} where the homology stabilizes for $n>iI+J$. However, by Theorem \ref{thm: cohen lenstra stab slope}, one can pick $I=2r+3$ (and $J=2r+1$), finishing the proof. 
\end{proof}

\subsection{Relating $H_1$ stability to Siegel threefolds}\label{sec: H1 conjectures}

In this section, we give a geometric heuristic for Conjecture \ref{conj: H1} from the perspective of Siegel modular threefolds. For convenience, instead of tackling the general case, we specialize to the case where \textbf{$k$ is prime}. We believe that a similar heuristic should also work in the general case. Let $c$ be the conjugacy class of reflections in $G=D_k$. We are concerned with the case of $n=5$ which is where our conjectured stability happens for $3\nmid k$. 

Let us first review some background on Siegel threefolds. Throughout this section we will only be working over $\C$, i.e. we will view the following constructions as complex manifolds/orbifolds. Recall that the Siegel upper half-space of genus $2$ is defined by $$\mathfrak H_2\coloneqq \left\{\tau\in M_2(\C)\mid \tau=\tau^t, \text{Im}(\tau) \text{ is positive definite}\right\},$$
and we define the quotient $\mathcal A_2\coloneq \Sp_4(\Z)\backslash \mathfrak H_2$, which parametrizes principally polarized abelian surfaces. Furthermore, for a congruence subgroup $\Gamma \subseteq \Sp_4(\Z)$, one can also define the quotient $\Gamma\backslash \mathfrak H_2$ which is a finite étale cover of $\mathcal A_2$ which parametrizes principally polarized abelian surfaces with some level structure. 

There is a Torelli map from the moduli space of genus $2$ curves $\mathcal M_2\hookrightarrow \mathcal A_2$ sending a curve to its Jacobian. This is an open immersion, and then we denote the complement as $\mathcal A_{1,1}\coloneqq \mathcal A_2\setminus \mathcal M_2$. Since every principally polarized abelian surface is either a Jacobian or a product of elliptic curves, we have that $\mathcal A_{1,1}\cong \mathcal M_{1,1}\times \mathcal M_{1,1}/S_2$ where we quotient by the swapping action. This also tells us that $\mathcal A_{1,1}$ is a connected divisor in $\mathcal A_2$. Furthermore, since every genus $2$ curve is hyperelliptic, we have $\mathcal H_2\cong \mathcal M_2$.

We now relate the Hurwitz space $\Hur^{D_k,c}_5$ for $n=5$ to these Siegel modular threefolds. As we assumed $k$ is prime for simplicity, from the parametrization in Proposition \ref{prop: param cohen lenstra R}, since $n=5\ge 3=2r+1$, we have $2k$ connected components of $\Hur^{D_k,c}_5$ in total. Here, $k$ of these have global monodromy $\cong \Z/2\Z$ and are thus isomorphic to configuration space $\Conf_5$, for which the homology is known, and the other $k$ components have full global monodromy (these are also isomorphic to one another by conjugation). It thus suffices to understand a connected component $Z\subseteq \Hur^{D_k,c}_5$ with full global monodromy.

Recall that $\Hur^{D_k,c}_5$ parametrizes covers $X\rightarrow C\rightarrow \mathbb P^1$ where $X\rightarrow C$ is an unramified $\Z/k\Z$ Galois cover and $C\rightarrow \mathbb P^1$ is a genus $2$ hyperelliptic cover with $5$ branch points over $\A^1$ and one branch point over $\infty$. The framing marks a point over the chosen section of the order-$2$ root stack at $\infty$, as in Definition \ref{def: algebraic hurwitz space}. Let $w$ be the point of $C$ lying above $\infty$. Also recall that $X\rightarrow C$ corresponds to a $k$-torsion element of the class group, which is an element $v\in \Pic(C)[k]= \Pic^0(C)[k]= \Jac(C)[k]$ of the Jacobian $\Jac(C)$. On the connected component $Z$, we also have that $v\neq 0$. 

Let us first consider the moduli space of maps $C\rightarrow \mathbb P^1$. Taking the quotient by the group of affine linear transformations of $\A^1$
$$\Aff=\{ax+b\mid a\in \G_m, b\in \A^1\},$$ 
we get rid of the choice of the $5$ branch points over $\A^1$, so this now parametrizes pairs $(C,w)$ consisting of a hyperelliptic curve with a choice of Weierstrass point. 

However, the action of $\Aff$ does not lift to $Z$, as a generator of $\pi_1(\Aff)$ is a full twist $\Delta_5^2$ in $B_5$, which acts via the loop $\gamma(t)\colon z\mapsto e^{2\pi it}z$ for $0\le t\le 1$. This acts on tuples by the boundary monodromy, and since $n=5$ is odd, the boundary monodromy is a reflection, so on those components that have full global monodromy this action is nontrivial. Hence, we need to consider the double cover $\widetilde{\Aff} \rightarrow \Aff$, which will now act on $Z$, and the quotient $Z/\widetilde{\Aff}$ of a fixed full monodromy component parametrizes triplets $(C,w,v)$ where $0\neq v\in \Jac(C)[k]$.

Using the open immersion $\mathcal H_2\cong\mathcal M_2\hookrightarrow \mathcal A_2$, we see that this is the data of a Jacobian with extra structure. Hence, the connected component $Z/\widetilde{\Aff}$ must be isomorphic to an open subset of a quotient $\Gamma_k\backslash \mathfrak H_2$ for the appropriate congruence group $\Gamma_k$. More precisely, the open subset is the preimage $\pi^{-1}(\mathcal M_2)$ along the projection $\pi\colon \Gamma_k\backslash \mathfrak H_2\rightarrow \mathcal A_2$, and the congruence subgroup $\Gamma_k\subseteq \Sp_4(\Z)$ is the stabilizer of the Weierstrass point $w$ and a primitive vector $0\neq v\in \Jac(C)[k]$. 

We can now write $\Gamma_k$ down explicitly. It is known (see for example \cite[Section 2]{Bergström_Faber_vdG_2008}) that $\Sp_4(\Z)$ acts on the six Weierstrass points via the quotient $\Sp_4(\Z/2\Z)\cong S_6$, and the stabilizer of a single Weierstrass point is $S_5$. Hence, we have
$$\Gamma_k\coloneqq \left\{\gamma \in \Sp_4(\Z) \mid \gamma \text{ (mod }2\text{)} \in S_5, \ \gamma v\equiv v \pmod k\right\}.$$

We summarize the above discussion in the following lemma for clarity.
\begin{lemma}\label{lem: summarize Hur_5} Let $k$ be an odd prime. Then, $\Hur^{D_k,c}_5$ has $2k$ components, $k$ of them are isomorphic to $\Conf_5$, and the other $k$ are isomorphic to $Z$ where the quotient $Z/\widetilde{\Aff}= \pi^{-1}(\mathcal M_2)$ is the preimage of $\mathcal M_2\subseteq \mathcal A_2$ along the projection $\pi\colon \Gamma_k\backslash \mathfrak H_2\rightarrow \mathcal A_2$.
\end{lemma}

One could then expect that the cohomology of $\Hur^{D_k,c}_5$ can be computed using properties of $\Gamma_k\backslash \mathfrak H_2$ and $\mathcal A_2$. This turns out to be true, as we state in the next proposition. 
\begin{proposition}\label{prop: H1 siegel}
Let $k$ be an odd prime. We have
$$\dim_\Q\left( H_1(\Hur^{D_k,c}_5,\Q)\right) = 2k+kr_k,$$
where $r_k$ is defined as follows. The divisor $D\coloneqq \pi^{-1}(\mathcal A_{1,1})\subseteq \Gamma_k\backslash \mathfrak H_2$ has three connected components $D=D_1\sqcup D_2\sqcup D_3$, and we define
$$r_k\coloneqq \dim_\Q \ker\left( \Q^3  \xrightarrow{e_i\mapsto [D_i]}\Pic_\Q(\Gamma_k\backslash \mathfrak H_2)\right)$$
where $\Pic_\Q(\Gamma_k\backslash \mathfrak H_2)\coloneqq \Pic(\Gamma_k\backslash \mathfrak H_2) \otimes \Q$ is the rational Picard group.
\end{proposition}

We remark that the stable dimension of both $H_1(\Hur^{D_k,c}_5,\Q)$ and $H_0(\Hur^{D_k,c}_5,\Q)$ is $2k$, since the homology of each connected component stabilizes to that of configuration space by Theorem \ref{thm: stable hom conf all directions}. This directly implies the following corollary on the stability of $H_1$.
\begin{corollary}\label{cor: H1 Siegel}
Let $k>3$ be a prime. The following are equivalent.
\begin{enumerate}[(\alph*)]
\item Conjecture \ref{conj: H1} is true for $G=D_k$, i.e. $H_1(\Hur^{D_k,c}_n,\Q)$ stabilizes for $n\ge 5$.
\item $H_1(\Hur^{D_k,c}_5,\Q)$ has dimension $2k$.
\item $[D_1],[D_2],[D_3]$ are linearly independent in $\Pic_\Q(\Gamma_k\backslash \mathfrak H_2)$.
\end{enumerate}
\end{corollary}
Indeed, the equivalence of $(b)$ and $(c)$ follows from the proposition, while the equivalence of $(a)$ and $(b)$ follows from the usual argument after Corollary \ref{cor: stab from R}. 

Before we prove Proposition \ref{prop: H1 siegel}, we discuss condition $(c)$ in the next few paragraphs. Unfortunately, we were not able to prove $(c)$ for any $k>3$, but this reformulation places the problem in the well-developed theory of Picard groups of Siegel modular threefolds with close relations to theta lifts and special cycles, and we believe that further investigation could potentially lead to a proof. 

We give a brief summary of the surrounding literature in this context. \cite{Weissauer_1988,Weissauer_1992} gave an automorphic description of the degree $2$ cohomology of such arithmetic quotients $\Gamma\backslash \mathfrak H_2$ for any finite-index subgroup $\Gamma\subseteq \Sp_4(\Z)$. As summarized in \cite[Theorem 4.1,4.3]{He_Hoffman_2012}, Weissauer showed that $H^2(\Gamma \backslash \mathfrak H_2,\C)$ has a pure Hodge structure of weight $2$, and the $(1,1)$ part $H^{1,1}(\Gamma\backslash \mathfrak H_2)=\Pic_\C(\Gamma\backslash \mathfrak H_2)\coloneqq \Pic(\Gamma\backslash \mathfrak H_2)\otimes \C$ is simply the Picard group which is the relevant part for our context. More importantly, we have a decomposition
$\Pic_\C(\Gamma\backslash \mathfrak H_2) = \C  [\mathcal L]\oplus \Pic_\C(\Gamma\backslash \mathfrak H_2)_0$ where $[\mathcal L]$ is the Lefschetz class corresponding to the trivial automorphic representation of $\Sp_4$, and the second term  $\Pic_\C(\Gamma\backslash \mathfrak H_2)_0$ is the non-Lefschetz automorphic contribution which is generated by theta lifts of weight $5/2$ cusp forms. Furthermore, $\Pic_\Q(\Gamma\backslash \mathfrak H_2)$ is spanned by special cycles \cite[Theorem 7.1]{He_Hoffman_2012}, and $\mathcal A_{1,1}$ is one such special cycle. In summary it seems plausible that there might be a useful automorphic or special cycle interpretation of $(c)$.

Although we were not able to prove $(c)$, we can present a possible heuristic for why it should be true for primes $k>3$ (note that it is false for $k=3$), which justifies Conjecture \ref{conj: H1} for these $D_k$. As $k$ gets larger, the degree of the map $\Gamma_k\backslash \mathfrak H_2\rightarrow \mathcal A_2$ grows, which suggests that the space $\Gamma_k\backslash \mathfrak H_2$ should become more complicated. Consequently, the rational Picard group $\Pic_\Q(\Gamma_k\backslash \mathfrak H_2)$ should also grow. This is empirically supported by \cite[Corollary 10]{Bruinier_2002}, which shows that in the case of principal Siegel congruence groups the dimension of the Picard group grows at least linearly. This makes it less likely for there to be a relation between the divisors $[D_1],[D_2],[D_3]$. The case of $k=3$ is small, so it is plausible that there is a coincidental linear relationship between $[D_1],[D_2],[D_3]$. However, as we checked for $k=5,7,11,13,17,19,23$, this is not the case, and our heuristic suggests that this pattern will continue for larger primes $k$ as well. 

We remark that in the case when $k$ is not prime, the preimage $\pi^{-1}(\mathcal A_{1,1})$ will split into more than three connected components, as we can see later from modifying the proof of Claim \ref{claim 2}. One can prove an analogue of Proposition \ref{prop: H1 siegel} and Corollary \ref{cor: H1 Siegel} with the correct stability dimension and number of connected components. We would again need to show that the connected components are linearly independent in the Picard group to prove Conjecture \ref{conj: H1} for that $D_k$. We believe that the heuristic should still apply, albeit weaker in this context, as the ``complexity'' of $\Gamma_k\backslash \mathfrak H_2$ and its Picard group should grow faster than the number of connected components, making a linear dependence less likely as $k$ grows larger.

The rest of this section is devoted to proving the main proposition.
\begin{proof}[Proof of Proposition \ref{prop: H1 siegel}] We will default to using $\Q$ coefficients in this proof. By Lemma \ref{lem: summarize Hur_5}, since we know the cohomology of $\Conf_5$, it suffices to show that $H_1(Z)$ has dimension $1+r_k$. The first step is to relate the cohomology of $Z$ to that of the quotient $Z/\widetilde{\Aff}$. Applying the Serre spectral sequence to the fibration $\widetilde{\Aff}\rightarrow Z\rightarrow Z/\widetilde{\Aff}$ yields $$E_2^{p,q}= H^p(Z/\widetilde{\Aff},H^q(\widetilde{\Aff})) \Rightarrow H^{p+q}(Z).$$ 
The corresponding low-degree exact sequence (with $\Q$ coefficients) is
$$0\rightarrow H^1(Z/\widetilde{\Aff},H^0(\widetilde{\Aff}))\xrightarrow{\text{inf}} H^1(Z)\xrightarrow{\text{res}}H^0(Z/\widetilde{\Aff},H^1(\widetilde{\Aff}))\xrightarrow{d_2} H^2(Z/\widetilde{\Aff})\rightarrow\cdots.$$
Note that both $H^0(\widetilde{\Aff},\Q)$ and $H^1(\widetilde{\Aff},\Q)$ are isomorphic to $\Q$ since $\widetilde{\Aff}\simeq \C^\times \simeq S^1$. Let us first show that the third term $H^0(Z/\widetilde{\Aff},H^1(\widetilde{\Aff}))\cong H^1(\widetilde{\Aff})\cong \Q$. This is because $Z\rightarrow Z/\widetilde \Aff$ is a principal $\widetilde \Aff$ bundle, so transporting around a loop in the base $Z/\widetilde \Aff$ is given by translation by an element of $\widetilde \Aff$. Since $\widetilde \Aff$ is connected, these translations on the fibers are homotopic to the identity and hence act trivially on the cohomology of the fibers $H^1(\widetilde \Aff,\Q)\cong \Q$.

We now show that the restriction map $\text{res}\colon H^1(Z)\rightarrow H^0(Z/\widetilde{\Aff},H^1(\widetilde{\Aff}))\cong H^1(\widetilde{\Aff})$ is surjective. Unwinding definitions tells us that this map is given by taking a fiber $\widetilde{\Aff}\rightarrow Z$ and considering the pullback $H^1(Z)\rightarrow H^1(\widetilde{\Aff})$. Hence, it suffices to show that the dual pushforward of homology $H_1(\widetilde{\Aff})\rightarrow H_1(Z)$ is injective. Indeed, let us consider $Z\rightarrow \Conf_5 \xrightarrow\Delta \C^\times$ given by taking the branch points $(z_1,\ldots, z_5)$ and then the discriminant map $\Delta(z_1,\ldots, z_5)=\prod_{i<j}(z_i-z_j)^2$. This gives the map $\Q \cong H_1(\widetilde{\Aff})\rightarrow H_1(Z)\rightarrow H_1(\C^\times)\cong \Q$, and it suffices to show that this is an isomorphism. However, this is simply given by multiplication by $40$, as a generator of $H_1(\widetilde{\Aff})$ is given by the loop $\gamma(t)\colon z\mapsto e^{4\pi it}z$ for $0\le t\le 1$, so the pushforward acts on the branch points of $z$ via $(z_1,\ldots, z_5)\mapsto (e^{4\pi it}z_1,\ldots, e^{4\pi it}z_5)$, and thus the discriminant $\Delta(z_1,\ldots, z_5)$ winds around the origin $40=4\binom 5 2$ times.

Hence, it suffices to show that the first term $H^1(Z/\widetilde{\Aff})$ has dimension $r_k$. Recall from Lemma \ref{lem: summarize Hur_5} that $Z/\widetilde{\Aff}$ is the preimage $\pi^{-1}(\mathcal M_2)$ along $\pi\colon \Gamma_k\backslash \mathfrak H_2 \rightarrow \mathcal A_2$. Alternatively, it is the complement of $\pi^{-1}(\mathcal A_{1,1})$ in $\Gamma_k\backslash \mathfrak H_2$. We will defer the proof of the following claim later, but for now, this lets us write $\pi^{-1}(\mathcal A_{1,1})=D=D_1\sqcup D_2\sqcup D_3$ as the union of three connected components. Note that this is the only place we use the fact that $k$ is an odd prime.
\begin{claim}\label{claim 2}
For an odd prime $k$, the preimage $\pi^{-1}(\mathcal A_{1,1})$ has three connected components.
\end{claim}

Continuing with our proof, we use the Gysin sequence associated to $Z/\widetilde{\Aff} \hookrightarrow \Gamma_k\backslash \mathfrak H_2\hookleftarrow D$
$$\cdots \rightarrow H^1(\Gamma_k\backslash\mathfrak H_2)\rightarrow H^1(Z/\widetilde{\Aff})\rightarrow H^0(D)(-1)\xrightarrow{i_*} H^2(\Gamma_k\backslash\mathfrak H_2)\rightarrow \dots.$$
We show that the first term $H^1(\Gamma_k\backslash \mathfrak H_2)=0$. Since $\Gamma_k$ is a finite-index subgroup of $\Sp_4(\Z)$, it has Kazhdan's property $(T)$, as $\Sp_4(\R)$ has property $(T)$ and property $(T)$ passes to lattices and their finite-index subgroups, see \cite[Theorem 1.6.1, 1.7.1]{Bekka_de_la_Harpe_Valette_2008}. Then, by Corollary 1.7.2 of loc. cit., the abelianization $\Gamma_k^{\text{ab}}$ is finite, and thus $H_1(\Gamma_k,\Q)=\Gamma_k^{\text{ab}}\otimes \Q =0$, so we have the desired result.

Hence, we have $H^1(Z/\widetilde{\Aff})\cong \ker(i_*)$, but as we discussed $H^2(\Gamma_k\backslash \mathfrak H_2,\C)$ is pure of weight $2$ and $H^{1,1}(\Gamma_k\backslash \mathfrak H_2)=\Pic_\C(\Gamma_k\backslash \mathfrak H_2)$, so the map $i_*$ (over $\Q$) is simply the map described in the statement of the proposition.

We are left with the proof of the claim.
\begin{proof}[Proof of Claim \ref{claim 2}] We first recall that the data of a marked Weierstrass point over a genus $2$ curve corresponds to a choice of one of the six odd theta characteristics on the Jacobian, see for example \cite[Section 6]{Theta_II}. Now, let us take a generic point $E_1\times E_2\in \mathcal A_{1,1}$ which is a product of elliptic curves equipped with the product polarization. The fiber of $\pi$ over this point is finite, with points parametrized by the data of an odd theta characteristic $\theta$ on $E_1\times E_2$ and a nonzero vector $0\neq v\in  (E_1\times E_2)[k]$. Also recall that 
$$\mathcal A_{1,1}\cong \mathcal M_{1,1}\times \mathcal M_{1,1}/S_2\cong (\mathfrak H_1\times \mathfrak H_1)/\left((\SL_2(\Z)\times \SL_2(\Z))\rtimes S_2\right)$$
where $\mathfrak H_1\subseteq \C$ is the usual upper half plane, $\mathfrak H_1\times \mathfrak H_1$ lies in $\mathfrak H_2$ via the diagonal embedding, and each copy of $\SL_2(\Z)$ acts on the respective $\mathfrak H_1$ with $S_2$ swapping the two $\mathfrak H_1$. Then, the number of components of $\pi^{-1}(\mathcal A_{1,1})$ will simply be the number of orbits of $(\theta,v)$ under the action of $(\SL_2(\Z)\times\SL_2(\Z))\rtimes S_2$.

It is easy to see that $v\in (E_1\times E_2)[k]$ decomposes as $(v_1,v_2)$ with $v_i\in E_i[k]$ where they are not both zero. Then, the $i$-th copy of $\SL_2(\Z)$ acts on $E_i[k]\cong \F_k^2$ via $\SL_2(\Z/k\Z)$. On the other hand, the odd theta characteristic $\theta$ decomposes according to the two elliptic factors, and by parity (see the introduction of \cite{Mumford_1971} for more background), exactly one factor has odd theta characteristic and the other one has even theta characteristic. Let us denote this decomposition as $\theta=(\theta_o,\theta_e)$ for the odd and even characteristics, coming from the elliptic curves $E_o$ and $E_e$ respectively.

The theta characteristics of an elliptic curve $E$ are by definition the line bundles which square to the canonical divisor which in this case is the trivial bundle, and these correspond to the four $2$-torsion points $\Pic^0(E)[2]\cong E[2]$. Hence, $\SL_2(\Z)$ acts on $E[2]$ via $\SL_2(\Z/2\Z)$. Furthermore, the three nontrivial theta characteristics are even and the trivial one is odd. Hence, going back to our situation, one copy of $\SL_2(\Z)$ acts trivially on $\theta_o\in E_o[2]$, and the other copy acts on the three possible even theta characteristics $\theta_e\in E_e[2]$ via the permutation group $\SL_2(\Z/2\Z)\cong S_3$.

We can now put everything together to show that there are three orbits of $(\theta,v)$. By swapping if necessary, we can assume that $E_o=E_1$ and $E_e=E_2$. Then, we claim that the three orbits are given by the cases $(v_1\neq0, v_2\neq 0)$, $(v_1=0, v_2\neq 0)$ and $(v_1\neq0, v_2= 0)$. Indeed, this follows because $\SL_2(\Z/2\Z)\times \SL_2(\Z/2\Z)$ acts transitively on each theta characteristic, and $\SL_2(\Z/k\Z)\times \SL_2(\Z/k\Z)$ acts transitively on the vectors $(v_1,v_2)$ in each separate case.
\end{proof}
This finishes the proof of the proposition.
\end{proof}
\section{Applications to $S_d$-fields with squarefree discriminant}\label{sec: S_d}

In this section, we will use our homological theorems to obtain arithmetic results in the setting of $S_d$-function fields with squarefree discriminant as well as the Picard rank conjecture. The key part of this is deducing a homological stability slope from only knowledge of $R$. 

However, there is one key difference with the previous section on Cohen-Lenstra which makes the analysis here more difficult. In the case of the Cohen-Lenstra heuristics, the relevant Hurwitz space satisfies usual homological stability with respect to only one element. On the other hand, the relevant Hurwitz space here only has homological stability with respect to multiple elements. This is also where working in the generality of multiple quotients in Sections \ref{sec: iterated quotients as koszul} and \ref{sec: filtrations} pays off.

The relevant Hurwitz space for this section is $G=S_d$ and $c$ the conjugacy class of transpositions in $S_d$, and we will use this choice of $G$ and $c$ throughout this section. In line with the previous section, we will first give a parametrization of the connected components of $\Hur^{G,c}_n$ in Section \ref{sec: S_d param}, then prove some lemmas relating $A_c$ and $MA_c$ in Section \ref{sec: comparing A_c and MA_c}, and use both of these to come up with a method to deduce a stability slope for any $d$ in Section \ref{sec: S_d slope}. We will then do explicit computations to find this stability slope in Section \ref{sec: explicit slope computations}, which will then allow us to deduce our arithmetic applications in Section \ref{sec: S_d arithmetic}.
\subsection{Parametrization of $R$}\label{sec: S_d param}
We start off by parametrizing the connected components of $\Hur^{G,c}_n$ for $G=S_d$ and $c$ the conjugacy class of transpositions. It turns out that this case is much simpler than the Cohen-Lenstra case presented in Section \ref{sec: cohen lenstra R} as the obvious invariants suffice and there is no need for an additional lifting invariant like $\bigwedge^2H'$. 

We will need the following invariants. The global monodromy $H$ is a subgroup of $S_d$ generated by transpositions, so it must be of the form $H= S_{B_1}\times \cdots \times S_{B_j}$ for disjoint subsets $B_i\subseteq [d]$ with $|B_i|=d_i\ge 2$. We have the boundary monodromy $\gamma \in H$ which can be decomposed as the product of $\gamma_i\in S_{B_i}$ over $1\le i\le j$. Finally, we have the number $n_i$ of elements of the tuple belonging in each $S_{B_i}$ for each $1\le i\le j$. There are also some compatibility conditions between these invariants, like having $\gamma_i\in A_{B_i}$ if and only if $n_i$ is even, and in order to ensure that $\gamma_i$ can be a product of $n_i$ elements that generate $S_{B_i}$, we need to have $n_i\ge d_i+c(\gamma_i)-2$ where $c(\gamma_i)$ is the number of cycles of $\gamma_i$, including fixed points as $1$-cycles.

The following proposition tells us that these are all the considerations needed to parametrize a connected component of $\Hur^{G,c}_n$.
\begin{proposition}\label{prop: S_d param}
For each $n\ge 1$, there is an injective map
$$\pi_0(\Hur^{G,c}_n)\hookrightarrow \left\{\left(H,(\gamma_i,n_i)\right)\middle| \ H=S_{B_1}\times \cdots \times S_{B_j},\  \gamma_i\in\begin{cases}A_{B_i}& \text{for } n_i \text{ even}\\
S_{B_i}\setminus A_{B_i}& \text{for } n_i \text{ odd}\end{cases},\  \sum_{i=1}^{j} n_i=n\right\}.$$
Furthermore, its image consists of exactly those tuples satisfying $$n_i\ge d_i + c(\gamma_i)-2$$ for all $1\le i\le j$.
\end{proposition}

The parametrization above already appears in \cite[Proposition 7.13]{Bianchi_2025}, which was brought to our attention after this section was written. We retain our proof and exposition for completeness.

Hurwitz orbits of transpositions have a long history -- \cite{Kluitmann_1988} treats the case of full global monodromy $S_d$, while \cite{Ben-Itzhak_Teicher_2003,Seguin_phd} study the case of identity boundary monodromy. These works use graph-theoretic techniques, viewing transpositions as edges of a graph. Related results also appear in the broader Coxeter-group setting of \cite{Douvropoulos_Lewis_2024}, which includes the symmetric group $S_d$ as type $A_{d-1}$. We give a proof of the proposition using \cite{Kluitmann_1988}.
\begin{proof}
Let us consider the tuples with global monodromy $H=S_{B_1}\times \cdots \times S_{B_j}$. Each element in the word $(g_1,\ldots, g_n)$ must be in some $S_{B_i}$, and since elements from different $S_{B_i}$ commute we can use Hurwitz moves to gather them into separate subtuples. Now, \cite[Theorem 1]{Kluitmann_1988} tells us that all subtuples generating $S_{B_i}$ with length $n_i$ and boundary monodromy $\gamma_i$ are in the same Hurwitz orbit, proving injectivity.

For the statement on the image, we first show that tuples in the image must have $n_i\ge d_i+c(\gamma_i)-2$. Again consider each subtuple, and let us multiply the $n_i$ transpositions that are in $S_{B_i}$ one at a time. At each point, multiplying by a transposition either joins two cycles into one or splits one cycle into two. Let there be $J$ joins and $S$ splits. Then, we have $$c(\gamma_i)=d_i-J+S,\qquad n_i=J+S.$$ 

It is useful to take a graph-theoretic view where transpositions are edges in the graph; then the final graph on the subset $B_i$ must be connected as the transpositions generate $S_{B_i}$. However, only joins can reduce the number of connected components, and we start with the empty graph which has $d_i$ connected components, so we must have $J\ge d_i-1$. Putting these equations together, we have
$$n_i=J+S=2J+c(\gamma_i)-d_i\ge d_i+c(\gamma_i)-2$$
as desired. The construction is simple -- just construct each cycle in $\gamma_i$ of length $\ell$ using $\ell-1$ transpositions for a total of $d_i-c(\gamma_i)$ transpositions, then choose $c(\gamma_i)-1$ transpositions between these cycles and append each one twice, which makes the global monodromy $S_{B_i}$ for this subtuple.
\end{proof}
\subsection{Comparing $A_c$ and $MA_c$}\label{sec: comparing A_c and MA_c} A crucial part of the argument for the stability slope will be going between homological stability of $A_{c'}$ and $MA_{c'}$ for various subracks $c'\subseteq c$, which we will discuss in this subsection for a general rack $c$. The key to letting us do this is the filtration $F_iA_c$ in Equation \eqref{eqn: filtration A_c and MA_c} which relates $MA_{c'}$ and $A_{c'}$.

It turns out that we will need to consider directional stability in some stages of the proof, so it will be convenient to work in this setting. Let $c$ be a finite rack, and suppose that we write $c=c_1\sqcup c_2$ where each $c_i$ is a union of connected components. From Section \ref{sec: stability in chosen direction}, this defines a bigrading $(n_1,n_2)$ on $A_c$, as well as on $MA_c$. Likewise, for any subrack $c'\subseteq c$, by taking $c_1'=c'\cap c_1$ and $c_2'=c'\cap c_2$, we can induce the same bigrading on $A_{c'}$ and $MA_{c'}$. Also recall that if we let $U_1,\ldots, U_k$ be $\E_2$-central elements of $A_c$ with homological degree $d_i$ and algebraic bigrading $(a_i,b_i)$ where $a_i>0$, we can define the iterated quotients as $A_c/(U_j)\coloneqq A_c/(U_1,\ldots, U_k)$ and $MA_c/(U_j)\coloneqq MA_c/(U_1,\ldots, U_k)$. We will not be changing the set of stability elements in this subsection, so we default to the above shorthand of $(U_j)\coloneqq (U_1,\ldots, U_k)$.

As we are only focused on the subrack $c_1$ of $c$, we need a directional version of both the filtration in Equation \eqref{eqn: filtration A_c and MA_c} as well as maximal (chains of) Hurwitz spaces. We first define maximal Hurwitz space in a certain direction, modifying Definition \ref{def: hurwitz space}. 

For this to be well defined, we need to assume that $c_1$ and $c_2$ act on each other trivially; otherwise the action of $c_2$ on $c_1$ could change the maximality in the subrack of $c_1$. Hence, for the rest of this section, we will assume this condition, which will also be satisfied in our use case later as well.
\begin{definition}
Suppose that $c_1$ and $c_2$ act on each other trivially. Let $\MHur^{c_1\mid c_2}_{n_1\mid n_2}$ denote the \textit{maximal Hurwitz space in the $n_1$ direction}, which is defined to be the finite étale cover of $\Conf_{n_1,n_2}$ with fibers $(c_1^{n_1}\setminus \bigcup_{c_1'\subsetneq c_1}c_1'^n)\times c_2^{n_2}$. Alternatively, it is the union of connected components of $\Hur^{c_1,c_2}_{n_1,n_2}$ that are not in $\Hur^{c_1',c_2}_{n_1,n_2}$ for any subrack $c_1'\subsetneq c_1$. Also note that $\MHur^{c_1,c_2}_{n_1,n_2}$ is an open and closed subset of $\MHur^{c_1\mid c_2}_{n_1\mid n_2}$ as a connected component which is maximal in the $n_1$ direction may not be maximal in the $n_2$ direction as well.

Likewise, if we write $c_1=c_{1,1}\sqcup \cdots \sqcup c_{1,k_1}$ and $c_2=c_{2,1}\sqcup \cdots \sqcup c_{2,k_2}$, we can define the multicolored maximal Hurwitz space in the $n_1$ direction $\MHur^{c_{1,1},\ldots, c_{1,k_1}\mid c_{2,1},\ldots, c_{2,k_2}}_{n_{1,1},\ldots, n_{1,k_1}\mid n_{2,1},\ldots,n_{2,k_2}}$ to be the connected components which are not in $\MHur^{c_{1,1}',\ldots, c_{1,k_1}'\mid c_{2,1},\ldots, c_{2,k_2}}_{n_{1,1},\ldots, n_{1,k_1}\mid n_{2,1},\ldots,n_{2,k_2}}$ for some $c_1'\subsetneq c_1$ with $c_{1,i}'\coloneqq c_{1,i}\cap c_1'$.

Then, we can denote $$MA_{c_1|c_2}\coloneqq C_*(\MHur^{c_1\mid c_2})=\bigoplus_{n_1,n_2=0}^{\infty} C_*(\MHur^{c_1\mid c_2}_{n_1\mid n_2})$$ to be the maximal chains in the $n_1$ direction.
\end{definition}

We then have the directional generalization of the filtration in Equation \eqref{eqn: filtration A_c and MA_c}.
\begin{equation}\label{eqn: directional filtration A_c, MA_c}
F_mA_{c}/F_{m+1}A_c\cong \bigoplus_{c_1'\subseteq c_1,|c'_1|=m} MA_{c_1'\mid c_2}   
\end{equation}
where $F_m\Hur^c$ is the union of components of $\Hur^{c_1,c_2}$ such that if $c_1'\subseteq c_1$ is the smallest subrack such that the component comes from $\Hur^{c_1',c_2}$, then $|c_1'|\ge m$, and $F_mA_c\coloneqq C_*(F_m\Hur^c,k)$ is the corresponding $A_c$-bimodule. Note that we have $F_0A_c=A_c$ and $F_{|c_1|}A_c=MA_{c_1\mid c_2}$. 

This filtration gives the cofiber long exact sequence
\begin{equation}\label{eqn: LES for filtration}
\begin{split}
\cdots \ \rightarrow  \ H_i(F_{m+1}A_c/(U_j))\ \ \rightarrow \  H_i(F_{m}A_c/(U_j))\ \ &\rightarrow \bigoplus_{c_1'\subseteq c_1,|c_1'|=m} \ H_i(MA_{c_1'\mid c_2}/(U_j)) \ \ \rightarrow\\
\rightarrow H_{i-1}(F_{m+1}A_c/(U_j))\rightarrow H_{i-1}(F_{m}A_c/(U_j))&\rightarrow \bigoplus_{c_1'\subseteq c_1,|c_1'|=m} H_{i-1}(MA_{c_1'\mid c_2}/(U_j))\rightarrow \ \cdots ,
\end{split}
\end{equation}
from which we will deduce the following lemma.

\begin{lemma}\label{lem: deduce A_c from MA_c and A_c'}
Fix an algebraic bigrading $n_1,n_2$ (or a further refinement $n_{1,1},\ldots, n_{1,k_1},n_{2,1},\ldots, n_{2,k_2}$) and an integer $I$. Suppose that for every \textit{proper} subrack $c_1'\subsetneq c_1$ we have
\begin{equation}\label{eqn: proper subrack vanishing condition}
H_i((MA_{c_1'\mid c_2}/(U_j))_{n_1,n_2})=0 \qquad \text{for }i\le I.
\end{equation}
Then, the natural map $$H_i((MA_{c_1\mid c_2}/(U_j))_{n_1,n_2})\rightarrow H_i((A_c/(U_j))_{n_1,n_2})$$
is an isomorphism for $i<I$ and a surjection for $i=I$. 
\end{lemma}
\begin{proof} The long exact sequence \eqref{eqn: LES for filtration} tells us that for every $m$, 
$$H_i((F_{m+1}A_c/(U_j))_{n_1,n_2})\rightarrow H_i((F_mA_c/(U_j))_{n_1,n_2})$$
is an isomorphism for $i<I$ and a surjection for $i=I$. Composing these maps gives the result.
\end{proof}

This immediately gives the following corollaries, which will be crucial for our argument later.
\begin{corollary}\label{cor: deduce A_c}
Suppose for all subracks $c_1'\subseteq c_1$ that $MA_{c_1'\mid c_2}$ stabilizes in an $(N,I)$ staircase with respect to $U_1,\ldots, U_k$ in the $n_1$ direction, then $A_c$ stabilizes in an $(N,I)$ staircase with respect to $U_1,\ldots, U_k$ in the $n_1$ direction.
\end{corollary}
\begin{corollary}\label{cor: deduce MA_c}
Suppose for all proper subracks $c_1'\subsetneq c_1$ that $MA_{c_1'\mid c_2}$ stabilizes (strongly) in an $(N,I+1)$ staircase with respect to $U_1,\ldots, U_k$ in the $n_1$ direction, and that $A_c$ stabilizes (strongly) in an $(N,I)$ staircase with respect to $U_1,\ldots, U_k$ in the $n_1$ direction. Then, $MA_{c_1\mid c_2}$ stabilizes (strongly) in an $(N,I)$ staircase with respect to $U_1,\ldots, U_k$ in the $n_1$ direction.
\end{corollary}

\subsection{Deducing a stability slope}\label{sec: S_d slope}  With the lemmas introduced in Sections \ref{sec: lem extra elements} and \ref{sec: comparing A_c and MA_c}, we are now ready to prove a stability slope for $G=S_d$ and $c_{S_d}$ the conjugacy class of transpositions.

It is easy to check that all possible subracks $c'\subseteq c_{S_d}$ are of the form $c_{S_d}\cap H$ for a subgroup $H=S_{B_1}\times \cdots \times S_{B_j}\le G$ (recalling the notation in Section \ref{sec: S_d param}). It will be convenient to denote this by $c_{B_1,\ldots, B_j}\coloneqq c_H\coloneqq c\cap H$, and we can also decompose $c_{B_1,\ldots, B_j}=c_{B_1}\sqcup  \cdots \sqcup c_{B_j}$ into its connected components, where $c_{B_i}\cong c_{S_{d_i}}$. Importantly, note that these subracks act trivially on one another.

Our end goal will be to deduce homological stability in explicit staircases of $\Hur^{c_{S_d}}_n$ and $\CHur^{c_{S_d}}_n=\MHur^{c_{S_d}}_n$ for each integer $d$. As our argument will involve splitting $c_{S_d}$ into subracks with multiple components and then looking at stability in a particular direction, it makes more sense to work in the directional setting from the start to begin with. Hence, our goal will be to deduce stability in the $n_1$ direction for $\Hur^c$ for $c=c_1\sqcup c_2$, where we have $c_1=c_{S_d}$, and $c_2$ can be taken to be any $c_H$ for $H\cong S_{d'_{1}}\times \cdots \times S_{d'_{j'}}$, with no restrictions on the $d'_{i}$. We can basically ignore $c_2$ for most of this section and we will be focusing on $c_1$. We will be carrying this out by induction on $d$.

Our main mechanism to extrapolate homological stability to larger staircases is our generalized main theorem, Theorem \ref{thm: general main theorem}, which we can apply to $A_c$ in the $n_1$ direction. By Theorem \ref{thm: LL A_c stability}, we have that $A_c$ is homologically stable in the $n_1$ direction with respect to the $\E_2$-central elements $U_i=s_i^2$ where $s_i$ are the $k\coloneqq \binom d 2$ transpositions in $c_1=c_{S_d}$. We will later need to take large powers of these, and by Lemma \ref{lem: homological stability for powers}, we have also that $A_c$ is homologically stable with respect to $U_1^{e_1},\ldots, U_k^{e_k}$ for all $e_1,\ldots, e_k$. 

The least we need to know for a nontrivial application of our main theorem is stability of $A_c$ in an $(N,0)$ staircase with respect to $U_1^{e_1},\ldots, U_k^{e_k}$ in the $n_1$ direction for some $N$. Since we are carrying out our argument by induction on $d$, we would eventually want to deduce this from the case of smaller $d$. By Corollary \ref{cor: deduce A_c}, we need to show that $MA_{c_1'\mid c_2}$ stabilizes in an $(N,0)$ staircase with respect to the same elements and same direction, for all subracks $c_1'\subseteq c_1$. 

We will first deal with the case where $c_1'=c_1$, where we are concerned with $MA_{c_1\mid c_2}$. We first show that this is homologically stable with respect to just one element.
\begin{lemma}\label{lem: MA_c_1|c_2 hom stable}
$MA_{c_1\mid c_2}$ is homologically stable with respect to $U_1$ in the  $n_1$ direction.
\end{lemma}
\begin{proof}
We can basically deduce this from Theorem \ref{thm: LL CA_c iterated quotients}, except that the theorem only tells us for the case of $MA_c$ rather than for $MA_{c_1\mid c_2}$. To deduce the $MA_{c_1\mid c_2}$ case, note that because $c_1$ and $c_2$ act trivially on one another, we can write $\MHur^{c_1\mid c_2}_{n_1,n_2} = \bigsqcup_{c_2'\subseteq c_2} \MHur^{c_1, c_2'}_{n_1,n_2}$. Since multiplication by $U_1$ preserves each component, and because Theorem \ref{thm: LL CA_c iterated quotients} tells us the case of $MA_{c_1\sqcup c_2'}$, we have the desired result.
\end{proof}

This also tells us that $MA_{c_1\mid c_2}$ is homologically stable with respect to some power $U_1^{e_1}$ in the $n_1$ direction by Lemma \ref{lem: homological stability for powers}. Then, by applying Lemma \ref{lem: add extra elements}, the condition that $MA_{c_1\mid c_2}$ stabilizes in an $(N,0)$ staircase with respect to $U_1^{e_1},\ldots, U_k^{e_k}$ in the $n_1$ direction for all sufficiently large $e_1,\ldots, e_k$ is equivalent to it stabilizing in an $(N,0)$ staircase with respect to $U_1^{e_1}$ for all sufficiently large $e_1$. 

Just like in the proof of Corollary \ref{cor: stab from R}, one can in fact deduce this statement from the fact that $H_0((MA_{c_1\mid c_2})_{n_1,n_2})$ stabilizes under $U_1$ for $n\ge N$, which we will explain in more detail in this directional case in the next proof. One can then understand $H_0$ from the parametrization of Proposition \ref{prop: S_d param}. We work this out in the following proposition.
\begin{proposition}\label{prop: MA (N,0)}
We have that 
$$U_1\colon H_0(\MHur^{c_1\mid c_2}_{n_1,n_2},\Q)\rightarrow H_0(\MHur^{c_1\mid c_2}_{n_1+2,n_2},\Q)$$
is an isomorphism for $n_1\ge 2d-3$. 

Consequently, for all $N\ge 2d-2$, we have for sufficiently large $e_1$ that $MA_{c_1\mid c_2}$ stabilizes (strongly) in an $(N,0)$ staircase with respect to $U_1^{e_1}$ in the $n_1$ direction.
\end{proposition}
\begin{proof}
We can view $c=c_1\sqcup c_2$ as the conjugacy class of transpositions in $G=S_d\times H$. Then, it is easy to see that $\MHur^{c_1\mid c_2}_{n_1,n_2}$ contains exactly those components whose global monodromy is of the form $S_d\times H'$ for a subgroup $H'\leq H$. On the other hand, note that multiplication by $U_1=s_1^2$ does not change the boundary monodromy $(\gamma_i)_i$, and it also does not change the global monodromy since it was already of the form $S_d\times H'$. For the tuple $(n_i)_i$, multiplying by $U_1$ only changes $n_1$ to $n_1+2$ and keeps the other $n_i$ tracking the multigrading in $H'$ the same. 

We now check that if $n_1\ge 2d-3$, $U_1$ induces a bijection on the level of $\pi_0$ and hence an isomorphism of $H_0$. Let us first consider the case where $n_1\ge 2d-2$ is even. Since $\gamma_1\in A_d$, we must have $c(\gamma_1)\le d$ and hence $n_1\ge d+c(\gamma_1)-2$ is automatically satisfied for all $\gamma_1$, so every possible $\gamma_1 \in A_d$ is a valid choice. Furthermore, the validity of the equations $n_i\ge d_i+c(\gamma_i)-2$ for $i\neq 1$ stays the same as $n_i$ and $\gamma_i$ do not change. Hence, multiplication by $U_1$ induces a bijection of valid tuples $(S_d\times H',(\gamma_i,n_i))$ in Proposition \ref{prop: S_d param}, and hence gives a bijection on $\pi_0$ as desired. 

For the odd case where $n_1\ge 2d-3$ is odd, since $\gamma_1\in S_d\setminus A_d$ cannot be the identity, we have $c(\gamma_1)\le d-1$ and we still have $n_1\ge d+c(\gamma_1)-2$ for all $\gamma_1\in S_d\setminus A_d$, and we reach the desired conclusion.

We are left to deduce the second part of the statement, in an analogous way to the proof of Corollary \ref{cor: stab from R}, but we are now in the directional setting. Recall that for stability in an $(N,0)$ staircase with respect to one element, we need a surjection starting from $(n_1,i)=(N-1,0),(N,1)$ and an isomorphism starting from $(n_1,i)=(N,0)$. Taking $N\ge 2d-2$ guarantees all except for the $(N,1)$ case. Recall again that $\MHur^{c_1\mid c_2}_{n_1,n_2}=\bigsqcup_{c_2'\subseteq c_2} \MHur^{c_1, c_2'}_{n_1,n_2}$ in our setting, and we also have that $\MHur^{c_1, c_2'}_{n_1,n_2}\subseteq \CHur^{c_1, c_2'}_{n_1,n_2}$ (in fact, one can see that they are equal here). Hence, any component $Z_{n_1}\subseteq \MHur^{c_1\mid c_2}_{n_1,n_2}$ is contained in some $\CHur^{c_1, c_2'}_{n_1,n_2}$. Then, Theorem \ref{thm: stable hom conf one direction} tells us that if $n_1$ is large enough, projecting $Z_{n_1}\subseteq \CHur^{c_1, c_2'}_{n_1,n_2}$ to the connected component $Z'_{n_1}\subseteq \CHur^{1, c_2'}_{n_1,n_2}$ induces an isomorphism of $H_0$ and $H_1$. We can then draw the following diagram for sufficiently large $e_1$. We need to prove that the top arrow is a surjection, and this will follow from the properties of other arrows in the diagram. It is clear that the left arrow is surjective due to it being a finite étale cover, and the right arrow is an isomorphism as explained above. The only thing we need to justify is why the bottom arrow is an isomorphism.
\[\begin{tikzcd}
	{H_1(\MHur^{c_1\mid c_2}_{n_1,n_2},\Q)} & {H_1(\MHur^{c_1\mid c_2}_{n_1+2e_1,n_2},\Q)} \\
	{H_1(\sqcup Z'_{n_1},\Q)} & {H_1(\sqcup Z'_{n_1+2e_1},\Q)}
	\arrow["{U_1^{e_1}}", from=1-1, to=1-2]
	\arrow[two heads, from=1-1, to=2-1]
	\arrow["\sim", from=1-2, to=2-2]
	\arrow["\sim", from=2-1, to=2-2]
\end{tikzcd}\]

To do this, recall that as long as $n_1\ge 2d-3$ we have that $U_1^{e_1}$ induces a bijection between the components $\MHur^{c_1\mid c_2}_{n_1,n_2}$ and $\MHur^{c_1\mid c_2}_{n_1+2e_1,n_2}$ upstairs, and hence also induces a bijection downstairs between the various $Z'_{n_1}$ and $Z'_{n_1+2e_1}$. Take one such pair of components, and view $Z'_{n_1}$ as a fiber bundle $\Conf_{n_1}(\C\setminus \{p_1,\ldots, p_{n_2}\})\rightarrow Z'_{n_1}\rightarrow Y$ for the corresponding connected component $Y\subseteq \CHur^{c_2'}_{n_2}$. The map $U_1^{e_1}$ simply adds $2e_1$ points on the right to each fiber $\Conf_{n_1}(\C\setminus \{p_1,\ldots, p_{n_2}\})$. Homological stability in this punctured case is well known, indeed, \cite[Proposition A.1]{Segal_1979} tells us that the map described here gives an isomorphism on the $H_1$ and $H_0$ of the fibers as long as $n_1\ge 2$ which we already assumed. By naturality of this map across fibers, one can upgrade this to an isomorphism of local systems for the higher pushforwards of the fiber bundle, and finally the Serre spectral sequence associated to this fiber bundle will tell us that the bottom arrow is indeed an isomorphism.
\end{proof}

We are now done with the case of $c_1'=c_1$, and we are left with the case for proper racks $c_1'\subsetneq c_1$, for which we need to show that $MA_{c_1'\mid c_2}$ stabilizes in an $(N,0)$ staircase with respect to $U_1^{e_1},\ldots , U_k^{e_k}$. We will actually attempt to deduce an $(N,I)$ staircase, as we will need it in this generality later. Let us fix such a proper subrack $c_1'\subsetneq c_1$ which must be of the form $c_1'=c_{1,1}\sqcup \cdots \sqcup c_{1,j}$ where $c_{1,i}\cong c_{S_{d_i}}$ for integers $d_i$ satisfying $d_1+\cdots +d_j\le d$. Consider the multigrading $(n_{1,1},\ldots, n_{1,j},n_2)$ where $n_{1,i}$ corresponds to the $i$-th connected component $c_{S_{d_i}}\subseteq c_1'$, so we have $n_1=n_{1,1}+\cdots +n_{1,j}$. The general strategy here is to use our induction hypothesis, as each of these $d_i$ is now strictly smaller than $d$. 

Now, rearrange the indices of $U_i$ such that $U_1=s_1^2$ comes from a transposition $s_1\in c_{1,1}$, and do this for $U_2,\ldots, U_j$ as well. We now have that $MA_{c_1'\mid c_2}$ is homologically stable with respect to $U_1,\ldots, U_j$ in the $n_1$ direction by the same argument presented in Lemma \ref{lem: MA_c_1|c_2 hom stable}, so applying Lemmas \ref{lem: homological stability for powers} and \ref{lem: add extra elements} tells us that $MA_{c_1'\mid c_2}$ stabilizes in an $(N,I)$ staircase with respect to $U_1^{e_1},\ldots, U_k^{e_k}$ in the $n_1$ direction for all sufficiently large $e_i$ if and only if it stabilizes in an $(N,I)$ staircase with respect to $U_1^{e_1},\ldots, U_j^{e_j}$ in the $n_1$ direction for sufficiently large $e_i$.

Notice now that if $n_1$ is large, we will necessarily have that some $n_{1,i}$ is large. Thus, it seems plausible that we can deduce a stabilization staircase in the $n_1$ direction from various stabilization staircases in each of the $n_{1,i}$ directions. This is indeed the case, as we will show in the following lemma. Note here that we will need strong stability, which in this case simply includes the blue rectangular region as shown on the left of Figure \ref{fig: strong vs usual staircase}.
\begin{lemma}\label{lem: one direction from many directions}
For each $1\le i\le j$, suppose that $MA_{c_{1,i}\mid (c_1'\setminus c_{1,i})
\sqcup c_2}$ stabilizes strongly in an $(N_i,I+j-1)$ staircase with respect to $U_i^{e_i}$ in the $n_{1,i}$ direction. Then, $MA_{c_1'\mid c_2}$ stabilizes strongly in an $(N,I)$ staircase with respect to $U_1^{e_1},\ldots, U_j^{e_j}$ in the $n_1$ direction, where we have
$$N= 1+ \sum_i (N_i-1).$$
\end{lemma}
\begin{proof} 
For each $1\le i\le j$, since $$\MHur^{c_{1,i}\mid (c_1'\setminus c_{1,i})\sqcup c_2}=\bigsqcup_{c_1''\subseteq c_1'\setminus c_{1,i}} \MHur^{c_{1,i}\sqcup c_1''\mid c_2},$$ and multiplication by $U_i$ preserves this decomposition, the condition in the statement actually tells us that $MA_{c_1'\mid c_2}$ stabilizes strongly in an $(N_i,I+j-1)$ staircase with respect to $U_i^{e_i}$ in the $n_{1,i}$ direction. 

We claim that this tells us that $MA_{c_1'\mid c_2}$ stabilizes strongly in an $(N_i,I)$ staircase with respect to $U_1^{e_1},\ldots, U_j^{e_j}$ in the $n_{1,i}$ direction. It is important to note that we cannot add the elements $U_1^{e_1},\ldots, \widehat{U_i^{e_i}},\ldots, U_j^{e_j}$ via Lemma \ref{lem: add extra elements} because in the directional setup of that lemma we assumed that the elements we are adding have positive algebraic grading in the relevant direction, which is not the case. This is why we lose a constant term of $j-1$. To prove this claim, we simply use the usual long exact sequence \eqref{eqn: iterated quotient LES}, noting that we lose one homological degree whenever we quotient by an element.

Now, we go back to the definition of strong stability in a staircase, where in our case the slopes are all zero we need to show for all $i'\le \min(I-N+n',I)$ that $$H_{i'+j}\left((MA_{c_1'\mid c_2}/(U_1^{e_1},\ldots, U_j^{e_j}))_{n'+2e_1+\cdots +2e_j}\right)=0.$$
where we only show the $n_1$ grading $n'+2e_1+\cdots +2e_j$ here. First suppose that $n'\ge N$. Then by the definition of $N$, we must have for some $i$ that the grading of the above term in the $n_{1,i}$ direction is $\ge N_i+2e_i$. In that case, the claim in the previous paragraph on the vanishing in the $n_{1,i}$ direction tells us what we need. In the case when $n'<N$, let $n'=N-c$. Then, it is also easy to see that there is some $n_{1,i}$ grading which is $\ge N_i+2e_i-c$, giving us what we need. 
\end{proof}

With this, since $d_i<d$, we can simply use the induction to obtain the conditions we need. Going back to the start of this subsection, our work so far tells us that $A_{c}$ stabilizes in an $(N,0)$ staircase with respect to $U_1^{e_1},\ldots, U_k^{e_k}$ in the $n_1$ direction for some $N$. Applying Theorem \ref{thm: general main theorem}, we have stability in an $(r(N+1)-1,r(I+1)-1)$ staircase (with respect to the same elements and direction), and in fact strong stability by Theorem \ref{thm: usual vs strong stability direction}. Then, we can deduce strong stability of $MA_{c_1\mid c_2}$ using Corollary \ref{cor: deduce MA_c}, and we once again remove the extra stability elements using Lemma \ref{lem: MA_c_1|c_2 hom stable} and Lemma \ref{lem: add extra elements}.

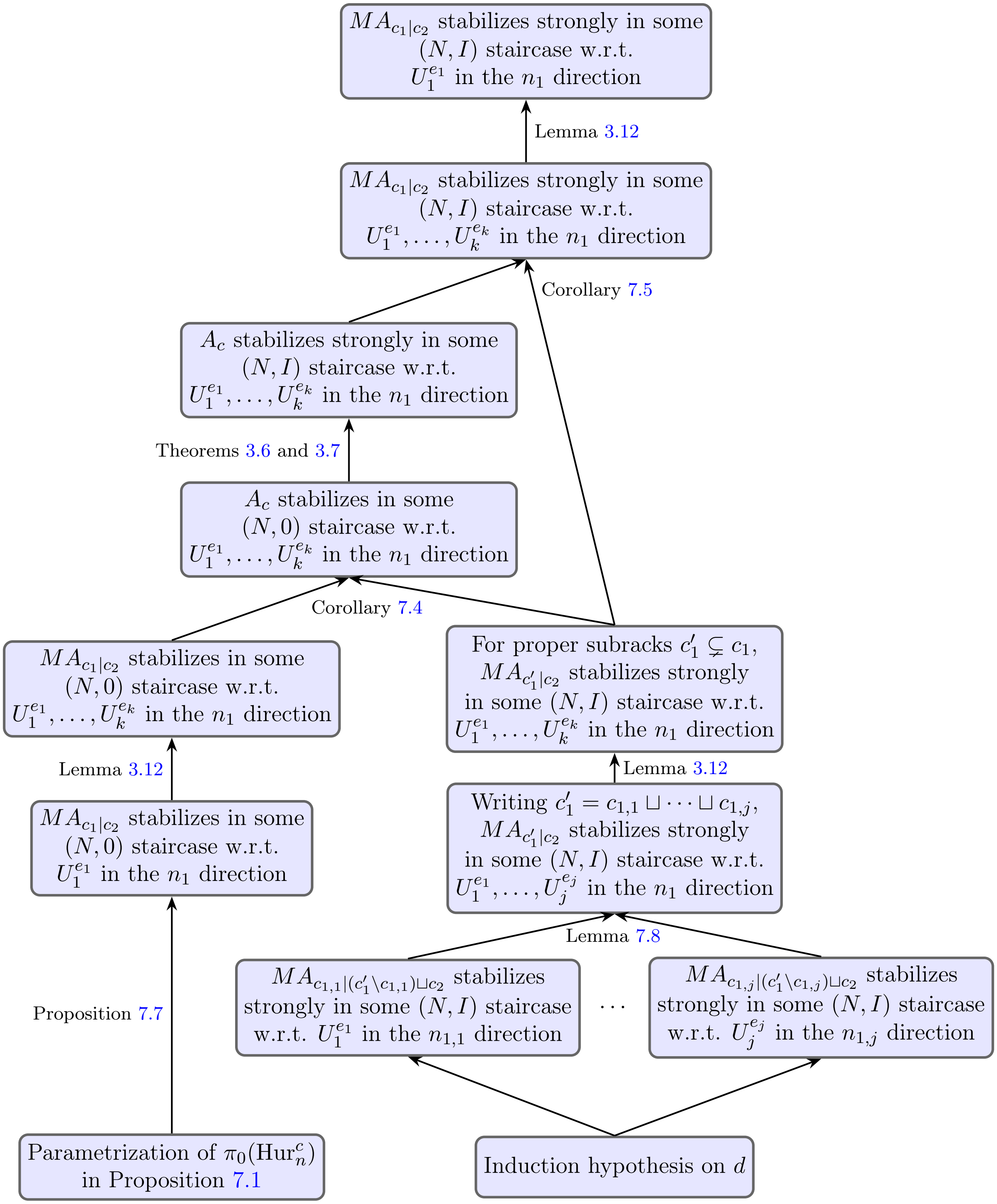
\begin{figure}[!htbp]
\centering
\begin{tikzpicture}[
    >=Stealth,
    box/.style={
        rectangle,
        draw=black!60,
        fill=blue!10,
        very thick,
        rounded corners,
        minimum width=28mm,
        minimum height=10mm,
        align=center
    },
    arrow/.style={->, thick},
    arrowlabel/.style={
        font=\footnotesize
    },
    scale = 0.98
]

\node[box] (L0) at (0,-2.7) {Parametrization of $\pi_0(\Hur^c_n)$\\in Proposition \ref{prop: S_d param}};
\node[box] (L1) at (0,2.7) {$MA_{c_1\mid c_2}$ stabilizes in some\\$(N,0)$ staircase w.r.t.\\$U_1^{e_1}$ in the $n_1$ direction};
\node[box] (L2) at (0,5.4) {$MA_{c_1\mid c_2}$ stabilizes in some\\$(N,0)$ staircase w.r.t.\\$U_1^{e_1},\ldots, U_k^{e_k}$ in the $n_1$ direction};
\node[box] (L3) at (3,8.1) {$A_c$ stabilizes in some\\$(N,0)$ staircase w.r.t.\\$U_1^{e_1},\ldots, U_k^{e_k}$ in the $n_1$ direction};
\node[box] (L4) at (3,10.8) {$A_{c}$ stabilizes strongly in some\\$(N,I)$ staircase w.r.t.\\$U_1^{e_1},\ldots, U_k^{e_k}$ in the $n_1$ direction};

\node[box] (T) at (6,13.5) {$MA_{c_1\mid c_2}$ stabilizes strongly in some\\$(N,I)$ staircase w.r.t.\\$U_1^{e_1},\ldots, U_k^{e_k}$ in the $n_1$ direction};
\node[box] (T2) at (6,16.2) {$MA_{c_1\mid c_2}$ stabilizes strongly in some\\$(N,I)$ staircase w.r.t.\\$U_1^{e_1}$ in the $n_1$ direction};

\node[box] (B0) at (7.5,-2.7) {Induction hypothesis on $d$};
\node[box] (B1) at (4,0) {$MA_{c_{1,1}\mid (c_1'\setminus c_{1,1})\sqcup c_2}$ stabilizes\\strongly in some $(N,I)$ staircase\\w.r.t. $U_1^{e_1}$ in the $n_{1,1}$ direction};
\node[box] (B2) at (7.5,2.7) {Writing $c_1'=c_{1,1}\sqcup \cdots\sqcup c_{1,j}$,\\$MA_{c_1'\mid c_2}$ stabilizes strongly\\in some $(N,I)$ staircase w.r.t.\\$U_1^{e_1},\ldots, U_j^{e_j}$ in the $n_1$ direction};
\node[box] (B3) at (11,0) {$MA_{c_{1,j}\mid (c_1'\setminus c_{1,j})\sqcup c_2}$ stabilizes\\strongly in some $(N,I)$ staircase\\w.r.t. $U_j^{e_j}$ in the $n_{1,j}$ direction};
\node[box] (B4) at (7.5,5.4) {For proper subracks $c_1'\subsetneq c_1$,\\$MA_{c_1'\mid c_2}$ stabilizes strongly\\in some $(N,I)$ staircase w.r.t.\\$U_1^{e_1},\ldots, U_k^{e_k}$ in the $n_1$ direction};
\node at (7.5,0) {$\cdots$};

\draw[arrow] (L0.north) -- node[arrowlabel, left] {Proposition \ref{prop: MA (N,0)}} (L1.south);
\draw[arrow] (L1.north) -- node[arrowlabel, left] {Lemma \ref{lem: add extra elements}} (L2.south);
\draw[arrow] (L2.north) -- node[arrowlabel, right=20pt] {Corollary \ref{cor: deduce A_c}} (L3.south);
\draw[arrow] (L3.north) -- node[arrowlabel, left] {Theorems \ref{thm: usual vs strong stability direction} and \ref{thm: general main theorem}} (L4.south);

\draw[arrow] (L4.north) -- node[arrowlabel, right=45pt] {Corollary \ref{cor: deduce MA_c}} (T.south);
\draw[arrow] (B4.north) -- node[arrowlabel, right] {$ $} (T.south);
\draw[arrow] (T.north) -- node[arrowlabel, right] {Lemma \ref{lem: add extra elements}} (T2.south);

\draw[arrow] (B0.north) -- node[arrowlabel, left] {$ $} (B1.south);
\draw[arrow] (B0.north) -- node[arrowlabel, below] {$ $} (B3.south);

\draw[arrow] (B1.north) -- node[arrowlabel, left] {$ $} (B2.south);
\draw[arrow] (B3.north) -- node[arrowlabel, left=23pt] {Lemma \ref{lem: one direction from many directions}} (B2.south);
\draw[arrow] (B2.north) -- node[arrowlabel, right] {Lemma \ref{lem: add extra elements}} (B4.south);

\draw[arrow] (B4.north) -- node[arrowlabel, above] {$ $} (L3.south);

\end{tikzpicture}
\caption{A flowchart summarizing Section \ref{sec: S_d slope}.}\label{fig: summary of S_d slope}
\end{figure}

This completes our argument on the inductive step on $d$, and we summarize this in the flowchart in Figure \ref{fig: summary of S_d slope}. All that is left to do is to explicitly find these stability staircases by tracing through the steps we took, which we will do in the following subsection.

\subsection{Explicit computations}\label{sec: explicit slope computations} Continuing from the previous subsection, given $d\ge 2$ and $I\ge 0$, we want to find constants $N(d,I)$ such that $MA_{c_1\mid c_2}$ (where $c_1\cong c_{S_d}$) stabilizes strongly in an $(N(d,I),I)$ staircase with respect to $U_1^{e_1}$ (where $U_1^{e_1}=s_1^{2e_1}$ and $s_1$ is any transposition in $c_1$) in the $n_1$ direction. Note that for the $I=0$ case we can immediately take $N(d,0)=2d-2$ by Proposition \ref{prop: MA (N,0)}. For $I\ge 1$, we will use the argument outlined previously to determine the constants $N(d,I)$ by induction on $d$.

Let us first work these out for the first few values of $d$ starting from $d=2$ by hand, to get a feel for the computations required. We will often omit which elements and direction the stability is with respect to, as it is usually clear from context and our previous discussion. It will be useful to refer to Figure \ref{fig: summary of S_d slope} here.
\begin{itemize}
\item $d=2$. In this case, $c_1\cong c_{S_2}$ has only one element and hence has no proper subracks, so $MA_{c_1\mid c_2}=A_c$ (for $n_1\ge 1$). Proposition \ref{prop: MA (N,0)} then tells us that $A_c$ stabilizes in a $(2,0)$ staircase with respect to $U_1^{e_1}$ in the $n_1$ direction, which then gives us strong stability in a $(2+3I,I)$ staircase by our main theorems. Hence, we have $N(2,I)=2+3I$.
\item $d=3$. Now, $c_1\cong c_{S_3}$ has $3$ possible subracks all of the form $c_1'\cong c_{S_2}$, and the previous case tells us that for each $I\ge 0$, $MA_{c_1'\mid c_2}$ stabilizes strongly in an $(N(2,I),I)=(2+3I,I)$ staircase with respect to $U_1^{e_1}$ in the $n_1$ direction. Meanwhile, $MA_{c_1\mid c_2}$ stabilizes in a $(4,0)$ staircase by Proposition \ref{prop: MA (N,0)}, so by Corollary \ref{cor: deduce A_c}, $A_c$ stabilizes in a $(\max(4,2),0)=(4,0)$ staircase with respect to $U_1^{e_1},U_2^{e_2},U_3^{e_3}$ in the $n_1$ direction. Our main theorems tell us that $A_c$ stabilizes strongly in a $(4+5I,I)$ staircase, so by Corollary \ref{cor: deduce MA_c}, $MA_{c_1\mid c_2}$ stabilizes strongly in a $(\max(4+5I,2+3(I+1)),I)=(4+5I,I)$ staircase. Hence, we have $N(3,I)=4+5I$.
\item $d=4$. There are three possible types of subracks of $c_1\cong c_{S_4}$ which are $c_1'\cong c_{S_3}$, $c_1'\cong c_{S_2}$ and $c_1'\cong c_{S_2}\sqcup c_{S_2}$. In the first case, we have $MA_{c_1'\mid c_2}$ stabilizes strongly in a $(4+5I,I)$ staircase for each $I\ge 0$. The second case is subsumed by either the first or the third case. For the third case, we apply Lemma \ref{lem: one direction from many directions} with $j=2$, and since $MA_{c_{S_2}\mid c_{S_2}\sqcup c_2}$ stabilizes strongly in a $(2+3(I+1),I+1)$ staircase, we have that $MA_{c_1'\mid c_2}$ stabilizes strongly in a $(3(2I+3),I)$ staircase for each $I\ge 0$. Combining with the fact that $MA_{c_1\mid c_2}$ also stabilizes in a $(9,0)$ staircase (since $9\ge 2d-2$), we have that $A_c$ stabilizes in a $(9,0)$ staircase and hence stabilizes strongly in a $(9+10I,I)$ staircase. Applying Corollary \ref{cor: deduce MA_c}, $MA_{c_1\mid c_2}$ stabilizes strongly in a $(\max(9+10I,3(2I+5)),I)$ staircase. Hence, we have $N(4,1)=21$, and $N(4,I)=9+10I$ for $I\ge 2$.
\end{itemize}

With these examples in mind, we now write a formal recurrence relation computing $N(d,I)$ for $I\ge 1$ given $N(d',I)$ for $d'<d$. First, we have by Lemma \ref{lem: one direction from many directions} that the smallest $N=\widetilde N(d,I)$ for which $MA_{c_1'\mid c_2}$ stabilizes strongly in an $(N,I)$ staircase for \textit{all} proper subracks $c_1'\subsetneq c_1$ is given by 
\begin{equation}\label{eqn: N tilde (d,I)}
\widetilde N(d,I)=\max_{\substack{d_1+\cdots +d_j\le d\\ d_i\ge 2\\(j,d_1)\neq (1,d)}} \left(-j+1+\sum_{i=1}^j N(d_i,I+j-1)\right).
\end{equation}
Then, it is clear for $d\ge 4$ that $\widetilde N(d,0)\ge 2d-2$ so $A_c$ stabilizes in a $(\widetilde N(d,0),0)$ staircase, from which we conclude that it stabilizes in a $((I+1)(\widetilde N(d,0)+1)-1,I)$ staircase. Then, we have for $I\ge 1$ that
\begin{equation}\label{eqn: N (d,I)}
N(d,I)=\max\left((I+1)(\widetilde N(d,0)+1)-1,\widetilde N(d,I+1)\right).
\end{equation}

Running this recurrence using code for a few more $d$, we can take the following values for $I\ge 1$:
\begin{align*}
N(5,I)&=\begin{cases}
    29&\text{if }I=1\\
    13+14I &\text{if } I\ge 2,
\end{cases}\\
N(6,I)&=25+26I,\\
N(7,I)&=\begin{cases}
    77&\text{if }I=1\\
    103&\text{if }I=2\\
    33+34I&\text{if }I\ge 3,
\end{cases}\\
N(8,I)&=\begin{cases}
    113&\text{if }I=1\\
    55+56I &\text{if } I\ge 2,
\end{cases}
\end{align*}

\begin{figure}[!htbp]
    \centering
    \includegraphics[width=0.48\textwidth]{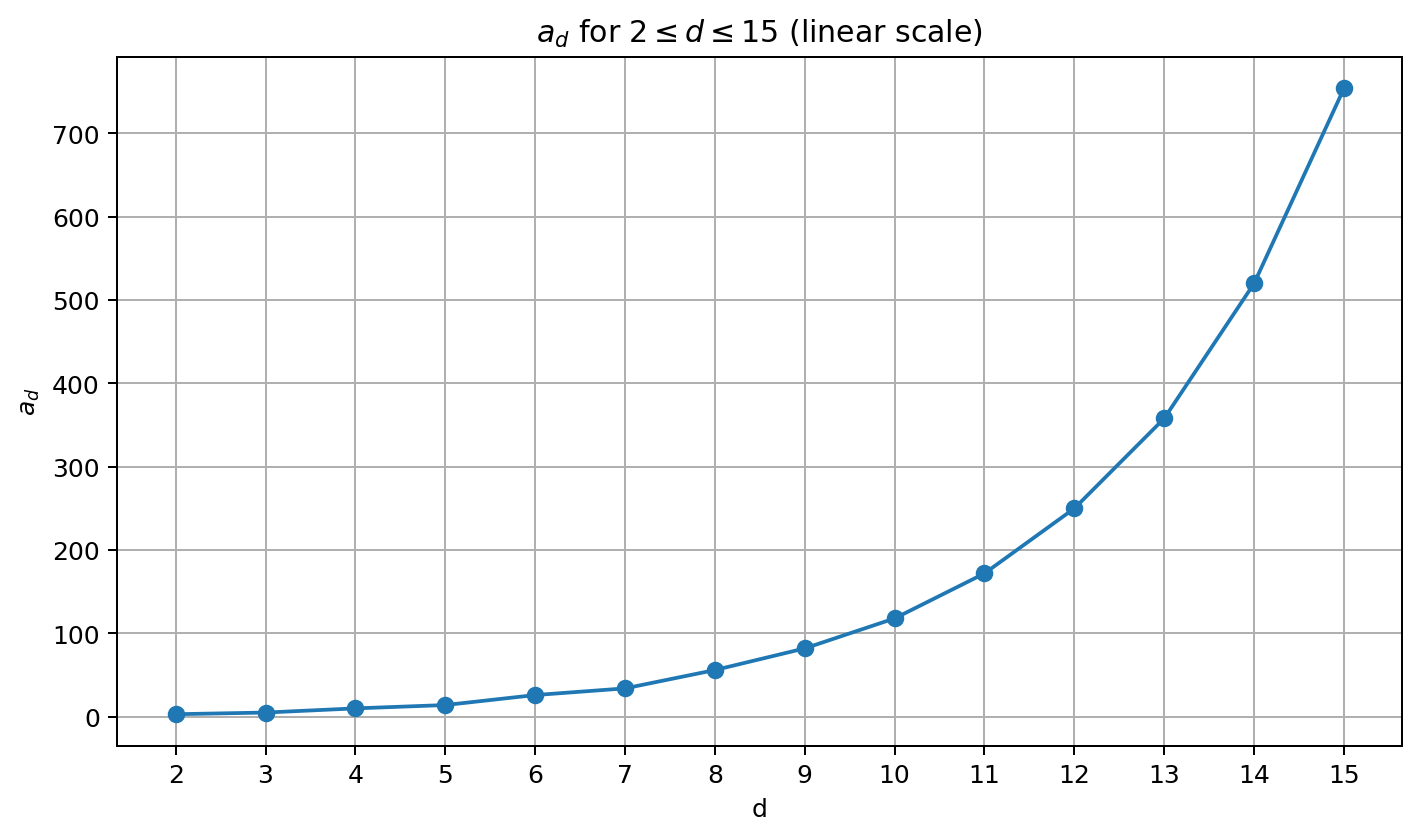}
    \hfill
    \includegraphics[width=0.48\textwidth]{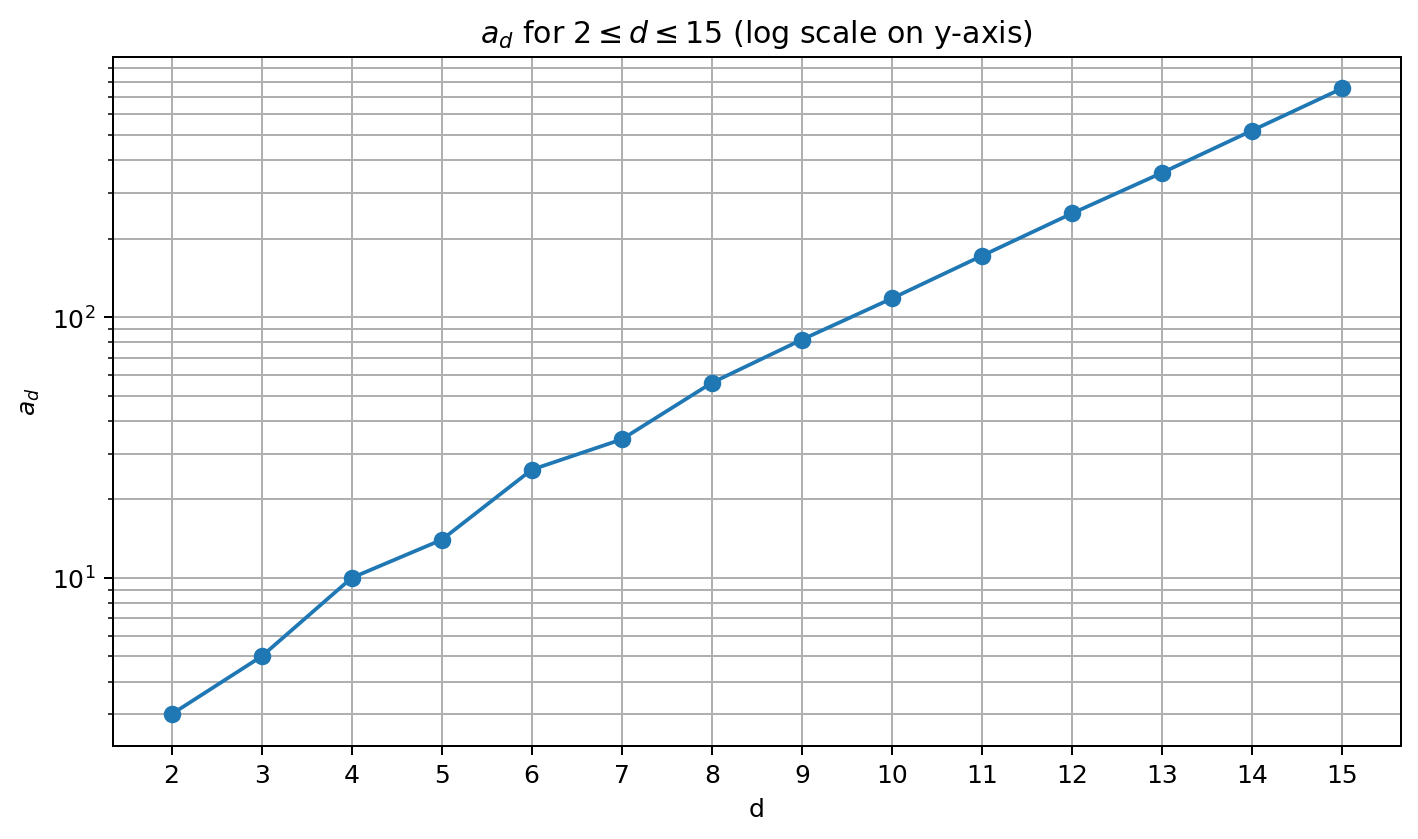}
    \caption{Growth of $a_d$ for $2\le d\le 15$ on linear and logarithmic scales.}
    \label{fig: ad-growth}
\end{figure}

We also notice the following pattern once $d\ge 9$.
\begin{proposition}\label{prop: explicit ad computations}
For every $d\ge 9$ and $I\ge 1$, we can take
$$N(d,I)=\begin{cases}
    3a_{d-1}-1&\text{if }I=1\\
    (I+1)a_d-1 &\text{if } I\ge 2
\end{cases}$$
where for $k\ge 0$ we have
\begin{equation}\label{eqn: a_d values}
\begin{split}
a_{8+3k}&=58\cdot 3^k-2,\\
a_{9+3k}&=84\cdot 3^k-2,\\
a_{10+3k}&=120\cdot 3^k-2.
\end{split}
\end{equation}
In particular, we have $N(d,I)=\Theta(3^{d/3}I)$.
\end{proposition}
Before we prove this, note that these values of $a_d$ for $d\ge 8$ are the inverse of the stability slopes. In line with this, we let 
\begin{equation}\label{eqn: small a_d values}
(a_2,a_3,a_4,a_5,a_6, a_7)=(3,5,10,14,26,34).
\end{equation}
We plot this in Figure \ref{fig: ad-growth}. We also have the explicit inequality
$a_d\le \frac{28}{9}3^{d/3}-2$ for all $d\ge 2$. 
\begin{proof}
The key is to find out which terms dominate when taking the maximum in Equations \eqref{eqn: N tilde (d,I)} and \eqref{eqn: N (d,I)} once $d$ is sufficiently large. We take say $d\ge 20$ to make our calculations simpler and proceed via strong induction on $d$. We verify the $d\le 19$ cases directly via code, which we omit. It will then suffice to prove the following three statements for $d\ge 19$.
\begin{enumerate}[(\alph*)]
\item The maximum defining $\widetilde N(d,0)$ in Equation \eqref{eqn: N tilde (d,I)} is attained by $(d_1,d_2)=(2,d-2)$, giving us
$$\widetilde N(d,0)=N(d-2,1)+4=a_d-1.$$
\item The maximum defining $\widetilde N(d,2)$ in Equation \eqref{eqn: N tilde (d,I)} is attained by $(d_1)=(d-1)$, giving us
$$\widetilde N(d,2)=N(d-1,2)=3a_{d-1}-1.$$
Combining $(a)$ and $(b)$, we then have
$$N(d,1) = \max(2a_d-1,3a_{d-1}-1)=3a_{d-1}-1.$$
\item For all $I\ge 2$, the maximum defining $N(d,I)$ in Equation \eqref{eqn: N (d,I)} is attained by the first term, giving us
$$N(d,I)=(I+1)(\widetilde N(d,0)+1)-1=(I+1)a_d-1.$$
\end{enumerate}

It remains to show $(a),(b)$ and $(c)$ via explicit bounding, which we carry out below. 
\begin{enumerate}[(\alph*)]
\item 
\begin{enumerate}[(\roman*)]
    \item If $j=1$, we have $N(d-1,0)=2d-4\ll a_d-1$.
    \item If $j=2$, then since $N(d,1)$ grows exponentially, by convexity\footnote{note that $N(d,1)$ is not convex for small values of $d$, but we can make it convex by changing some values slightly.}, $N(d_1,1)+N(d_2,1)-1$ is maximized at $(d_1,d_2)=(2,d-2)$ giving $a_d-1$ as discussed.
    \item If $j\ge 3$, by convexity, we have $\le N(d-2j+2,j-1)+(j-1)N(2,j-1) <ja_{d-2j+2}+10j^2\ll a_d-1$.
\end{enumerate}
\item 
\begin{enumerate}[(\roman*)]
\item If $j=1$, we have $N(d-1,2)=3a_{d-1}-1$.
\item If $j\ge 2$, by convexity, we have $\le (j+2)a_{d-2j+2}+10j^2\ll 3a_{d-1}-1$.
\end{enumerate}
\item We want to show that $\widetilde N(d,I+1)\le (I+1)a_d$ for $I\ge 2$.
\begin{enumerate}[(\roman*)]
\item If $j=1$, we have $N(d-1,I+1)\le (I+2)a_{d-1}<(I+1)a_d$.
\item If $j\ge 2$, by convexity, we have $\le (I+j+1)a_{d-2j+2}+10Ij^2\ll (I+1)a_d$.
\end{enumerate}
\end{enumerate}
This finishes the proof.
\end{proof}
\subsection{Arithmetic applications}\label{sec: S_d arithmetic}

We will now deduce the two arithmetic applications in Section \ref{sec: intro S_d fields}, Theorems \ref{thm: S_d count} and \ref{thm: picard rank}. This will basically follow from the explicit stability staircases that we computed in Section \ref{sec: explicit slope computations}. 

Let us specialize back to the case where $G=S_d$ and $c=c_{S_d}$ is the rack of transpositions for $d\ge 3$. Then, from the previous section, we have that $H_i(\CHur^{G,c}_n,\Q)$ stabilizes for $n\ge N(d,i)$, where we take $c_2$ to be the empty rack and use the fact that $\CHur^{G,c}_n=\MHur^{G,c}_n$ as $c$ is connected. The stability slopes are then given by $1/a_d$ as defined in Equations \eqref{eqn: a_d values} and \eqref{eqn: small a_d values}.

\begin{proof}[Proof of Theorem \ref{thm: S_d count}]
We deduce the asymptotic count of $S_d$-extensions of $\mathbb F_q(t)$ with squarefree discriminant by following the point-counting argument of \cite[Proposition 5.3.1]{landesman2025cohenlenstramomentsfunctionfields}. Such extensions correspond to $\mathbb F_q$-points coming from either the components of $[\CHur^{G,c}_{2n}/S_d]$ with zero boundary monodromy or the components of $[\CHur^{G,c}_{2n-1}/S_d]$ with boundary monodromy a transposition. By our parametrization in Proposition \ref{prop: S_d param}, both strata are geometrically connected for $2n\ge 2d-2$. Theorem \ref{thm: stable hom conf all directions} identifies the stable homology of each full monodromy component with that of configuration space, and this identification descends to the $S_d$-quotients. One can apply the comparison theorem with étale cohomology since $\gcd(q,d!)=1$. The point-counting argument of \cite[Proposition 5.3.1]{landesman2025cohenlenstramomentsfunctionfields} then tells us that the two cases have a main term of $(1-q^{-1})q^{2n}$ and $(1-q^{-1})q^{2n-1}$ and this sums to $(1-q^{-2})q^{2n}$ as desired. Finally, since we have a stability slope of $1/a_d$ from the previous section, this gives the desired error term by Deligne's weight bounds and the Betti number bounds for sufficiently large $q$.
\end{proof}

\begin{proof}[Proof of Theorem \ref{thm: picard rank}]
We work throughout with rational coefficients $\Q$ and obtain the explicit bound by tracking the number $n$ of branch points through the proof of \cite[Theorem 7.1.1]{landesman2025homologicalstabilityhurwitzspaces}. The only restrictions on $n$ arise in Lemmas 7.3.1 and 7.3.2 of loc. cit.. Proposition \ref{prop: S_d param} gives the required component counts in gradings $n$ and $n-1$ whenever $n\ge 2d-2$ is even, which is automatic since $n=2g+2d-2$. In the proof of Lemma 7.3.2, the Gysin sequence relates $H^2$ of the projective-line Hurwitz space to $H^1,H^2$ of the affine stratum with $n$ branch points and $H^0,H^1$ of the affine stratum with $n-1$ branch points. The five lemma therefore gives the required comparison whenever $n\ge N(d,2)$ and $n-1\ge N(d,1)$. Consequently, the Picard rank conjecture holds for $n\ge\max\{N(d,2),N(d,1)+1\}$. The computations in Section \ref{sec: explicit slope computations} give $N(d,2)\ge N(d,1)+1$ and $N(d,2)\le 10\cdot 3^{d/3}$, proving the first assertion. 

Under Conjecture \ref{conj: S_d}, both $H_1$ and $H_2$ stabilize for $n\ge 2d+3$, so the same argument applies for even $n\ge 2d+4$, equivalently $g\ge 3$. Combining this with \cite{Deopurkar_Patel_2015,Mullane_2023} proves the second assertion.
\end{proof}

\section{Homological computations for slope and exponent bounds}\label{sec: computations}

In this section, we report the results of our homological computations for various Hurwitz spaces. These homological computations were done in order to obtain the improved bounds on stability slopes and secondary exponents which we presented in the introduction. We will obtain these bounds from the homological computations by applying Theorem \ref{thm: main theorem} in the non-splitting case, and Theorem \ref{thm: general main theorem} in the general case. We present our computational results concerning Cohen--Lenstra moments in Section \ref{sec: compute cohen lenstra}, squarefree degree-$d$ extensions with Galois group $S_d$ in Section \ref{sec: compute S_d}, and other Hurwitz spaces in Section \ref{sec: compute other hurwitz spaces}.

\textit{All homological computations in this section were performed using code generated by ChatGPT, which also devised the overall computational strategy.} The AI-generated code is available in the file \href{https://drive.google.com/drive/folders/1V3231gmKYnkIrTjt3vzWn-ZfoA8ppqbO?usp=drive_link}{here}, with the code for computing homology of Hurwitz spaces in \textit{hurwitz\_cpp}. The method used here is a big improvement over the one used in \cite{prequel}. While the previous algorithm for computing homology had complexity exponential in $n$ with no dependence on $i$, the new algorithm is exponential in homological degree $i$, and for fixed $i$ it is polynomial in $n$ (and empirically, almost linear). This allows us to compute the homology for large $n$ provided that $i$ is reasonably small, letting us see where each $H_i$ stabilizes. 

We briefly explain the main idea here; a full AI-generated writeup is available at the same link in \textit{Hurwitz\_Computational\_Methods}. The previous algorithm computes homology using the bar complex, whose terms have dimensions exponential in $n$. In contrast, the new algorithm constructs a minimal free resolution, introducing only the generators needed to kill homology classes, which allows us to compute low-degree homology using a much smaller complex. We start from the identification $H_i(\Hur_n,\Q)^*\cong \Tor_{n-i,n}^{\mathfrak A}(\Q,\Q)$, where $\mathfrak A$ is the associated quantum shuffle algebra (the dual of the tensor algebra, see \cite[Equation (2.1)]{prequel}), and the second grading on $\Tor$ is the algebraic grading. We construct a minimal free resolution
$$
\cdots \rightarrow P_2 \rightarrow P_1 \rightarrow P_0 \rightarrow \Q \rightarrow 0
$$
of $\Q$ over $\mathfrak A$ with $P_0=\mathfrak A$. Recall that a free resolution is minimal if $d_p(P_p)\subseteq \mathfrak A_{>0}P_{p-1}$ where $\mathfrak A_{>0}$ is the positive-degree part of $\mathfrak A$. In particular, since $P_0=\mathfrak A$ is generated in algebraic grading $0$, minimality implies inductively that
$$
P_p=\bigoplus_{n\geq p}\mathfrak A(-n)\otimes_{\Q}E_{p,n},
$$
where $\mathfrak A(-n)$ denotes the algebraic grading shift by $n$, and $E_{p,n}$ is the vector space of free generators in resolution degree $p$ and algebraic grading $n$. Moreover, since $\Q=\mathfrak A/\mathfrak A_{>0}$, all differentials vanish after tensoring the resolution with $\Q$ over $\mathfrak A$. Thus $\Tor_{p,n}^{\mathfrak A}(\Q,\Q)\cong E_{p,n}$, and consequently $\dim H_i(\Hur_n,\Q)=\dim E_{n-i,n}$, so it suffices to construct the minimal resolution to obtain the Betti numbers.

The resolution is constructed inductively in the algebraic grading. Having constructed all generators and differentials in gradings less than $n$, we compute the homology of the degree-$n$ part of the partial resolution and introduce new generators whose differentials represent and kill these homology classes. The computational advantage is that, to compute homology in degrees at most $i$, we only need generators $E_{p,n}$ satisfying $p\leq n\leq p+i$. Thus, although the resolution must be constructed in increasingly large resolution degrees $p$ as $n$ grows, we only need a strip of width $i+1$ in the $(p,n)$-plane. This avoids the exponentially large chain groups of the original bar complex and makes computations for large $n$ feasible.

To speed up the computations, we work over finite fields. The resulting Betti numbers give upper bounds for the rational Betti numbers, which suffice for the stability checks below. In particular, vanishing over a finite field implies rational vanishing. In every case we checked in two different characteristics, the dimensions agreed, providing evidence that these bounds are equalities.

In the figures below, blue shading marks the stable range. To save space, we show only enough of this range to indicate where stability begins. When a stability staircase required checks beyond the plotted range, we computed those additional gradings as well.

\subsection{Cohen-Lenstra moments}\label{sec: compute cohen lenstra}

In this subsection, we compute the homology of $\Hur^{G,c}_n$ for $G=H\rtimes \Z/2\Z$ and $c=c^{\text{rfl}}_G$ for various choices of $H$, so as to justify the bounds on slopes and exponents in Figures \ref{fig: table H without 3-torsion} and \ref{fig: table H with 3-torsion} as well as the surrounding exposition there. We illustrate some of our computations in Figures \ref{fig: homology computation non 3 torsion} and \ref{fig: homology computation 3-torsion}. 
\begin{figure}[!htbp]
    \centering
 \begin{tikzpicture}[x=1cm,y=1cm,scale=0.55]
\drawnaxis{13}
\drawiaxis{6}
\drawstaircase{1}{6}{0}{black!10}{0}
\drawblock{7}{13}{0}{6}{blue!20}{0}

\foreach \n/\i/\txt in {
1/0/5,
2/0/9,2/1/9,
6/2/2,
12/5/4
 }{
    \cell{\n}{\i}{\txt}{black!10}
  }

\foreach \n/\i/\txt in {
3/0/10,3/1/25,3/2/15,
4/1/21,4/2/39,4/3/28,
5/0/10,5/1/10,5/2/25,
6/0/10,6/1/10,6/3/47,
7/0/10,7/1/10,7/3/30,
8/0/10,8/1/10,8/3/10,8/4/58,
9/0/10,9/1/10,9/5/90,
10/0/10,10/1/10,10/5/45,
11/0/10,11/1/10,11/5/20,
12/0/10,12/1/10,12/6/13,
13/0/10,13/1/10
 }{
    \cell{\n}{\i}{\txt}{black!10}
  }

\foreach \n/\i/\txt in {
5/3/125,5/4/100,
6/4/307,6/5/262,
7/4/170,
8/5/406,
10/6/273,
11/6/110
 }{
    \smallcell{\n}{\i}{\txt}{black!10}
  }

  \foreach \n/\i/\txt in {
7/5/1205,7/6/1065,
8/6/3161,
9/6/1095
 }{
    \tinycell{\n}{\i}{\txt}{black!10}
  }

\drawblock{3}{13}{0}{0}{blue!20}{10}
\drawblock{5}{13}{1}{1}{blue!20}{10}
\end{tikzpicture}
\begin{tikzpicture}[x=1cm,y=1cm,scale=0.55]
\drawnaxis{9}
\drawiaxis{4}
\drawstaircase{1}{4}{0}{black!10}{0}
\drawblock{5}{9}{0}{4}{blue!20}{0}
\foreach \n/\i/\txt in {
1/0/7,
8/4/2
 }{
    \cell{\n}{\i}{\txt}{black!10}
  }

\foreach \n/\i/\txt in {
2/0/13,2/1/13,
3/0/14,3/1/49,3/2/35,
4/0/14,4/1/31,4/2/85,4/3/68,
5/0/14,5/1/14,
6/0/14,6/1/14,6/3/45,
7/0/14,7/1/14,
8/0/14,8/1/14,
9/0/14,9/1/14
 }{
    \cell{\n}{\i}{\txt}{black!10}
  }

\foreach \n/\i/\txt in {
5/3/427,5/4/427,
7/4/147
 }{
    \smallcell{\n}{\i}{\txt}{black!10}
  }

  \foreach \n/\i/\txt in {
6/4/1233
 }{
    \tinycell{\n}{\i}{\txt}{black!10}
  }
\drawblock{3}{9}{0}{0}{blue!20}{14}
\drawblock{5}{9}{1}{1}{blue!20}{14}
\end{tikzpicture}
\begin{tikzpicture}[x=1cm,y=1cm,scale=0.55]
\drawnaxis{7}
\drawiaxis{3}
\drawstaircase{1}{4}{0}{black!10}{0}
\drawblock{5}{7}{0}{3}{blue!20}{0}
\foreach \n/\i/\txt in {
6/2/3
 }{
    \cell{\n}{\i}{\txt}{black!10}
  }

\foreach \n/\i/\txt in {
1/0/11,
2/0/21,2/1/21,
3/0/22,
4/0/22,4/1/53,
5/0/22,5/1/22,5/2/33,
6/0/22,6/1/22,
7/0/22,7/1/22
 }{
    \cell{\n}{\i}{\txt}{black!10}
  }

\foreach \n/\i/\txt in {
3/1/143,3/2/121,
4/2/271,4/3/240,
6/3/158
 }{
    \smallcell{\n}{\i}{\txt}{black!10}
  }

  \foreach \n/\i/\txt in {
5/3/4103
 }{
    \tinycell{\n}{\i}{\txt}{black!10}
  }
\drawblock{3}{7}{0}{0}{blue!20}{22}
\drawblock{5}{7}{1}{1}{blue!20}{22}
\end{tikzpicture}
\begin{tikzpicture}[x=1cm,y=1cm,scale=0.55]
\drawnaxis{7}
\drawiaxis{2}
\drawstaircase{1}{3}{0}{black!10}{0}
\drawblock{4}{7}{0}{2}{blue!20}{0}

\foreach \n/\i/\txt in {
1/0/25
 }{
    \cell{\n}{\i}{\txt}{black!10}
  }
  
\foreach \n/\i/\txt in {
2/0/145,2/1/145,
3/0/275,
4/0/299,
6/2/456
 }{
    \smallcell{\n}{\i}{\txt}{black!10}
  }

\foreach \n/\i/\txt in {
5/0/300,5/1/300,
6/0/300,6/1/300,
7/0/300,7/1/300
 }{
    \smallcell{\n}{\i}{\txt}{blue!20}
  }

\foreach \n/\i/\txt in {
3/1/1825,3/2/1550,
4/1/2089,
5/2/6875
 }{
    \tinycell{\n}{\i}{\txt}{black!10}
  }

\foreach \n/\i/\txt in {
4/2/10966
 }{
    \verytinycell{\n}{\i}{\txt}{black!10}
  }
\end{tikzpicture}
\begin{tikzpicture}[x=1cm,y=1cm,scale=0.55]
\drawnaxis{5}
\drawiaxis{1}
\drawstaircase{1}{2}{0}{black!10}{0}
\drawblock{3}{5}{0}{1}{blue!20}{0}

\foreach \n/\i/\txt in {
1/0/49
 }{
    \cell{\n}{\i}{\txt}{black!10}
  }
  
\foreach \n/\i/\txt in {
2/0/385,2/1/385,
3/0/735,
4/0/783
 }{
    \smallcell{\n}{\i}{\txt}{black!10}
  }

\foreach \n/\i/\txt in {
5/0/784,5/1/784
 }{
    \smallcell{\n}{\i}{\txt}{blue!20}
  }

\foreach \n/\i/\txt in {
4/1/7909
 }{
    \tinycell{\n}{\i}{\txt}{black!10}
  }

\foreach \n/\i/\txt in {
3/1/11221
 }{
    \verytinycell{\n}{\i}{\txt}{black!10}
  }
\end{tikzpicture}
  \caption{Homology of $\Hur^{H\rtimes \Z/2\Z,c^{\text{rfl}}}_n$ in the non-$3$-torsion case, for $H=\Z/5\Z, \Z/7\Z,\Z/11\Z,(\Z/5\Z)^2,(\Z/7\Z)^2$ respectively.}\label{fig: homology computation non 3 torsion}
\end{figure}

\begin{figure}[!htbp]
    \centering
    \begin{tikzpicture}[x=1cm,y=1cm,scale=0.55]
\drawnaxis{15}
\drawiaxis{3}
\drawstaircase{1}{4}{0}{black!10}{0}
\drawblock{5}{15}{0}{3}{blue!20}{0}
\foreach \n/\i/\txt in {
1/0/9,
8/2/6,
14/3/6
 }{
    \cell{\n}{\i}{\txt}{black!10}
  }

\foreach \n/\i/\txt in {
2/0/21,2/1/21,
3/0/27,3/1/99,3/2/72,
4/0/27,4/1/79,
5/0/27,5/1/63,5/2/81,
6/0/27,6/1/57,6/2/61,
7/0/27,7/1/45,7/2/36,7/3/72,
8/0/27,8/1/33,8/3/30,
9/0/27,9/1/27,9/3/18,
10/0/27,10/1/27,10/3/30,
11/0/27,11/1/27,11/3/36,
12/0/27,12/1/27,12/3/30,
13/0/27,13/1/27,13/3/18,
14/0/27,14/1/27,
15/0/27,15/1/27
 }{
    \cell{\n}{\i}{\txt}{black!10}
  }

\foreach \n/\i/\txt in {
4/2/230,4/3/178,
6/3/173
 }{
    \smallcell{\n}{\i}{\txt}{black!10}
  }

  \foreach \n/\i/\txt in {
5/3/1548
 }{
    \tinycell{\n}{\i}{\txt}{black!10}
  }
\drawblock{3}{15}{0}{0}{blue!20}{27}
\drawblock{9}{15}{1}{1}{blue!20}{27}
\end{tikzpicture}
\begin{tikzpicture}[x=1cm,y=1cm,scale=0.55]
\drawnaxis{9}
\drawiaxis{2}
\drawstaircase{1}{3}{0}{black!10}{0}
\drawblock{4}{9}{0}{2}{blue!20}{0}

\foreach \n/\i/\txt in {
1/0/15,
2/0/45,2/1/45,
3/0/60,
4/0/60,
5/0/60,
6/0/60,
7/0/60,7/1/90,7/2/60,
8/0/60,8/1/70,8/2/10,
9/0/60,9/1/60
 }{
    \cell{\n}{\i}{\txt}{black!10}
  }

\foreach \n/\i/\txt in {
3/1/375,3/2/315,
4/1/221,
5/1/120,5/2/315,
6/1/110,6/2/115
 }{
    \smallcell{\n}{\i}{\txt}{black!10}
  }

  \foreach \n/\i/\txt in {
4/2/1135
 }{
    \tinycell{\n}{\i}{\txt}{black!10}
  }
\drawblock{3}{9}{0}{0}{blue!20}{60}
\drawblock{9}{9}{1}{1}{blue!20}{60}
\end{tikzpicture}
\begin{tikzpicture}[x=1cm,y=1cm,scale=0.55]
\drawnaxis{15}
\drawiaxis{3}
\drawstaircase{1}{4}{0}{black!10}{0}
\drawblock{5}{15}{0}{3}{blue!20}{0}
\foreach \n/\i/\txt in {
1/0/9
 }{
    \cell{\n}{\i}{\txt}{black!10}
  }

\foreach \n/\i/\txt in {
2/0/33,2/1/33,
3/0/63,3/2/90,
4/0/71,
5/0/72,
6/0/72,
7/0/72,
8/0/72,8/2/75,
9/0/72,9/1/72,
10/0/72,10/1/72,
11/0/72,11/1/72,
12/0/72,12/1/72,
13/0/72,13/1/72,
14/0/72,14/1/72,14/3/48,
15/0/72,15/1/72
 }{
    \cell{\n}{\i}{\txt}{black!10}
  }

\foreach \n/\i/\txt in {
3/1/153,
4/1/273,4/2/570,4/3/368,
5/1/351,
6/1/312,
7/1/216,7/2/612,
8/1/120,8/3/840,
9/3/180,
10/3/244,
11/3/288,
12/3/240,
13/3/144
 }{
    \smallcell{\n}{\i}{\txt}{black!10}
  }

  \foreach \n/\i/\txt in {
5/2/1206,5/3/3114,
6/2/1114,6/3/4818,
7/3/3402
 }{
    \tinycell{\n}{\i}{\txt}{black!10}
  }
\drawblock{5}{15}{0}{0}{blue!20}{72}
\drawblock{9}{15}{1}{1}{blue!20}{72}
\end{tikzpicture}
\begin{tikzpicture}[x=1cm,y=1cm,scale=0.55]
\drawnaxis{9}
\drawiaxis{2}
\drawstaircase{1}{3}{0}{black!10}{0}
\drawblock{4}{9}{0}{2}{blue!20}{0}
\foreach \n/\i/\txt in {
1/0/27
 }{
    \cell{\n}{\i}{\txt}{black!10}
  }

\foreach \n/\i/\txt in {
2/0/153,2/1/153,
3/0/324,
4/0/372,
8/1/657,8/2/441,
 }{
    \smallcell{\n}{\i}{\txt}{black!10}
  }

\foreach \n/\i/\txt in {
5/0/378,
6/0/378,
7/0/378,
8/0/378,
9/0/378,9/1/378
 }{
    \smallcell{\n}{\i}{\txt}{blue!20}
  }

  \foreach \n/\i/\txt in {
3/1/2187,3/2/1863,
4/1/2373,
5/1/1998,5/2/8505,
6/1/1773,6/2/6687,
7/1/1215,7/2/3618
 }{
    \tinycell{\n}{\i}{\txt}{black!10}
  }

   \foreach \n/\i/\txt in {
4/2/12981
 }{
    \verytinycell{\n}{\i}{\txt}{black!10}
  }
\end{tikzpicture}




\caption{Homology of $\Hur^{H\rtimes \Z/2\Z,c^{\text{rfl}}}_n$ in the $3$-torsion case, for $H=\Z/9\Z, \Z/15\Z,(\Z/3\Z)^2,\Z/3\Z\times \Z/9\Z$ respectively.}\label{fig: homology computation 3-torsion}
\end{figure}

We also performed computations for additional choices of $H$, which we omit from the figures for brevity but describe here. In the case without $3$-torsion, we considered $H=\Z/k\Z$ with $3\nmid k$. For $k=13$, we verified that $H_0,H_1,H_2,H_3$ stabilize at $n=3,5,7,7$, respectively, matching the pattern observed for $k=11$ in Figure \ref{fig: homology computation non 3 torsion}. For $17\le k\le 25$, we observed the same pattern for $H_0,H_1,H_2$, and for $29\le k\le 49$, we verified it for $H_0,H_1$. In the $3$-torsion case, we verified that for $H=\Z/21\Z$ and $\Z/27\Z$, the groups $H_0,H_1,H_2$ stabilize at $n=3,9,9$, respectively, matching the pattern observed for $H=\Z/15\Z$ in Figure \ref{fig: homology computation 3-torsion}.

In all of these cases, $c$ is non-splitting, so $\Hur^{G,c}_n$ is homologically stable with respect to $U=\sum_{g\in c} g^2$. To get the bounds in the introduction, we first verify that we have stability in an $(N,I)$-staircase with respect to a sufficiently large power $U^e$, just like in the Davenport-Heilbronn case discussed in the introduction. This can be done by simply reading off the dimensions of the homology that we computed, with the $H_0$ stabilization being a consequence of Theorem \ref{thm: cohen lenstra stab R}. Then, applying Theorem \ref{thm: main theorem} gives us the desired stability slope, and hence the exponent bound via the Grothendieck-Lefschetz trace formula like in the proof of Theorem \ref{thm: cohen-lenstra H}.

\begin{figure}[!htbp]
    \centering
    \begin{tikzpicture}[x=1cm,y=1cm,scale=0.55]
\drawnaxis{7}
\drawiaxis{3}
\drawstaircase{1}{4}{0}{black!10}{0}
\drawblock{5}{7}{0}{3}{blue!20}{0}
\foreach \n/\i/\txt in {
1/0/9,
5/2/9,
6/2/1
 }{
    \cell{\n}{\i}{\txt}{black!10}
  }

\foreach \n/\i/\txt in {
2/0/21,2/1/21,
3/0/27,3/1/81,3/2/54,
4/0/27,4/1/49,
5/0/27,5/1/27,
6/0/27,6/1/27,6/3/71,
7/0/27,7/1/27
 }{
    \cell{\n}{\i}{\txt}{black!10}
  }

\foreach \n/\i/\txt in {
4/2/164,4/3/142
 }{
    \smallcell{\n}{\i}{\txt}{black!10}
  }

  \foreach \n/\i/\txt in {
5/3/1458
 }{
    \tinycell{\n}{\i}{\txt}{black!10}
  }
\drawblock{3}{7}{0}{0}{blue!20}{27}
\drawblock{5}{7}{1}{1}{blue!20}{27}
\end{tikzpicture}
\begin{tikzpicture}[x=1cm,y=1cm,scale=0.55]
\drawnaxis{7}
\drawiaxis{2}
\drawstaircase{1}{3}{0}{black!10}{0}
\drawblock{4}{7}{0}{2}{blue!20}{0}

\foreach \n/\i/\txt in {
1/0/15,
2/0/45,2/1/45,
3/0/60,
4/0/60,
5/0/60,5/1/60,
6/0/60,6/1/60,6/2/15,
7/0/60,7/1/60
 }{
    \cell{\n}{\i}{\txt}{black!10}
  }

\foreach \n/\i/\txt in {
3/1/345,3/2/285,
4/1/171,
5/2/195
 }{
    \smallcell{\n}{\i}{\txt}{black!10}
  }

  \foreach \n/\i/\txt in {
4/2/1025
 }{
    \tinycell{\n}{\i}{\txt}{black!10}
  }
\drawblock{3}{7}{0}{0}{blue!20}{60}
\drawblock{5}{7}{1}{1}{blue!20}{60}
\end{tikzpicture}
\begin{tikzpicture}[x=1cm,y=1cm,scale=0.55]
\drawnaxis{11}
\drawiaxis{3}
\drawstaircase{1}{4}{0}{black!10}{0}
\drawblock{5}{11}{0}{3}{blue!20}{0}
\foreach \n/\i/\txt in {
1/0/9,
10/3/4
 }{
    \cell{\n}{\i}{\txt}{black!10}
  }

\foreach \n/\i/\txt in {
2/0/33,2/1/33,
3/0/63,3/1/81,3/2/18,
4/0/71,4/1/93,
5/0/72,5/1/81,
6/0/72,6/1/77,
7/0/72,7/1/72,
8/0/72,8/1/72,8/2/27,
9/0/72,9/1/72,9/3/36,
10/0/72,10/1/72,
11/0/72,11/1/72
 }{
    \cell{\n}{\i}{\txt}{black!10}
  }

\foreach \n/\i/\txt in {
4/2/174,4/3/152,
5/2/666,
6/2/644,
7/2/324,
8/3/600
 }{
    \smallcell{\n}{\i}{\txt}{black!10}
  }

  \foreach \n/\i/\txt in {
5/3/2439,
6/3/4019,
7/3/2826
 }{
    \tinycell{\n}{\i}{\txt}{black!10}
  }
\drawblock{5}{11}{0}{0}{blue!20}{72}
\drawblock{7}{11}{1}{1}{blue!20}{72}
\end{tikzpicture}
\begin{tikzpicture}[x=1cm,y=1cm,scale=0.55]
\drawnaxis{9}
\drawiaxis{2}
\drawstaircase{1}{3}{0}{black!10}{0}
\drawblock{4}{9}{0}{2}{blue!20}{0}
\foreach \n/\i/\txt in {
1/0/27
 }{
    \cell{\n}{\i}{\txt}{black!10}
  }

\foreach \n/\i/\txt in {
2/0/153,2/1/153,
3/0/324,
4/0/372,
5/1/432,
6/1/408,
8/2/162,
 }{
    \smallcell{\n}{\i}{\txt}{black!10}
  }

\foreach \n/\i/\txt in {
5/0/378,
6/0/378,
7/0/378,7/1/378,
8/0/378,8/1/378,
9/0/378,9/1/378
 }{
    \smallcell{\n}{\i}{\txt}{blue!20}
  }

  \foreach \n/\i/\txt in {
3/1/1782,3/2/1458,
4/1/1338,
5/2/5373,
6/2/3957,
7/2/1944
 }{
    \tinycell{\n}{\i}{\txt}{black!10}
  }

   \foreach \n/\i/\txt in {
4/2/10704
 }{
    \verytinycell{\n}{\i}{\txt}{black!10}
  }
\end{tikzpicture}




\caption{Homology of $\Hur^{H\rtimes \Z/2\Z,c^{\text{rfl}}}_n$ in the $3$-torsion case, with contribution from the reduced homology of $\CHur^{S_3,c^{\text{trn}}}_n$ removed, for $H=\Z/9\Z, \Z/15\Z,\Z/21\Z,(\Z/3\Z)^2,\Z/3\Z\times \Z/9\Z$ respectively.}\label{fig: homology computation 3-torsion with S3 removed}

\end{figure}

Next, we briefly sketch how we remove the contribution of the $S_3$ transposition homology in Figure \ref{fig: homology computation 3-torsion with S3 removed}. The code is available in \textit{hurwitz\_s3\_removal}, with a more detailed explanation in the accompanying writeup. Recall that the Hurwitz spaces $\Hur^{G,c}_n$ and $\Hur^{S_3,c^{\text{trn}}_{S_3}}_n$ define finite covers of $\Conf_n$, and hence give rise to finite-dimensional permutation representations of $B_n$, which we denote by $P_n$ and $Q_n$, respectively. We define the reduced $S_3$ transposition representation $\widetilde Q_n$ by removing all trivial subrepresentations from $Q_n$, and call the corresponding local system the reduced $S_3$ local system.

We then remove the maximal number $m_n$ of copies of $\widetilde Q_n$ from $P_n$, writing $P_n\cong \widetilde Q_n^{\oplus m_n}\oplus\widetilde P_n$ where $\widetilde P_n$ contains no further copy of $\widetilde Q_n$. This can be computed using finite-group representation theory, since the braid actions factor through finite groups. More specifically, using the Burau representation discussed in Section \ref{sec: cohen lenstra}, we reduce the computation to finite symplectic groups and their stabilizers. The multiplicities of the relevant irreducible constituents can then be recovered from traces on the center of a small endomorphism algebra, without explicitly decomposing the exponentially large permutation representations. We also observe from our computations that the multiplicities $m_n$ eventually stabilize as $n$ increases.

\subsection{Squarefree $S_d$-extensions}\label{sec: compute S_d} 

In this subsection, we present the homological computations in Figures \ref{fig: homology computation S_d}, \ref{fig: compute MA quotient 1}, and \ref{fig: compute MA quotient 2}, and explain how they establish the bounds in Figure \ref{fig: improved c_d} and support the surrounding discussion.

\begin{figure}[!htbp]
    \centering
 \begin{tikzpicture}[x=1cm,y=1cm,scale=0.45]
\drawnaxis{33}
\drawiaxis{9}
\drawstaircase{1}{10}{0}{black!10}{0}
\drawblock{11}{33}{0}{9}{black!10}{0}
\drawblock{3}{33}{0}{0}{blue!20}{3}
\drawblock{9}{33}{1}{1}{blue!20}{3}
\drawblock{9}{33}{2}{2}{blue!20}{0}
\drawblock{15}{33}{3}{4}{blue!20}{0}
\drawblock{21}{33}{5}{6}{blue!20}{0}
\drawblock{27}{33}{7}{8}{blue!20}{0}
\drawblock{33}{33}{9}{9}{blue!20}{0}
\foreach \n/\i/\txt in {
2/0/2,2/1/2,
3/1/6,3/2/3,
4/1/8,4/3/6,
5/1/9,5/4/9,
6/1/8,
7/1/6,7/2/6,
8/1/4,8/2/1,8/3/5,
9/3/3,9/4/9,
10/3/5,
11/3/6,
12/3/5,
13/3/3,13/4/6,
14/3/1,14/4/1,14/5/5,
15/5/3,15/6/9,
16/5/5,
17/5/6,
18/5/5,
19/5/3,19/6/6,
20/5/1,20/6/1,20/7/5,
21/7/3,21/8/9,
22/7/5,
23/7/6,
24/7/5,
25/7/3,25/8/6,
26/7/1,26/8/1,26/9/5,
27/9/3,
28/9/5,
29/9/6,
30/9/5,
31/9/3,
32/9/1
 }{
    \cell{\n}{\i}{\txt}{black!10}
  }

\foreach \n/\i/\txt in {
4/2/11,
5/2/12,5/3/15,
6/2/10,6/3/17,6/4/41,6/5/29,
7/3/12,7/4/33,7/5/72,7/6/48,
8/4/21,8/5/75,
9/5/30,
10/4/10,10/5/16,10/6/43,
11/4/12,11/5/15,11/6/27,11/7/63,
12/4/10,12/5/15,12/6/33,12/7/73,
13/5/12,13/6/33,13/7/78,
14/6/21,14/7/58,
15/7/30,15/8/90,
16/6/10,16/7/16,16/8/40,
17/6/12,17/7/15,17/8/27,17/9/63,
18/6/10,18/7/15,18/8/33,18/9/73,
19/7/12,19/8/33,19/9/78,
20/8/21,20/9/58,
21/9/30,
22/8/10,22/9/16,
23/8/12,23/9/15,
24/8/10,24/9/15,
25/9/12
 }{
    \smallcell{\n}{\i}{\txt}{black!10}
  }

\foreach \n/\i/\txt in {
8/6/167,8/7/108,
9/6/120,9/7/300,9/8/204,
10/7/198,10/8/718,10/9/552,
11/8/333,
12/8/165,12/9/748,
13/8/165,13/9/354,
14/8/138,14/9/300,
15/9/222,
16/9/109
 }{
    \tinycell{\n}{\i}{\txt}{black!10}
  }

  \foreach \n/\i/\txt in {
11/9/1353
 }{
    \verytinycell{\n}{\i}{\txt}{black!10}
  }

\end{tikzpicture}
 \begin{tikzpicture}[x=1cm,y=1cm,scale=0.55]
\drawnaxis{23}
\drawiaxis{5}
\drawstaircase{1}{6}{0}{black!10}{0}
\drawblock{7}{23}{0}{5}{black!10}{0}
\drawblock{5}{23}{0}{0}{blue!20}{12}
\drawblock{11}{23}{1}{1}{blue!20}{12}
\drawblock{11}{23}{2}{2}{blue!20}{0}
\drawblock{17}{23}{3}{4}{blue!20}{0}
\drawblock{23}{23}{5}{5}{blue!20}{0}
\foreach \n/\i/\txt in {
3/0/6,3/2/6,
9/2/6,
10/2/1,
15/3/6,15/4/6,
16/3/1,16/4/1,
21/5/6,
22/5/1
 }{
    \cell{\n}{\i}{\txt}{black!10}
  }

\foreach \n/\i/\txt in {
3/1/12,
4/0/11,4/1/36,4/2/67,4/3/42,
5/1/54,
6/1/58,
7/1/48,
8/1/31,8/2/39,
9/1/18,9/3/36,
10/1/13,10/3/70,
11/3/78,
12/3/64,
13/3/42,
14/3/20,14/4/23,14/5/61,
15/5/60,
16/5/77,
17/5/78,
18/5/64,
19/5/42,
20/5/20
 }{
    \cell{\n}{\i}{\txt}{black!10}
  }

\foreach \n/\i/\txt in {
5/2/132,5/3/270,5/4/180,
6/2/151,6/3/397,6/5/857,
7/2/114,7/3/324,
8/3/122,8/4/659,
9/4/204,
10/4/268,10/5/739,
11/4/288,11/5/744,
12/4/210,12/5/522,
13/4/102,13/5/204
 }{
    \smallcell{\n}{\i}{\txt}{black!10}
  }

  \foreach \n/\i/\txt in {
6/4/1149,
7/4/1134,7/5/3984,
8/5/3177,
9/5/1122
 }{
    \tinycell{\n}{\i}{\txt}{black!10}
  }

\end{tikzpicture}

 \begin{tikzpicture}[x=1cm,y=1cm,scale=0.55]
\drawnaxis{19}
\drawiaxis{4}
\drawstaircase{1}{5}{0}{black!10}{0}
\drawblock{6}{19}{0}{4}{black!10}{0}
\drawblock{7}{19}{0}{0}{blue!20}{60}
\drawblock{13}{19}{1}{1}{blue!20}{60}
\drawblock{13}{19}{2}{2}{blue!20}{0}
\drawblock{19}{19}{3}{4}{blue!20}{0}
\foreach \n/\i/\txt in {
12/2/1,
18/3/1,18/4/1
 }{
    \cell{\n}{\i}{\txt}{black!10}
  }

\foreach \n/\i/\txt in {
4/0/24,4/1/48,4/2/72,4/3/48,
5/0/50,
6/0/59,
10/2/57,
11/1/70,11/2/10,
12/1/61,
16/3/41,16/4/42,
17/3/10,17/4/10
 }{
    \cell{\n}{\i}{\txt}{black!10}
  }

\foreach \n/\i/\txt in {
5/1/150,5/2/350,5/4/760,
6/1/247,6/2/719,
7/1/280,7/2/800,
8/1/235,8/2/531,
9/1/160,9/2/210,9/3/490,
10/1/101,10/3/409,
11/3/440,
12/3/368,
13/3/270,13/4/590,
14/3/181,14/4/240,
15/3/100,15/4/110
 }{
    \smallcell{\n}{\i}{\txt}{black!10}
  }

  \foreach \n/\i/\txt in {
5/3/1010,
6/3/2173,
7/3/2510,7/4/9040,
8/3/1477,8/4/6082,
9/4/2370,
10/4/2054,
11/4/2050,
12/4/1372
 }{
    \tinycell{\n}{\i}{\txt}{black!10}
  }

 \foreach \n/\i/\txt in {
6/4/10490
 }{
    \verytinycell{\n}{\i}{\txt}{black!10}
  }

\end{tikzpicture}
 \begin{tikzpicture}[x=1cm,y=1cm,scale=0.55]
\drawnaxis{21}
\drawiaxis{3}
\drawstaircase{1}{4}{0}{black!10}{0}
\drawblock{5}{21}{0}{3}{blue!20}{0}
\foreach \n/\i/\txt in {
14/2/1
 }{
    \cell{\n}{\i}{\txt}{black!10}
  }

\foreach \n/\i/\txt in {
13/2/15,
18/3/92,
19/3/15,
20/3/1
 }{
    \cell{\n}{\i}{\txt}{black!10}
  }

\foreach \n/\i/\txt in {
5/0/120,5/1/240,5/2/360,
6/0/274,6/1/822,
7/0/345,
8/0/359,
11/1/675,11/2/555,
12/1/452,12/2/114,
13/1/375,
14/1/361,
16/3/722,
17/3/315
 }{
    \smallcell{\n}{\i}{\txt}{black!10}
  }

\foreach \n/\i/\txt in {
9/0/360,
10/0/360,
11/0/360,
12/0/360,
13/0/360,
14/0/360,
15/0/360,15/1/360,
16/0/360,16/1/360,
17/0/360,17/1/360,
18/0/360,18/1/360,
19/0/360,19/1/360,
20/0/360,20/1/360,
21/0/360,21/1/360
 }{
    \smallcell{\n}{\i}{\txt}{blue!20}
  }
  
  \foreach \n/\i/\txt in {
5/3/1680,
6/2/1918,6/3/6184,
7/1/1485,7/2/4080,
8/1/1776,8/2/4893,
9/1/1545,9/2/3675,9/3/9495,
10/1/1075,10/2/1772,10/3/4854,
11/3/3270,
12/3/2865,
13/3/2445,
14/3/1884,
15/3/1275
 }{
    \tinycell{\n}{\i}{\txt}{black!10}
  }

 \foreach \n/\i/\txt in {
7/3/13155,
8/3/14848
 }{
    \verytinycell{\n}{\i}{\txt}{black!10}
  }
  
\end{tikzpicture}
\begin{tikzpicture}[x=1cm,y=1cm,scale=0.55]
\drawnaxis{17}
\drawiaxis{2}
\drawstaircase{1}{3}{0}{black!10}{0}
\drawblock{4}{17}{0}{2}{black!10}{0}
\drawblock{17}{17}{2}{2}{blue!20}{0}
\foreach \n/\i/\txt in {
16/2/1
 }{
    \cell{\n}{\i}{\txt}{black!10}
  }

\foreach \n/\i/\txt in {
15/2/21
 }{
    \cell{\n}{\i}{\txt}{black!10}
  }

\foreach \n/\i/\txt in {
6/0/720,
14/2/212
 }{
    \smallcell{\n}{\i}{\txt}{black!10}
  }

  \foreach \n/\i/\txt in {
6/1/1440,6/2/2160,
7/0/1764,7/1/5292,
8/0/2344,
9/0/2499,
10/0/2519,
11/1/8295,
12/1/5216,12/2/5721,
13/1/3402,13/2/1344,
14/1/2703,
15/1/2541,
16/1/2521
 }{
    \tinycell{\n}{\i}{\txt}{black!10}
  }

 \foreach \n/\i/\txt in {
11/0/2520,
12/0/2520,
13/0/2520,
14/0/2520,
15/0/2520,
16/0/2520,
17/0/2520,17/1/2520
 }{
    \tinycell{\n}{\i}{\txt}{blue!20}
  }

\foreach \n/\i/\txt in {
7/2/12348,
8/1/10280,8/2/28864,
9/1/12936,9/2/38073,
10/1/11660,10/2/31184,
11/2/16422
 }{
    \verytinycell{\n}{\i}{\txt}{black!10}
  }
  
\end{tikzpicture}
    \caption{Homology of $\CHur^{S_d,c^{\text{trn}}_{S_d}}_n$ for $d=3,4,5,6,7$ respectively.}\label{fig: homology computation S_d}
\end{figure}

\begin{figure}[!htbp]
    \centering
 \begin{tikzpicture}[x=1cm,y=1cm,scale=0.45]
\drawnaxis{13}
\drawiaxis{6}
\drawstaircase{1}{7}{0}{black!10}{0}
\drawblock{8}{13}{0}{6}{black!10}{0}
\drawblock{5}{13}{0}{0}{blue!20}{0}
\drawblock{6}{13}{1}{1}{blue!20}{0}
\drawblock{7}{13}{2}{2}{blue!20}{0}
\drawblock{9}{13}{3}{3}{blue!20}{0}
\drawblock{10}{13}{4}{4}{blue!20}{0}
\drawblock{12}{13}{5}{5}{blue!20}{0}
\drawblock{13}{13}{6}{6}{blue!20}{0}
\foreach \n/\i/\txt in {
 2/0/1,2/1/1,
3/0/2,3/1/4,3/2/2,
4/0/1,4/1/5,4/2/7,4/3/3,
5/1/2,5/2/8,5/4/6,
6/2/3,
7/3/8,
8/3/1,
9/4/6,
11/5/2
 }{
    \cell{\n}{\i}{\txt}{black!10}
  }

\foreach \n/\i/\txt in {
5/3/12,
6/3/16,6/4/24,6/5/11,
7/4/36,7/5/46,7/6/18,
8/4/24,8/5/77,8/6/84,
9/5/62,
10/5/22,
11/6/64,
12/6/12
 }{
    \smallcell{\n}{\i}{\txt}{black!10}
  }

\foreach \n/\i/\txt in {
9/6/154,
10/6/140
 }{
    \tinycell{\n}{\i}{\txt}{black!10}
  }

\end{tikzpicture}
\begin{tikzpicture}[x=1cm,y=1cm,scale=0.45]
\drawnaxis{14}
\drawiaxis{5}
\drawstaircase{1}{6}{0}{black!10}{0}
\drawblock{7}{14}{0}{5}{black!10}{0}
\drawblock{7}{14}{0}{0}{blue!20}{0}
\drawblock{8}{14}{1}{1}{blue!20}{0}
\drawblock{9}{14}{2}{2}{blue!20}{0}
\drawblock{11}{14}{3}{3}{blue!20}{0}
\drawblock{12}{14}{4}{4}{blue!20}{0}
\drawblock{14}{14}{5}{5}{blue!20}{0}

\foreach \n/\i/\txt in {
3/0/1,3/1/3,3/2/2,
4/0/3,4/3/9,
5/0/3,
6/0/1,
7/1/3,
8/2/9,
10/3/3
 }{
    \cell{\n}{\i}{\txt}{black!10}
  }

\foreach \n/\i/\txt in {
4/1/12,4/2/18,
5/1/18,5/2/51,5/3/66,5/4/30,
6/1/12,6/2/65,6/5/91,
7/2/39,
9/3/44,
11/4/33,
13/5/12
 }{
    \smallcell{\n}{\i}{\txt}{black!10}
  }

\foreach \n/\i/\txt in {
6/3/179,6/4/216,
7/3/234,7/4/606,7/5/666,
8/3/153,8/4/855,
9/4/639,
10/4/237,
12/5/244
 }{
    \tinycell{\n}{\i}{\txt}{black!10}
  }

  \foreach \n/\i/\txt in {
8/5/1956,
9/5/2990,
10/5/2541,
11/5/1164
 }{
    \verytinycell{\n}{\i}{\txt}{black!10}
  }

\end{tikzpicture}
\begin{tikzpicture}[x=1cm,y=1cm,scale=0.55]
\drawnaxis{25}
\drawiaxis{6}
\drawstaircase{1}{6}{0}{black!10}{0}
\drawblock{7}{25}{0}{6}{blue!20}{0}

\foreach \n/\i/\txt in {
 3/0/2,3/1/4,3/2/2,
4/0/5,
5/0/4,
6/0/1,
7/1/5,
8/1/1,
9/2/8,
10/2/7,
11/2/4,
12/2/1,
13/3/5,
14/3/1,
15/4/8,
16/4/7,
17/4/4,
18/4/1,
19/5/5,
20/5/1,
21/6/8,
22/6/7,
23/6/4,
24/6/1
 }{
    \cell{\n}{\i}{\txt}{black!10}
  }

\foreach \n/\i/\txt in {
4/1/15,4/2/21,4/3/11,
5/1/21,5/2/54,5/3/79,5/4/42,
6/1/14,6/2/57,
7/2/31,
8/2/13,8/3/83,
9/3/42,
10/3/35,
11/3/26,
12/3/14,12/4/65,
13/4/31,
14/4/13,14/5/75,
15/5/42,
16/5/35,
17/5/26,
18/5/14,18/6/65,
19/6/31,
20/6/13
 }{
    \cell{\n}{\i}{\txt}{black!10}
  }

\foreach \n/\i/\txt in {
6/3/167,6/4/256,6/5/133,
7/3/159,7/4/539,7/5/832,7/6/426,
8/4/493,
9/4/247,
10/4/138,10/5/777,
11/4/102,11/5/419,
12/5/284,
13/5/160,13/6/771,
14/6/410,
15/6/205,
16/6/135,
17/6/102
 }{
    \smallcell{\n}{\i}{\txt}{black!10}
  }

\foreach \n/\i/\txt in {
8/5/1630,8/6/2480,
9/5/1465,9/6/4773,
10/6/4230,
11/6/2288,
12/6/1267
 }{
    \tinycell{\n}{\i}{\txt}{black!10}
  }

\end{tikzpicture}

\begin{tikzpicture}[x=1cm,y=1cm,scale=0.55]
\drawnaxis{21}
\drawiaxis{4}
\drawstaircase{1}{5}{0}{black!10}{0}
\drawblock{6}{21}{0}{4}{blue!20}{0}

\foreach \n/\i/\txt in {
4/0/6,
7/0/7,
8/0/1,
9/1/8,
10/1/1,
13/2/7,
14/2/1,
15/3/7,
16/3/1,
19/4/7,
20/4/1
 }{
    \cell{\n}{\i}{\txt}{black!10}
  }

\foreach \n/\i/\txt in {
4/1/18,4/2/30,4/3/18,
5/0/17,5/1/71,
6/0/17,
7/1/89,
8/1/37,
10/2/73,
11/2/49,
12/2/24,
13/3/96,
14/3/29,
16/4/90,
17/4/56,
18/4/25
 }{
    \cell{\n}{\i}{\txt}{black!10}
  }

\foreach \n/\i/\txt in {
5/2/187,5/3/333,5/4/200,
6/1/112,6/2/408,
7/2/434,
8/2/262,
9/2/121,9/3/785,
10/3/504,
11/3/380,
12/3/231,
13/4/581,
14/4/241,
15/4/128
 }{
    \smallcell{\n}{\i}{\txt}{black!10}
  }

\foreach \n/\i/\txt in {
6/3/1192,6/4/2183,
7/3/1818,7/4/6312,
8/3/1444,8/4/7561,
9/4/4945,
10/4/2834,
11/4/1927,
12/4/1221
 }{
    \tinycell{\n}{\i}{\txt}{black!10}
  }

\end{tikzpicture}

\begin{tikzpicture}[x=1cm,y=1cm,scale=0.55]
\drawnaxis{21}
\drawiaxis{4}
\drawstaircase{1}{5}{0}{black!10}{0}
\drawblock{6}{21}{0}{4}{blue!20}{0}

\foreach \n/\i/\txt in {
4/0/4,4/3/4,
7/0/6,
8/0/1,
10/1/2,
14/2/2,
16/3/4,
20/4/3
 }{
    \cell{\n}{\i}{\txt}{black!10}
  }

\foreach \n/\i/\txt in {
4/1/12,4/2/12,
5/0/12,5/1/48,5/2/96,5/4/60,
6/0/13,6/1/77,
7/1/66,
8/1/35,
9/1/12,
10/2/65,
11/2/48,
12/2/31,
13/2/12,
14/3/74,
15/3/24,
17/4/84,
18/4/49,
19/4/18
 }{
    \cell{\n}{\i}{\txt}{black!10}
  }

\foreach \n/\i/\txt in {
5/3/120,
6/2/250,6/3/558,6/4/957,
7/2/306,
8/2/216,
9/2/114,9/3/674,
10/3/415,
11/3/310,
12/3/245,
13/3/156,13/4/856,
14/4/501,
15/4/246,
16/4/130
 }{
    \smallcell{\n}{\i}{\txt}{black!10}
  }

\foreach \n/\i/\txt in {
7/3/1038,7/4/2616,
8/3/1023,8/4/4089,
9/4/3542,
10/4/2207,
11/4/1500,
12/4/1186
 }{
    \tinycell{\n}{\i}{\txt}{black!10}
  }

\end{tikzpicture}

    \caption{Homology of $MA_{c_{S_2}\sqcup c_{S_2}}/(s_1^2,s_2^2)$,  $MA_{c_{S_2}\sqcup c_{S_2}\sqcup c_{S_2}}/(s_1^2,s_2^2,s_3^2)$, $MA_{c_{S_2}\sqcup c_{S_3}}/(s_1^2,s_2^2)$, $MA_{c_{S_2}\sqcup c_{S_4}}/(s_1^2,s_2^2)$, $MA_{c_{S_3}\sqcup c_{S_3}}/(s_1^2,s_2^2)$ respectively.}\label{fig: compute MA quotient 1}
\end{figure}

\begin{figure}[!htbp]
    \centering
    \begin{tikzpicture}[x=1cm,y=1cm,scale=0.55]
\drawnaxis{21}
\drawiaxis{4}
\drawstaircase{1}{5}{0}{black!10}{0}
\drawblock{6}{21}{0}{4}{blue!20}{0}

\foreach \n/\i/\txt in {
4/0/2,4/1/8,4/3/6,
5/0/7,
6/0/9,
7/0/5,
8/0/1,
9/1/7,
10/1/1,
13/2/5,
14/2/1,
15/3/7,
16/3/1,
19/4/5,
20/4/1
 }{
    \cell{\n}{\i}{\txt}{black!10}
  }

\foreach \n/\i/\txt in {
4/2/12,
5/1/35,5/2/86,5/4/57,
6/1/61,
7/1/54,
8/1/26,
9/2/94,
10/2/32,
11/2/16,
12/2/11,
13/3/61,
14/3/26,
15/4/94,
16/4/32,
17/4/16,
18/4/11
 }{
    \cell{\n}{\i}{\txt}{black!10}
  }

\foreach \n/\i/\txt in {
5/3/115,
6/2/228,6/3/545,6/4/705,
7/2/299,
8/2/216,
9/3/859,
10/3/381,
11/3/169,
12/3/103,12/4/808,
13/4/458,
14/4/236
 }{
    \smallcell{\n}{\i}{\txt}{black!10}
  }

\foreach \n/\i/\txt in {
7/3/1151,7/4/2793,
8/3/1300,8/4/5356,
9/4/5687,
10/4/3629,
11/4/1668
 }{
    \tinycell{\n}{\i}{\txt}{black!10}
  }
\end{tikzpicture}
\begin{tikzpicture}[x=1cm,y=1cm,scale=0.55]
\drawnaxis{19}
\drawiaxis{3}
\drawstaircase{1}{4}{0}{black!10}{0}
\drawblock{5}{19}{0}{3}{blue!20}{0}

\foreach \n/\i/\txt in {
10/0/1,
12/1/1,
16/2/1,
18/3/1
 }{
    \cell{\n}{\i}{\txt}{black!10}
  }

\foreach \n/\i/\txt in {
5/0/24,5/1/72,
6/0/74,
7/0/85,
8/0/45,
9/0/11,
10/1/69,
11/1/12,
14/2/50,
15/2/11,
16/3/53,
17/3/11
 }{
    \cell{\n}{\i}{\txt}{black!10}
  }

\foreach \n/\i/\txt in {
5/2/168,5/3/384,
6/1/296,6/2/890,
7/1/508,
8/1/459,
9/1/233,
10/2/747,
11/2/401,
12/2/244,
13/2/131,
14/3/421,
15/3/164
 }{
    \smallcell{\n}{\i}{\txt}{black!10}
  }

\foreach \n/\i/\txt in {
6/3/2860,
7/2/1932,7/3/7399,
8/2/2180,8/3/9549,
9/2/1469,9/3/7479,
10/3/4622,
11/3/3041,
12/3/2015,
13/3/1043
 }{
    \tinycell{\n}{\i}{\txt}{black!10}
  }
\end{tikzpicture}
\begin{tikzpicture}[x=1cm,y=1cm,scale=0.55]
\drawnaxis{19}
\drawiaxis{3}
\drawstaircase{1}{4}{0}{black!10}{0}
\drawblock{5}{19}{0}{3}{blue!20}{0}

\foreach \n/\i/\txt in {
9/0/9,
10/0/1,
12/1/2,
16/2/2,
18/3/4
 }{
    \cell{\n}{\i}{\txt}{black!10}
  }

\foreach \n/\i/\txt in {
5/0/12,5/1/48,5/2/96,
6/0/40,
7/0/51,
8/0/31,
10/1/79,
11/1/18,
14/2/72,
15/2/18,
17/3/35
 }{
    \cell{\n}{\i}{\txt}{black!10}
  }

\foreach \n/\i/\txt in {
5/3/120,
6/1/206,6/2/624,
7/1/369,
8/1/361,
9/1/213,
10/2/848,
11/2/466,
12/2/293,
13/2/171,
15/3/473,
16/3/157
 }{
    \smallcell{\n}{\i}{\txt}{black!10}
  }

\foreach \n/\i/\txt in {
6/3/1436,
7/2/1536,7/3/5052,
8/2/1959,8/3/8465,
9/2/1511,9/3/8171,
10/3/5501,
11/3/3484,
12/3/2509,
13/3/1793,
14/3/1055
 }{
    \tinycell{\n}{\i}{\txt}{black!10}
  }
\end{tikzpicture}

    \caption{Homology of $MA_{c_{S_2}\sqcup c_{S_2}\sqcup c_{S_3}}/(s_1^2,s_2^2,s_3^2)$, $MA_{c_{S_2}\sqcup c_{S_5}}/(s_1^2,s_2^2)$, $MA_{c_{S_3}\sqcup c_{S_4}}/(s_1^2,s_2^2)$ respectively.}\label{fig: compute MA quotient 2}
\end{figure}

The main difficulty here is that the relevant rack $c=c^{\text{trn}}_{S_d}$ is splitting for $d>3$, and so we no longer have homological stability with respect to just a single element. One set of stability elements that we can use comes from Theorem \ref{thm: LL A_c stability}, which tells us that $A_c$ is stable with respect to $U_1\coloneqq s_1^2,\ldots, U_{|c|}\coloneqq s_{|c|}^2$ where $s_1,\ldots, s_{|c|}$ are the transpositions in $c$. If we were to verify stability in an $(N,I)$-staircase directly using Definition \ref{def: staircase iterated version}, we would need to compute the homology of a $|c|$-fold iterated quotient. This would take a long time, as the homotopy coherent hypercube corresponding to the iterated quotient has $2^{|c|}$ terms.

We present a better way to verify that $A_c$ stabilizes in an $(N,I)$-staircase with respect to $U_1^{e_1},\ldots, U_{|c|}^{e_{|c|}}$ for all sufficiently large $e_1,\ldots, e_{|c|}$ (note that $A_c$ is homologically stable with respect to these powers by Lemma \ref{lem: homological stability for powers}). We make use of some ideas presented in Section \ref{sec: S_d}, in particular, the filtration \eqref{eqn: filtration A_c and MA_c} of $A_c$ into maximal chains $MA_{c'}$ for subracks $c'\subseteq c$. This leads to Corollary \ref{cor: deduce A_c}, which when applied in the non-directional case with $c_2=\emptyset$ tells us that it suffices to verify the $(N,I)$-staircase for each $MA_{c'}$ with respect to the elements $U_1^{e_1},\ldots, U_{|c|}^{e_{|c|}}$ for each subrack $c'\subseteq c$. 

Let us look at a single subrack which we decompose into its connected components $c'=c_1'\sqcup c_2'\sqcup \cdots \sqcup c_k'$. Then by Theorem \ref{thm: LL CA_c iterated quotients}, $MA_{c'}$ is already stable with respect to only $k$ elements $U_1^{e_1},\ldots, U_k^{e_k}$ where we rearrange the indices so that $s_i$ is a transposition in the component $c_i'$. Then, by Lemma \ref{lem: add extra elements}, it suffices to show that $MA_{c'}$ stabilizes in an $(N,I)$-staircase with respect to $U_1^{e_1},\ldots, U_k^{e_k}$ for all sufficiently large $e_1,\ldots, e_k$. 

At first sight, it seems impossible to show this via computation as we need to verify it for all sufficiently large $e_1,\ldots, e_k$. We first claim that it suffices to show that $MA_{c'}$ stabilizes \textit{strongly} in an $(N,I)$-staircase with respect to just $U_1,\ldots, U_k$. Indeed, one can see from the long exact sequence argument in the proof of Lemma \ref{lem: homological stability for powers} that we can raise each $U_i$ to an arbitrary power. This is still an infinite amount of computation, but thankfully our results in Section \ref{sec: S_d} already give an explicit stability slope for each $MA_{c'}$, which reduces this to a \textit{finite} computation. Although we need to compute to a rather large $n$ here, our algorithm takes time polynomial in $n$ (for fixed homological degree), which makes this feasible.

Hence, to summarize, we first need to check that $\MHur^{S_d,c^{\text{trn}}_{S_d}}_n=\CHur^{S_d,c^{\text{trn}}_{S_d}}_n$ stabilizes with respect to $U_1^{e_1}$ for sufficiently large powers $e_1$. This is straightforward from the homological computations in Figure \ref{fig: homology computation S_d}. Then, for each subrack $c'\subseteq c^{\text{trn}}_{S_d}$, which we saw in Section \ref{sec: S_d} is necessarily of the form $c'=c_{S_{d_1}}\sqcup \cdots \sqcup c_{S_{d_k}}$ for $d_1+\cdots +d_k\le d$, we need to check that $MA_{c'}$ stabilizes strongly in an $(N,I)$-staircase with respect to $U_1,\ldots, U_k$ (again with rearranged indices). This is given in Figures \ref{fig: compute MA quotient 1} and \ref{fig: compute MA quotient 2}. 

Now that we have stability in an $(N,I)$-staircase, we apply Theorem \ref{thm: general main theorem} to deduce a stability slope of $s=\frac{I+1}{N+1}$ for $A_c$. However, our arithmetic application in Figure \ref{fig: improved c_d} requires the stability slope for $MA_c$ instead of $A_c$. We can obtain this by doing the usual yoga with the filtration \eqref{eqn: filtration A_c and MA_c}, which we describe here. Applying Corollary \ref{cor: deduce A_c} and Theorem \ref{thm: general main theorem} in the same way gives a stability slope $s$ for $A_{c'}$ for each subrack $c'\subseteq c$. We then induct on $|c'|$, applying Corollary \ref{cor: deduce MA_c} to staircases just below a common stability line for $A_{c'}$ and the proper maximal pieces. The extra degree required on the proper pieces shifts the line down by one, so each $MA_{c'}$ (including $MA_c$) still has slope $s$. We can then remove the excess stability elements for $MA_{c}$ using Lemma \ref{lem: remove extra}.

Finally, we briefly discuss how we compute the iterated quotient. The code is available in \textit{hurwitz\_quotients\_cpp}, with the full detailed explanation in the writeup. The main idea is to reuse the minimal resolution constructed above to obtain a much smaller chain complex. Tensoring the minimal resolution with $\Q$ over $\mathfrak A$ gives a complex with zero differential, whose graded pieces are precisely the spaces $E_{p,n}$. To compute the derived quotient by the stabilization classes $s_1^2,\ldots,s_k^2$, we use the differentials of the original resolution to compute the corresponding stabilization maps on this smaller model, together with higher homotopies encoding their commutativity. We then assemble these maps into a homotopy-coherent iterated mapping cone, whose homology computes the desired quotient. The main computational advantage is that we again only need a bounded-width portion of the minimal resolution, avoiding the exponentially large bar complex.

\subsection{Other Hurwitz spaces}\label{sec: compute other hurwitz spaces}
Finally, we compute the homology of some other Hurwitz spaces in Figure \ref{fig: computations other hurwitz}, which proves the bounds in Figures \ref{fig: table G A_4} and \ref{fig: table G semidirect}. All of these conjugacy classes are non-splitting, so it is easy to deduce the bounds from the dimensions of homology via the same arguments used earlier.

\begin{figure}[!htbp]
    \centering
    \begin{tikzpicture}[x=1cm,y=1cm,scale=0.54]
\drawnaxis{28}
\drawiaxis{7}
\drawstaircase{1}{8}{0}{black!10}{0}
\drawblock{9}{28}{0}{7}{blue!20}{0}
\drawblock{4}{28}{0}{0}{blue!20}{12}
\drawblock{10}{28}{1}{1}{blue!20}{12}
\foreach \n/\i/\txt in {
1/0/4,
2/0/8,2/1/8,
3/0/11,3/2/9,
4/3/8,
8/2/4,8/3/4,
9/2/1,
13/3/8,
14/3/4,14/4/4,14/5/4,14/6/8,
15/3/1,15/4/1,
19/5/8,
20/5/4,20/6/4,20/7/4,
21/5/1,21/6/1,
25/7/8,
26/7/4,
27/7/1
 }{
    \cell{\n}{\i}{\txt}{black!10}
  }

\foreach \n/\i/\txt in {
3/1/20,
4/1/24,4/2/20,
5/1/24,5/2/24,5/3/40,5/4/28,
6/1/23,6/2/21,6/3/40,
7/1/20,7/2/12,7/3/16,7/4/68,
8/1/16,8/4/12,8/5/88,
9/1/13,9/3/10,9/4/21,9/5/44,
10/3/12,10/4/24,10/5/40,
11/3/12,11/4/24,11/5/36,11/6/84,
12/3/11,12/4/21,12/5/31,12/6/74,
13/4/12,13/5/16,13/6/44,
14/7/12,
15/5/10,15/6/20,15/7/30,
16/5/12,16/6/24,16/7/36,
17/5/12,17/6/24,17/7/36,
18/5/11,18/6/21,18/7/31,
19/6/12,19/7/16,
21/7/10,
22/7/12,
23/7/12,
24/7/11
 }{
    \cell{\n}{\i}{\txt}{black!10}
  }

\foreach \n/\i/\txt in {
6/4/141,6/5/111,
7/5/220,7/6/164,
8/6/500,8/7/420,
9/6/208,
10/6/120,10/7/528,
11/7/288,
12/7/191,
13/7/108
 }{
    \smallcell{\n}{\i}{\txt}{black!10}
  }

  \foreach \n/\i/\txt in {
9/7/2623
 }{
    \tinycell{\n}{\i}{\txt}{black!10}
  }

\end{tikzpicture}
    \begin{tikzpicture}[x=1cm,y=1cm,scale=0.45]
\drawnaxis{17}
\drawiaxis{6}
\drawstaircase{1}{7}{0}{black!10}{0}
\drawblock{8}{17}{0}{6}{blue!20}{0}

\foreach \n/\i/\txt in {
1/0/5,
5/2/5,
8/2/1,
9/4/5,
12/4/3,
16/6/5
 }{
    \cell{\n}{\i}{\txt}{black!10}
  }

\foreach \n/\i/\txt in {
2/1/10,
3/1/20,3/2/10,
4/1/21,4/2/30,4/3/19,
5/3/50,5/4/45,
6/2/15,6/3/60,
7/2/10,7/3/60,
8/3/25,8/4/91,
9/5/35,
10/4/40,
11/4/25,
12/5/46,
13/6/20,
14/6/85,
15/6/50
 }{
    \smallcell{\n}{\i}{\txt}{black!10}
  }

\foreach \n/\i/\txt in {
2/0/10,
3/0/10,
4/0/10,
5/0/10,5/1/10,
6/0/10,6/1/10,
7/0/10,7/1/10,
8/0/10,8/1/10,
9/0/10,9/1/10,
10/0/10,10/1/10,
11/0/10,11/1/10,
12/0/10,12/1/10,
13/0/10,13/1/10,
14/0/10,14/1/10,
15/0/10,15/1/10,
16/0/10,16/1/10,
17/0/10,17/1/10
 }{
    \smallcell{\n}{\i}{\txt}{blue!20}
  }

\foreach \n/\i/\txt in {
6/4/335,6/5/290,
7/4/145,7/5/630,7/6/535,
8/5/289,
9/6/285,
10/5/215,10/6/650,
11/5/160,11/6/545,
12/6/188
 }{
    \tinycell{\n}{\i}{\txt}{black!10}
  }

  \foreach \n/\i/\txt in {
8/6/2874
 }{
    \verytinycell{\n}{\i}{\txt}{black!10}
  }

\end{tikzpicture}
\begin{tikzpicture}[x=1cm,y=1cm,scale=0.45]
\drawnaxis{10}
\drawiaxis{4}
\drawstaircase{1}{5}{0}{black!10}{0}
\drawblock{6}{10}{0}{4}{blue!20}{0}

\foreach \n/\i/\txt in {
1/0/7,
6/2/1,
9/4/7
 }{
    \cell{\n}{\i}{\txt}{black!10}
  }

\foreach \n/\i/\txt in {
2/1/14,
3/1/28,3/2/14,
4/2/56,4/3/56,
5/2/28,
6/3/83,
7/4/77,
 }{
    \smallcell{\n}{\i}{\txt}{black!10}
  }

\foreach \n/\i/\txt in {
2/0/14,
3/0/14,
4/0/14,4/1/14,
5/0/14,5/1/14,
6/0/14,6/1/14,
7/0/14,7/1/14,
8/0/14,8/1/14,
9/0/14,9/1/14,
10/0/14,10/1/14
 }{
    \smallcell{\n}{\i}{\txt}{blue!20}
  }

\foreach \n/\i/\txt in {
5/3/364,5/4/336,
8/4/133
 }{
    \tinycell{\n}{\i}{\txt}{black!10}
  }

  \foreach \n/\i/\txt in {
6/4/1189
 }{
    \verytinycell{\n}{\i}{\txt}{black!10}
  }

\end{tikzpicture}
\begin{tikzpicture}[x=1cm,y=1cm,scale=0.45]
\drawnaxis{19}
\drawiaxis{4}
\drawstaircase{1}{5}{0}{black!10}{0}
\drawblock{6}{19}{0}{4}{blue!20}{0}

\foreach \n/\i/\txt in {
1/0/7,
7/2/7,
12/2/3,
18/4/5
 }{
    \cell{\n}{\i}{\txt}{black!10}
  }

\foreach \n/\i/\txt in {
2/0/21,2/1/21,
3/0/35,3/1/77,3/2/42,
4/0/35,
5/0/21,5/1/63,
6/0/15,6/1/34,6/2/59,
7/3/28,7/4/63,
8/2/49,
9/2/77,
10/2/63,
11/2/21,11/3/91,
12/3/28,12/4/60,
13/4/21,
14/4/91,
17/4/35
 }{
    \smallcell{\n}{\i}{\txt}{black!10}
  }

\foreach \n/\i/\txt in {
7/0/14,7/1/14,
8/0/14,8/1/14,
9/0/14,9/1/14,
10/0/14,10/1/14,
11/0/14,11/1/14,
12/0/14,12/1/14,
13/0/14,13/1/14,
14/0/14,14/1/14,
15/0/14,15/1/14,
16/0/14,16/1/14,
17/0/14,17/1/14,
18/0/14,18/1/14,
19/0/14,19/1/14
 }{
    \smallcell{\n}{\i}{\txt}{blue!20}
  }

\foreach \n/\i/\txt in {
4/1/105,4/2/182,4/3/112,
5/2/140,5/3/308,5/4/210,
6/3/200,
8/3/175,8/4/441,
9/3/273,9/4/588,
10/3/231,10/4/441,
11/4/175,
15/4/133,
16/4/105
 }{
    \tinycell{\n}{\i}{\txt}{black!10}
  }

  \foreach \n/\i/\txt in {
6/4/1415
 }{
    \verytinycell{\n}{\i}{\txt}{black!10}
  }

\end{tikzpicture}
\begin{tikzpicture}[x=1cm,y=1cm,scale=0.45]
\drawnaxis{6}
\drawiaxis{3}
\drawstaircase{1}{4}{0}{black!10}{0}
\drawblock{5}{6}{0}{3}{blue!20}{0}

\foreach \n/\i/\txt in {
1/0/8
 }{
    \cell{\n}{\i}{\txt}{black!10}
  }

\foreach \n/\i/\txt in {
2/1/16,
3/1/32,3/2/16,
4/2/32,4/3/32
 }{
    \smallcell{\n}{\i}{\txt}{black!10}
  }

\foreach \n/\i/\txt in {
2/0/16,
3/0/16,
4/0/16,4/1/16,
5/0/16,5/1/16,
6/0/16,6/1/16
 }{
    \smallcell{\n}{\i}{\txt}{blue!20}
  }

\foreach \n/\i/\txt in {
5/3/256
 }{
    \tinycell{\n}{\i}{\txt}{black!10}
  }

\end{tikzpicture}
\begin{tikzpicture}[x=1cm,y=1cm,scale=0.45]
\drawnaxis{15}
\drawiaxis{3}
\drawstaircase{1}{4}{0}{black!10}{0}
\drawblock{5}{15}{0}{3}{blue!20}{0}

\foreach \n/\i/\txt in {
1/0/9,
8/2/6,
14/3/6
 }{
    \cell{\n}{\i}{\txt}{black!10}
  }

\foreach \n/\i/\txt in {
2/0/21,2/1/21,
3/1/99,3/2/72,
4/1/75,
5/1/63,5/2/81,
6/1/57,6/2/61,
7/1/45,7/2/36,7/3/72,
8/1/33,8/3/30,
9/3/18,
10/3/30,
11/3/36,
12/3/30,
13/3/18
 }{
    \smallcell{\n}{\i}{\txt}{black!10}
  }

\foreach \n/\i/\txt in {
3/0/27,
4/0/27,
5/0/27,
6/0/27,
7/0/27,
8/0/27,
9/0/27,9/1/27,
10/0/27,10/1/27,
11/0/27,11/1/27,
12/0/27,12/1/27,
13/0/27,13/1/27,
14/0/27,14/1/27,
15/0/27,15/1/27
 }{
    \smallcell{\n}{\i}{\txt}{blue!20}
  }

\foreach \n/\i/\txt in {
4/2/186,4/3/138,
5/3/576,
6/3/135
 }{
    \tinycell{\n}{\i}{\txt}{black!10}
  }

\end{tikzpicture}
\begin{tikzpicture}[x=1cm,y=1cm,scale=0.45]
\drawnaxis{9}
\drawiaxis{3}
\drawstaircase{1}{4}{0}{black!10}{0}
\drawblock{5}{9}{0}{3}{blue!20}{0}

\foreach \n/\i/\txt in {
1/0/9,
8/2/3
 }{
    \cell{\n}{\i}{\txt}{black!10}
  }

\foreach \n/\i/\txt in {
2/0/27,2/1/27,
3/1/81,3/2/45,
5/1/63,5/2/81,
6/2/54,
7/2/45
 }{
    \smallcell{\n}{\i}{\txt}{black!10}
  }

\foreach \n/\i/\txt in {
3/0/36,
4/0/36,
5/0/36,
6/0/36,6/1/36,
7/0/36,7/1/36,
8/0/36,8/1/36,
9/0/36,9/1/36
 }{
    \smallcell{\n}{\i}{\txt}{blue!20}
  }

\foreach \n/\i/\txt in {
4/1/106,4/2/275,4/3/205,
5/3/846,
6/3/414,
7/3/396,
8/3/159
 }{
    \tinycell{\n}{\i}{\txt}{black!10}
  }

\end{tikzpicture}
\begin{tikzpicture}[x=1cm,y=1cm,scale=0.45]
\drawnaxis{15}
\drawiaxis{3}
\drawstaircase{1}{4}{0}{black!10}{0}
\drawblock{5}{15}{0}{3}{blue!20}{0}

\foreach \n/\i/\txt in {
1/0/9,
12/2/3
 }{
    \cell{\n}{\i}{\txt}{black!10}
  }

\foreach \n/\i/\txt in {
2/0/33,2/1/33,
3/0/63,3/2/90,
4/0/69,
5/0/54,
6/0/46,
7/1/90,
8/1/57,8/2/60,
9/2/90,
10/2/78,
11/2/27,
13/3/36,
14/3/12
 }{
    \smallcell{\n}{\i}{\txt}{black!10}
  }

\foreach \n/\i/\txt in {
7/0/45,
8/0/45,
9/0/45,9/1/45,
10/0/45,10/1/45,
11/0/45,11/1/45,
12/0/45,12/1/45,
13/0/45,13/1/45,
14/0/45,14/1/45,
15/0/45,15/1/45
 }{
    \smallcell{\n}{\i}{\txt}{blue!20}
  }

\foreach \n/\i/\txt in {
3/1/153,
4/1/249,4/2/492,4/3/312,
5/1/189,5/2/486,
6/1/146,6/2/370,
7/2/117,7/3/405,
8/3/309,
9/3/414,
10/3/384,
11/3/207,
12/3/110
 }{
    \tinycell{\n}{\i}{\txt}{black!10}
  }

\foreach \n/\i/\txt in {
5/3/1116,
6/3/1474
 }{
    \verytinycell{\n}{\i}{\txt}{black!10}
  }

\end{tikzpicture}
\begin{tikzpicture}[x=1cm,y=1cm,scale=0.45]
\drawnaxis{7}
\drawiaxis{3}
\drawstaircase{1}{4}{0}{black!10}{0}
\drawblock{5}{7}{0}{3}{blue!20}{0}

\foreach \n/\i/\txt in {
1/0/9
 }{
    \cell{\n}{\i}{\txt}{black!10}
  }

\foreach \n/\i/\txt in {
2/1/18,
3/1/36,3/2/18,
5/2/27,
6/3/18
 }{
    \smallcell{\n}{\i}{\txt}{black!10}
  }

\foreach \n/\i/\txt in {
2/0/18,
3/0/18,
4/0/18,4/1/18,
5/0/18,5/1/18,
6/0/18,6/1/18,
7/0/18,7/1/18
 }{
    \smallcell{\n}{\i}{\txt}{blue!20}
  }

\foreach \n/\i/\txt in {
4/2/288,4/3/288,
5/3/432
 }{
    \tinycell{\n}{\i}{\txt}{black!10}
  }

\end{tikzpicture}
    \caption{Homology of $\Hur^{G,c}_n$ for $G,c$ in Figures \ref{fig: table G A_4} and \ref{fig: table G semidirect} in order.}
    \label{fig: computations other hurwitz}
\end{figure}

\clearpage
\newpage

\printbibliography
\end{document}